\documentclass[11pt]{nwsd} 
\mathchardef\varSigma="0106  
\usepackage{amssymb,latexsym,amsthm} 
\font\medit=ptmri at 11pt 
\font\smallit=ptmri at 10pt 
\font\medsmallit=ptmri at 9pt

\font\medbf=ptmbo at 12pt

\begin{document} 
\theoremstyle{plain} 
\newtheorem{theorem}{\bf Theorem}[section] 
\newtheorem{lemma}[theorem]{\bf Lemma} 
\newtheorem{corollary}[theorem]{\bf Corollary} 
\newtheorem{proposition}[theorem]{\bf Proposition} 
\theoremstyle{definition} 
\newtheorem{definition}[theorem]{\sc Definition} 
\newtheorem{remark}[theorem]{\medbf Remark} 
\newtheorem{example}[theorem]{\medbf Example} 
\renewcommand{\proofname}{\medbf Proof.} 
\renewcommand{\theequation}{\thesection.\arabic{equation}} 
\mathchardef\varGamma="0100 
\mathchardef\varTheta="0102 
\mathchardef\varLambda="0103 
\mathchardef\varXi="0104 
\mathchardef\varPi="0105 
\mathchardef\varSigma="0106 
\mathchardef\varPhi="0108 
\mathchardef\varPsi="0109 
\mathchardef\varOmega="010A 
\def\dpi{\hbox{$(\hskip-.6ptd\pi\hskip-1pt)$}} 
\def\nft 
{\hbox{$n$\hskip3pt$\equiv$\hskip4pt$5$\hskip4.4pt$($mod\hskip2pt$3)$}} 
\def\bbC{{\mathchoice {\setbox0=\hbox{$\displaystyle\mathrm{C}$}\hbox{\hbox 
to0pt{\kern0.4\wd0\vrule height0.9\ht0\hss}\box0}} 
{\setbox0=\hbox{$\textstyle\mathrm{C}$}\hbox{\hbox 
to0pt{\kern0.4\wd0\vrule height0.9\ht0\hss}\box0}} 
{\setbox0=\hbox{$\scriptstyle\mathrm{C}$}\hbox{\hbox 
to0pt{\kern0.4\wd0\vrule height0.9\ht0\hss}\box0}} 
{\setbox0=\hbox{$\scriptscriptstyle\mathrm{C}$}\hbox{\hbox 
to0pt{\kern0.4\wd0\vrule height0.9\ht0\hss}\box0}}}} 
\def\bbP{\mathrm{I\!P}} 
\def\bbR{\mathrm{I\!R}} 
\def\rto{\bbR\hskip-.5pt^2} 
\def\rtr{\bbR\hskip-.7pt^3} 
\def\rn{{\bbR}^{\hskip-.6ptn}} 
\def\rk{{\bbR}^{\hskip-.6ptk}} 
\def\bbZ{\mathbf{Z}} 
\def\hyp{\hskip.5pt\vbox 
{\hbox{\vrule width3ptheight0.5ptdepth0pt}\vskip2.2pt}\hskip.5pt} 
\def\er{k} 
\def\df{d\hskip-.8ptf} 
\def\dfh{d\hskip-.8pt\fh} 
\def\trf{[\hskip-1.7pt[\hskip-.4ptF\hskip.8pt]\hskip-1.7pt]} 
\def\trfk{[\hskip-1.7pt[\hskip-.4ptF^K\hskip0pt]\hskip-1.7pt]} 
\def\trfo{[\hskip-2pt[\hskip-.4ptF_0]\hskip-2pt]} 
\def\trz{[\hskip-2pt[\hskip.4pt\zeta\hskip.7pt]\hskip-2pt]} 
\def\trfz{[\hskip-2pt[\hskip.8ptf\zeta\hskip.7pt]\hskip-2pt]} 
\def\tb{T\hskip.3pt\bs} 
\def\tab{{T\hskip.2pt^*\hskip-.9pt\bs}} 
\def\tayb{{T_y^*\hskip-1pt\bs}} 
\def\tyb{{T\hskip-2pt_y\hskip-.2pt\bs}} 
\def\txm{{T\hskip-2.3pt_x\hskip-.5ptM}} 
\def\tym{{T\hskip-2pt_y\hskip-.9ptM}} 
\def\txn{{T\hskip-2pt_x\hskip-.9ptN}} 
\def\txtm{{T\hskip-2pt_{x(t)}\hskip-.9ptM}} 
\def\txom{{T\hskip-2pt_{x(0)}\hskip-.9ptM}} 
\def\txtsm{{T\hskip-2pt_{x(t,s)}\hskip-.9ptM}} 
\def\olm{\hskip2.1pt\overline{\hskip-2.1ptM}} 
\def\bw{\hskip1.8pt\overline{\hskip-1.8ptw\hskip-.8pt}\hskip.8pt} 
\def\bu{\hskip1pt\overline{\hskip-1ptu\hskip-.4pt}\hskip.4pt} 
\def\bv{\hskip.8pt\overline{\hskip-.8ptv\hskip-.4pt}\hskip.4pt} 
\def\bc{\hskip1pt\overline{\hskip-1pt\mathrm{D}\hskip-1.8pt}\hskip1.8pt} 
\def\br{\hskip3.3pt\overline{\hskip-3.3ptR\hskip-.7pt}\hskip.7pt} 
\def\hd{\widehat{\mathrm{D}}} 
\def\tm{{T\hskip-.3ptM}} 
\def\tn{{T\hskip-.3ptN}} 
\def\tu{{T\hskip-.3ptU}} 
\def\tam{{T^*\!M}} 
\def\so{\mathfrak{so}\hh} 
\def\gi{\mathfrak{g}} 
\def\hi{\mathfrak{h}} 
\def\fh{f} 
\def\hf{\varLambda} 
\def\line{\Lambda} 
\def\plane{\Pi} 
\def\cst{^{\mathrm{cst}}} 
\def\lin{^{\mathrm{lin}}} 
\def\qdr{^{\mathrm{qdr}}} 
\def\cub{^{\mathrm{cub}}} 
\def\qrt{^{\mathrm{qrt}}} 
\def\qnt{^{\mathrm{qnt}}} 
\def\rc{c} 
\def\rd{X} 
\def\kri{\mathrm{Ker}\hskip1pt\ri} 
\def\kr{\mathrm{Ker}\hskip1ptR} 
\def\kw{\mathrm{Ker}\hskip2ptW} 
\def\kb{\mathrm{Ker}\hskip1.5ptB} 
\def\xc{\mathcal{X}_c} 
\def\fd{F\hskip-2pt.} 
\def\yw{Y\hskip-2.7pt_w} 
\def\ta{\hskip.7pt\widetilde{\hskip-.7pt\alpha\hskip-.4pt}\hskip.4pt} 
\def\tga{\hskip1pt\widetilde{\hskip-1pt\gamma\hskip-.8pt}\hskip.8pt} 
\def\tet{\hskip.7pt\widetilde{\hskip-.7pt\eta\hskip-.4pt}\hskip.4pt} 
\def\tth{\hskip.7pt\widetilde{\hskip-.7pt\theta\hskip-.4pt}\hskip.4pt} 
\def\tg{\hskip2pt\widetilde{\hskip-2ptg\hskip-.4pt}\hskip.4pt} 
\def\tts{\hskip.7pt\widetilde{\hskip-.7pt\ts\hskip-.4pt}\hskip.4pt} 
\def\thor{\hskip1pt\widetilde{\hskip-1pt\mathcal{H}\hskip1pt}\hskip-1pt} 
\def\tna{\hskip1pt\widetilde{\hskip-1pt\nabla\hskip-1pt}\hskip1pt} 
\def\tB{\hskip2.3pt\widetilde{\hskip-2.3ptB\hskip.7pt}\hskip-.7pt} 
\def\tR{\hskip.7pt\widetilde{\hskip-.7ptR\hskip-.4pt}\hskip.4pt} 
\def\tX{\hskip.7pt\widetilde{\hskip-.7pt\varXi\hskip1pt}\hskip-1pt} 
\def\tyw{\hskip.7pt\widetilde{\hskip-.7ptY\hskip-.4pt}\hskip-2.3pt_w} 
\def\tZ{\hskip.7pt\widetilde{\hskip-.7ptZ\hskip-.4pt}\hskip.4pt} 
\def\tvt{\hskip2pt\widetilde{\hskip-2pt\vt\hskip.6pt}\hskip-.6pt} 
\def\tw{{\hskip.7pt\widetilde{\hskip-.7ptw\hskip-.4pt}\hskip.4pt}} 
\def\hw{{\hskip.7pt\widehat{\hskip-.7ptw\hskip-.4pt}\hskip.4pt}} 
\def\trd{\hskip.7pt\widetilde{\hskip-.7pt\rd\hskip-.4pt}\hskip.4pt} 
\def\xe{\mathcal{E}} 
\def\as{\mathcal{S}} 
\def\xz{\mathcal{X}} 
\def\yz{\mathcal{Y}} 
\def\zz{\mathcal{Z}} 
\def\tim{\hskip1.5pt\widetilde{\hskip-1.5ptM\hskip-.5pt}\hskip.5pt} 
\def\hm{\hskip1.9pt\widehat{\hskip-1.9ptM\hskip-.2pt}\hskip.2pt} 
\def\hmt{\hskip1.9pt\widehat{\hskip-1.9ptM\hskip-.5pt}_t} 
\def\hmz{\hskip1.9pt\widehat{\hskip-1.9ptM\hskip-.5pt}_0} 
\def\hmp{\hskip1.9pt\widehat{\hskip-1.9ptM\hskip-.5pt}_p} 
\def\hg{\hskip1.2pt\widehat{\hskip-1.2ptg\hskip-.4pt}\hskip.4pt} 
\def\nao{\hbox{$\nabla\!\!^{^{^{_{\!\!\circ}}}}$}} 
\def\ro{\hbox{$R\hskip-4.5pt^{^{^{_{\circ}}}}$}{}} 
\def\mppp{\hbox{$-$\hskip1pt$+$\hskip1pt$+$\hskip1pt$+$}} 
\def\mpdp{\hbox{$-$\hskip1pt$+$\hskip1pt$\dots$\hskip1pt$+$}} 
\def\mmpp{\hbox{$-$\hskip1pt$-$\hskip1pt$+$\hskip1pt$+$}} 
\def\mmp{\hbox{$-$\hskip1pt$-$\hskip1pt$+$}} 
\def\mmmp{\hbox{$-$\hskip1pt$-$\hskip1pt$-$\hskip1pt$+$}} 
\def\pppp{\hbox{$+$\hskip1pt$+$\hskip1pt$+$\hskip1pt$+$}} 
\def\mpmp{\hbox{$-$\hskip1pt$\pm$\hskip1pt$+$}} 
\def\mpmpp{\hbox{$-$\hskip1pt$\pm$\hskip1pt$+$\hskip1pt$+$}} 
\def\mmpmp{\hbox{$-$\hskip1pt$-$\hskip1pt$\pm$\hskip1pt$+$}} 
\def\q{q} 
\def\bq{\hat q} 
\def\p{p} 
\def\x{v} 
\def\y{y} 
\def\vdx{\vd{}\hskip-4.5pt_x} 
\def\bz{b\hh} 
\def\cy{{y}} 
\def\rkwo{\,\hs\mathrm{rank}\hskip2.7ptW\hskip-2.7pt=\hskip-1.2pt1} 
\def\rkwho{\,\hs\mathrm{rank}\hskip2.2ptW^h\hskip-2.2pt=\hskip-1pt1} 
\def\rkw{\,\hs\mathrm{rank}\hskip2.4ptW\hskip-1.5pt} 
\def\rw{\varPsi} 
\def\nd{\Phi} 
\def\ft{\Psi} 
\def\js{J} 
\def\ism{H} 
\def\fe{F} 
\def\fy{f} 
\def\dfc{dF\hskip-2.3pt_\cy\hskip.4pt} 
\def\dfct{dF\hskip-2.3pt_\cy(t)\hskip.4pt} 
\def\dic{d\im\hskip-1.4pt_\cy\hskip.4pt} 
\def\vl{\Lambda} 
\def\qt{\mathcal{E}} 
\def\tqt{\tilde{\qt}} 
\def\vh{h} 
\def\mv{V} 
\def\vy{\mathcal{V}} 
\def\xv{\mathcal{X}} 
\def\yv{\mathcal{Y}} 
\def\iv{\mathcal{I}} 
\def\gkp{\Sigma} 
\def\bs{\Sigma} 
\def\dbs{\dot\bs} 
\def\das{\hskip3pt\dot{\hskip-3pt\as}} 
\def\hs{\hskip.7pt} 
\def\hh{\hskip.4pt} 
\def\nh{\hskip-.7pt} 
\def\hn{\hskip-.4pt} 
\def\nnh{\hskip-1pt} 
\def\hrz{^{\hskip.5pt\mathrm{hrz}}} 
\def\vrt{^{\hskip.2pt\mathrm{vrt}}} 
\def\vt{\varTheta} 
\def\zh{\zeta} 
\def\vg{\varGamma} 
\def\gm{\gamma} 
\def\gp{\mathrm{G}} 
\def\hp{\mathrm{H}} 
\def\kp{\mathrm{K}} 
\def\Gm{\Gamma} 
\def\Lm{\Lambda} 
\def\Dt{\Delta} 
\def\sj{\sigma} 
\def\lg{\langle} 
\def\rg{\rangle} 
\def\lr{\lg\,,\rg} 
\def\uv{\underline{v\hskip-.8pt}\hskip.8pt} 
\def\uvp{\underline{v\hh'\hskip-.8pt}\hskip.8pt} 
\def\uw{\underline{w\hskip-.8pt}\hskip.8pt} 
\def\uxs{\underline{x_s\hskip-.8pt}\hskip.8pt} 
\def\vs{vector space} 
\def\vf{{q}} 
\def\tf{tensor field} 
\def\tvn{the vertical distribution} 
\def\dn{distribution} 
\def\od{Ol\-szak distribution} 
\def\pt{point} 
\def\tc{tor\-sion\-free connection} 
\def\rt{Ric\-ci tensor} 
\def\pde{partial differential equation} 
\def\pf{projectively flat} 
\def\pfs{projectively flat surface} 
\def\pfc{projectively flat connection} 
\def\pftc{projectively flat tor\-sion\-free connection} 
\def\su{surface} 
\def\sco{simply connected} 
\def\psr{pseu\-d\hbox{o\hs-}\hskip0ptRiem\-ann\-i\-an} 
\def\inv{-in\-var\-i\-ant} 
\def\trinv{trans\-la\-tion\inv} 
\def\feo{dif\-feo\-mor\-phism} 
\def\feic{dif\-feo\-mor\-phic} 
\def\feicly{dif\-feo\-mor\-phi\-cal\-ly} 
\def\Feicly{Dif\-feo\-mor\-phi\-cal\-ly} 
\def\nw{\hbox{non\hh-\hskip-1.5pt}\hskip0ptWalk\-er} 
\def\itnw{\hbox{non\hh-\hskip-1.8pt}\hskip0ptWalk\-er} 
\def\diml{-di\-men\-sion\-al} 
\def\prl{-par\-al\-lel} 
\def\skc{skew-sym\-met\-ric} 
\def\sky{skew-sym\-me\-try} 
\def\Sky{Skew-sym\-me\-try} 
\def\dbly{-dif\-fer\-en\-ti\-a\-bly} 
\def\cs{con\-for\-mal\-ly symmetric} 
\def\cf{con\-for\-mal\-ly flat} 
\def\ls{locally symmetric} 
\def\kf{Killing field} 
\def\om{\omega} 
\def\vol{\varOmega} 
\def\dv{\delta} 
\def\ve{\varepsilon} 
\def\zt{\zeta} 
\def\kx{\kappa} 
\def\mf{manifold} 
\def\mfd{-man\-i\-fold} 
\def\bmf{base manifold} 
\def\bd{bundle} 
\def\ed{\Phi} 
\def\bna{\hs\overline{\nh\nabla\nh}\hs} 
\def\gs{s} 
\def\hj{t} 
\def\ea{A} 
\def\qx{Q} 
\def\lx{L} 
\def\lp{\lx\nnh^+} 
\def\lm{\lx\nnh^-} 
\def\lmo{\lx_1^{\hskip-1.3pt-}} 
\def\lpt{\lx_2^{\hskip-1.3pt+}} 
\def\ex{E} 
\def\ty{T} 
\def\gy{G} 
\def\ny{N\nh} 
\def\oy{{I\hskip-3.3ptP\hskip-.8pt}} 
\def\apb{{\mathcal{A}}} 
\def\alb{{\mathcal{J}}} 
\def\pb{{\mathcal{P}}} 
\def\lb{{\mathcal{L}}} 
\def\tl{{\mathcal{D}}} 
\def\fb{{\mathcal{W}}} 
\def\tlp{\tl^\perp} 
\def\tbd{tangent bundle} 
\def\ctb{cotangent bundle} 
\def\bp{bundle projection} 
\def\pr{pseu\-d\hbox{o\hs-}\hskip0ptRiem\-ann\-i\-an} 
\def\prc{pseu\-d\hbox{o\hs-}\hskip0ptRiem\-ann\-i\-an metric} 
\def\prd{pseu\-d\hbox{o\hs-}\hskip0ptRiem\-ann\-i\-an manifold} 
\def\Prd{pseu\-d\hbox{o\hs-}\hskip0ptRiem\-ann\-i\-an manifold} 
\def\npd{null parallel distribution} 
\def\pj{-pro\-ject\-a\-ble} 
\def\pd{-pro\-ject\-ed} 
\def\lcc{Le\-vi-Ci\-vi\-ta connection} 
\def\vb{vector bundle} 
\def\vbm{vec\-tor-bun\-dle morphism} 
\def\kerd{\mathrm{Ker}\hskip2.7ptd} 
\def\ro{\rho} 
\def\sy{\sigma} 
\def\ts{\tau} 
\def\pmb{\pi} 
\def\ri{{\rho}}

\title[{{\smallit Self\hskip.4pt-\hskip.4ptDu\-al Neutral Ein\-stein 
Four\hskip.3pt-\hskip.3ptMan\-i\-folds}}]{Non\hs-\hskip-2.7ptWalk\-er 
Self-\hskip-.4ptDu\-al Neutral Ein\-stein Four-Man\-i\-folds of Pe\-trov Type 
III} 
\author[{{\smallit A. Derdzinski}}]{\medit Andrzej Derdzinski} 
\address{Department of Mathematics, The Ohio State University, 
Columbus, OH 43210, USA} 
\email{andrzej@math.ohio-state.edu} 
\begin{abstract} 
{\smallit ABSTRACT. \ \ The local structure of the manifolds named in the 
title is described. Although curvature homogeneous, they are not, in general, 
locally homogeneous. Not all of them are Ric\-\hbox{ci\hh-}\hskip0ptflat, 
which answers an existence question about type III Jor\-dan-Osser\-man 
metrics, raised by \hbox{D\'\i az\hh-}\hskip0ptRamos, 
Gar\-\hbox{c\'\i a\hh-}\hskip0ptR\'\i o and 
V\'az\-\hbox{quez\hh-}\hskip0ptLo\-ren\-zo (2006).} 
\end{abstract}

\subjclass{53B30 (Primary); 53B05 (Secondary)} 
\keywords{ 
Cur\-va\-ture\hskip.8pt-ho\-mo\-ge\-ne\-ous neutral Ein\-stein metric, type 
III Osser\-man metric} 
\acknow{Work begun during the author's visit to the University of 
San\-tia\-go de Com\-pos\-tela, supported by Grant MTM2006-01432 (Spain)}

\maketitle

\setcounter{theorem}{0}

\voffset=0pt\hoffset=0pt 

\baselineskip14pt 
\section{Introduction}\label{intr} 
\setcounter{equation}{0} 
The main result of this paper, Theorem~\ref{lcstr}, describes the local 
structure of all \nw\ self-du\-al oriented Ein\-stein 
\hbox{four\hh-}\hskip0ptman\-i\-folds $\,(M,g)\,$ of the neutral metric 
signature $\,\mmpp\,$ which are of Pe\-trov type III, in the sense that so 
is the self-dual Weyl tensor at every point of $\,M$. Such $\,(M,g)$, also 
referred to as the {\it \itnw, type\/} III {\it 
\hbox{four\hh-}\hskip0ptdi\-men\-sion\-al Jor\-\hbox{dan\hs-}Osser\-man 
manifolds\/} \cite[Remark 2.1]{diaz-ramos-garcia-rio-vazquez-lorenzo}, are 
known to be curvature homogeneous, cf.\ Remark~\ref{cvhmg}.

According to Theorem~\ref{lcstr}, all $\,(M,g)\,$ with the stated properties 
are, locally, pa\-ram\-e\-trized by arbitrary solutions to equations 
(\ref{eqn}), which, by Remark~\ref{unkfc}, are equivalent to the system 
(\ref{loc}) of four first-or\-der qua\-si-lin\-e\-ar partial differential 
equations imposed on eight unknown real-val\-ued functions of two real 
variables.

A description of all solutions to (\ref{eqn}) is given in Section~\ref{stsy}. 
It applies, however, only to the dense open subset of $\,\rto$ formed by 
points which are in {\it general position}, in the sense that, on some 
neighborhood of the point in question, each of several specific 
vec\-tor-val\-ued functions associated with the solution is either identically 
zero, or nonzero everywhere.

The precise meaning of the adjective `\nw' used above is that the metrics 
$\,g\,$ are assumed to represent the {\it strictly \itnw\ case}, in which a 
certain natural $\,1$-form $\,\beta\,$ is nonzero everywhere. By contrast, the 
{\it Walker case}, defined in Section~\ref{wnwc} by requiring $\,\beta\,$ to 
vanish identically, is equivalent to the existence of a 
\hbox{two\hh-}\hskip0ptdi\-men\-sion\-al null parallel distribution compatible 
with the orientation. In the Walker case, \hbox{D\'\i az\hh-}\hskip0ptRamos, 
Gar\-\hbox{c\'\i a\hs-}\hskip0ptR\'\i o and 
V\nnh\'az\-\hbox{quez\hh-}\hskip0ptLo\-ren\-zo have already found a canonical 
coordinate form of such metrics 
\cite[Theorem~3.1(ii.3)]{diaz-ramos-garcia-rio-vazquez-lorenzo}. For a 
co\-or\-di\-nate-free version of their result, see 
\cite[Theorem~13.1]{derdzinski-08}.

The Walker case implies Ric\-\hbox{ci\hh-}\hskip0ptflat\-ness of $\,(M,g)$, 
which is why \hbox{D\'\i az\hh-}\hskip0ptRamos, 
Gar\-\hbox{c\'\i a\hs-}\hskip0ptR\'\i o and 
V\nnh\'az\-\hbox{quez\hh-}\hskip0ptLo\-ren\-zo asked in 
\cite[Remark~3.5]{diaz-ramos-garcia-rio-vazquez-lorenzo} whether a type III 
self-du\-al neutral Ein\-stein \hbox{four\hh-}\hskip0ptman\-i\-fold can have 
nonzero scalar curvature. Theorem~\ref{rfnrf} of this paper shows, by means of 
explicit examples, that the answer is `yes' while, at the same time, there 
also exist Ric\-\hbox{ci\hh-}\hskip0ptflat strictly \nw\ self-du\-al neutral 
\hbox{four\hh-}\hskip0ptman\-i\-folds, and, whether 
Ric\-\hbox{ci\hh-}\hskip0ptflat or not, such manifolds need not, in general, 
be locally homogeneous.

Besides \cite{diaz-ramos-garcia-rio-vazquez-lorenzo}, a few other papers 
contain results in this direction. Bla\-\v zi\'c, Bo\-kan and Raki\'c 
\cite{blazic-bokan-rakic} found a characterization of type III self-du\-al 
neutral Ein\-stein \hbox{four\hh-}\hskip0ptman\-i\-folds in terms of a system 
of first-or\-der differential equations imposed on the Le\-vi-Ci\-vi\-ta 
connection forms in a suitable local trivialization of the tangent bundle. 
Brans \cite{brans} described all {\it Lo\-rentz\-i\-an\/} Ein\-stein metrics 
of Pe\-trov type III in dimension four. 
Cur\-va\-ture\hskip1.3pt-ho\-mo\-ge\-ne\-ous Ein\-stein 
\hbox{four\hh-}\hskip0ptman\-i\-folds of any metric signature, with a 
curvature operator that is com\-plex-di\-ag\-o\-nal\-izable, are known to be 
locally homogeneous, and have been fully classified \cite{derdzinski-03}. In 
particular, cur\-va\-\hbox{ture\hskip1.3pt-}\hskip0ptho\-mo\-ge\-ne\-ous 
\hbox{four\hh-}\hskip0ptdi\-men\-sion\-al {\it Riemannian\/} Ein\-stein 
manifolds, which obviously satisfy the di\-ag\-o\-nal\-iza\-bil\-i\-ty 
condition, are all locally symmetric 
\cite[p.\ 476, Corollary~7.2]{dillen-verstraelen}.

For more on Osser\-man metrics and 
cur\-va\-ture\hskip1.1pt-ho\-mo\-ge\-ne\-i\-ty, see 
\cite{garcia-rio-kupeli-vazquez-lorenzo}, \cite{gilkey} and 
\cite{calvino-louzao-garcia-rio-gilkey-vazquez-lorenzo}.

I wish to thank Eduardo Garc\'\i a\hh-R\'\i o for his hospitality during my 
visit to the University of San\-tia\-go de Com\-pos\-tela in April 2007, for 
bringing to my attention the question raised in 
\cite[Remark~3.5]{diaz-ramos-garcia-rio-vazquez-lorenzo}, and for sharing with 
me his insights about it. I am also grateful to Saleh Tanveer for helpful 
comments about equations of type (\ref{rzo}).

\section{An outline of the argument}\label{outl} 
\setcounter{equation}{0} 
For any type III self-du\-al neutral oriented Ein\-stein 
\hbox{four\hh-}\hskip0ptman\-i\-fold $\,(M,g)$, it is shown in 
Lemmas~\ref{cptqt}(a)\hs--\hs(b) and~\ref{sdept}(i) 
that, locally, $\,M\,$ admits a natural identification with the total space of 
an af\-fine plane bundle over a surface, endowed with a distinguished 
``nonlinear connection'' in the form of a horizontal distribution 
$\,\mathcal{H}$, transverse to the vertical distribution $\,\mathcal{V}\,$ of 
the bundle. Both $\,\mathcal{V}\,$ and $\,\mathcal{H}\,$ consist of $\,g$-null 
vectors.

In the next step, $\,\mathcal{H}\,$ is ignored. What is kept in the picture 
consists of some natural differential forms (such as $\,\beta$, mentioned in 
the Introduction) along with the vertical distribution $\,\mathcal{V}\,$ on 
$\,M$, the family $\,\mathrm{D}\,$ of the standard flat tor\-sion\-free 
connections of the af\-fine leaves of $\,\mathcal{V}\nh$, and a partial 
version $\,h\,$ of the original metric $\,g$. Specifically, $\,h\,$ 
``remembers'' only how to evaluate inner products in which one of the vectors 
is vertical. The data just listed form a {\it basic octuple\/} defined, as an 
abstract object, in Section~\ref{boct}.

The original metric $\,g\,$ may be reconstructed from its associated basic 
octuple through a choice of a suitable horizontal distribution 
$\,\mathcal{H}$, declared to be $\,g$-null. The desired properties of $\,g\,$ 
(being a self-du\-al Ein\-stein metric of Pe\-trov type III such that the 
given $\,\mathcal{V}\,$ and $\,\mathcal{H}\,$ correspond to it as in 
Lemma~\ref{cptqt}(a)) can be rephrased as a system of four curvature 
conditions, appearing in Theorem~\ref{crvtc}, which are differential equations 
with the unknown $\,\mathcal{H}$.

What makes basic octuples convenient to use is the fact that they all 
represent a unique local \feic\ type (Theorem~\ref{unimo}), and so choosing to 
work with just one of them leads to no loss of generality. Secondly, the 
horizontal distributions for a fixed basic octuple form an af\-fine space 
(and, in fact, constitute arbitrary sections of a certain af\-fine bundle). 
The discussion of the four curvature conditions may thus be simplified by 
selecting one horizontal distribution $\,\mathcal{H}\,$ to serve as the 
origin, and expressing other horizontal distributions as sums 
$\,\thor=\mathcal{H}+F\nnh$, where $\,F\,$ is a section of a specific vector 
bundle over $\,M$.

A simplification of this kind is provided by {\it tw\hbox{o\hh-}plane 
systems}, introduced in Section~\ref{afsy}. Specifically, $\,M\,$ is 
replaced with the product $\,\bs\times\plane_+$ of an af\-fine plane $\,\bs\,$ 
and a (vector) half-plane $\,\plane_+$, so that the $\,\plane_+$ and $\,\bs\,$ 
factor distributions serve as $\,\mathcal{V}\,$ and the ``origin'' 
$\,\mathcal{H}$. Of the four curvature conditions in Theorem~\ref{crvtc}, 
imposed on $\,\thor=\mathcal{H}+F\nnh$, with $\,F\,$ as above, three then turn 
out to be (nonhomogeneous) linear, of first or second order, and involve only 
derivatives in $\,\plane_+$ directions. Therefore, they can be 
solved explicitly in each fibre $\,\{y\}\times\plane_+$, and the solutions 
form a \hbox{nine\hs-}\hskip0ptdi\-men\-sion\-al af\-fine space, which is the 
same for all $\,y$.

The sections $\,F\,$ just mentioned may thus be viewed as functions on 
nonempty open subsets $\,\,U\,$ of $\,\bs$, valued in the af\-fine 
\hbox{$\,9$\hs-}\hskip0ptspace. For such functions $\,F\nnh$, the fourth 
curvature condition in Theorem~\ref{crvtc} is equivalent to the system 
(\ref{eqn}) of qua\-si-lin\-e\-ar first-or\-der differential equations with 
the unknown function $\,(\vf,\lambda,\mu):U\nh\to V\nnh$, related to 
$\,F\nnh$, and taking values in a specific 
\hbox{eight-}\hskip0ptdi\-men\-sion\-al vector space $\,V\nnh$. That solutions 
to (\ref{eqn}) exist is obvious (Example~\ref{lccne}). A description of all 
solutions is achieved by interpreting a function $\,\,U\to V\hs$ as a pair 
$\,(\vf,\bna\hh)$ consisting of a section $\,\vf\,$ of a (vector) plane bundle 
$\,\pb\,$ over the surface $\,\,U\,$ and a unimodular connection $\,\bna\,$ in 
$\,\pb\nh$. The system (\ref{eqn}) then amounts to the algebraic condition 
(\ref{brd}.i) on the curvature of $\,\bna\,$ coupled with the 
co\-var\-i\-ant-de\-riv\-a\-tive equation (\ref{brd}.ii) on $\,\vf$. As a 
result, (\ref{eqn}) is easily solved with the aid of gauge transformations and 
the method of characteristics.

Solutions to (\ref{eqn}) play a central role in Theorem~\ref{lcstr}, which 
presents a construction giving rise, locally and up to isometries, to all 
strictly \nw\ type III self-du\-al neutral Ein\-stein metrics in dimension 
four, and only to such metrics. The construction is explicit enough to yield 
easy answers to the most obvious questions: for instance, metrics with the 
properties just listed need not be Ric\-\hbox{ci\hh-}\hskip0ptflat or locally 
homogeneous (Theorem~\ref{rfnrf}). Theorem~\ref{lcstr} does not, however, 
provide a complete local classification of the these metrics, since it fails 
to describe a local moduli space, that is, to determine when two choices of 
the parameters $\,(\vf,\bna\hh)\,$ mentioned above lead to two 
\hbox{four\hh-}\hskip0ptman\-i\-folds which are locally isometric to each 
other.

\section{Preliminaries}\label{prel} 
\setcounter{equation}{0} 
All manifolds, \bd s, their sections and sub\bd s, as well as connections and 
mappings, including \bd\ morphisms, are assumed to be of class 
$\,C^\infty\nnh$. A manifold is by definition connected; a \bd\ morphism may 
operate only between two \bd s with the same \bmf, and acts by identity on the 
base.

We treat the co\-var\-i\-ant derivatives of a vector field $\,v\,$ and of a 
$\,1$-form $\,\xi$, relative to any fixed connection $\,\nabla\,$ on a 
manifold $\,M$, as a morphism $\,\nabla v:\tm\to\tm\,$ and, respectively, a 
twice-co\-var\-i\-ant \tf, acting on vector fields $\,u,w\,$ by 
$\,(\nabla v)u=\nabla_{\!u}v\,$ and $\,(\nabla\xi)(u,w)=(\nabla_{\!u}\xi)(w)$. 
For the tensor and exterior products of $\,1$-forms $\,\beta,\alpha\,$ on a 
manifold, the exterior derivative of $\,\beta$, and any tangent vector fields 
$\,u,v$, we have 
\begin{equation}\label{bwa} 
\begin{array}{rl} 
\mathrm{i)}\hskip7pt& 
(\beta\otimes\alpha)(\nh u,v)\,=\,\beta(u)\alpha(v)\hs,\hskip31pt 
\mathrm{ii)}\hskip7pt\beta\wedge\alpha=\beta\otimes\alpha 
-\alpha\otimes\beta\hs,\\ 
\mathrm{iii)}\hskip7pt& 
(d\beta)(\nh u,v)\,=\,d_u[\beta(v)]\,-\,d_v[\beta(u)]\,-\,\beta([u,v])\hs. 
\end{array} 
\end{equation} 
We use the metric $\,g\,$ of a \psr\ manifold $\,(M,g)\,$ to identify any 
vector field $\,u\,$ on $\,M\,$ with the $\,1$-form $\,g(u,\,\cdot\,)$. 
Similarly, we identify a vec\-tor-bun\-dle morphism $\,C:\tm\to\tm$ with the 
twice-co\-var\-i\-ant \tf\ $\,\bz\,$ such that $\,\bz(u,w)=g(Cu,w)\,$ for all 
vector fields $\,u,w$. In other words, $\,C\,$ is the result of raising the 
second index in $\,\bz$. A twice-co\-var\-i\-ant \tf\ $\,\bz\,$ thus 
associates with a vector field $\,u\,$ a new vector field $\,\bz u$, and 
with a $\,1$-form $\,\xi\,$ a new $\,1$-form $\,\bz\xi$, characterized by 
\begin{equation}\label{bue} 
\mathrm{i)}\hskip7pt\bz u\,=\,C\hh u\hs,\hskip39pt\mathrm{ii)}\hskip7pt 
\bz\hs\xi\,=\,\bz\hh(v,\,\cdot\,)\hs, 
\end{equation} 
where $\,C\,$ corresponds to $\,\bz\,$ as above, and $\,v\,$ is the vector 
field such that $\,\xi=g(v,\,\cdot\,)$.

By the Leib\-niz rule, when $\,\nabla\,$ is the \lcc\ of a \prc\ $\,g\,$ and 
$\,u,v,w\,$ are tangent vector fields, $\,2\hs g(\nabla_{\!v}u,w)\,$ equals 
\begin{equation}\label{lcc} 
d_v[\hs g(u,w)\hh]+d_u[\hs g(v,w)\hh]-d_w[\hs g(v,u)\hh]+g(u,[w,v]) 
+g(w,[v,u])-g(v,[u,w])\hh, 
\end{equation} 
cf.\ \cite[p.\ 160]{kobayashi-nomizu}, where $\,d_v$ is the directional 
derivative. Our sign convention about the curvature tensor $\,R\,$ of any 
connection $\,\nabla\,$ in a real vector bundle $\,\mathcal{E}$ over a 
manifold $\,M\,$ is 
\begin{equation}\label{cur} 
R(u,v)\hs w\hskip7pt=\hskip7pt\nabla_{\!v}\nabla_{\!u}w\, 
-\,\nabla_{\!u}\nabla_{\!v}w\,+\,\nabla_{[u,v]}w 
\end{equation} 
for vector fields $\,u,v\,$ tangent to $\,M\,$ and a section $\,w\,$ of 
$\,\mathcal{E}\nh$. Then 
\begin{equation}\label{ope} 
R(u,v):\mathcal{E}\nh\to\mathcal{E} 
\end{equation} 
denotes the bundle morphism sending any section $\,w\,$ to 
$\,R(u,v)\hs w$. We use the symbol $\,R\,$ also for the four-times covariant 
curvature tensor of a \prd\ $\,(M,g)$, with 
\begin{equation}\label{rww} 
R(w,w\hh'\nnh,u,v)\,=\,g(R(w,w\hh'\hh)\hs u,v)\hs. 
\end{equation} 
Aa a consequence of (\ref{cur}) and the Leib\-niz rule, any vec\-tor-bun\-dle 
morphism $\,C:\tm\to\tm$ in a pseu\-d\hbox{o\hs-}\hskip0ptRiem\-ann\-i\-an 
manifold $\,(M,g)\,$ satisfies the {\it Ric\-ci identity} 
\begin{equation}\label{rid} 
\nabla_{\!w}\nabla_{\!v}\hh C\,-\,\nabla_{\!v}\nabla_{\!w}\hh C\,+ 
\,\nabla_{[v,w]}\hh C\,=\,[R(v,w),C] 
\end{equation} 
for tangent vector fields $\,v,w$, with $\,[\hskip2.5pt,\hskip1pt]\,$ on the 
right-hand side denoting the commutator of bundle morphisms $\,\tm\to\tm$. (In 
coordinates, this reads 
$\,C^p_{\hh j,kl}-C^p_{\hh j,lk}=R_{lkj}{}^s\nh C^p_s-R_{lks}{}^p\nh C^s_j$.)

The four-times co\-var\-i\-ant curvature and Weyl tensors 
$\,R\,$ and $\,W$ of a \pr\ Ein\-stein \hbox{four\hh-}\hskip0ptman\-i\-fold $\,(M,g)\,$ with the 
the scalar curvature $\,12\hs K\,$ are related by 
\begin{equation}\label{tre} 
R\,=\,W\hs+\,Kg\wedge g\hs,\hskip8pt\mathrm{or,\ in\ coordinates,}\hskip6pt 
R_{jklp}=\hs W_{jklp}+\hs K(g_{jl}g_{kp}\nh-g_{kl}g_{jp})\hh. 
\end{equation} 
Our conventions about cur\-va\-ture\hs-like tensors $\,R\,$ acting on 
$\,2$-forms $\,\zeta\,$ and the inner product $\,\langle\,,\rangle\,$ of 
$\,2$-forms, are, in local coordinates, 
\begin{equation}\label{rom} 
\mathrm{a)}\hskip7pt2(R\hskip1.1pt\zeta)_{jk} 
=R_{jklp}\hs\zeta\hh^{lp}, 
\hskip23pt\mathrm{b)}\hskip7pt2\hs\lg\zeta\hh,\hs\eta\rg 
=-\hs\mathrm{tr}\hskip3pt\zeta\eta\hs, 
\end{equation} 
$\zeta\eta\,$ being the composite of the bun\-dle morphisms $\,\tm\to\tm\,$ 
corresponding to $\,\zeta\,$ and $\,\eta\,$ as in (\ref{bue}.i). The 
coordinate versions $\,2\lg\zeta\hh,\eta\rg 
=\hs\zeta_{jk}\hs\eta\hs^{jk}$ of (\ref{rom}.b) and 
$\,(\beta\wedge\alpha)_{jk}=\beta_j\hh\alpha_k-\beta_k\hh\alpha_j$ of 
(\ref{bwa}.ii) now give, for tangent vector fields $\,u\,$ and $\,v$, 
\begin{equation}\label{zuv} 
\lg\zeta\hh,\hs\eta\rg\,=\,\zeta(u,v)\hskip16pt\mathrm{if}\hskip13pt 
\eta\,=\,g(u,\,\cdot\,)\wedge g(v,\,\cdot\,)\hs. 
\end{equation} 
The Hodge star $\,*\,$ of an oriented \psr\ 
\hbox{four\hh-}\hskip0ptman\-i\-fold $\,(M,g)\,$ of the neutral signature 
$\,\mmpp\hs$, acting in the bundle $\,[\tam]^{\wedge2}$ of $\,2$-forms is an 
involution, and so 
$\,[\tam]^{\wedge2}\nh= 
\varLambda^+\hskip-1.8ptM\oplus\varLambda^-\hskip-1.8ptM$, for the 
$\,\pm1$-eigen\-space bundles $\,\varLambda^\pm\hskip-1.5ptM\,$ of $\,*\hh$, 
both of fibre dimension $\,3$, known as the bundles of {\it self-du\-al\/} and 
{\it anti-self-du\-al\/} $\,2$-forms \cite[p.\ 641, formula 
37.26]{dillen-verstraelen}. According to \cite[p.\ 643, formula 
(37.31)]{dillen-verstraelen}, with notation as in (\ref{rom}.b), 
\begin{equation}\label{ant} 
\zeta\hh\eta\,+\,\eta\hs\zeta\,\, 
=\,\,-\hs\lg\zeta\hh,\hs\eta\rg\hs\mathrm{Id}\hskip12pt 
\mathrm{whenever}\hskip7pt\zeta\hh,\hs\eta\hskip12pt\mathrm{are\ 
sections\ of}\hskip7pt\varLambda^\pm\hskip-1.5ptM\hs. 
\end{equation} 
\begin{remark}\label{liebr}A vector field $\,v\,$ on the total space of a 
bundle projection $\,\pi\,$ is called {\it vertical\/} if it is a section of 
the vertical distribution $\,\mathcal{V}=\kerd\pi$. As one easily verifies in 
suitable local coordinates, a vector field $\,w\,$ on the total space is 
$\,\pi$-pro\-ject\-a\-ble onto some vector field on the base manifold if and 
only if, for every vertical vector field $\,v$, the Lie bracket $\,[w,v]$ is 
also vertical. More generally, given an integrable distribution 
$\,\mathcal{V}\,$ on a manifold $\,M$, by a 
$\,\mathcal{V}\nh${\it-pro\-ject\-a\-ble local vector field\/} in $\,M\,$ we 
will mean any vector field $\,w\,$ defined on an open set $\,\,U\subset M\,$ 
and such that, whenever $\,v\,$ is a section of $\,\mathcal{V}\,$ defined on 
$\,\,U\nh$, so is $\,[w,v]$. 
\end{remark} 
\begin{remark}\label{orien}For any oriented \psr\ 
\hbox{four\hh-}\hskip0ptman\-i\-fold $\,(M,g)\,$ of the neutral signature 
$\,\mmpp\hs$, the bundles $\,\varLambda^\pm\hskip-1.5ptM\,$ 
of self-du\-al and anti-self-du\-al $\,2$-forms are canonically oriented. In 
fact, a pos\-i\-tive-oriented $\,(\mmpp)$-or\-tho\-nor\-mal basis of 
$\,\txm\,$ naturally gives rise to a basis of $\,\varLambda_x^\pm\nnh M$, 
which is also $\,(\mmp)$-or\-tho\-nor\-mal 
\cite[formula (37.25) on p.\ 641]{dillen-verstraelen}, in such a way that 
changing the signs of the first and third vectors in the basis of $\,\txm\,$ 
leaves the resulting orientations of $\,\varLambda_x^\pm\nnh M\,$ unchanged. 
Thus, both connected components of the set of pos\-i\-tive-oriented 
$\,(\mmpp)$-or\-tho\-nor\-mal bases of $\,\txm\,$ produce the same orientation 
in $\,\varLambda_x^\pm\nnh M$. 
\end{remark} 
\begin{remark}\label{trvol}For an $\,n$-form $\,\zeta\in[V^*]^{\wedge n}$ in a 
real vector space $\,V\hs$ of dimension $\,n$, any en\-do\-mor\-phism 
$\,\ed\,$ of $\,V\nh$, and any $\,w_1,\dots,w_n\in V\nh$, we have 
$\,\zeta(\ed w_1,w_2,\dots,w_n)+\zeta(w_1,\ed w_2,w_3,\dots,w_n)+\ldots 
+\zeta(w_1,w_2,\dots,w_{n-1},\ed w_n) 
=(\mathrm{tr}\,\ed)\hs\zeta(w_1,w_2,\dots,w_n)$, as one sees using the matrix 
of $\,\ed\,$ in the basis $\,w_1,\dots,w_n$, if $\,w_1,\dots,w_n$ are linearly 
independent, and noting that both sides vanish due to their skew-sym\-me\-try 
in $\,w_1,\dots,w_n$, if $\,w_1,\dots,w_n$ are linearly dependent. 
\end{remark} 
\begin{remark}\label{cclsm}If a tri\-lin\-e\-ar mapping 
$\,(w,w\hh'\nnh,w\hh''\hh)\mapsto\delta(w,w\hh'\nnh,w\hh'')\,$ from a 
\hbox{two\hh-}\hskip0ptdi\-men\-sion\-al vector space $\,\plane\,$ into any 
vector space is skew-sym\-met\-ric $\,w\hh'\nnh,w\hh''\nnh$, then 
$\,\delta(w,w\hh'\nnh,w\hh'')\,$ summed cyclically over $\,w,w\hh'\nnh,w\hh''$ 
yields $\,0$. In fact, the cyclic sum then depends on $\,w,w\hh'$ and 
$\,w\hh''$ skew-sym\-met\-ri\-cal\-ly, so that it vanishes as 
$\,\dim\plane=2$. 
\end{remark}

\section{Rank versus Pe\-trov type}\label{rvpt} 
\setcounter{equation}{0} 
We say that a traceless en\-do\-mor\-phism of a 
\hbox{three\hh-}\hskip0ptdi\-men\-sion\-al 
pseu\-\hbox{do\hs-}\hskip0ptEuclid\-e\-an space is of {\it Pe\-trov type\/} II 
(or, III) if it is self-ad\-joint and sends some ordered basis 
$\,(\mathbf{x},\mathbf{y},\mathbf{z})\,$ to $\,(0,0,\mathbf{x})$ (or, 
respectively, to $\,(0,\mathbf{x},\mathbf{y})$).

The symbol $\,W^\pm\nnh$ always denotes the self-du\-al and anti-self-du\-al 
parts of the Weyl tensor of a given oriented \psr\ \hbox{four\hh-}\hskip0ptman\-i\-fold $\,(M,g)\,$ of 
the neutral signature $\,\mmpp\hs$. Thus, $\,W^\pm\nnh$ may be viewed as an 
en\-do\-mor\-phism of the bundle $\,\varLambda^\pm\hskip-1.5ptM\,$ of 
(anti)self-du\-al $\,2$-forms. As usual, we call $\,(M,g)\,$ {\it 
self-du\-al\/} if $\,W^-\nnh=0$, and refer to $\,g\,$ as a {\it self-du\-al 
metric of Pe\-trov type\/} II (or, III) if $\,W^-\nnh=0\,$ and $\,W^+\nnh$ is 
of Pe\-trov type II (or III) at every point. 
\begin{remark}\label{satrz}Note that $\,\hs\mathrm{tr}\hskip3ptW^\pm\nnh=0$. 
See, for instance, \cite[formula (38.15) on p.\ 650]{dillen-verstraelen}. 
\end{remark} 
\begin{remark}\label{trtpt}A trace\-less self-ad\-joint en\-do\-mor\-phism 
$\,\nd\,$ of rank $\,1\,$ in a pseu\-\hbox{do\hs-}\hskip0ptEuclid\-e\-an 
\hbox{$\,3\hh$-}\hskip0ptspace $\,V\hs$ is necessarily of Pe\-trov type II, 
and its image $\,I\,$ is a null line.

In fact, if $\,I\,$ were not null, we would have $\,V\nh=I\oplus I^\perp$ and 
$\,I^\perp\nh=\hs\mathrm{Ker}\,\nd$, and hence the restriction 
$\,\nd:I\to I\,$ would be trace\-less, contrary to its surjectivity. Now we 
can use $\,\mathbf{x},\mathbf{y},\mathbf{z}$ with 
$\,\mathbf{x}\in I\subset I^\perp\nh=\hs\mathrm{Ker}\,\nd\,$ and 
$\,\mathbf{y}\in\hs\mathrm{Ker}\,\nd$, normalizing $\,\mathbf{z}\,$ so that 
$\,\nd\mathbf{z}=\mathbf{x}$. 
\end{remark} 
\begin{lemma}\label{rkvpt}Let\/ $\,(M,g)\,$ be a \psr\ 
\hbox{four\hh-}\hskip0ptman\-i\-fold of the neutral signature $\,\mmpp\,$ such 
that the Weyl tensor $\,W\nh$ acting in the bundle $\,[\tam]^{\wedge2}$ of 
$\,2$-forms has constant rank\/ $\,\er$, and the restriction of the fibre metric 
$\,\lr\,$ in $\,[\tam]^{\wedge2}$ to the sub\-bun\-dle\/ $\,\mathcal{E}\,$ 
forming the image of\/ $\,W\nnh$ is degenerate at every point. 
\begin{enumerate} 
  \def\theenumi{{\rm\roman{enumi}}} 
\item[{\rm(i)}] If\/ $\,\er=1$, then\/ $\,M\,$ admits an orientation for which\/ 
$\,g\,$ is self-du\-al of Pe\-trov type\/ {\rm II}. 
\item[{\rm(ii)}] If\/ $\,M\,$ is oriented, $\,\er=2$, and\/ $\,\mathcal{E}$ is 
$\,\lr$-null at some point, then, at every point, $\,\mathcal{E}$ is 
$\,\lr$-null and\/ $\,W^\pm\nnh$ are both of Pe\-trov type\/ {\rm II}. 
\item[{\rm(iii)}] If\/ $\,\er=2\,$ and, at some point, $\,\mathcal{E}$ is 
non-null relative to $\,\lr$, then\/ $\,\mathcal{E}$ is $\,\lr$-non-null at 
every point, while\/ $\,M\hs$ admits an orientation for which\/ $\,g\,$ is 
self-du\-al of Pe\-trov type\/ {\rm III}. 
\end{enumerate} 
\end{lemma} 
\begin{proof}Let $\,\er^\pm$ be the function assigning to each $\,x\in M\,$ 
the rank of $\,W^\pm\nh$ at $\,x$. Since $\,\er^\pm$ is lower 
sem\-i\-con\-tin\-u\-ous, constancy of $\,\er=\er^+\hskip-1.8pt+\hs \er^-$ 
implies that both $\,\er^\pm$ are constant. Thus, locally, 
$\,\mathcal{E}=\hskip1.8pt\mathcal{E}^+\hskip-2.4pt\oplus\mathcal{E}^-$ for 
some sub\-bun\-dles $\,\mathcal{E}^\pm$ of $\,\varLambda^\pm\hskip-1.5ptM$.

If $\,\er=1$, the requirement that $\,\er^+\nnh=1\,$ and $\,\er^-\nnh=0\,$ 
uniquely defines an orientation of $\,M$, and (i) is obvious from 
Remarks~\ref{satrz} and ~\ref{trtpt}.

If $\,M\,$ is oriented, $\,\er=2$, and $\,\er^+\hskip-1.8pt=\er^-\nnh=1$, 
Remarks~\ref{satrz} and ~\ref{trtpt} yield two conclusions. First, both 
$\,W^\pm\nh$ are of Pe\-trov type II at each point. Secondly, both 
$\,\mathcal{E}^\pm$ are $\,\lr$-null at every point, and hence so is 
$\,\mathcal{E}=\hskip1pt\mathcal{E}^+\hskip-2.5pt\oplus\mathcal{E}^-\nnh$.

Finally, let $\,\er=2\,$ and $\,\er^+\hskip-1.8pt\ne\hs \er^-\nnh$. For the 
orientation defined 
by requiring that $\,\er^+\nnh=2\,$ and $\,\er^-\nnh=0$, we have $\,W^-\nnh=0$, 
while $\,\mathcal{E}\nh=\mathcal{E}^+$ cannot be $\,\lr$-null at any point 
$\,x$. (If it were, $\,\mathcal{E}_x$ would be a $\,2$\diml\ subspace of 
$\,\varLambda_x^+\nnh M\,$ contained in its own $\,1$\diml\ orthogonal 
complement.) Since the restriction of $\,\lr\,$ to the plane bundle 
$\,\mathcal{E}\,$ (the image of $\,W^+$) is assumed degenerate, Pe\-trov's 
classification \cite[Proposition 39.2 on p.\ 652]{dillen-verstraelen} implies 
that, at each point, $\,W^+\nh$ is of Pe\-trov type III. Combined with the 
preceding paragraph, this proves (ii) and (iii). 
\end{proof}

\section{The vertical-horizontal decomposition}\label{tvhd} 
\setcounter{equation}{0} 
In the following two lemmas, the meaning of $\,\varLambda^+\hskip-1.8ptM\,$ 
and $\,W^+\nnh$ is the same as in the second paragraph of Section~\ref{rvpt}. An 
en\-do\-mor\-phism $\,\rw\,$ of $\,\varLambda^+\hskip-1.8ptM$, such as 
$\,W^+\nnh$, is identified with a four-times co\-var\-i\-ant tensor field 
(namely, a morphism $\,[\tam]^{\wedge2}\nh\to[\tam]^{\wedge2}$ vanishing on 
$\,\varLambda^-\hskip-1.8ptM$), and $\,\hs\mathrm{div}\,\rw\,$ is defined, in 
local coordinates, by $\,(\mathrm{div}\,\rw)_{klm}=\rw^j{}_{klm,j}$. 
\begin{lemma}\label{cptqt}Given an oriented 
\hbox{four\hh-}\hskip0ptdi\-men\-sion\-al \prd\/ $\,(M,g)\,$ of 
the neutral metric signature $\,\mmpp\hs$, let\/ $\,\rw\,$ be a bundle 
en\-do\-mor\-phism of\/ $\,\varLambda^+\hskip-1.8ptM\,$ with\/ 
$\,\mathrm{tr}\,\rw=0$ and\/ $\,\mathrm{div}\,\rw=0$, which is of Pe\-trov 
type\/ {\rm III} at every point, as defined in Section\/~{\rm\ref{rvpt}}.

Then\/ $\,\varLambda^+\hskip-1.8ptM\,$ has a $\,C^\infty\nnh$ global 
trivialization\/ $\,(\zeta,\eta,\theta)\,$ with 
\begin{equation}\label{wze} 
\begin{array}{rl} 
\mathrm{i)}\hskip7pt&\rw\zeta=0,\,\,\rw\eta=-\hs\zeta,\,\,\rw\theta=\eta\hs,\\ 
\mathrm{ii)}\hskip7pt&\lg\zeta,\theta\rg=2=-\hs\lg\eta,\eta\rg\hs,\hskip12pt 
\lg\zeta,\zeta\rg=\lg\zeta,\eta\rg=\lg\eta,\theta\rg=\lg\theta,\theta\rg=0\hs, 
\end{array} 
\end{equation} 
unique up to replacement by $\,(-\hs\zeta,-\hs\eta,-\hs\theta)$, and\/ 
$\,\mathrm{rank}\hskip3pt\zeta=\hs\mathrm{rank}\hskip3pt\theta=2\,$ at each 
point, while 
\begin{enumerate} 
  \def\theenumi{{\rm\alph{enumi}}} 
\item[{\rm(a)}] $\mathcal{V}=\mathrm{Ker}\,\zeta\,$ and\/ 
$\,\mathcal{H}=\mathrm{Ker}\,\theta\,$ are null\/ $\,2$\diml\ distributions 
on\/ $\,(M,g)$. 
\end{enumerate} 
The trivialization\/ $\,(\zeta,\eta,\theta)\,$ becomes unique if one requires, 
in addition, that 
\begin{equation}\label{eig} 
\eta\hs v=v\,\ \ \mathrm{and\ \ }\,\eta\hs w=-\hs w\,\mathrm{\ \ for\ 
sections\ }\,\,v\,\,\mathrm{\ of\ }\,\,\mathcal{V}\,\ \mathrm{and\ } 
\,\,w\,\,\mathrm{\ of\ }\,\,\mathcal{H}, 
\end{equation} 
with $\,\eta\hs v\,$ as in\/ {\rm(\ref{bue}.i)}. Whether\/ {\rm(\ref{eig})} is 
assumed or not, we have the following conclusions. 
\begin{enumerate} 
  \def\theenumi{{\rm\alph{enumi}}} 
\item[{\rm(b)}] $\mathcal{V}\,$ is integrable, its leaves are totally 
geodesic, and\/ $\,\tm\,=\,\mathcal{H}\oplus\mathcal{V}\nh$. 
\item[{\rm(c)}] $\nabla\zeta=2\hh\alpha\otimes\zeta+2\beta\otimes\eta$, 
$\,\,\nabla\eta=2\gamma\otimes\zeta+2\beta\otimes\theta\,$ and\/ 
$\,\nabla\theta=2\gamma\otimes\eta-2\hh\alpha\otimes\theta\,$ for some unique 
$\,1$-forms $\,\alpha,\beta,\gamma\,$ on $\,M$. 
\item[{\rm(d)}] $\beta(v)=0\,$ for every section $\,v\,$ of\/ $\,\mathcal{V}$. 
\item[{\rm(e)}] $2\hs \rw\hs=\,\zeta\otimes\eta\,+\hs\eta\otimes\zeta\,$ for 
$\,\rw\hs$ treated as a four-times co\-var\-i\-ant tensor field. 
\item[{\rm(f)}] $2\hs\zeta\gamma+\eta\alpha+\theta\beta=2\eta\beta+\zeta\alpha 
=\zeta\beta=0$, in the notation of\/ {\rm(\ref{bue}.ii)}. 
\end{enumerate} 
\end{lemma} 
\begin{proof}The assumption made about $\,\rw$, combined with Pe\-trov's 
classification \cite[Proposition 39.2 on p.\ 652]{dillen-verstraelen}, gives 
(\ref{wze}) at each point $\,x$, for some basis $\,\zeta_x,\eta_x,\theta_x$ of 
$\,\varLambda_x^\pm\nnh M\,$ which is unique up to an overall sign change 
\cite[p.\ 656, Remark 39.3(iv-c)]{dillen-verstraelen}. Our trivialization 
$\,(\zeta,\eta,\theta)\,$ is therefore unique, locally, up to a change of 
sign, and so it has a global sin\-gle\hs-val\-ued $\,C^\infty$ branch (as the 
bundle $\,\varLambda^+\hskip-1.8ptM\,$ is orientable, cf.\ 
Remark~\ref{orien}). According to 
\cite[p.\ 645, Lemma 37.8]{dillen-verstraelen}, at each point $\,x$, the 
$\,2$-forms $\,\zeta_x$ and $\,\theta_x$ (written in the next five 
lines without the subscript $\,x$), being nonzero, self-du\-al and 
null, can be decomposed as $\,\zeta=\xi_1\wedge\,\xi_2$ and 
$\,\theta=\xi_3\wedge\,\xi_4$, with $\,\xi_j=g_x(e_j,\,\cdot\,)$, 
$\,j=1,2,3,4$, for some $\,e_j\in\txm\,$ such that 
$\,\mathcal{V}_{\nh x}=\,\mathrm{span}\hs\{e_1,e_2\}\,$ and 
$\,\mathcal{H}_x=\,\mathrm{span}\hs\{e_3,e_4\}\,$ are null 
planes. This proves (a). As 
$\,\xi_1\wedge\,\ldots\,\wedge\,\xi_4=\zeta\wedge\theta\,$ equals 
$\,\lg\zeta,*\hh\theta\rg$ times the volume form, while 
$\,\lg\zeta,*\hh\theta\rg=\lg\zeta,\theta\rg=2\ne0$, the vectors 
$\,e_1,e_2,e_3,e_4$ must form a basis of $\,\txm$. Hence 
$\,\tm\hs=\hh\mathcal{H}\oplus\mathcal{V}$.

By (\ref{ant}) and (\ref{wze}.ii), $\,\eta\eta=\mathrm{Id}\,$ and $\,\eta\,$ 
an\-ti\-com\-mutes with $\,\zeta\,$ and $\,\theta$.

We will now verify that (\ref{eig}) holds, possibly after $\,\eta\,$ has been 
replaced by $\,-\hs\eta\,$ (which will clearly imply the uniqueness claim 
concerning (\ref{eig})). To this end, let us first note that, if $\,\eta_x$, 
at any point $\,x\in M$, has eigenvectors $\,v\in\mathcal{V}_x$ and 
$\,u\in\mathcal{H}_x$, then they cannot correspond to the same eigenvalue. In 
fact, if they did, it would follow that $\,g_x(v,u)=0$, since each 
eigen\-space of $\,\eta_x$ is null due to skew-sym\-me\-try of $\,\eta\,$ and 
its nondegeneracy. On the other hand, the relations 
$\,\mathcal{H}=\mathrm{Ker}\,\theta\,$ and 
$\,\tm=\mathcal{H}\oplus\mathcal{V}\,$ imply injectivity of $\,\theta_x$ 
restricted to $\,\mathcal{V}_x$, and so, since $\,\eta\,$ and $\,\theta$ 
an\-ti\-com\-mute, $\,\theta_xv\in\mathcal{H}_x^\perp\nnh=\mathcal{H}_x$ would 
be an eigenvector of $\,\eta_x$ for the opposite of the original 
$\,v$-eigen\-value, and, consequently, $\,u\,$ and $\,\theta_xv$, being 
linearly independent, would span $\,\mathcal{H}_x$, while 
$\,g_x(v,\theta_xv)=0\,$ as $\,\theta\,$ is skew-sym\-met\-ric. Thus, $\,v\,$ 
would be orthogonal to both $\,u\,$ and $\,\theta_xv$, so that it would lie in 
$\,\mathcal{H}_x^\perp\nnh=\mathcal{H}_x$, contradicting the relation 
$\,\mathcal{H}_x\cap\hh\mathcal{V}_x=\{0\}$.

Since $\,\eta\,$ an\-ti\-com\-mutes with $\,\zeta\,$ and $\,\theta$, it leaves 
the distributions $\,\mathcal{V}=\mathrm{Ker}\,\zeta\,$ and 
$\,\mathcal{H}=\mathrm{Ker}\,\theta$ invariant, so that, as 
$\,\eta\eta\,$ equals the identity, $\,\txm\,$ is, at each point $\,x\in M$, 
spanned by eigenvectors of $\,\eta_x$ for the eigenvalues $\,\pm1$, lying in 
$\,\mathcal{V}_x$ and $\,\mathcal{H}_x$. Combined with the last paragraph, 
this proves (\ref{eig}) up to a change of sign.

The existence and uniqueness of $\,1$-forms $\,\alpha,\beta,\gamma\,$ with 
(c) is immediate from (\ref{wze}.ii) and invariance of 
$\,\varLambda^+\hskip-1.8ptM\,$ under parallel transports. Also, (e) 
follows since, by (\ref{wze}.ii) and (\ref{rom}), both sides act on 
$\,\varLambda^+\hskip-1.8ptM\,$ as described in (\ref{wze}.i). In view of (e) 
and (c), 
$\,\nabla \rw=2\gamma\otimes\zeta\otimes\zeta+2\beta\otimes\eta\otimes\eta 
+\alpha\otimes(\zeta\otimes\eta+\eta\otimes\zeta) 
+\beta\otimes(\zeta\otimes\theta+\theta\otimes\zeta)$. Contraction now gives 
$\,\hs\mathrm{div}\,\rw 
=(2\hs\zeta\gamma+\eta\alpha+\theta\beta)\otimes\zeta 
+(2\eta\beta+\zeta\alpha)\otimes\eta+\zeta\beta\otimes\theta$. As 
$\,\mathrm{div}\,\rw=0$, this implies (f). The vector field associated by 
$\,g\,$ with $\,\beta\,$ thus is a section of 
$\,\mathcal{V}=\mathrm{Ker}\,\zeta\,$ and, as $\,\mathcal{V}\,$ is null, we 
get (d).

Finally, let $\,v,v\hh'$ be any sections of $\,\mathcal{V}$. The Leib\-niz 
rule and (c) with $\,\zeta v=\zeta v\hh'\nh=0\,$ show that 
$\,\zeta\nabla_{\hskip-2ptv}v\hh'\nh 
=-\hs(\nabla_{\hskip-2ptv}\zeta)v\hh'\nh 
=-2\beta(v)\hs\eta v\hh'\nnh$. However, $\,\beta(v)=0\,$ in view of (e). 
Consequently, $\,\nabla_{\hskip-2ptv}v\hh'$ is a section of 
$\,\mathcal{V}=\mathrm{Ker}\,\zeta$, which yields (b). 
\end{proof} 
\begin{lemma}\label{sdept}Suppose that an oriented \psr\ Ein\-stein 
\hbox{four\hh-}\hskip0ptman\-i\-fold\/ $\,(M,g)\,$ with the neutral metric 
signature $\,\mmpp\,$ is self-du\-al of Pe\-trov type\/ {\rm III}. The 
assumptions of Lemma\/~{\rm\ref{cptqt}} then are satisfied by $\,\rw=W^+\nnh$. 
For $\,\zeta,\eta,\alpha,\beta,\gamma,\theta,\mathcal{V}\,$ and 
$\,\mathcal{H}\,$ uniquely defined as in Lemma\/~{\rm\ref{cptqt}} with 
$\,\rw=W^+\nnh$, and with $\,\nabla\nnh,R\,$ and\/ $\,K\,$ denoting the \lcc, 
four-times co\-var\-i\-ant curvature tensor, and\/ $\,1/12\,$ of the scalar 
curvature of\/ $\,g$, 
\begin{enumerate} 
  \def\theenumi{{\rm\roman{enumi}}} 
\item[{\rm(i)}] the connection induced by\/ $\,\nabla\hh$ on each leaf of\/ 
$\,\mathcal{V}\,$ is flat, 
\item[{\rm(ii)}] $2R(v,w)=\hs\zeta(v,w)\hs\eta 
+\eta(v,w)\hs\zeta+2\hh K\hs\xi\wedge\hs\xi\hs'\nnh$, where $\,v,w\,$ are any 
vector fields, $\,\xi=g(v,\,\cdot\,)$ and\/ 
$\,\xi\hs'\nh=g(w,\,\cdot\,)$, with\/ $\,\xi\wedge\hs\xi\hs'$ as in\/ 
{\rm(\ref{bwa}.ii)}, 
\item[{\rm(iii)}] $d_u[\hs g(v,w)]\,=\,g(v,\nabla_{\hskip-2ptu}w)\, 
=\,g(\nabla_{\!w}u,v)\,=\,\beta(w)\hs\theta(u,v)\,$ whenever $\,u,v\,$ are 
sections of\/ $\,\mathcal{V}$ parallel in the direction of\/ 
$\,\mathcal{V}\nh$, while $\,w\,$ is a $\,\mathcal{V}$-pro\-ject\-a\-ble local 
section of\/ $\,\mathcal{H}$, 
\item[{\rm(iv)}] $2\hskip1.4ptd\beta+4\beta\wedge\alpha=-\hs K\zeta\hs, 
\hskip16pt 
2\hskip1.4ptd\gamma+4\hh\alpha\wedge\gamma=K\theta+\eta$, 
\item[{\rm(v)}] $\mathcal{V}=\mathcal{V}^\perp\nnh=\mathrm{Ker}\,\zeta 
=\mathrm{Im}\,\zeta=\mathrm{Ker}\,(\eta\hh-\nh\mathrm{Id})\,$\hskip3.5ptand 
\hskip3pt$\,\mathcal{H}=\mathcal{H}^\perp\nnh=\mathrm{Ker}\,\theta 
=\mathrm{Im}\,\theta=\mathrm{Ker}\,(\eta\hh+\nh\mathrm{Id})$, with\/ 
$\,\mathrm{Im}\,$ meaning `\hskip.8ptimage\nnh' and\/ $\,\zeta,\eta,\theta\,$ 
treated as morphisms $\,\tm\to\tm$, cf.\ {\rm(\ref{bue}.i)}, 
\item[{\rm(vi)}] $[\nabla_{\!w}u]^{\mathcal{H}}\nh=\beta(w)\hs\theta\hh u\,$ 
for any vector field\/ $\,w\,$ and any section $\,u\,$ of\/ 
$\,\mathcal{V}\nh$, where $\,[\hskip3pt]^{\mathcal{H}}$ denotes the 
$\,\mathcal{H}$-com\-po\-nent projection in\/ 
$\,\tm\,=\,\mathcal{H}\oplus\mathcal{V}\nh$. 
\end{enumerate} 
\end{lemma} 
\begin{proof}The assumptions of Lemma~\ref{cptqt} are satisfied by 
$\,\rw=W^+\nnh$. Namely, it is the well known, 
cf.\ \cite[p.\ 460, Lemma 5.2]{dillen-verstraelen}, that 
$\,\mathrm{div}\hskip3ptW\nnh=0\,$ for any 
pseu\-d\hbox{o\hs-}\hskip0ptRiem\-ann\-i\-an Ein\-stein metric. Now 
(\ref{tre}) and Lemma~\ref{cptqt}(e) with $\,\rw=W^+\nnh=W\hs$ yield (ii). 
Also, by (ii), $\,R(v,v\hh'\hs)\hs w=0\,$ for any sections 
$\,v,v\hh'\nnh,w\hs$ of the null distribution 
$\,\mathcal{V}=\mathrm{Ker}\,\zeta$, which proves (i). 
Next, (v) is obvious from Lemma~\ref{cptqt}(a) and (\ref{eig}), since 
$\,\zeta,\theta:\tm\to\tm\,$ are 
skew-ad\-joint at each point, while the distributions $\,\mathcal{V}\,$ and 
$\,\mathcal{H}\,$ are null and $\,2$\diml. On the other hand, (iii) follows 
since, by (\ref{eig}), Lemma~\ref{cptqt}(c) and (v), $\,g(\nabla_{\!w}u,v) 
=g(\nabla_{\!w}(\eta u),v)=g((\nabla_{\!w}\eta)u,v) 
+g(\eta\hs[\nabla_{\!w}u]^{\mathcal{H}}\nnh,v) 
=2g(\beta(w)\hs\theta\hh u,v)-g(\nabla_{\!w}u,v)$, with 
$\,[\hskip3pt]^{\mathcal{H}}$ as in (vi), while 
$\,d_u[\hs g(w,v)]=g(\nabla_{\hskip-2ptu}w,v)\,$ from the Leib\-niz rule, and 
$\,g(\nabla_{\hskip-2ptu}w,v)=g(\nabla_{\!w}u,v)\,$ in view of 
Remark~\ref{liebr}, as $\,[w,v]=\nabla_{\!w}v-\nabla_{\hskip-2ptv}w\,$ and 
$\,\mathcal{V}\,$ is null.

For any fixed vector fields $\,v,w$, let $\,[R(v,w),\eta\hh]\,$ denote the 
$\,2$-form given by 
$\,[R(v,w),\eta\hh]=\nabla_{\!w}\nabla_{\!v}\hh\eta 
-\nabla_{\!v}\nabla_{\!w}\hh\eta+\nabla_{[v,w]}\hh\eta$. Evaluating 
$\,[R(v,w),\eta\hh]\,$ with the aid of Lemma~\ref{cptqt}(c), and using 
(\ref{wze}.ii) along with (\ref{bwa}), we see that 
\begin{equation}\label{thr} 
\lg\theta,[R(v,w),\eta\hh]\rg 
=2\hh[\hh2\hskip1.4ptd\gamma+4\hh\alpha\wedge\gamma](w,v)\hs,\hskip10pt 
\lg\zeta,[R(v,w),\eta\hh]\rg 
=2\hh[\hh2\hskip1.4ptd\beta+4\beta\wedge\alpha](w,v)\hs. 
\end{equation} 
Identifying $\,\zeta,\eta,\theta\,$ with bun\-dle morphisms $\,\tm\to\tm\,$ 
as in (\ref{bue}.i), and using the multiplicative notation for their 
composites, we have 
\begin{equation}\label{tee} 
\theta\eta=-\hs\eta\hh\theta=\theta\hs,\hskip22pt 
\zeta\eta=-\hs\eta\hs\zeta=-\hs\zeta\hs. 
\end{equation} 
In fact, by (v), both sides in each equality agree separately on 
$\,\mathcal{V}=\mathrm{Ker}\,\zeta\,$ and on 
$\,\mathcal{H}=\mathrm{Ker}\,\theta$.

Defining commutators of $\,2$-forms, as usual, in terms of their composites, 
we see that $\,[R(v,w),\eta\hh]\,$ introduced above becomes such a commutator 
if we identify $\,R(v,w)\,$ with the $\,2$-form 
$\,R(v,w,\,\cdot\,,\,\cdot\,)$, cf.\ (\ref{bue}.i). In fact, this is immediate 
from the Ric\-ci identity (\ref{rid}) applied to $\,C\,$ which corresponds to 
$\,\bz=\eta\,$ as in (\ref{bue}.i). By (ii), 
$\,[R(v,w),\eta\hh]=-\hs\eta(v,w)\hs\zeta 
+K[\hs\xi\wedge\hs\xi\hs'\nnh,\eta\hh]$, 
with $\,\xi,\hh\xi\hs'$ as in (ii), since $\,[\eta,\eta\hh]=0$, while 
(\ref{tee}) gives $\,[\hh\zeta,\eta\hh]=-2\hs\zeta$. Consequently, 
\begin{equation}\label{ltr} 
\lg\theta,[R(v,w),\eta\hh]\rg=2\hh\eta(w,v)+2K\theta(w,v)\hs,\hskip22pt 
\lg\zeta,[R(v,w),\eta\hh]\rg=-2K\zeta(w,v)\hs. 
\end{equation} 
To verify (\ref{ltr}), note that $\,\lg\zeta,\theta\rg=2\,$ and 
$\,\lg\zeta,\zeta\rg=0\,$ by (\ref{wze}.ii), while 
$\,-2\lg\theta,[\hs\xi\wedge\hs\xi\hs'\nnh,\eta\hh]\rg 
=\mathrm{tr}\hskip3pt\theta[\hs\xi\wedge\hs\xi\hs'\nnh,\eta\hh] 
=\mathrm{tr}\hskip3pt(\xi\wedge\hs\xi\hs'\hh)\eta\theta 
-\mathrm{tr}\hskip3pt\theta\eta(\xi\wedge\hs\xi\hs'\hh) 
=-2\hskip1.2pt\mathrm{tr}\hskip3pt(\xi\wedge\hs\xi\hs'\hh)\theta 
=4\lg\theta,\xi\wedge\hs\xi\hs'\hh\rg=4\hh\theta(v,w)\,$ by (\ref{rom}.b), 
(\ref{zuv}) and (\ref{tee}), so that 
$\,\lg\theta,[\hs\xi\wedge\hs\xi\hs'\nnh,\eta\hh]\rg 
=-2\hh\theta(v,w)\,$ and, similarly, 
$\,\lg\zeta,[\hs\xi\wedge\hs\xi\hs'\nnh,\eta\hh]\rg 
=2\hs\zeta(v,w)$.

Combining (\ref{thr}) with (\ref{ltr}), we obtain (iv).

Finally, both sides in (vi) are, by (v), sections of the null distribution 
$\,\mathcal{H}$. Their inner products with any section of $\,\mathcal{H}\,$ 
(or, respectively, any section $\,v\,$ of $\,\mathcal{V}$) thus are both zero, 
(or, respectively, are equal in view of the last equality in (iii)). This 
yields (vi). 
\end{proof} 
Given $\,(M,g)\,$ as in Lemma~\ref{sdept}, with the corresponding objects 
$\,\zeta,\eta,\theta,\mathcal{V}\,$ and $\,\mathcal{H}$, we may choose, 
locally, sections $\,w,w\hh'$ of $\,\mathcal{H}\,$ such that 
$\,\zeta(w,w\hh'\hh)=1$. Setting $\,v=-\hs\zeta w\hh'$ and 
$\,v\hh'\nnh=\zeta w$, we obtain a local trivialization $\,w,w\hh'\nnh,v,v\hh'$ 
of the tangent bundle $\,\tm$, in which the only nonzero components of 
$\,g,\zeta,\eta\,$ and $\,\theta\,$ are $\,g(v,w)=g(v\hh'\nnh,w\hh'\hh)=1$, 
$\,\zeta(w,w\hh'\hh)=1$, $\,\eta(v,w)=\eta(v\hh'\nnh,w\hh'\hh)=1$, 
$\,\theta(v,v\hh'\hh)=2$, and those arising from them due to symmetry of $\,g\,$ 
and skew-sym\-me\-try of $\,\zeta,\eta,\theta$.

In fact, such $\,w\,$ and $\,w\hh'$ exist since 
$\,\tm\,=\,\mathcal{H}\oplus\mathcal{V}\hs$ and 
$\,\mathcal{V}=\mathrm{Ker}\,\zeta$, which also implies 
injectivity of $\,\zeta_x:\mathcal{H}_x\to\mathcal{V}_x$, at every 
point $\,x$. Hence $\,w,w\hh'\nnh,v,v\hh'$ form a local trivialization 
of $\,\tm$. As $\,\mathcal{H}=\mathrm{Ker}\,\theta$, applying to 
$\,w\,$ and $\,w\hh'$ the relation 
$\,(\zeta\theta+\theta\zeta)/2=-\hs\mathrm{Id}$, 
immediate from (\ref{wze}.ii) for $\,\rw=W^+$ and (\ref{ant}), we get 
$\,\theta\hh v\hh'\nnh=-2w\,$ and $\,\theta\hh v=2w\hh'\nnh$. As 
$\,\mathcal{V}\,$ and $\,\mathcal{H}\,$ are $\,g$-null, our claim is now 
obvious from (\ref{eig}). 
\begin{remark}\label{cvhmg}All self-du\-al oriented Ein\-stein 
\hbox{four\hh-}\hskip0ptman\-i\-folds of the neutral metric signature 
$\,\mmpp\hs$, which are of Pe\-trov type III, are curvature homogeneous. This 
well-known fact is an obvious consequence of the last paragraph: by 
Lemma~\ref{sdept}(ii), $\,w,w\hh'\nnh,v,v\hh'$ chosen as above form, at any 
point, a basis of the tangent space providing a canonical expression for both 
the metric and the curvature tensor. 
\end{remark}

\section{The Walker and strictly \nw\ cases}\label{wnwc} 
\setcounter{equation}{0} 
Suppose that $\,\mathcal{T}\hs$ is a pseu\-\hbox{do\hs-}\hskip0ptEuclid\-e\-an 
\hbox{$\,4\hh$-}\hskip0ptspace of the neutral signature $\,(\mmpp)$. By a {\it 
null plane\/} in $\,\mathcal{T}\hs$ we mean a null 
\hbox{two\hh-}\hskip0ptdi\-men\-sion\-al vector subspace of 
$\,\mathcal{T}\nnh$. We also use the inner product of $\,\mathcal{T}\nnh$, 
denoted by $\,\lr$, to identify the space $\,[\mathcal{T}^*]^{\wedge2}$ of 
$\,2$-forms with the space $\,\mathcal{T}^{\wedge2}$ of bi\-vec\-tors. Thus, 
if $\,\mathcal{T}\hs$ is oriented, we can treat the Hodge star $\,*\,$ as an 
involution of $\,\mathcal{T}^{\wedge2}\nnh$, and speak of self-du\-al or 
an\-ti-self-du\-al bi\-vec\-tors in $\,\mathcal{T}\nnh$. 
\begin{lemma}\label{orntt}Let\/ $\,\mathcal{T}\hs$ be a 
pseu\-\hbox{do\hs-}\hskip0ptEuclid\-e\-an \hbox{$\,4\hh$-}\hskip0ptspace of 
the neutral metric signature. 
\begin{enumerate} 
  \def\theenumi{{\rm\alph{enumi}}} 
\item[{\rm(a)}] Any null plane\/ $\,\mathcal{N}\,$ in\/ $\,\mathcal{T}\hs$ 
naturally distinguishes an orientation of\/ $\,\mathcal{T}\nnh$, namely, the one 
which, for some\hh/any basis\/ $\,u,v\,$ of\/ $\,\mathcal{N}$, makes the 
bi\-vec\-tor\/ $\,u\wedge v\,$ self-du\-al. 
\item[{\rm(b)}] If\/ $\,\mathcal{N}\hs$ and\/ $\,\mathcal{N}\hh'\nh$ are null 
planes in\/ $\,\mathcal{T}\hs$ and\/ 
$\,\dim\hskip2pt(\hh\mathcal{N}\nnh\cap\hskip.2pt\mathcal{N}\hh')=1$, 
then\/ $\,\mathcal{N}\hh$ and\/ $\,\mathcal{N}\hh'\nh$ distinguish, in the 
sense of\/ {\rm(a)}, two opposite orientations of\/ $\,\mathcal{T}\nnh$. 
\end{enumerate} 
\end{lemma} 
\begin{proof}A bi\-vec\-tor in $\,\mathcal{T}\hs$ equals $\,u\wedge v\,$ for 
some basis $\,u,v\,$ of some null plane if and only if it is nonzero, null, 
and self-du\-al or an\-ti-self-du\-al 
\cite[p.\ 645, Lemma~37.8]{dillen-verstraelen}. This yields (a), and at the 
same time shows that, under the assumptions of (b), if 
$\,u,v,w\in\mathcal{T}\hs$ and $\,u,v\,$ span $\,\mathcal{N}\nh$, while 
$\,v,w\,$ span $\,\mathcal{N}\hh'\nnh$, then $\,u\wedge v\,$ and 
$\,v\wedge w\,$ are linearly independent null bi\-vec\-tors, and each of them 
is self-du\-al or an\-ti-self-du\-al. If both $\,\mathcal{N}\,$ and 
$\,\mathcal{N}\hh'$ distinguished the same orientation of $\,\mathcal{T}\nnh$, 
this orientation would make both bi\-vec\-tors self-du\-al, and so they would 
have to be orthogonal, since, for any self-du\-al bi\-vec\-tors 
$\,\zeta,\eta$, the $\,4$-form $\,\zeta\wedge\eta\,$ is the product of 
$\,\lg\zeta,*\hh\eta\rg=\lg\zeta,\eta\rg\,$ and the volume form. Thus, 
$\,u\wedge v\,$ and $\,v\wedge w\,$ would span a null plane in the 
pseu\-\hbox{do\hs-}\hskip0ptEuclid\-e\-an \hbox{$\,3\hh$-}\hskip0ptspace of 
self-du\-al bi\-vec\-tors, which cannot exist for dimensional reasons, as it 
would be transverse to a space\-like or time\-like plane. 
\end{proof} 
For a null plane $\,\mathcal{N}\,$ as in Lemma~\ref{orntt}(a), we will say 
that $\,\mathcal{N}\,$ is {\it compatible\/} with the orientation of 
$\,\mathcal{T}\hs$ distinguished by it. Thus, a 
\hbox{two\hh-}\hskip0ptdi\-men\-sion\-al null distribution on an oriented 
pseu\-d\hbox{o\hs-}\hskip0ptRiem\-ann\-i\-an 
\hbox{four\hh-}\hskip0ptman\-i\-fold $\,(M,g)\,$ of the neutral metric 
signature is either {\it compatible with the orientation\/} (at every point), 
or not compatible with it at any point.

The orientation distinguished by $\,\mathcal{N}\,$ also has a description 
that does not invoke self-du\-al\-i\-ty. Consequently, it can be generalized 
to all even dimensions $\,n\ge4\,$ (see Remark~\ref{hidim}).

Part (iii) of the next theorem is due to \hbox{D\'\i az\hh-}\hskip0ptRamos, 
Gar\-\hbox{c\'\i a\hs-}\hskip0ptR\'\i o and 
V\nnh\'az\-\hbox{quez\hh-}\hskip0ptLo\-ren\-zo 
\cite{diaz-ramos-garcia-rio-vazquez-lorenzo}. 
\begin{theorem}\label{psiez}Let\/ $\,\mathcal{V}\hs$ and\/ $\,\beta\,$ be 
the \hbox{two\hh-}\hskip0ptdi\-men\-sion\-al null distribution and the 
$\,1$-form, defined as in Lemma\/~{\rm\ref{sdept}} for a given 
neu\-tral-sig\-na\-ture oriented self-du\-al Ein\-stein 
\hbox{four\hh-}\hskip0ptman\-i\-fold\/ $\,(M,g)$ of Pe\-trov type\/ 
{\rm III}. 
\begin{enumerate} 
  \def\theenumi{{\rm\roman{enumi}}} 
\item[{\rm(i)}] $\mathcal{V}\hs$ is compatible with the orientation. 
\item[{\rm(ii)}] $\mathcal{V}\hs$ is parallel if and only if\/ $\,\beta=0\,$ 
identically. 
\item[{\rm(iii)}] If\/ $\,\mathcal{V}\hs$ is parallel, the metric\/ 
$\,g\,$ must be Ric\-\hbox{ci\hh-}\hskip0ptflat. 
\item[{\rm(iv)}] If\/ $\,\beta\ne0\,$ everywhere, $\,(M,g)\,$ does not admit 
a \hbox{two\hh-}\hskip0ptdi\-men\-sion\-al null parallel distribution 
compatible with the orientation. 
\end{enumerate} 
\end{theorem} 
\begin{proof}Choosing $\,w,w\hh'\nnh,v,v\hh'$ as in the lines preceding 
Remark~\ref{cvhmg}, we see that $\,v\wedge v\hh'$ is the bi\-vec\-tor 
corresponding to $\,\zeta\,$ via $\,g$. Since $\,v\,$ and $\,v\hh'$ span 
$\,\mathcal{V}\nh$, while $\,\zeta\,$ is self-du\-al, (i) follows.

As $\,\mathcal{V}\,$ has totally geodesic leaves (Lemma~\ref{cptqt}(b)), 
$\,\nabla_{\hskip-2ptv}u\,$ is a section of $\,\mathcal{V}\,$ if so are 
$\,u\,$ and $\,v$. Thus, by Lemma~\ref{sdept}(i), for $\,\mathcal{V}\hs$ to be 
parallel, it is necessary and sufficient that $\,\nabla_{\!w}u\,$ be a section 
of $\,\mathcal{V}\,$ whenever $\,u\,$ is a section of $\,\mathcal{V}\nh$, 
parallel in the direction of $\,\mathcal{V}\nh$, and $\,w\,$ is any 
$\,\mathcal{V}$-pro\-ject\-a\-ble local vector field. Since 
$\,\mathcal{V}=\mathcal{V}^\perp$ and $\,\theta\ne0\,$ everywhere (cf.\ 
Lemma~\ref{sdept}(v)), the last equality in Lemma~\ref{sdept}(iii) shows that 
$\,\mathcal{V}\,$ has the property just stated if and only if 
$\,\beta\,$ is identically zero, which yields (ii). Next, if $\,\beta=0$, the 
first equality in Lemma~\ref{sdept}(iv) gives $\,K\nh=0\,$ (since $\,\zeta\ne0\,$ 
by Lemma~\ref{sdept}(v)), proving (iii).

Suppose now that $\,\beta\ne0\,$ everywhere and $\,\tl\,$ is any 
\hbox{two\hh-}\hskip0ptdi\-men\-sion\-al null parallel distribution on 
$\,(M,g)$. If $\,v,w\,$ are sections of $\,\tl$, setting 
$\,\xi=g(v,\,\cdot\,)\,$ and $\,\xi\hs'\nh=g(w,\,\cdot\,)$, we have 
\begin{equation}\label{zvv} 
\zeta(v,w)\hs\eta+\eta(v,w)\hs\zeta+2\hh K\hs\xi\wedge\hs\xi\hs'\nh=\hh0\hs. 
\end{equation} 
This is immediate from Lemma~\ref{sdept}(ii), as $\,R(v,w,u,w\hh'\hh)=0\,$ 
for arbitrary vector fields $\,u,w\hh'\nnh$, due to the fact that, by 
(\ref{cur}), $\,R(u,w\hh'\hh)\hh v$, being a section of $\,\tl$, must be 
orthogonal to $\,w$.

Therefore, $\,\zeta(v,w)=0\,$ whenever $\,v,w\,$ are sections of $\,\tl$. In 
fact, evaluating (\ref{zvv}) on $\,(v,w)$, we get 
$\,\zeta(v,w)\hs\eta(v,w)=0$. Thus, at points where $\,\zeta(v,w)\ne0\,$ we 
would have $\,\eta(v,w)=0\,$ and, by (\ref{zvv}), $\,\eta\,$ would be a 
multiple of $\,\xi\wedge\hs\xi\hs'\nnh$, contrary to its nondegeneracy (cf.\ 
(\ref{eig})).

As $\,\beta\ne0\,$ everywhere, there exists $\,x\in M\,$ with 
$\,\tl_{\nh x}\ne\mathcal{V}_{\nh x}$ (or else $\,\mathcal{V}=\tl\,$ would be 
parallel, contradicting (ii)), and we may choose $\,v\in\tl_{\nh x}$ such that 
$\,v\notin\mathcal{V}_{\nh x}=\mathrm{Ker}\,\zeta_x$. Thus, 
$\,\zeta_xv\in\mathcal{V}_{\nh x}\smallsetminus\{0\}\,$ (by 
Lemma~\ref{sdept}(v)), while, according to the last paragraph, 
$\,\zeta_xv\in\tl^\perp_{\nh x}=\tl_{\nh x}$. Hence 
$\,\dim\hskip2pt(\hh\mathcal{V}_{\nh x}\cap\tl_{\nh x})=1$. 
Lemma~\ref{orntt}(b) and (i) now show that $\,\tl\,$ is not compatible with 
the orientation at $\,x$, or, equivalently, at any point, and (iv) follows. 
\end{proof} 
Suppose that $\,(M,g)\,$ is a type III self-du\-al oriented Ein\-stein 
\hbox{four\hh-}\hskip0ptman\-i\-fold\/ of the neutral metric signature. We 
call a point $\,x\in M\,$ {\it generic\/} if $\,\beta_x\ne0\,$ or 
$\,\beta=0\,$ at all points of some neighborhood of $\,x$. Generic points 
obviously form a dense open subset of $\,M$. Each connected component of this 
set represents either the {\it Walker case\/} ($\hh\beta=0\,$ identically), or 
the {\it strictly \itnw\ case\/} ($\hh\beta\ne0\,$ everywhere). Our 
terminology is motivated by Theorem~\ref{psiez} and the fact that null 
parallel distributions on pseu\-d\hbox{o\hs-}\hskip0ptRiem\-ann\-i\-an 
manifolds are described by Walker's classical theorem \cite{walker}.

It should be pointed out that, in any neu\-tral-sig\-na\-ture oriented 
self-du\-al \hbox{four\hh-}\hskip0ptman\-i\-fold $\,(M,g)\,$ which is 
Ric\-\hbox{ci\hh-}\hskip0ptflat but not flat, there exists a whole family, 
dif\-feo\-mor\-phic to the circle, of 
tw\hbox{o\hh-}\hskip0pt-di\-men\-sion\-al null parallel distributions which 
are {\it not\/} compatible with the orientation. In fact, for such $\,(M,g)$, 
the Le\-vi-Ci\-vi\-ta connection in the bundle $\,\varLambda^-\hskip-1.8ptM\,$ 
is well-known to be flat. The distributions in question now arise from 
nonzero null parallel local sections of $\,\varLambda^-\hskip-1.8ptM$, treated 
as an\-ti-self-du\-al bi\-vec\-tor fields. See 
\cite[p.\ 638, Proposition~37.1(i) and p.\ 645, 
Lemma~37.8]{dillen-verstraelen}.

On the other hand, not all manifolds with the stated properties admit 
tw\hbox{o\hh-}\hskip0pt-di\-men\-sion\-al null parallel distributions 
compatible with the orientation (cf.\ Section~\ref{lsth}).

\begin{remark}\label{hidim}In a $\,2m$\diml\ 
pseu\-\hbox{do\hs-}\hskip0ptEuclid\-e\-an space $\,\mathcal{T}\hs$ of the 
neutral metric signature ($m\,$ pluses, $\,m\,$ minuses), any $\,m$\diml\ null 
sub\-space $\,\mathcal{N}\hs$ naturally distinguishes an orientation of 
$\,\mathcal{T}\nnh$. Specifically, this is the orientation represented by the 
basis $\,v_1,\dots,v_m,w_1,\dots,w_m$, where $\,v_1,\dots,v_m$ is any basis of 
$\,\mathcal{N}\hs$ and $\,w_1,\dots,w_m\in\mathcal{T}\hs$ are any vectors 
satisfying the in\-ner-prod\-uct relations 
$\,\langle v_j,w_k\rangle=\delta_{jk}$. Note that, given $\,v_1,\dots,v_m$, 
such $\,w_1,\dots,w_m$ exist since the functionals $\,\mathcal{N}\to\bbR\,$ 
forming the basis dual to $\,v_1,\dots,v_m$ may be extended to 
$\,\mathcal{T}\hs$ and then represented as $\,\langle\,\cdot\,,w_k\rangle$, 
$\,k=1,\dots,m$. Also, the $\,2m$-tuple $\,v_1,\dots,v_m,w_1,\dots,w_m$ is a 
basis since the corresponding Gram matrix of inner products is nonsingular. 
Finally, the transition matrix between two such bases is easily seen to have a 
positive determinant.

When $\,m=2$, the orientation distinguished by $\,\mathcal{N}\hs$ in the 
manner just described coincides with that of Lemma~\ref{orntt}(a). See 
\cite[p.\ 638, Proposition~37.1(i)]{dillen-verstraelen}. 
\end{remark}

\section{Partial metrics and af\-fine foliations}\label{pmaf} 
\setcounter{equation}{0} 
Let $\,\mathcal{E}$ and $\,\mathcal{F}\,$ be real vector bundles over a 
manifold $\,\bs$. By a {\it pairing\/} of $\,\mathcal{E}$ and 
$\,\mathcal{F}\,$ we mean any $\,C^\infty$ section of 
$\,(\mathcal{E}\nh\otimes\mathcal{F})^*\nnh$. In other words, such a pairing 
is a $\,C^\infty$ assignment of a bi\-lin\-e\-ar mapping 
$\,\mathcal{E}_{\nh y}\times\mathcal{F}_{\nnh y}\to\bbR\,$ to every 
$\,y\in\bs$, and may also be regarded as a vec\-tor-bun\-dle morphism 
$\,\mathcal{E}\nh\to\mathcal{F}\hh^*\nnh$, or 
$\,\mathcal{F}\to\mathcal{E}^*\nnh$.

Given an $\,m$\diml\ distribution $\,\mathcal{V}\,$ on a manifold $\,M\,$ of 
dimension $\,2m$, we define a {\it partial metric\/} for $\,(M,\mathcal{V})\,$ 
to be any pairing $\,h\,$ of the vector bundles $\,\mathcal{V}\,$ and 
$\,\tm\,$ over $\,M$ which, treated as a morphism 
$\,\tm\to\mathcal{V}^*\nnh$, has the kernel $\,\mathcal{V}\nh$. (Cf.\ 
\cite[Sec.\ IV]{derdzinski-roter-06}.) Clearly, 
\begin{equation}\label{vbi} 
\mathrm{such\ }\,h\,\mathrm{\ amount\ to\ (arbitrary)\ 
vec\-tor}\hyp\mathrm{bun\-dle\ 
isomorphisms\ }\,\mathcal{V}\to[(\tm)/\hs\mathcal{V}\hs]^*\nh. 
\end{equation} 
An obvious example of a partial metric for $\,(M,\mathcal{V})\,$ is 
the restriction $\,h\,$ to $\,\mathcal{V}\,$ and $\,\tm\,$ of the pairing of 
$\,\tm\,$ and $\,\tm\,$ provided by a \psr\ metric $\,g\,$ on $\,M\,$ such 
that $\,\mathcal{V}\,$ is $\,g$-null. In this case we refer to $\,h\,$ (or, 
$\,g$) as the {\it restriction\/} of $\,g\,$ to $\,\mathcal{V}\,$ and 
$\,\tm\,$ (or, respectively, a {\it to\-tal-met\-ric extension\/} of $\,h$).

By an {\it af\-fine foliation\/} on a manifold $\,M\,$ we mean a pair 
$\,(\mathcal{V}\nh,\mathrm{D}\nh)\,$ consisting of an integrable distribution 
$\,\mathcal{V}\,$ on $\,M\,$ along with a fixed choice of a flat 
tor\-sion\-free connection $\,\mathrm{D}\,$ on each leaf of $\,\mathcal{V}\,$ 
such that, in an obvious sense, the connection depends 
$\,C^\infty$-dif\-fer\-en\-ti\-a\-bly on the leaf. Our notation ignores 
the dependence of $\,\mathrm{D}\,$ on the leaf, and, instead, treats 
$\,\mathrm{D}\,$ as a mapping that sends sections $\,u,v\,$ of 
$\,\mathcal{V}\,$ to a section $\,\mathrm{D}_uv\,$ of $\,\mathcal{V}\nh$. 
Obviously, for $\,(\mathcal{V}\nh,\mathrm{D}\nh)\,$ as above, 
\begin{equation}\label{trv} 
\begin{array}{l} 
\mathrm{the\ vector\ sub\-bundle\ }\,\,\mathcal{V}\hs\,\mathrm{\ of\ 
}\,\hs\tm\hs\,\mathrm{\ is\ locally\ trivialized\ by\ sections\ of\ 
}\,\,\mathcal{V}\\ 
\mathrm{that\ are\ }\,\mathcal{V}\hyp\mathit{parallel\hskip1.3pt}\mathrm{\ in\ 
the\ sense\ of\ being\ }\,\mathrm{D}\hh\hyp\mathrm{par\-al\-lel\ 
along\ each\ leaf\ of\ }\,\mathcal{V}\nh. 
\end{array} 
\end{equation} 
Let $\,(\mathcal{V}\nh,\mathrm{D}\nh)\,$ be an af\-fine foliation of dimension 
$\,m\,$ on a $\,2m$\diml\ manifold $\,M$, and let $\,h\,$ be a partial metric 
for $\,(\mathcal{V}\nh,\mathrm{D})\,$ (that is, for $\,(M,\mathcal{V})$). We 
will say that $\,h\,$ is 
\begin{enumerate} 
  \def\theenumi{{\rm\roman{enumi}}} 
\item {\it af\-fine\/} if, for any $\,\mathcal{V}$-par\-al\-lel section 
$\,v\,$ of $\,\mathcal{V}\,$ and any $\,\mathcal{V}$-pro\-ject\-a\-ble vector 
field $\,w\,$ (cf.\ Remark~\ref{liebr}), both defined on any open subset of 
$\,M$, the function $\,h(v,w)\,$ restricted to each leaf $\,N\hs$ of 
$\,\mathcal{V}\,$ is locally af\-fine or, equivalently, the $\,1$-form on 
$\,N\hs$ obtained by restricting $\,d\hs[\hs h(v,w)]\,$ to $\,N\hs$ is 
$\,\mathrm{D}\hs$-par\-al\-lel, 
\item {\it skew-af\-fine\/} if $\,d_v[\hs h(v,w)]=0\,$ for any $\,v,w\,$ as in 
(i), 
\item {\it trivial\/} if, for any $\,v,w\,$ as in (i), $\,h(v,w)\,$ is locally 
constant along every leaf of $\,\mathcal{V}\nh$. 
\end{enumerate} 
For $\,(M,\mathcal{V})\,$ and a partial metric $\,h\,$ as above, being trivial 
obviously implies being skew-af\-fine, while being skew-af\-fine implies being 
af\-fine: the last claim is clear since, whenever $\,u,v\hh'$ are sections of 
$\,\mathcal{V}\,$ parallel along $\,\mathcal{V}\,$ and $\,v,w\,$ are as in 
(i), $\,d_ud_v[\hs h(v\hh'\nnh,w)]\,$ must vanish due to its simultaneous 
skew-sym\-me\-try in $\,v,v\hh'$ and symmetry in $\,u,v$.

An af\-fine foliation $\,(\mathcal{V}\nh,\mathrm{D}\nh)\,$ on $\,M\,$ 
obviously arises when $\,M\,$ is an open submanifold of the total space of an 
af\-fine bundle over a manifold, $\,\mathcal{V}\,$ is the restriction to 
$\,M\,$ of the vertical distribution $\,\kerd\pi$, where $\,\pi\,$ is the 
bundle projection, and $\,\mathrm{D}\,$ is the standard flat tor\-sion\-free 
connection of each fibre. Locally, there are no other examples: any af\-fine 
foliation $\hs(\mathcal{V}\nh,\mathrm{D}\nh)\hs$ of dimension $\,k\,$ on an 
$\,n$\diml\ manifold $\,M\,$ is, locally, obtained in the manner just 
described. In fact, let us fix an $\,(n-k)$\diml\ submanifold $\,\bs\,$ of 
$\,M$, transverse to $\,\mathcal{V}\nh$, and treat it as the zero section 
$\,\bs\subset\pb\,$ in the total space $\,\pb\,$ of the vector bundle over 
$\,\bs\,$ which is the restriction of $\,\mathcal{V}\,$ to $\,\bs$. Then, at 
any point $\,y\,$ of the zero section $\,\bs$, the exponential mapping of 
$\,\mathrm{D}\,$ sends a neighborhood of $\,y\,$ in $\,\pb\hs$ \feicly\ onto 
an open set in $\,M$, in such a way that the vertical distribution in 
$\,\pb\hs$ corresponds to $\,\mathcal{V}\nh$.

The total space $\,M=\tab\,$ of the cotangent bundle of any manifold $\,\bs\,$ 
carries a trivial partial metric $\,h$, for 
$\,(\mathcal{V}\nh,\mathrm{D}\nh)\,$ defined as in the last paragraph, 
obtained by setting $\,h_x(\xi,w)=\xi(d\pi_xw)\,$ for any $\,x\in\tab=M$, any 
vertical vector $\,\xi\in\kerd\pi_x=\tayb$, with $\,y=\pi(x)$, and any 
$\,w\in\txm$, where $\,\pi:M\to\bs\,$ is the bundle projection, cf.\ 
\cite{patterson-walker}. Again, these are, locally, the only examples: 
for any trivial partial metric $\,h\,$ for 
$\,(\mathcal{V}\nh,\mathrm{D}\nh)\,$ on a manifold $\,M$, treating the leaves 
of $\,\mathcal{V}\nh$, locally, as the fibres of a bundle projection 
$\,\pi:M\to\bs$, we obtain a natural bijective correspondence between 
$\,\mathcal{V}$-par\-al\-lel sections $\,v\,$ of $\,\mathcal{V}\,$ and 
$\,1$-forms $\,\xi\,$ on $\,\bs$, given by 
$\,\xi(d\pi\hs w)=h(v,w)$, where $\,w\,$ is any 
$\,\mathcal{V}$-pro\-ject\-a\-ble local vector field in $\,M$, and 
$\,d\pi\hs w\,$ denotes its $\,\pi$-im\-age in $\,\bs$. Thus, $\,M\,$ can be 
identified, locally, with the total space of an af\-fine bundle over $\,\bs$, 
the associated vector bundle of which is $\,\tab$. The required local 
identification of $\,M$ with $\,\tab\,$ may now be obtained by choosing an 
$\,(n-k)$\diml\ submanifold $\,\bs\,$ of $\,M$, transverse to 
$\,\mathcal{V}\nh$, as in the last paragraph.

The partial metrics that naturally appear in the geometric situation discussed 
in this paper are skew-af\-fine, though not trivial.

\section{Basic octuples}\label{boct} 
\setcounter{equation}{0} 
By a {\it basic octuple\/} we mean a system 
$\,(M,\mathcal{V}\nh,\hh\mathrm{D},h,\alpha,\beta,\theta,\zeta)\,$ formed by 
a skew-af\-fine partial metric $\,h\,$ for a 
\hbox{two\hh-}\hskip0ptdi\-men\-sion\-al af\-fine foliation 
$\,(\mathcal{V}\nh,\hh\mathrm{D})\,$ on a manifold $\,M\,$ of dimension 
four, along with sections $\,\alpha,\beta,\zeta,\theta\,$ of 
$\,\mathcal{V}^*\nnh,\hs\tam,\hs[\tam]^{\wedge2}$ and, respectively, 
$\,[\mathcal{V}^*]^{\wedge2}\nnh$, such that $\,\beta\ne0$ everywhere, 
$\,\mathrm{rank}\hskip3pt\zeta=2\,$ everywhere, and 
\begin{equation}\label{dub} 
\begin{array}{rlrl} 
\mathrm{a)}\hskip-1pt&d_u[\hs h(v,w)]=\beta(w)\hs\theta(u,v)\hs, 
\hskip-2pt&\mathrm{b)}\hskip-1pt&d_u\hs[\hh\alpha(v)] 
=\alpha(u)\hs\alpha(v)\hs,\\ 
\mathrm{c)}\hskip-1pt&d_u\hs[\hh\zeta(w,w\hh'\hh)] 
=\alpha(u)\hs\zeta(w,w\hh'\hh)\hskip-2pt 
&\mathrm{d)}\hskip-1pt&d_u\hs[\hh\beta(w)]=2\hh\alpha(u)\hs\beta(w)\hs,\\ 
\mathrm{e)}\hskip-1pt&d_u\hs[\hs\theta(v,v\hh'\hh)] 
=-2\hh\alpha(u)\hs\theta(v,v\hh'\hh)\hs,\hskip-2pt 
&\mathrm{f)}\hskip-1pt&\theta(v,\zeta w)=2\hh h(v,w)\hs,\\ 
\mathrm{g)}\hskip-1pt&\alpha(\zeta w)=2\beta(w)\hs,\hskip-2pt 
&\mathrm{h)}\hskip-1pt&\mathcal{V}=\hs\mathrm{Ker}\,\zeta 
\end{array} 
\end{equation} 
for any $\,\mathcal{V}$-par\-al\-lel sections $\,u,v,v\hh'$ of 
$\,\mathcal{V}\,$ (see Section~\ref{pmaf}) and 
$\,\mathcal{V}$-pro\-ject\-a\-ble 
local vector fields $\,w,w\hh'$ in $\,M\,$ (cf.\ Remark~\ref{liebr}). Unlike 
the $\,1$-form $\,\beta\,$ and the $\,2$-form $\,\zeta\,$ on $\,M$, the 
objects $\,\alpha\,$ and $\,\theta\,$ in a basic octuple are only ``partial'' 
differential forms: $\,\alpha(v)\,$ and $\,\theta(v,v\hh'\hh)\,$ are not 
defined unless $\,v\,$ and $\,v\hh'$ are sections of $\,\mathcal{V}\nh$.

In (\ref{dub}.f), $\,\zeta w\,$ denotes the unique section of 
$\,\mathcal{V}\hs$ with $\,h(\zeta w,w\hh'\hh)=\zeta(w,w\hh'\hh)\,$ for all 
vector fields $\,w\hh'\nnh$, the existence and uniqueness of $\,\zeta w\,$ 
being clear from (\ref{vbi}). By (\ref{dub}.h) and (\ref{dub}.g), 
\begin{equation}\label{vim} 
\mathrm{i)}\hskip6pt\mathcal{V}=\hs\mathrm{Im}\,\zeta\hs,\hskip34pt 
\mathrm{ii)}\hskip6pt\mathcal{V}=\hs\mathrm{Ker}\,\zeta 
\subset\hs\mathrm{Ker}\,\beta\hs, 
\end{equation} 
where (i) expresses surjectivity of $\,\zeta\,$ treated as a morphism 
$\,\tm\to\mathcal{V}\,$ acting by $\,w\mapsto\zeta w$. For sections $\,u,v\,$ 
of $\,\mathcal{V}\,$ and vector fields $\,w,w\hh'\nnh$, 
\begin{equation}\label{bwh} 
\begin{array}{rl} 
\mathrm{a)}\hskip7pt&2\beta(w)\hs h(u,w\hh'\hh)-2\beta(w\hh'\hh)\hs h(u,w) 
=\alpha(u)\hs\zeta(w,w\hh'\hh)\hs,\\ 
\mathrm{b)}\hskip7pt&\theta(u,v)\hskip1.2pt\zeta w=2\hs[\hh h(u,w)\hs v 
-h(v,w)\hs u\hh]\hs.\\ 
\end{array} 
\end{equation} 
Namely, we obtain (\ref{bwh}.a) (or, (\ref{bwh}.b)) by first selecting a 
vector field $\,w\hh''$ with $\,u=\zeta w\hh''$ (or, vector fields 
$\,w,w\hh'$ with $\,u=\zeta w$, $\hs v=\zeta w\hh'$), then using 
(\ref{dub}.g) (or, (\ref{dub}.f)) to replace 
$\,\alpha(u)=\alpha(\zeta w\hh''\hh)\,$ by $\,2\beta(w\hh''\hh)\,$ 
(or, $\,\theta(u,v)=\theta(\zeta w,\hh\zeta w\hh'\hh)\,$ by 
$\,2\hh h(\zeta w,w\hh'\hh)=2\hs\zeta(w,w\hh'\hh)$), and, finally, 
applying Remark~\ref{cclsm}, at any $\,x\in M$, to 
$\,\plane=\txm\nnh/\hh\mathcal{V}_{\nh x}$, cf.\ (\ref{vim}.ii). 
\begin{remark}\label{duzwe}For a basic octuple 
$\,(M,\mathcal{V}\nh,\hh\mathrm{D},h,\alpha,\beta,\theta,\zeta)$, a section 
$\,u\,$ of $\,\mathcal{V}\nh$, a $\,\mathcal{V}$-pro\-ject\-a\-ble local 
vector field $\,w\,$ in $\,M$, and $\,\zeta w\,$ as above, we have 
$\,\mathrm{D}_u(\zeta w)\nh=\nh2\hh\alpha(u)\hs\zeta w\nh-\nh2\beta(w)\hs u$.

In fact, let $\,v=\mathrm{D}_u(\zeta w)-2\hh\alpha(u)\hs\zeta w\nh 
+\nh2\beta(w)\hs u$. That $\,v=0\,$ will clearly follow from (\ref{vbi}) once 
we show that $\,h(v,w\hh'\hh)=0\,$ for every $\,\mathcal{V}$-pro\-ject\-a\-ble 
local vector field $\,w\hh'$ in $\,M$. To this end, note that (\ref{dub}.c), 
(\ref{dub}.a) and the Leib\-niz rule give 
$\,\alpha(u)\hs\zeta(w,w\hh'\hh)=d_u\hs[\hs h(\zeta w,w\hh'\hh)] 
=\beta(w\hh'\hh)\hs\theta(u,\zeta w)+h(\mathrm{D}_u(\zeta w),w\hh'\hh)$, and 
so, by (\ref{dub}.f), 
$\,h(\mathrm{D}_u(\zeta w),w\hh'\hh)=\alpha(u)\hs\zeta(w,w\hh'\hh) 
-2\beta(w\hh'\hh)\hs h(u,w)$. Thus, 
$\,h(v,w\hh'\hh)=-\hs\alpha(u)\hs\zeta(w,w\hh'\hh)-2\beta(w\hh'\hh)\hs h(u,w) 
+2\beta(w)\hs h(u,w\hh'\hh)$, which vanishes in view of (\ref{bwh}.a). 
\end{remark} 
\begin{remark}\label{notus}The conclusion in Remark~\ref{duzwe} was obtained 
without using condition (\ref{dub}.b). 
\end{remark} 
\begin{remark}\label{aetbu}If 
$\,(M,\mathcal{V}\nh,\hh\mathrm{D},h,\alpha,\beta,\theta,\zeta)\,$ is 
a basic octuple, then $\,\alpha=\theta(\bu,\,\cdot\,)\,$ on 
$\,\mathcal{V}\nh$, where $\,\bu$ is the unique section of $\,\mathcal{V}\,$ 
with $\,h(\bu,\,\cdot\,)=\beta$. (Its existence and uniqueness are immediate 
from (\ref{vbi}), since (\ref{vim}.ii) allows us to treat $\,\beta\,$ as a 
section of $\,[(\tm)/\hs\mathcal{V}\hs]^*\nnh$.) Namely, writing an 
arbitrary section of $\,\mathcal{V}\,$ as $\,\zeta w$, which is allowed in 
view of (\ref{vim}.i), we see that, by (\ref{dub}.g) and (\ref{dub}.f), 
$\,\alpha(\zeta w)=2\beta(w)=2\hh h(\bu,w)=\theta(\bu,\zeta w)$. 
\end{remark} 
Our interest in basic octuples is due to the fact that they naturally arise in 
the strictly \nw\ case of our geometric situation. Specifically, we have the 
following result. 
\begin{theorem}\label{chboc}Given a neu\-tral-sig\-na\-ture oriented 
self-du\-al Ein\-stein \hbox{four\hh-}\hskip0ptman\-i\-fold\/ $\,(M,g)\,$ of 
Pe\-trov type\/ {\rm III}, let us define\/ 
$\,\alpha,\beta,\zeta,\theta,\mathcal{V}\,$ as in Lemma\/~{\rm\ref{sdept}}, 
denote by\/ $\,\mathrm{D}\,$ the restriction of the \lcc\ of\/ $\,g\,$ to the 
leaves of\/ $\,\mathcal{V}\nh$, and declare\/ $\,h\,$ to be the partial metric 
for\/ $\,(M,\mathcal{V})\,$ obtained by restricting $\,g\,$ to\/ 
$\,\mathcal{V}\,$ and\/ $\,\tm$. If\/ $\,\beta\ne0\,$ everywhere in\/ $\,M$, 
then\/ $\,(M,\mathcal{V}\nh,\hh\mathrm{D},h,\alpha,\beta,\theta,\zeta)\,$ is a 
basic octuple. 
\end{theorem} 
\begin{proof}According to Lemma~\ref{sdept}(i), 
$\,(\mathcal{V}\nh,\mathrm{D}\nh)\,$ is an af\-fine foliation, while the 
partial metric $\,h\,$ for $\,(\mathcal{V}\nh,\hh\mathrm{D})\,$ is 
skew-af\-fine in view of skew-sym\-me\-try of $\,\theta\,$ and the equality 
$\,d_u[\hs g(v,w)]=\beta(w)\hs\theta(u,v)$ in Lemma~\ref{sdept}(iii), which 
also yields (\ref{dub}.a). Next, (\ref{dub}.e) follows from the third 
equality in Lemma~\ref{cptqt}(c) with $\,\rw=W^+\nnh$, as 
$\,\nabla_{\hskip-2ptu}v=\nabla_{\hskip-2ptu}v\hh'\nnh=0$. (Note that 
(\ref{eig}) gives $\,\eta\hs v=v\,$ and $\,\eta(v,v\hh'\hh)=g(v,v\hh'\hh)=0$, 
since $\,\mathcal{V}\,$ is null by Lemma~\ref{sdept}(v).) Conditions 
(\ref{dub}.h) and (\ref{dub}.g) are in turn immediate from 
Lemma~\ref{sdept}(v) and, respectively, the equality 
$\,2\eta\beta+\zeta\alpha=0\,$ in Lemma~\ref{cptqt}(f) with $\,\rw=W^+\nnh$, 
as the vector field corresponding to $\,\beta\,$ via $\,g\,$ is a section of 
$\,\mathcal{V}$ (see Lemmas~\ref{cptqt}(e) and~\ref{sdept}(v)), so that 
$\,\eta\beta=\beta\,$ by (\ref{eig}).

On the other hand, (\ref{wze}.ii) for $\,\rw=W^+$ and (\ref{ant}) imply that 
$\,(\zeta\theta+\theta\zeta)/2=-\hs\mathrm{Id}$. Also, by 
Lemma~\ref{sdept}(v), $\,\zeta\,$ vanishes on $\,\mathcal{V}\nh$, while 
$\,\theta\,$ vanishes on $\,\mathcal{H}\,$ and maps $\,\mathcal{V}\,$ onto 
$\,\mathcal{H}$. The last equality now shows that the composite 
$\,\theta\zeta$, treated as a morphism $\,\tm\to\tm$, equals $\,-2$ times 
$\,\mathrm{Id}\,$ on $\,\mathcal{H}$, and $\,0\,$ on $\,\mathcal{V}\nh$, 
which, combined with skew-sym\-me\-try of $\,\theta$, yields (\ref{dub}.f).

Furthermore, using the first equality in Lemma~\ref{cptqt}(c) with $\,\rw=W^+$ 
and the Leib\-niz rule, we see that the left-hand side in (\ref{dub}.c) equals 
$\,\alpha(u)\hs\zeta(w,w\hh'\hh)\,$ plus 
$\,\alpha(u)\hs\zeta(w,w\hh'\hh)+\zeta(\nabla_{\hskip-2ptu}w,w\hh'\hh) 
+\zeta(w,\nabla_{\hskip-2ptu}w\hh'\hh)$, since $\,\beta(u)=0\,$ by 
Lemma~\ref{cptqt}(d). Replacing $\,\nabla_{\hskip-2ptu}w\,$ with 
$\,[\nabla_{\!u}w]^{\mathcal{H}}\nh=[\nabla_{\!w}u]^{\mathcal{H}}$ (cf.\ 
(\ref{dub}.h) and Remark~\ref{liebr}), and using Lemma~\ref{sdept}(vi), and 
setting $\,w\hh''\nnh=\theta\hh u$, we can rewrite this last expression as 
$\,\alpha(u)\hs\zeta(w,w\hh'\hh)+\beta(w)\hs\zeta(w\hh''\nnh,w\hh'\hh) 
\beta(w\hh'\hh)\hs\zeta(w,w\hh''\nnh)$, which vanishes in view of 
Remark~\ref{cclsm}. (As in the last paragraph, we see that the morphism 
$\,\zeta\theta:\tm\to\tm\,$ equals $\,-2$ times $\,\mathrm{Id}\,$ on 
$\,\mathcal{V}\nh$, and so, by (\ref{dub}.g), 
$\,\alpha(u)=-\hs\alpha(\zeta w\hh''\hh)/2=-\hs\beta(w\hh''\hh)$.)

Next, for $\,u,w\,$ as in (\ref{dub}.d), $\,\beta(u)=\beta([u,w])=0\,$ and 
$\,\zeta u=0\,$ by Lemma~\ref{cptqt}(d), Remark~\ref{liebr} and 
(\ref{dub}.h), so that, evaluating $\,(d\beta+2\beta\wedge\alpha)(u,w)\,$ from 
(\ref{bwa}) and then using the first equality in Lemma~\ref{sdept}(iv), we get 
(\ref{dub}.d).

Finally, by (\ref{trv}), equality (\ref{dub}.b) amounts to the relation 
$\,[\mathrm{D}_u\alpha](v)=\alpha(u)\hs\alpha(v)\,$ for {\it all\/} sections 
$\,u,v\,$ of $\,\mathcal{V}\nh$. Thus, since 
$\,\mathcal{V}=\hs\mathrm{Ker}\,\zeta=\hs\mathrm{Im}\,\zeta\,$ (see 
Lemma~\ref{sdept}(v)), (\ref{dub}.b) will follow if we prove the latter 
relation for $\,v=\zeta w$, where $\,w\,$ is any 
$\,\mathcal{V}$-pro\-ject\-a\-ble local vector field. However, as 
$\,[\mathrm{D}_u\alpha](\zeta w) 
=d_u\hs[\hh\alpha(\zeta w)]-\alpha(\mathrm{D}_u(\zeta w))$, while 
$\,\mathrm{D}_u(\zeta w)\,$ may be replaced by 
$\,2\hh\alpha(u)\hs\zeta w\nh-\nh2\beta(w)\hs u$ (which is allowed according 
to Remark~\ref{notus}), (\ref{dub}.b) is immediate from (\ref{dub}.d) and 
(\ref{dub}.g). 
\end{proof} 
Our discussion of basic octuples can be simplified as follows. Let us assume, 
for the remainder of this section, that $\,h\,$ is a skew-af\-fine partial 
metric for a \hbox{two\hh-}\hskip0ptdi\-men\-sion\-al af\-fine foliation 
$\,(\mathcal{V}\nh,\hh\mathrm{D})\,$ on a manifold $\,M\,$ of dimension four.

Given $\,\alpha,\beta,\zeta,\theta\,$ such that 
$\,(M,\mathcal{V}\nh,\hh\mathrm{D},h,\alpha,\beta,\theta,\zeta)\,$ is a basic 
octuple, we may choose, locally in $\,M$, a positive function $\,\phi\,$ with 
$\,\alpha(v)=-\hskip1.7ptd_v\log\hs\phi\,$ for all sections $\,v\,$ of 
$\,\mathcal{V}\nh$. (In fact, by (\ref{bwa}) and (\ref{dub}.b), the 
restriction of $\,\alpha\,$ to every leaf of $\,\mathcal{V}\,$ is a closed 
$\,1$-form.) Rephrased in terms of $\,\phi$, (\ref{dub}.b) states that 
$\,d_ud_v\phi=0$, and hence 
\begin{enumerate} 
  \def\theenumi{{\rm\roman{enumi}}} 
\item[{\rm(i)}] the restriction of $\,\phi\,$ to each leaf of 
$\,\mathcal{V}\,$ is a nonconstant positive af\-fine function. 
\end{enumerate} 
Note that $\,\phi\,$ is nonconstant since $\,\beta\ne0\,$ everywhere, and so 
$\,\alpha\ne0\,$ everywhere by (\ref{dub}.g).

If one now sets $\,\hat\alpha=\phi\hs\alpha$, $\,\hat\beta=\phi^2\beta$, 
$\,\hat\zeta=\phi\hs\zeta\,$ and $\,\hat\theta=\phi^{-2}\theta$, then, 
according to (\ref{dub}), 
\begin{enumerate} 
  \def\theenumi{{\rm\roman{enumi}}} 
\item[{\rm(ii)}] $\hat\beta\ne0\,$ everywhere, 
$\,\mathrm{rank}\hskip3pt\hat\zeta=2\,$ everywhere, and 
$\,\mathrm{Ker}\,\hat\zeta=\mathcal{V}\nh$, 
\item[{\rm(iii)}] 
$\hat\alpha(u)=-\hskip1.7ptd_u\phi$, 
$\,\hat\theta(v,\hat\zeta w)=2\hh\phi^{-1} h(v,w)$, 
$\,\hat\alpha(\hat\zeta w)=2\hat\beta(w)\,$ and 
$\,d_u[\hs h(v,w)]=\hat\beta(w)\hs\hat\theta(u,v)$, while $\,\hat\beta(w)$, 
$\hat\theta(v,v\hh'\hh)\,$ and $\,\hat\zeta(w,w\hh'\hh)\,$ are constant along 
$\,\mathcal{V}\nh$, for any $\,\mathrm{D}\hs$-par\-al\-lel sections 
$\,u,v,v\hh'$ of $\,\mathcal{V}\,$ and $\,\mathcal{V}$-pro\-ject\-a\-ble local 
vector fields $\,w,w\hh'$ in $\,M$. 
\end{enumerate} 
Conversely, if sections $\,\hat\alpha,\hat\beta,\hat\zeta,\hat\theta\,$ of 
$\,\mathcal{V}^*\nnh,\hs\tam,\hs[\tam]^{\wedge2}\nnh, 
[\mathcal{V}^*]^{\wedge2}\,$ and a function $\,\phi:M\to\bbR\,$ 
satisfy (i) -- (iii), a basic octuple 
$\,(M,\mathcal{V}\nh,\hh\mathrm{D},h,\alpha,\beta,\theta,\zeta)\,$ can clearly 
be defined by setting 
\begin{equation}\label{aeh} 
\alpha=\phi^{-1}\hat\alpha\hs,\hskip9pt\beta=\phi^{-2}\hat\beta\hs,\hskip9pt 
\zeta=\phi^{-1}\hat\zeta\hs,\hskip9pt\theta=\phi^2\hat\theta\hs. 
\end{equation}

\section{Tw\hbox{o\hh-\nh}\hskip0ptplane systems}\label{afsy} 
\setcounter{equation}{0} 
By a {\it tw\hbox{o\hh-}plane system\/} we mean a sextuple 
$\,(\bs,\xi,\ts,\plane,c,\varOmega)\,$ consisting of a real af\-fine space 
$\,\bs\,$ and a real vector space $\,\plane\,$ with $\,\dim\bs=\dim\plane=2$, 
two linearly independent constant $\,1$-forms $\,\xi,\ts\,$ on $\,\bs$, a 
nonzero constant vector field $\,c\,$ on $\,\plane$, and a nonzero constant 
$\,2$-form $\,\varOmega$ on $\,\plane$. (In other words, 
$\,c\in\plane\smallsetminus\{0\}\,$ and 
$\,\varOmega\in[\plane^*]^{\wedge2}\nnh\smallsetminus\{0\}$.)

Any tw\hbox{o\hh-\nh}\hskip0ptplane system 
$\,(\bs,\xi,\ts,\plane,c,\varOmega)\,$ gives rise to a basic octuple 
$\,(M,\mathcal{V}\nh,\hh\mathrm{D},h,\alpha,\beta,\theta,\zeta)$ defined as 
follows. Let $\,\plane_+\subset\plane\,$ be the open set on which 
$\,\varOmega(\,\cdot\,,c)>0$. Thus, $\,\plane_+$ is a connected component of 
$\,\plane\smallsetminus\line$, for the line $\,\line=\bbR\hs c\,$ spanned by 
$\,c\,$ in $\,\plane$. On the \hbox{four\hh-}\hskip0ptdi\-men\-sion\-al 
product manifold $\,M=\bs\times\plane_+$ one has the 
\hbox{two\hh-}\hskip0ptdi\-men\-sion\-al 
af\-fine foliation $\,(\mathcal{V}\nh,\mathrm{D}\nh)$ formed by the 
distribution $\,\mathcal{V}\,$ tangent to the $\,\plane_+$ factor and the 
standard flat tor\-sion\-free connection $\,\mathrm{D}\,$ on each leaf of 
$\,\mathcal{V}\nh$, the leaf being identified with the open set $\,\plane_+$ 
in the plane $\,\plane$. Next, we denote by $\,\rd\,$ the {\it radial vector 
field\/} on $\,\plane$, that is, the identity mapping $\,\plane\to\plane$ 
treated as a vector field on $\,\plane$. Vector fields on the factor manifolds 
$\,\bs\,$ and $\,\plane_+$, including constant fields $\,v\,$ (such as 
$\,v=c$) and the radial field $\,\rd\,$ on $\,\plane\,$ (restricted to 
$\,\plane_+$), and all vector fields $\,w\,$ on $\,\bs$, will also be treated 
as vector fields on $\,M=\bs\times\plane_+$, tangent to the factor 
distributions. Similarly, we will use the same symbols for differential forms 
on $\,\bs\,$ and $\,\plane_+$ as for their pull\-backs to $\,M$. Thus, 
$\,\xi\,$ and $\,\ts\,$ are now $\,1$-forms on $\,M$, and $\,\varOmega\,$ 
is a $\,2$-form on $\,M$. Using these conventions, we declare $\,h\,$ to be 
the partial metric for $\,(M,\mathcal{V})\,$ such that, for all vector fields 
$\,v,u\,$ on $\,\plane_+$ and $\,w\,$ on $\,\bs$, treated as vector fields 
on $\,M$, 
\begin{equation}\label{hvv} 
\mathrm{a)}\hskip7pth(v,u)=0\hs,\hskip17pt\mathrm{b)}\hskip7pth(v,w) 
=\varOmega(\yw,v)\hs,\hskip9pt\mathrm{where} 
\hskip7pt\mathrm{c)}\hskip7pt\yw=\xi(w)\hh\rd+\ts(w)\hh c\hs, 
\end{equation} 
$\rd\,$ being the radial field on $\,\plane$. Thus, for $\,u,v,w\,$ as above, 
$\,d_u[\hs h(v,w)]=\xi(w)\hs\varOmega(u,v)$. Skew-sym\-me\-try of 
$\,\varOmega\,$ now implies that $\,h\,$ is a skew-af\-fine partial metric for 
$\,(\mathcal{V}\nh,\hh\mathrm{D})$.

Conditions (i) -- (iii) at the end of Section~\ref{boct} are in turn satisfied 
if one sets 
\begin{equation}\label{peo} 
\phi=\varOmega(\rd,c)\hs,\hskip9pt\hat\alpha 
=-\hskip1.7ptd\hs\phi\hskip7pt(\mathrm{on\ }\,\mathcal{V})\hs,\hskip9pt 
\hat\beta=\xi\hs,\hskip9pt\hat\zeta 
=2\hs\xi\wedge\ts\hs,\hskip9pt\hat\theta=\varOmega\hskip7pt(\mathrm{on\ } 
\,\mathcal{V})\hs. 
\end{equation} 
In fact, (i) -- (iii) follow since, for vector fields $\,v\,$ on 
$\,\plane_+$ and $\,w\,$ on $\,\bs$, one clearly has 
\begin{equation}\label{hcw} 
h(c,w)=\xi(w)\hs\phi\hs,\hskip12pt 
h(\rd,w)=-\hs\ts(w)\hs\phi\hs,\hskip12pt 
\hat\zeta w=-2\phi^{-1}\yw\hs,\hskip12pt 
\hat\alpha(v)=\varOmega(c,v)\hs, 
\end{equation} 
with $\,\hat\zeta w\,$ defined as in the lines following (\ref{dub}).

Consequently, formula (\ref{aeh}) defines a basic octuple 
$\,(M,\mathcal{V}\nh,\hh\mathrm{D},h,\alpha,\beta,\theta,\zeta)$, naturally 
associated with $\,(\bs,\xi,\ts,\plane,c,\varOmega)$.

We will use the following well-known lemma to prove Theorem~\ref{unimo}, 
stating that all basic octuples represent just one local \feic\ type. 
\begin{lemma}\label{volfm}Given a manifold\/ $\,\bs\,$ of dimension\/ $\,m$, 
a point $\,y\in\bs$, a differential\/ $\,m$-form\/ $\,\nu\,$ on\/ $\,\bs$, 
and closed\/ $\,1$-forms\/ $\,\xi^1\nnh,\dots,\xi^{m-1}$ which are linearly 
independent at\/ $\,y$, there exists a closed\/ $\,1$-form\/ $\,\ts\hs$ on 
a neighborhood\/ $\,\,U\,$ of\/ $\,y\,$ such that\/ 
$\,\nu=\xi^1\nnh\wedge\ldots\wedge\xi^{m-1}\nnh\wedge\ts\,$ on\/ 
$\,\,U\nh$. 
\end{lemma} 
Namely, choosing a closed $\,1$-form $\,\xi^m$ on a neighborhood of $\,y\,$ 
so that $\,\xi^1\nnh,\dots,\xi^m$ are linearly independent at $\,y$, we have 
$\,\xi^j\nh=d\hskip.2pty^{\hs j}$ for some local coordinates 
$\,y^{\hs j}$ at $\,y\,$ and $\,j=1,\dots,m$, so that we may set 
$\,\ts=d\chi$, where $\,\chi\,$ is a function defined near $\,y\,$ with the 
partial derivative $\,\psi=\hh\partial\chi/\partial y^m$ characterized by 
$\,\nu=\psi\hs\xi^1\nnh\wedge\ldots\wedge\xi^m\nnh$. 
\begin{theorem}\label{unimo}All basic octuples, at all points in their 
underlying \hbox{four\hh-}\hskip0ptman\-i\-folds, represent one single type of 
local \feic\ equivalence.

In other words, if\/ $\,(\bs,\xi,\ts,\plane,c,\varOmega)\,$ is a fixed 
tw\hbox{o\hh-}plane system, then every basic octuple is locally \feicly\ 
equivalent to the basic octuple associated with 
$\,(\bs,\xi,\ts,\plane,c,\varOmega)$. 
\end{theorem} 
\begin{proof}Given a basic octuple 
$\,(M,\mathcal{V}\nh,\hh\mathrm{D},h,\alpha,\beta,\theta,\zeta)$, let us 
choose $\,\phi\,$ and the corresponding objects 
$\,\hat\alpha,\hat\beta,\hat\zeta,\hat\theta\,$ as at the end of 
Section~\ref{boct}. A neighborhood of any given point of $\,M\,$ may, clearly, 
be \feicly\ identified with an open subset of the total space of a real 
af\-\hbox{fine\hh-}\hskip0ptplane bundle $\,\apb\,$ over a surface $\,\bs\,$ 
in such a way that $\,\mathcal{V}\,$ and $\,\mathrm{D}\,$ become the vertical 
distribution and the standard flat tor\-sion\-free connection in each fibre of 
$\,\apb\,$ (treated as an open set in an af\-fine plane). Since $\,\phi$, 
restricted to each fibre of $\,\apb$, is a nonconstant af\-fine function, its 
zero set is the total space of a real af\-\hbox{fine\hh-}\hskip0ptline 
sub\-bun\-dle $\,\alb\,$ of $\,\apb$. Conditions (i) -- (iii) at the end of 
Section~\ref{boct} will remain unaffected if we multiply $\,\phi\,$ by any 
positive function $\,\bs\to\bbR\,$ (pulled back to the total space $\,\apb$), 
at the same time multiplying $\,\hat\alpha,\hat\beta,\hat\zeta,\hat\theta\,$ 
by its appropriate powers. A suitable choice of such a positive function 
allows us to assume that $\,\hat\beta\,$ is closed. (Note that $\,\hat\beta\,$ 
and $\,\hat\zeta\,$ are the pull\-backs to $\,\apb\,$ of a $\,1$-form and a 
$\,2$-form without zeros on the surface $\,\bs$.) Setting $\,\xi=\hat\beta$, 
we may now use Lemma~\ref{volfm} (with $\,m=2$, $\,\nu=\hat\zeta/2\,$ and 
$\,\xi^1\nh=\xi$) to select, locally, a closed $\,1$-form $\,\ts\,$ in 
$\,\bs\,$ with $\,\hat\zeta=2\hs\xi\wedge\ts$.

Let us also set $\,\varOmega=\hat\theta\,$ and $\,c=\phi^3\bu$, for the 
section $\,\bu\,$ of $\,\mathcal{V}\,$ defined in Remark~\ref{aetbu}. Then, 
on $\,\mathcal{V}\nh$, we have $\,d\hs\phi=\varOmega(\,\cdot\,,c)\,$ (in other 
words, $\,d_v\phi=\varOmega(v,c)\,$ for every section $\,v\,$ of 
$\,\mathcal{V}$). This is clear since, on $\,\mathcal{V}\nh$, (iii) in 
Section~\ref{boct} and Remark~\ref{aetbu} give 
$\,d\hs\phi=-\hs\hat\alpha=-\hs\phi\hs\alpha=\phi^{-2}\theta(\,\cdot\,,c) 
=\hat\theta(\,\cdot\,,c)=\varOmega(\,\cdot\,,c)$. On the other hand, as 
$\,h(\bu,\,\cdot\,)=\beta\,$ (see Remark~\ref{aetbu}), it follows that 
$\,h(c,\,\cdot\,)=\phi\hs\xi$, and, consequently, 
$\,d_v[\hs h(c,w)]=\xi(w)\hskip1ptd_v\phi=\xi(w)\hskip1pt\varOmega(v,c)\,$ 
for every $\,\mathcal{V}$-pro\-ject\-a\-ble local vector field $\,w\,$ in 
$\,\apb\,$ and every section $\,v\,$ of $\,\mathcal{V}\nh$. At the same time, 
(\ref{dub}.a) and the Leib\-niz rule give 
$\,d_v[\hs h(c,w)]-h(\mathrm{D}_vc,w)=\beta(w)\hs\theta(v,c) 
=\xi(w)\hskip1pt\varOmega(v,c)$, so that $\,\mathrm{D}_vc=0$. Hence $\,c$, 
restricted to each fibre of $\,\apb\,$ (which is an open subset of an af\-fine 
plane), is a nonzero constant vector field.

Any fixed local section $\,z\,$ of 
$\,\alb\,$ gives rise to the section $\,\trd\,$ of the vertical distribution 
$\,\mathcal{V}\,$ on $\,\apb$, with the value at $\,x\in\apb\,$ equal to 
$\,x-z_{\hs\pi(x)}$, where $\,\pi\,$ is the bundle projection. (Thus, 
$\,\trd\,$ restricted to the fibre of $\,\apb\,$ containing $\,x\,$ is the 
radial vector field relative to the origin $\,z_{\hs\pi(x)}$.) Then 
$\,\phi=\varOmega(\trd,c)$. In fact, as we saw above, 
$\,d\hs\phi=\varOmega(\,\cdot\,,c)\,$ on $\,\mathcal{V}\nh$, so that 
$\,\phi\,$ and $\,\varOmega(\trd,c)\,$ have the same 
$\,d_v$-de\-riv\-a\-tive for any section $\,v\,$ of $\,\mathcal{V}\nh$, and, 
consequently, differ in each fibre of $\,\apb\,$ by a constant, while, due to 
our choice of $\,z$, they both vanish at the origin $\,z_{\hs\pi(x)}$ in the 
fibre containing $\,x$. Defining a $\,1$-form $\,\tts\,$ on the total space 
$\,\apb$ by $\,\tts=-\hs\phi^{-1}h(\trd,\,\cdot\,)$, we in turn obtain 
$\,h(v,w)=\varOmega(\tyw,v)$, with $\,\tyw=\xi(w)\hh\trd\nh+\tts(w)\hh c$, for 
all sections $\,v\,$ of $\,\mathcal{V}\,$ and all vector fields $\,w$. Namely, 
since $\,c\,$ and $\,\trd\,$ span $\,\mathcal{V}\,$ away from $\,\alb$, it 
suffices to consider the cases $\,v=c\,$ and $\,v=\trd$, in which the required 
equality follows since $\,h(c,w)=\xi(w)\hs\phi=\xi(w)\hs\varOmega(\trd,c)\,$ 
(as we saw earlier) and 
$\,h(\trd,w)=-\hs\tts(w)\hs\phi=\tts(w)\hs\varOmega(c,\trd)\,$ (by the 
definition of $\,\tts$). Suppose now that $\,w\,$ is a 
$\,\mathcal{V}$-pro\-ject\-a\-ble local vector field in $\,\apb\,$ and $\,v\,$ 
is a section of $\,\mathcal{V}\nh$. By (\ref{dub}.a), 
$\,d_v[\hs h(u,w)]=\xi(w)\hskip1pt\varOmega(v,u)\,$ if $\,u\,$ is a 
$\,\mathrm{D}\hs$-par\-al\-lel section of $\,\mathcal{V}\nh$, while, as shown 
above, $\,d\hs\phi=\varOmega(\,\cdot\,,c)\,$ on $\,\mathcal{V}\nh$. Thus, the 
Leib\-niz rule yields $\,d_v[\hs\tts(w)]=-\hs d_v[\hs\phi^{-1}h(\trd,w)] 
=-\hs\phi^{-1}\xi(w)\hskip1pt\varOmega(v,\trd)-\phi^{-1}h(v,w) 
+\phi^{-2}\varOmega(v,c)\hs h(\trd,w)$. Since $\,h(v,w)=\varOmega(\tyw,v)$, 
we get $\,d_v[\hs\tts(w)]=0$, and so $\,\tts\,$ is the pull\-back to $\,\apb\,$ 
of a $\,1$-form in $\,\bs$. Furthermore, 
$\,\hat\zeta=2\hs\xi\wedge\tts\,$ as a consequence of (\ref{bwh}.a) with 
$\,\xi=\hat\beta=\phi^2\beta$, the definition of $\,\tts$, (\ref{bwa}.ii), 
and the relation $\,\alpha(\trd)=-\hs1\,$ (immediate since, on 
$\,\mathcal{V}\nh$, we have $\,d\hs\phi=-\hs\hat\alpha=-\hs\phi\hs\alpha\,$ 
and $\,d\hs\phi=\varOmega(\,\cdot\,,c)$, while $\,\phi=\varOmega(\trd,c)$). As 
$\,\hat\zeta=2\hs\xi\wedge\ts\,$ for the closed $\,1$-form $\,\ts\,$ 
selected earlier, there exists a function $\,\psi\,$ in $\,\bs\,$ with 
$\,\ts=\tts+\psi\hs\xi$. Replacing $\,z\,$ by $\,z-\psi c\hs\,$ causes 
$\,\trd\,$ and $\,\tts\,$ to be replaced by $\,\trd+\psi c\,$ 
and, respectively, by $\,\ts$. With $\,z-\psi c\,$ (the new choice of $\,z$) 
declared the zero section, our af\-\hbox{fine\hh-}\hskip0ptplane bundle 
$\,\apb\,$ may be treated as a vector bundle $\,\pb\nnh$, in such a way that 
$\,c\,$ and $\,\varOmega\,$ are sections, both without zeros, of $\,\pb\,$ and 
$\,[\pb^*]^{\wedge2}\nnh$. Choosing, locally, a section $\,a\,$ of $\,\pb\,$ 
with $\,\varOmega(a,c)=1$, we obtain a local trivialization $\,a,c\,$ of 
$\,\pb$, which allows us to view $\,\pb\,$ as a product bundle of the form 
$\,\bs\times\plane$. Sections of $\,\pb\,$ now become functions on $\,\bs\,$ 
valued in the vector space $\,\plane$, with $\,c\,$ corresponding in this way 
to a constant function (an element of $\,\hs\plane$). Finally, $\,\bs\,$ may 
be identified, locally, with the space $\,\rto$ so as to make $\,\hs\xi\,$ 
and $\,\ts\,$ correspond to $\,\nnh d\hskip.2pty^{\nh1}\nnh$ and 
$d\hskip.2pty^2\nnh$, for the standard coordinates $\,y^{\hs j}$ in 
$\hh\rto\nnh$. The resulting 
tw\hbox{o\hh-}plane system $\,(\bs,\xi,\ts,\plane,c,\varOmega)\,$ clearly has 
$\,(M,\mathcal{V}\nh,\hh\mathrm{D},h,\alpha,\beta,\theta,\zeta)\,$ as its 
associated basic octuple. 
\end{proof}

\section{Horizontal distributions}\label{hrzd} 
\setcounter{equation}{0} 
By a {\it horizontal distribution\/} for a basic octuple 
$\,(M,\mathcal{V}\nh,\hh\mathrm{D},h,\alpha,\beta,\theta,\zeta)$, cf.\ 
Section~\ref{boct}, we mean a vector sub\-bun\-dle $\,\mathcal{H}\,$ of 
$\,\tm\,$ with $\,\tm=\mathcal{H}\oplus\mathcal{V}\nh$.

Any such $\,\mathcal{H}\,$ gives rise to a neu\-tral-sig\-na\-ture \prc\ 
$\,g\,$ on $\,M$. Namely, 
\begin{equation}\label{gxh} 
g\,\mathrm{\ is\ the\ unique\ to\-tal}\hyp\mathrm{met\-ric\ extension\ of\ }\,h 
\,\mathrm{\ such\ that\ }\,\mathcal{H}\,\mathrm{\ is\ }\,g\hyp\mathrm{null.} 
\end{equation} 
We denote by $\,\nabla\,$ the \lcc\ of $\,g$, and by $\,\gamma\,$ the 
$\,1$-form on $\,M\,$ with 
\begin{equation}\label{gam} 
g(\nabla_{\!w}w\hh'\nnh,w\hh''\hh) 
=-\hs\gamma(w)\hskip1pt\zeta(w\hh'\nnh,w\hh''\hh)\hskip4pt\mathrm{\ for\ all\ 
vector\ fields\ }\,w\,\mathrm{\ and\ sections\ }\,w\hh'\nnh,w\hh''\mathrm{\ 
of\ }\,\mathcal{H}\hs, 
\end{equation} 
$\gamma\,$ being well defined since $\,\zeta\,$ trivializes 
$\,\mathcal{H}^{\wedge2}\nnh$, while skew-sym\-me\-try of 
$\,g(\nabla_{\!w}w\hh'\nnh,w\hh''\hh)\,$ in $\,w\hh'\nnh,w\hh''$ implies its 
{\it val\-ue\-wise\/} dependence on $\,w,w\hh'$ and $\,w\hh''\nnh$. If $\,v\,$ 
is a section of $\,\mathcal{V}\nh$, we let $\,\theta\hh v\,$ stand for the 
unique section of $\,\mathcal{H}\,$ such that 
$\,g(\theta\hh v,u)=\theta(v,u)\,$ whenever $\,u\,$ is a section of 
$\,\mathcal{V}\nh$.

Next, we denote by $\,R\,$ the curvature tensor of $\,g$, and by $\,\eta\,$ 
the $\,2$-form satisfying (\ref{eig}) with our 
$\,\mathcal{V}\nh,\hs\mathcal{H}\,$ and $\,g$, so that, for sections 
$\,v,v\hh'$ of $\,\mathcal{V}\hs$ and $\,w,w\hh'$ of $\,\mathcal{H}$, 
\begin{equation}\label{eta} 
\eta(v,v\hh'\hh)=\eta(w,w\hh'\hh)=0\hs,\hskip18pt\eta(v,w)=-\hs\eta(w,v) 
=h(v,w)\hs. 
\end{equation} 
The symbol $\,\bw\,$ will be used for the unique section of $\,\mathcal{H}\,$ 
with 
\begin{equation}\label{brw} 
h(v,\bw)\,=\,\alpha(v)\hs\mathrm{\ for\ all\ sections\ }\,v\,\mathrm{\ 
of\ }\,\mathcal{V}. 
\end{equation} 
That (\ref{brw}) defines $\,\bw\,$ uniquely is clear from (\ref{gxh}), since 
$\,\tm\,=\,\mathcal{H}\oplus\mathcal{V}\nh$. For this $\,\bw$, 
\begin{equation}\label{bwz} 
\mathrm{i)}\hskip6pt\zeta(w,\bw)=\alpha(\zeta w)=2\beta(w)\hs,\hskip14pt 
\mathrm{ii)}\hskip6pt\beta(\bw)\,\,=\,\,0\hs, 
\end{equation} 
where $\,w\,$ in (i) is an arbitrary vector field. In fact, (\ref{brw}) and 
(\ref{dub}.g) yield (\ref{bwz}.i), while (\ref{bwz}.ii) follows from 
(\ref{bwz}.i) along with skew-sym\-me\-try of $\,\zeta$.

Three further objects associated with a horizontal distribution 
$\,\mathcal{H}\,$ for 
$\,(M,\mathcal{V}\nh,\hh\mathrm{D},h,\alpha,\beta,\theta,\zeta)$ are 
extensions of the ``partial'' differential forms $\,\alpha,\theta\,$ 
(see Section~\ref{boct}) to differential forms on $\,M$, still written as 
$\,\alpha,\theta$, which are given by 
$\,\alpha(w)=2\gamma(\zeta w)\,$ and $\,\theta(w,\,\cdot\,)=0\,$ for sections 
$\,w\,$ of $\,\mathcal{H}$, with $\,\gamma\,$ and $\,\zeta w\,$ as in 
(\ref{gam}) and (\ref{dub}.f), and a $\,1$-form $\,Z\,$ on $\,M\,$ 
characterized by $\,Z(w)\hs\zeta=\nabla_{\!w}\hs\zeta-2\hh\alpha(w)\hs\zeta\,$ 
on $\,\mathcal{H}$, that is, $\,Z(w)\hs\zeta(w\hh'\nnh,w\hh''\hh) 
=2\hh[\nabla_{\!w}\hs\zeta]\hh(w\hh'\nnh,w\hh''\hh) 
-4\hh\alpha(w)\hs\zeta(w\hh'\nnh,w\hh''\hh)\,$ for any sections 
$\,w\hh'\nnh,w\hh''$ of $\,\mathcal{H}$, and any vector field $\,w$. Note that 
we thus have $\,\mathcal{H}=\mathrm{Ker}\,\theta$, and $\,Z$ is well defined, 
since $\,\zeta\,$ trivializes $\,\mathcal{H}^{\wedge2}\nnh$.

For simplicity, our notation ignores the dependence of 
$\,g,\nabla\nnh,\gamma,\theta\hh v,R,\eta,\bw,\alpha,\theta\,$ and $\,Z\,$ on 
$\,\mathcal{H}$. 
\begin{lemma}\label{pctbl}If\/ $\,\mathcal{H}\,$ is a horizontal distribution 
for a basic octuple 
$\,(M,\mathcal{V}\nh,\hh\mathrm{D},h,\alpha,\beta,\theta,\zeta)$, and\/ 
$\,\bw\,$ is the section of\/ $\,\mathcal{H}\,$ defined by\/ {\rm(\ref{brw})}, 
while\/ $\,\phi\,$ is chosen as at the end of Section\/~{\rm\ref{boct}}, then 
the vector field\/ $\,\phi\hh\bw\,$ is $\,\mathcal{V}$-pro\-ject\-a\-ble. 
\end{lemma} 
\begin{proof}Let $\,w\,$ be a $\,\mathcal{V}$-pro\-ject\-a\-ble local section 
of $\,\mathcal{H}$, chosen so as to agree with $\,\phi\hh\bw\,$ at all points of 
a given surface $\,\bs\hh'$ embedded in $\,M\,$ and transverse to 
$\,\mathcal{V}\nh$. Since $\,\beta(w)=0\,$ on $\,\bs\hh'$ by (\ref{bwz}.ii), 
relation (\ref{dub}.d) combined with uniqueness of solutions for 
first-or\-der linear ordinary differential equations gives 
$\,\beta(w)=0\,$ on the union $\,\,U\,$ of all leaves of $\,\mathcal{V}\,$ 
that intersect $\,\bs\hh'\nnh$. However, $\,\beta\ne0\,$ everywhere, and 
so $\,\mathcal{H}\cap\mathrm{Ker}\,\beta\,$ is spanned by $\,\bw$, cf.\ 
(\ref{bwz}.ii) and (\ref{vim}.ii). Thus, $\,w=\chi\phi\hh\bw\,$ on $\,\,U\nh$, 
where $\,\chi\,$ is some function without zeros. For any 
$\,\mathcal{V}$-pro\-ject\-a\-ble local section $\,w\hh'$ of $\,\mathcal{H}$, 
(iii) in Section~\ref{boct} implies that $\,\phi\hs\zeta(w,w\hh'\hh)\,$ and 
$\,\phi^2\beta(w\hh'\hh)\,$ are constant along $\,\mathcal{V}\nh$, while, by 
(\ref{bwz}.i), $\,\phi\hs\zeta(w,w\hh'\hh)=\chi\phi^2\hn\zeta(\bw,w\hh'\hh) 
=-2\chi\phi^2\hn\beta(w\hh'\hh)$. Hence $\,\chi\,$ is constant along 
$\,\mathcal{V}\nh$, and $\,\phi\hh\bw=\chi^{-1}\nh w\,$ is 
$\,\mathcal{V}$-pro\-ject\-a\-ble. 
\end{proof}

\section{Properties of the associated metric}\label{pram} 
\setcounter{equation}{0} 
Let $\,\mathcal{H}\,$ be a horizontal distribution for a basic octuple 
$\,(M,\mathcal{V}\nh,\hh\mathrm{D},h,\alpha,\beta,\theta,\zeta)$. We have 
\begin{equation}\label{ned} 
\mathrm{a)}\hskip10pt 
\nabla_{\hskip-2ptv}u=\mathrm{D}_vu\hs,\hskip36pt\mathrm{b)}\hskip10pt 
\nabla_{\hskip-2ptv}w=\beta(w)\hs\theta\hh v-\gamma(v)\hskip1pt\zeta w\hs, 
\end{equation} 
whenever $\,w\,$ is a $\,\mathcal{V}$-pro\-ject\-a\-ble local section of 
$\,\mathcal{H}\,$ and $\,u,v\,$ are sections of $\,\mathcal{V}\nh$. 
(As before, $\,\zeta w$ denotes the section of $\,\mathcal{V}\,$ 
appearing in (\ref{dub}.f), and $\,g,\nabla,\theta\hh v,\gamma\,$ are defined 
as in Section~\ref{hrzd}.) In fact, (\ref{ned}.a) follows since, by 
(\ref{lcc}), $\,\nabla_{\hskip-2ptv}u=0\,$ for sections $\,u,v\,$ of 
$\,\mathcal{V}\,$ which are $\,\mathcal{V}$-par\-al\-lel. Namely, 
$\,g(\nabla_{\hskip-2ptv}u,w)=0\,$ both when $\,w\,$ is a section of 
$\,\mathcal{V}\,$ (all six terms resulting from (\ref{lcc}) then vanish as 
$\,\mathcal{V}\,$ is integrable and $\,g$-null), and when $\,w\,$ is a 
$\,\mathcal{V}$-pro\-ject\-a\-ble section of $\,\mathcal{H}\,$ (the last four 
terms in (\ref{lcc}) vanish, again, according to Remark~\ref{liebr}, while the 
sum of the first two is zero in view of (\ref{dub}.a) and skew-sym\-me\-try of 
$\,\theta$). Similarly, to obtain (\ref{ned}.b), we take the $\,g$-in\-ner 
product of both sides with any section of $\,\mathcal{H}$, or, respectively, 
with any $\,\mathcal{V}$-par\-al\-lel section $\,u$ of $\,\mathcal{V}\,$ 
(assuming $\,v\,$ to be $\,\mathcal{V}$-par\-al\-lel as well): in the former 
case the equality is obvious from (\ref{gam}); in the latter, as 
$\,\mathcal{V}\,$ is $\,g$-null, the right-hand side yields 
$\,\beta(w)\hs\theta(v,u)$, which, by (\ref{dub}.a), is the same as 
$\,g(\nabla_{\hskip-2ptv}w,u)=d_v[\hs g(w,u)]=d_v[\hs h(u,w)]$.

Also, for sections $\,v\,$ of $\,\mathcal{V}\,$ and $\,w,w\hh'$ of 
$\,\mathcal{H}$, with $\,\theta\hh v\,$ and $\,Z\,$ as in Section~\ref{hrzd}, 
\begin{equation}\label{ztz} 
\mathrm{i)}\hskip6pt\zeta\theta\hh v=-2v\hs,\hskip12pt 
\mathrm{ii)}\hskip6pt\theta\zeta w=-2w\hs,\hskip12pt 
\mathrm{iii)}\hskip6ptZ(v)=0\hs,\hskip12pt 
\mathrm{iv)}\hskip6pt[\nabla_{\!w}w\hh'\hh]^{\mathcal{V}}\nh 
=-\hs\gamma(w)\hs\zeta w\hh'\nh, 
\end{equation} 
$[\hskip3pt]^{\mathcal{V}}$ denoting the $\,\mathcal{V}$-com\-po\-nent 
relative to the decomposition $\,\tm=\mathcal{H}\oplus\mathcal{V}\nh$. Namely, 
for such $\,v\,$ and $\,w$, (\ref{dub}.f) gives 
$\,-\hs g(\zeta\theta\hh v,w)=g(\zeta w,\theta\hh v) 
=\theta(v,\zeta w)=2g(v,w)\,$ and, similarly, 
$\,-\hs g(\theta\zeta\hh w,v)=\theta(v,\zeta w)=2g(w,v)$, so that 
(\ref{ztz}.i) and (\ref{ztz}.ii) follow as $\,\mathcal{V}\,$ and 
$\,\mathcal{H}\,$ are $\,g$-null. On the other hand, in view of (\ref{dub}.c) 
and the Leib\-niz rule, the definition of $\,Z\,$ in Section~\ref{hrzd} gives 
$\,-\hs Z(v)\hs\zeta(w,w\hh'\hh)/2= 
\zeta(\nabla_{\hskip-2ptv}w,w\hh'\hh)+\zeta(w,\nabla_{\hskip-2ptv}w\hh'\hh) 
+\alpha(v)\hs\zeta(w,w\hh'\hh)\,$ for any 
$\,\mathcal{V}$-pro\-ject\-a\-ble local sections $\,w,w\hh'$ of 
$\,\mathcal{H}$. Since $\,\zeta(\nabla_{\hskip-2ptv}w,w\hh'\hh) 
=2\beta(w)\hs h(u,w\hh'\hh)\,$ by (\ref{ned}.b), (\ref{vim}) and 
(\ref{ztz}.i), the relation $\,Z(v)=0\,$ is now immediate from (\ref{bwh}.a). 
Finally, (\ref{ztz}.iv) is an obvious consequence of (\ref{gam}) and 
(\ref{gxh}).

Furthermore, for any $\,\mathcal{V}$-par\-al\-lel sections $\,u,v\,$ of 
$\,\mathcal{V}\nh$, 
\begin{equation}\label{nut} 
\nabla_{\hskip-2ptu}(\theta\hh v)\, 
=\,2\gamma(u)\hs v\,-\,2\hh\alpha(u)\hs\theta\hh v\hs. 
\end{equation} 
To verify (\ref{nut}), we will show that both sides have equal 
$\,g$-in\-ner products with any $\,\mathcal{V}$-par\-al\-lel section 
$\,v\hh'$ of $\,\mathcal{V}\nh$, and with any 
$\,\mathcal{V}$-pro\-ject\-a\-ble section $\,w\,$ of $\,\mathcal{H}$. For 
$\,v\hh'\nnh$, this is clear as $\,\theta(v,v\hh'\hh)\,$ equals 
$\,h(\theta\hh v,v\hh'\hh)$, that is, $\,g(\theta\hh v,v\hh'\hh)$, and so 
applying $\,d_u$ we get, from (\ref{dub}.e), 
$\,g(\nabla_{\hskip-2ptu}(\theta\hh v),v\hh'\hh) 
=-2\hh\alpha(u)\hs\theta(v,v\hh'\hh)$, as required. (By (\ref{ned}.a), 
$\,\nabla_{\hskip-2ptu}v\hh'\nh=0$, while $\,g(v,v\hh'\hh)=0\,$ since 
$\,\mathcal{V}\,$ is $\,g$-null.) For $\,w$, (\ref{gxh}) allows us to 
differentiate by parts, obtaining $\,g(\nabla_{\hskip-2ptu}(\theta\hh v),w) 
=-\hs g(\nabla_{\hskip-2ptu}w,\theta\hh v)$. In view of (\ref{gam}), the last 
expression equals 
$\,\gamma(u)\hs\zeta(w,\theta\hh v)=-\hs\gamma(u)\hs g(\zeta\theta\hh v,w)$, 
which, by (\ref{ztz}.i), coincides with $\,2\gamma(u)\hs g(v,w)$. 
\begin{lemma}\label{gwzww}Given a horizontal distribution\/ $\,\mathcal{H}\,$ 
for a basic octuple\/ 
$\,(M,\mathcal{V}\nh,\hh\mathrm{D},h,\alpha,\beta,\theta,\zeta)$, let\/ 
$\,g,\nabla,\gamma,Z\,$ and\/ $\,R\,$ be as in\/ Section\/~{\rm\ref{hrzd}}. 
Then, with\/ $\,R(w,w\hh'\nnh,u,v)\,$ given by\/ {\rm(\ref{rww})}, 
\begin{enumerate} 
  \def\theenumi{{\rm\roman{enumi}}} 
\item[{\rm(i)}] $R(v,v\hh'\hh)\hh v\hh''\nh=0\,$ for all sections 
$\,v,v\hh'\nnh,v\hh''$ of\/ $\,\mathcal{V}\nh$, 
\item[{\rm(ii)}] $\gamma(w)\hskip1pt\zeta(w\hh'\nnh,w\hh''\hh) 
=g(w,[w\hh'\nnh,w\hh''\hh])\,$ for all sections $\,w,w\hh'\nnh,w\hh''$ of\/ 
$\,\mathcal{H}$, 
\item[{\rm(iii)}] $R(w,u)\hh v=\nabla_{\hskip-2ptu}\nabla_{\!w}v 
=\mathrm{D}_u[w,v]+2\beta(w)\hs\gamma(u)\hs v 
-\mathrm{D}_u[\gamma(v)\hs\zeta w]\,$ whenever $\,u,v\,$ are 
$\,\mathcal{V}$-par\-al\-lel sections of\/ $\,\mathcal{V}\,$ and $\,w\,$ is 
a $\,\mathcal{V}$-pro\-ject\-a\-ble local vector field in\/ $\,M$, 
\item[{\rm(iv)}] $R(u,v,w,w\hh'\hh) 
=[(d\gamma)(v,u)]\hs\zeta(w,w\hh'\hh) 
+\gamma(v)\hs[\nabla_{\hskip-2ptu}\zeta](w,w\hh'\hh) 
-\gamma(u)\hs[\nabla_{\hskip-2ptv}\zeta](w,w\hh'\hh)\,$ for all sections 
$\,w,w\hh'$ of\/ $\,\mathcal{H}\,$ and all vector fields $\,u,v$, 
\item[{\rm(v)}] $R(w,w\hh'\nnh,\,\cdot\,,\,\cdot\,) 
=\zeta(w,w\hh'\hh)\hs[\varGamma-\gamma\wedge Z+\eta/2\hh] 
+K\hs\xi\wedge\hs\xi\hs'$ for any real constant\/ $\,K\,$ and sections\/ 
$\,w,w\hh'$ of\/ $\,\mathcal{H}$, with\/ $\,\xi=g(w,\,\cdot\,)$, 
$\,\xi\hs'\nh=g(w\hh'\nnh,\,\cdot\,)\,$ and\/ 
$\,\varGamma\hs=\,d\gamma+2\hh\alpha\wedge\gamma-(K\theta+\eta)/2$. 
\end{enumerate} 
\end{lemma} 
\begin{proof}Flatness of $\,\mathrm{D}$, (\ref{cur}) and (\ref{ned}.a) 
yield (i). Next, (\ref{gxh}) implies (ii): by (\ref{gam}) and (\ref{lcc}), 
$\,2\gamma(w)\hskip1pt\zeta(w\hh'\nnh,w\hh''\hh) 
=-2\hh g(\nabla_{\!w}w\hh'\nnh,w\hh''\hh) 
=-\hs g(w\hh'\nnh,[w\hh''\nnh,w])-g(w\hh''\nnh,[w,w\hh'\hh]) 
+g(w,[w\hh'\nnh,w\hh''\hh])=2\hh g(w,[w\hh'\nnh,w\hh''\hh])$. The last 
equality follows here from Remark~\ref{cclsm}, since the dependence of 
$\,g(w,[w\hh'\nnh,w\hh''\hh])\,$ on $\,w\hh'\nnh,w\hh''$ is both 
skew-sym\-met\-ric and valuewise (in view of (\ref{gxh}), one may replace 
$\,[w\hh'\nnh,w\hh''\hh]\,$ with its $\,\mathcal{V}$-com\-po\-nent 
$\,[w\hh'\nnh,w\hh''\hh]^{\mathcal{V}}$).

In (iii), $\,R(w,u)\hh v=\nabla_{\hskip-2ptu}\nabla_{\!w}v\,$ by 
(\ref{cur}) and Remark~\ref{liebr}. The other equality is now immediate from 
(\ref{ned}.b), since $\,\nabla_{\!w}v=[w,v]+\nabla_{\hskip-2ptv}w$, and 
(\ref{ned}.a) combined with Remark~\ref{liebr} imply that 
$\,\nabla_{\hskip-2ptu}[w,v]=\mathrm{D}_u[w,v]$, while 
$\,\nabla_{\hskip-2ptu}[\beta(w)\hs\theta\hh v]=2\beta(w)\hs\gamma(u)\hs v\,$ 
in view of (\ref{dub}.d) and (\ref{nut}).

For $\,u,v,w,w\hh'$ as in (iv), the Leib\-niz rule and (\ref{gam}) give 
$\,g(\nabla_{\!u}\hskip-1pt\nabla_{\!v}w,\nh w\hh'\hh) 
=-\hs d_u[\gamma(v)\hskip1pt\zeta(w,\nh w\hh'\hh)] 
-g(\nabla_{\!v}w,\nabla_{\!u}w\hh'\hh)$. As $\,\mathcal{V}\,$ and 
$\,\mathcal{H}\,$ are $\,g$-null, 
$\,g(\nabla_{\!v}w,\nabla_{\!u}w\hh'\hh) 
=g([\nabla_{\!v}w]^{\mathcal{H}}\nnh,\nabla_{\!u}w\hh'\hh) 
+g(\nabla_{\!v}w,[\nabla_{\!u}w\hh'\hh]^{\mathcal{H}}\hh)$, where 
$\,[\hskip3pt]^{\mathcal{H}}$ denotes the $\,\mathcal{H}$-component. Using 
(\ref{gam}) and (\ref{dub}.h), we now obtain 
$\,g(\nabla_{\!v}w,\nabla_{\!u}w\hh'\hh) 
=-\hs\gamma(u)\hskip1pt\zeta(w\hh'\nnh,\nabla_{\!v}w) 
-\gamma(v)\hskip1pt\zeta(w,\nabla_{\!u}w\hh'\hh)$, and (iv) easily 
follows from the above equalities combined with 
(\ref{rww}), (\ref{cur}), (\ref{bwa}.iii), (\ref{gam}) and the Leib\-niz 
rule.

Finally, (v) is immediate from (iv), the definitions of $\,Z\,$ and 
$\,\theta\,$ in Section~\ref{hrzd}, and (\ref{bwa}.ii), since (\ref{bwh}.b) 
gives $\,\zeta(w,w\hh'\hh)\hskip1.2pt\theta(u,v) 
=2\hs[\hh g(w,u)\hs g(w\hh'\nnh,v)-g(w\hh'\nnh,u)\hs g(w,v)\hh]\,$ 
for any sections $\,w,w\hh'$ of $\,\mathcal{H}\,$ and any vector fields 
$\,u,v$. 
\end{proof} 
\begin{lemma}\label{gflat}Let\/ 
$\,(M,\mathcal{V}\nh,\hh\mathrm{D},h,\alpha,\beta,\theta,\zeta)\,$ be the 
basic octuple obtained as in Section\/~{\rm\ref{afsy}} from a given 
tw\hbox{o\hh-}plane system\/ $\,(\bs,\xi,\ts,\plane,c,\varOmega)$. Then the 
distribution $\,\mathcal{H}\,$ on\/ $\,M=\bs\times\plane_+$ tangent to the 
factor plane\/ $\,\bs\,$ is a horizontal distribution for\/ 
$\,(M,\mathcal{V}\nh,\hh\mathrm{D},h,\alpha,\beta,\theta,\zeta)$, the 
metric\/ $\,g\,$ on\/ $\,M\,$ is flat, $\,\mathcal{H}\,$ is 
$\,g$-par\-al\-lel, and\/ $\,\gamma=0$, where\/ $\,g\,$ and the $\,1$-form 
$\,\gamma\,$ are associated with $\,\mathcal{H}\,$ as in\/ 
Section\/~{\rm\ref{hrzd}}. 
\end{lemma} 
\begin{proof}We fix a function $\,f:\bs\to\bbR\,$ with $\,\df=\xi$. For 
$\,u,v\in\plane$, let $\,\chi^{u,v}$ be the $\,1$-form on $\,M\,$ equal to 
$\,h(v,\,\cdot\,)\,$ on $\,\mathcal{H}\,$ and to 
$\,\varOmega(u-fv,\,\cdot\,)\,$ on $\,\mathcal{V}\nh$. As the $\,1$-forms 
$\,\xi,\ts\,$ on $\,\bs\,$ are constant, and hence closed, using 
(\ref{hvv}) and (\ref{bwa}.iii) we easily verify that 
$\,d\chi^{u,v}\nh=0\,$ for all $\,u,v\in\plane$.

The assignment $\,(u,v)\mapsto\chi^{u,v}$ is a linear operator, with the 
domain $\,\plane\times\plane$, and so its image $\,\mathcal{X}\,$ is a vector 
space. The $\,g$-in\-ner product $\,g(\chi,\chi\hh'\hh)\,$ of any 
$\,\chi,\chi\hh'\nh\in\mathcal{X}\,$ is constant on $\,M$. In fact, we may 
assume that $\,\chi=\chi\hh'\nh=\chi^{u,v}\nnh$. Now, as $\,v\,$ is the 
$\,\mathcal{V}$-com\-po\-nent $\,w^{\mathcal{V}}$ of the vector field $\,w\,$ 
such that $\,\chi=g(w,\,\cdot\,)$, while $\,\mathcal{V}\,$ and 
$\,\mathcal{H}\,$ are $\,g$-null, we get 
$\,g(\chi,\chi)/\hh2=g(w^{\mathcal{H}}\nnh,w^{\mathcal{V}}) 
=g(w,w^{\mathcal{V}})=\chi(w^{\mathcal{V}})=\chi(v) 
=\varOmega(u,v)\,$ due to skew-sym\-me\-try of $\,\varOmega$, as required.

Any fixed basis of $\,\mathcal{X}\,$ thus consists of forms which, locally, 
are the differentials of functions forming a coordinate system in $\,M$. 
According to the last paragraph, the components of $\,g$ in such 
coordinates are constant, so that $\,g\,$ is flat, and all 
$\,\chi^{u,v}$ are $\,g$-par\-al\-lel. Hence $\,\mathcal{H}\,$ is 
$\,g$-par\-al\-lel, being the simultaneous kernel of all $\,\chi^{u,v}$ with 
$\,v=0$. Finally, as $\,\mathcal{H}\,$ is $\,g$-par\-al\-lel and $\,g$-null, 
(\ref{gam}) gives $\,\gamma=0$. 
\end{proof} 
\begin{remark}\label{tsmgm}For $\,(M,g)\,$ satisfying the assumptions of 
Theorem~\ref{chboc} and such that $\,\beta\ne0$ everywhere in $\,M$, let 
$\,(M,\mathcal{V}\nh,\hh\mathrm{D},h,\alpha,\beta,\theta,\zeta)\,$ be the 
corresponding basic octuple. Then $\,\gamma$ defined by (\ref{gam}) is the 
same as in Lemma~\ref{cptqt}(c) with $\,\rw=W^+\nnh$. In fact, letting 
$\,\gamma\,$ stand for the latter, we have, by (\ref{eig}), 
$\,-\hs g(\nabla_{\!w}w\hh'\nnh,w\hh''\hh) 
=g(\nabla_{\!w}(\eta\hs w\hh'\hh),w\hh''\hh) 
=g(\eta(\nabla_{\!w}w\hh'\hh),w\hh''\hh) 
-g([\nabla_{\!w}\eta]\hh w\hh'\nnh,w\hh''\hh)$ whenever $\,w\hh'\nnh,w\hh''$ 
are sections of $\,\mathcal{H}\,$ and $\,w\,$ is any vector field. On the 
other hand, (\ref{eig}) gives 
$\,g(\eta(\nabla_{\!w}w\hh'\hh),w\hh''\hh) 
=-\hs g(\nabla_{\!w}w\hh'\nnh,\eta\hs w\hh''\hh) 
=g(\nabla_{\!w}w\hh'\nnh,w\hh''\hh)$, and, as $\,\theta\hh w\hh'\nh=0\,$ 
(see Lemma~\ref{cptqt}(a)), using Lemma~\ref{cptqt}(c) we get 
$\,g([\nabla_{\!w}\eta]\hh w\hh'\nnh,w\hh''\hh) 
=2\gamma(w)\hskip1pt\zeta(w\hh'\nnh,w\hh''\hh)$.

Furthermore, with $\,\bu\,$ denoting the section of $\,\mathcal{V}\,$ defined 
in Remark~\ref{aetbu}, the function $\,\gamma(\bu)$ is a local geometric 
invariant of $\,\tg$, since so are $\,\mathcal{V}\nh,\alpha,\beta,\gamma\,$ 
(due to the uniqueness assertions in Lemma~\ref{cptqt}), and, consequently, 
$\,\bu$. 
\end{remark}

\section{Curvature conditions}\label{crvc} 
\setcounter{equation}{0} 
Our next goal is to determine which horizontal distributions $\,\mathcal{H}\,$ 
for a given basic octuple 
$\,(M,\mathcal{V}\nh,\hh\mathrm{D},h,\alpha,\beta,\theta,\zeta)\,$ lead to 
metrics $\,g\,$ that are Ein\-stein and, at the same time, self-du\-al of 
Pe\-trov type III. Rather than approach this property of $\,g\,$ directly, we 
begin by describing some conditions, namely, (a) -- (d) in 
Theorem~\ref{crvtc}, which are equivalent to it, yet easier to verify. We 
refer to them as {\it curvature conditions}, since the curvature tensor 
$\,R\,$ explicitly appears in (a), while (b) and (d) involve the curvature 
forms of the Le\-vi-Ci\-vi\-ta connection in the bundle 
$\,\varLambda^+\hskip-1.8ptM$, expressed in terms of the connection forms 
$\,\alpha,\beta,\gamma$. 
\begin{lemma}\label{ptiii}If\/ $\,\mathcal{H}\,$ is a horizontal distribution 
for a basic octuple 
$\,(M,\mathcal{V}\nh,\hh\mathrm{D},h,\alpha,\beta,\theta,\zeta)$, and\/ 
$\,K\,$ is a real constant, while\/ $\,g,R\,$ and\/ $\,\eta\,$ correspond to\/ 
$\,\mathcal{H}\,$ as in\/ Section\/~{\rm\ref{hrzd}}, then the following two 
conditions are equivalent\/{\rm:} 
\begin{enumerate} 
  \def\theenumi{{\rm\alph{enumi}}} 
\item[{\rm(i)}] $M\,$ is orientable and, for a suitable orientation, 
$\,(M,g)\,$ is a self-du\-al Ein\-stein \hbox{four\hh-}\hskip0ptman\-i\-fold 
of Pe\-trov type\/ {\rm III}, its scalar curvature equals $\,12\hh K$, while 
our$\,(M,\mathcal{V}\nh,\hh\mathrm{D},h,\alpha,\beta,\theta,\zeta)\,$ and\/ 
$\,\mathcal{H}\,$ coincide with those determined by\/ $\,g\,$ as in 
Theorem\/{\rm~\ref{chboc}} and Lemma\/~{\rm\ref{sdept}}, 
\item[{\rm(ii)}] $2R\,=\,\hs\zeta\otimes\eta\,+\,\eta\otimes\zeta\, 
+\,2Kg\wedge g\hs$, where the notation of\/ {\rm(\ref{tre})} is used. 
\end{enumerate} 
\end{lemma} 
\begin{proof}That (i) implies (ii) is obvious from Lemma~\ref{sdept}(ii) and 
(\ref{bwa}.i). Next, let us assume (ii). The Ric\-ci tensor of $\,g\,$ then 
equals $\,3Kg$. (The Ric\-ci contraction applied to 
$\,\zeta\otimes\eta+\eta\otimes\zeta\,$ yields $\,0$, as 
$\,\mathcal{V}=\mathrm{Ker}\,\zeta=\mathrm{Im}\,\zeta\,$ according to 
(\ref{vim}), so that $\,\zeta\,$ sends 
$\,\mathrm{Ker}\,(\eta\hh\pm\nh\mathrm{Id})\,$ into 
$\,\mathrm{Ker}\,(\eta\hh\mp\nh\mathrm{Id})$, and hence an\-ti\-com\-mutes 
with $\,\eta$.) Therefore, by (\ref{tre}), the Weyl tensor $\,W\hs$ of $\,g\,$ 
is equal to $\,(\zeta\otimes\eta+\eta\otimes\zeta)/2$. The hypotheses of 
Lemma~\ref{rkvpt} are thus satisfied by our $\,(M,g)\,$ and $\,\er=2$, since 
the $\,2$-forms $\,\zeta\,$ and $\,\eta$, spanning the image of 
$\,W\nh$, are linearly independent at each point by (\ref{dub}.h) and 
(\ref{eig}), while the relations $\,\lg\zeta,\eta\rg=\lg\zeta,\zeta\rg=0\,$ 
and $\,\lg\eta,\eta\rg=-2\,$ (immediate from (\ref{rom}.b), as 
$\,\zeta,\eta\,$ an\-ti\-com\-mute, 
$\,\mathrm{Ker}\,\zeta=\mathrm{Im}\,\zeta$, while $\,\eta\eta=\mathrm{Id}$) 
show that the image of $\,W\hs$ is $\,\lr$-de\-gen\-er\-ate, but not 
$\,\lr$-null. Now Lemma~\ref{rkvpt}(iii) yields (i). 
\end{proof} 
\begin{theorem}\label{crvtc}Under the assumptions of Lemma~{\rm\ref{ptiii}}, 
condition\/ {\rm(i)} in Lemma~{\rm\ref{ptiii}} holds if and only if, for all 
sections $\,u,v\nh\,$ of\/ $\,\mathcal{V}\nh$, and\/ $\,w\,$ of\/ 
$\,\mathcal{H}$, with $\,g,\nabla\nnh,R,\alpha,\theta,\gamma\,$ as in\/ 
Section\/~{\rm\ref{hrzd}}, 
\begin{enumerate} 
  \def\theenumi{{\rm\roman{enumi}}} 
\item[{\rm(a)}] $R(w,u)\hh v\,=Kh(v,w)\hh u$,\hskip79.6pt{\rm(b)} 
\hskip2.6pt$d\beta+2\beta\wedge\alpha=-\hs K\zeta/\hh2$, 
\item[{\rm(c)}] $[\nabla_{\!w}\theta](u,v)=-2\hh\alpha(w)\hh\theta(u,v)$, 
\hskip49.4pt{\rm(d)}\hskip5pt$(d\gamma+2\hh\alpha\wedge\gamma)(\,\cdot\,,w) 
=g(\,\cdot\,,w)/2$. 
\end{enumerate} 
\end{theorem} 
\begin{proof}Assuming condition (i) in Lemma~\ref{ptiii}, we obtain (c) (or, 
respectively, (b) and (d)) from Lemma~\ref{cptqt}(c) with $\,\rw=W^+$ or, 
respectively, from Lemma~\ref{sdept}(iv), cf.\ Remark~\ref{tsmgm}. Note that 
$\,\eta\,$ satisfies (\ref{eig}), and hence (\ref{eta}), cf.\ 
Lemma~\ref{cptqt}(a), while $\,\theta(\,\cdot\,,w)=0\,$ by 
Lemma~\ref{sdept}(v). Also, $\,\alpha(w)\,$ defined in Section~\ref{hrzd} is 
the same as in Lemma~\ref{cptqt} with $\,\rw=W^+\nnh$, as $\,\theta\hh w=0\,$ 
and $\,\eta\hs w=-\hs w\,$ by Lemma~\ref{sdept}(v) and (\ref{eig}), and so the 
equality $\,2\hs\zeta\gamma+\eta\alpha+\theta\beta=0\,$ in 
Lemma~\ref{cptqt}(f), evaluated on $\,w$, gives 
$\,0=-\hs\beta(\theta\hh w)-\alpha(\eta\hs w) 
-2\gamma(\zeta w)=\alpha(w)-2\gamma(\zeta w)$.

Next, if (i) in Lemma~\ref{ptiii} holds, so does (ii). Since, for sections 
$\,u,v\,$ of $\,\mathcal{V}\nh$, Lemma~\ref{sdept}(v) gives 
$\,\zeta(\,\cdot\,,u)=\zeta(v,\,\cdot\,)=0\,$ and $\,g(u,v)=0$, while 
$\,g(u,\,\cdot\,)=h(u,\,\cdot\,)$, this yields (a).

Conversely, suppose that (a) -- (d) are satisfied.

For sections $\,w\,$ of $\,\mathcal{H}\,$ and $\,u\,$ of $\,\mathcal{V}\nh$, 
using the notation of (\ref{bue}.i), we now have 
\begin{equation}\label{nwt} 
[\nabla_{\!w}\theta]\hh u=-2\hh\alpha(w)\hh\theta\hh u+2\gamma(w)\hh u\hs. 
\end{equation} 
This is verified by taking the $\,g$-in\-ner products of both sides in 
(\ref{nwt}) with sections $\,v\,$ of $\,\mathcal{V}$ and $\,w\hh'$ of 
$\,\mathcal{H}$. In the former case, the agreement is obvious from (c), as 
$\,\mathcal{V}\,$ is $\,g$-null. In the latter, the Leib\-niz rule implies 
that $\,[\nabla_{\!w}\theta](u,w\hh'\hh)=-\hs\theta(u,\nabla_{\!w}w\hh'\hh) 
=-\hs\theta(u,[\nabla_{\!w}w\hh'\hh]^{\mathcal{V}}\hh)$, with 
$\,[\hskip3pt]^{\mathcal{V}}$ denoting the $\,\mathcal{V}$-com\-po\-nent 
(since $\,\theta(\,\cdot\,,w\hh'\hh)=0$, cf.\ Section~\ref{hrzd}). The 
required equality 
$\,[\nabla_{\!w}\theta](u,w\hh'\hh)=2\gamma(w)\hh h(u,w\hh'\hh)\,$ now 
follows from (\ref{ztz}.iv) and (\ref{dub}.f).

If $\,w,w\hh'$ are sections of $\,\mathcal{H}$, using the Leib\-niz rule we 
obtain $\,\theta\hh[(\nabla_{\!w}\hs\zeta)\hh w\hh'\hh] 
=\theta\hh[\nabla_{\!w}(\zeta w\hh'\hh)]-\theta\zeta\nabla_{\!w}w\hh'\nh 
=\nabla_{\!w}(\theta\zeta w\hh'\hh)-[\nabla_{\!w}\theta]\hs\zeta w\hh'\nh 
-\theta\zeta\nabla_{\!w}w\hh'\nnh$. By (\ref{ztz}.ii) and (\ref{dub}.h), 
$\,\theta\zeta w\hh'\nh=-2w\hh'$ and 
$\,\theta\zeta\nabla_{\!w}w\hh'\nh=[\nabla_{\!w}w\hh'\hh]^{\mathcal{H}}$ 
(the $\,\mathcal{H}$-com\-po\-nent of $\,\nabla_{\!w}w\hh'$), so that 
(\ref{nwt}) applied to $\,u=\zeta w\hh'$ yields 
$\,\theta\hh[(\nabla_{\!w}\hs\zeta)\hh w\hh'\hh] 
=-2\hs[\nabla_{\!w}w\hh'\hh]^{\mathcal{V}}\nh 
+2\hh\alpha(w)\hh\theta\zeta w\hh'\nh-2\gamma(w)\hs\zeta w\hh'\nh$. Thus, 
(\ref{ztz}.iv) gives $\,Z(w)=0\,$ for all sections $\,w\,$ of $\,\mathcal{H}$, 
where $\,Z\,$ is the $\,1$-form defined in Section~\ref{hrzd}. (We have just 
shown that $\,Z(w)\hs\zeta(w\hh'\nnh,w\hh''\hh)=0$ for sections 
$\,w,w\hh'\nnh,w\hh''$ of $\,\mathcal{H}\,$ with 
$\,w\hh''\nh=\theta\hh u\,$ for some vector field $\,u$, while such $\,w\hh''$ 
range over all sections of $\,\mathcal{H}\,$ due to (\ref{ztz}.ii).)

Combined with (\ref{ztz}.iii), the conclusion of the last paragraph yields 
$\,Z=0$. In addition, $\,\varGamma(w,\,\cdot\,)=0$, for the $\,2$-form 
$\,\varGamma\,$ appearing in Lemma~\ref{gwzww}(v), and all sections $\,w\,$ of 
$\,\mathcal{H}$. Namely, (\ref{eta}) and the definition of $\,\theta\,$ in 
Section~\ref{hrzd} give $\,\eta(\,\cdot\,,w)=h(\,\cdot\,,w)=g(\,\cdot\,,w)$, 
cf.\ (\ref{gxh}), and $\,\theta(\,\cdot\,,w)=0$, so that 
$\,\varGamma(\,\cdot\,,w)=0\,$ by (d).

In view of Lemma~\ref{ptiii}, it now suffices to verify that both sides in 
Lemma~\ref{ptiii}(ii) yield the same value when applied to any quadruple of 
vector fields, each of which is a section of $\,\mathcal{H}$ or 
$\,\mathcal{V}\nh$. In the following discussion of the possible cases, we will 
evaluate the right-hand side in Lemma~\ref{ptiii}(ii) on the four vector 
fields using, without further explanation, relations (\ref{eta}) and 
(\ref{dub}.h) along with the fact that $\,\mathcal{V}\,$ and $\,\mathcal{H}\,$ 
are both $\,g$-null. Due to well-known symmetries of $\,R$, only four cases 
need to be considered.

When three or four of the vector fields are sections of $\,\mathcal{V}\nh$, 
both sides vanish (Lemma~\ref{gwzww}(i)).

When the first vector field is a section of 
$\,\mathcal{H}$, while the second and third ones are sections of 
$\,\mathcal{V}\nh$, both sides yield the same value in view of (a).

When the first two vector fields are sections $\,v,u\,$ of $\,\mathcal{V}\,$ 
and the third one is a section $\,w$ of $\,\mathcal{V}\nh$, the first 
Bianchi identity gives $\,R(v,u)\hh w=R(w,u)\hh v-R(w,v)\hh u$, and our 
equality is an obvious consequence of (a).

Finally, when the first three vector fields are sections of $\,\mathcal{H}$, 
the required equality is immediate from Lemma~\ref{gwzww}(v), since, as we 
saw, $\,\varGamma(w,\,\cdot\,)=0$. 
\end{proof} 
\begin{remark}\label{rednd}The reader may have noticed that relation (b) in 
Theorem~\ref{crvtc} was not used in the second (sufficiency) part of the 
proof. In other words, (b) is a consequence of (a), (c) and (d). It is 
nevertheless convenient, due to the structure of our argument, to list (b) as 
a separate condition. See the proof of Lemma~\ref{dcfeq}. 
\end{remark}

\section{Deformations of horizontal distributions}\label{dfhd} 
\setcounter{equation}{0} 
Horizontal distributions for a fixed basic octuple 
$\,(M,\mathcal{V}\nh,\hh\mathrm{D},h,\alpha,\beta,\theta,\zeta)\,$ may be 
thought of as arbitrary sections of a specific locally trivial bundle 
$\,\mathcal{C}\hs$ over $\,M$. Its fibre $\,\mathcal{C}\nh_x$ at $\,x\in M\,$ 
consists of all vector subspaces $\,\mathcal{H}_x\subset\txm\,$ with 
$\,\txm=\mathcal{H}_x\oplus\mathcal{V}_x$. One can turn $\,\mathcal{C}\hs$ 
into an af\-fine bundle over $\,M$, having as its associated vector bundle the 
sub\-bun\-dle $\,\mathcal{F}\,$ of $\,\hs\mathrm{Hom}\hs(\tm,\mathcal{V})\,$ 
with the fibre $\,\mathcal{F}\hskip-2pt_x$ at any $\,x\in M\,$ formed by all 
operators $\,\txm\to\mathcal{V}_x$ sending $\,\mathcal{V}_x$ to $\,\{0\}$. 
Thus, 
\begin{equation}\label{mor} 
\mathrm{\ sections\ }\,F\,\mathrm{\ of\ }\,\mathcal{F}\hs\mathrm{\ are\ 
morphisms\ }\,\tm\to\tm\,\mathrm{\ valued\ in\ }\,\mathcal{V}\,\mathrm{\ and\ 
vanishing\ on\ }\,\mathcal{V}\nh. 
\end{equation} 
Specifically, given a horizontal distribution $\,\mathcal{H}\,$ and a section 
$\,F\,$ of $\,\mathcal{F}\nnh$, we declare the sum $\,\thor=\mathcal{H}+F\,$ to be 
a new horizontal distribution for 
$\,(M,\mathcal{V}\nh,\hh\mathrm{D},h,\alpha,\beta,\theta,\zeta)$, the sections 
of which have the form $\,\tw=w+Fw$, with $\,w\,$ ranging over all sections of 
$\,\mathcal{H}$.

A section $\,F\,$ of $\,\mathcal{F}\,$ associates with any 
twice-co\-var\-i\-ant \tf\ $\,\bz\,$ on $\,M\,$ two further tensor fields, 
$\,\fd\hs\bz\,$ and $\,F^*\bz$, defined by 
$\,(\fd\hs\bz)(v,v\hh'\hh)=\bz\hh(Fv,v\hh'\hh)+\bz\hh(v,Fv\hh'\hh)\,$ and 
$\,(F^*\bz)(v,v\hh'\hh)=\bz\hh(Fv,Fv\hh'\hh)\,$ for arbitrary vector fields 
$\,v,v\hh'\nnh$. Next, we denote by $\,\trf\,$ the function $\,M\to\bbR$ 
equal to $\,-1/4\,$ times the (pointwise) trace of the bundle morphism 
$\,\mathcal{H}\to\mathcal{H}\,$ sending a section $\,w\,$ of $\,\mathcal{H}\,$ 
to the section $\,\theta\hh v\,$ of $\,\mathcal{H}\,$ 
defined as in Section~\ref{hrzd} with $\,v=Fw$. For any vector fields 
$\,w,w\hh'$ on $\,M$, and any morphism $\,\varPhi\,$ of $\,\tm\,$ into any 
vector bundle over $\,M$, vanishing on $\,\mathcal{V}\nh$, 
\begin{equation}\label{trc} 
\begin{array}{rl} 
\mathrm{i)}\hskip0pt&h(Fw,w\hh'\hh)\,-\,h(Fw\hh'\nnh,w)\,\hs 
=\hs\,2\hs\trf\hs\zeta(w,w\hh'\hh)\hs,\\ 
\mathrm{ii)}\hskip0pt&2\beta(w)\hh\varPhi w\hh'\nh-2\beta(w\hh'\hh)\hh\varPhi 
w\,=\,\zeta(w,w\hh'\hh)\hh\varPhi\bw\hs, 
\end{array} 
\end{equation} 
with $\,\bw\,$ as in (\ref{brw}). Namely, since $\,\mathcal{V}\,$ is 
$\,h$-null, both sides in (\ref{trc}.i) and (\ref{trc}.ii) equal $\,0$ due 
to (\ref{mor}) and (\ref{vim}.ii) when one of $\,w,w\hh'$ is a section of 
$\,\mathcal{V}$. We may therefore assume that $\,w\,$ and $\,w\hh'$ are 
sections of $\,\mathcal{H}$. Remark~\ref{trvol} gives 
$\,-\hs4\hs\trf\hs\zeta(w,w\hh'\hh)=\zeta(\theta Fw,w\hh'\hh) 
+\zeta(w,\theta Fw\hh'\hh)=h(\zeta\theta Fw,w\hh'\hh) 
-h(\zeta\theta Fw\hh'\nnh,w)$, and so (\ref{ztz}.i) now implies (\ref{trc}.i), 
while, by (\ref{bwz}.i), $\,2\beta(w)=\zeta(w,\bw)\,$ and 
$\,-2\beta(w\hh'\hh)=\zeta(\bw,w\hh'\hh)$, and hence (\ref{trc}.ii) follows as 
$\,\zeta(w,w\hh'\hh)\hh\varPhi\bw\,$ summed cyclically over 
$\,w,w\hh'\nnh,\bw\,$ yields $\,0\,$ in view of Remark~\ref{cclsm}. 
\begin{lemma}\label{hpehf}Let\/ $\,\mathcal{H}\,$ be a horizontal 
distribution for a basic octuple\/ 
$\,(M,\mathcal{V}\nh,\hh\mathrm{D},h,\alpha,\beta,\theta,\zeta)$. If\/ 
$\,g,\gamma,\bw\,$ and the $\,1$-form $\,\alpha\,$ on $\,M\,$ are associated 
with $\,\mathcal{H}\,$ as in Section\/~{\rm\ref{hrzd}}, and\/ $\,F\nh\,$ is a 
section of\/ $\,\,\mathcal{F}\nnh$, then\/ $\,g,\alpha,\gamma\,$ and their 
analogues $\,\tg,\ta,\tga\,$ corresponding to the horizontal distribution\/ 
$\,\thor=\mathcal{H}\hh+F\,$ are related by 
\begin{equation}\label{ape} 
\begin{array}{rl} 
\mathrm{a)}\hskip0pt&\tg\hs=g-\nnh\fd g\hs,\\ 
\mathrm{b)}\hskip0pt&\ta\hs=\hs\alpha-F^*\nnh\alpha- 
2\hskip1.4ptd_{\hs\zeta(\hs\cdot\hs)}\trf-4\hh\trf\hh\beta 
+2\hh h(F\bw,\,\cdot\,)\hs,\\ 
\mathrm{c)}\hskip0pt&\tga(v) 
=\gamma(v)-d_v\trf-\trf\hs\alpha(v)+\theta(F\bw,v)/\hh2\hs,\\ 
\mathrm{d)}\hskip0pt&[\hs\tga(\tw)-\gamma(w)\hh]\hs\zeta(w\hh'\nnh,w\hh''\hh) 
=h([w\hh'\nnh,Fw\hh''\hh],w)+h([Fw\hh'\nnh,w\hh''\hh],w)\\ 
&\phantom{[\hs\tga(\tw)-\gamma(w)\hh]\hs\zeta(w\hh'\nnh,w\hh''\hh)} 
+\hskip2pth([Fw\hh'\nnh,Fw\hh''\hh],w)-h(F\hh[w\hh'\nnh,w\hh''\hh],w)\hs, 
\end{array} 
\end{equation} 
for\/ $\,\trf\,$ defined above, whenever\/ $\,v\,$ is a section of\/ 
$\,\mathcal{V}\,$ and\/ $\,w,w\hh'\nnh,w\hh''$ are\/ 
$\,\mathcal{V}$-pro\-ject\-a\-ble local sections\/ of $\,\mathcal{H}$, 
the section\/ $\,\tw\,$ of\/ $\,\thor=\mathcal{H}+F\,$ is given by\/ 
$\,\tw=w+Fw$, while\/ $\,F^*\nnh\alpha\,$ and\/ 
$\,d_{\hs\zeta(\hs\cdot\hs)}\hh f$, for any function $\,f$, denote the 
$\,1$-forms such that\/ $\,[F^*\nnh\alpha](u)=\alpha(Fu)\,$ and\/ 
$\,[\hh d_{\hs\zeta(\hs\cdot\hs)}\hh f\hh](u)=d_{\hs\zeta u}f$ for 
all vector fields\/ $\,u$. 
\end{lemma} 
\begin{proof}Let $\,g\hs'$ be the right-hand side of (\ref{ape}.a). We thus 
have $\,g\hs'\nnh(v,\,\cdot\,)=g(v,\,\cdot\,)=h(v,\,\cdot\,)\,$ for sections 
$\,v\,$ of $\,\mathcal{V}\,$ (in view of (\ref{mor}), since $\,\mathcal{V}\,$ 
is $\,g$-null), and, for the same reason, $\,g\hs'\nnh(\tw,\tw)=0$ if 
$\,\tw\,$ is a section of $\,\thor\,$ (that is, $\,\tw=w+Fw$ for some section 
$\,w\,$ of $\,\mathcal{H}$), which proves (\ref{ape}.a).

For $\,\mathcal{V}$-pro\-ject\-a\-ble sections $\,w,w\hh'$ of 
$\,\mathcal{H}\,$ and a $\,\mathcal{V}$-par\-al\-lel section $\,v\,$ of 
$\,\mathcal{V}\nh$, one has 
\begin{equation}\label{gvz} 
\begin{array}{rcl} 
-2\gamma(v)\hskip1pt\zeta(w,w\hh'\hh)& 
=&d_w[\hs h(v,w\hh'\hh)\hh]\,\,-\,\,\hs d_{w\hh'}[\hs h(v,w)\hh]\\ 
&+&h([w\hh'\nnh,v],w)\,\,+\,\,h([v,w],w\hh'\hh)\,\,-\,\,h(v,[w,w\hh'\hh])\hh, 
\end{array} 
\end{equation} 
where we write $\,h\,$ rather than $\,g\,$ since, in each inner product, one 
of the vector fields involved is a section of $\,\mathcal{V}$, cf.\ 
Remark~\ref{liebr}. This is immediate from (\ref{gam}) combined with 
(\ref{lcc}); the first of the six terms provided by (\ref{lcc}) vanishes here 
in view of (\ref{gxh}).

As (\ref{gvz}) holds for any horizontal distribution, including $\,\thor$, it 
remains valid if one replaces $\,\gamma(v),w\,$ and $\,w\hh'$ with 
$\,\tga(v),\tw=w+Fw\,$ and $\,\tw\hh'\nnh=w\hh'\nnh+Fw\hh'\nnh$. Since 
$\,Fw\,$ and $\,Fw\hh'$ are sections of the $\,h$-null distribution 
$\,\mathcal{V}$, Remark~\ref{liebr} implies that the right-hand side of the 
analogue of (\ref{gvz}) corresponding to the triple 
$\,(\thor,\tw,\tw\hh'\hh)\,$ equals its original version for 
$\,(\mathcal{H},w,w\hh'\hh)\,$ plus 
\begin{equation}\label{plu} 
d_{Fw}[\hs h(v,w\hh'\hh)\hh]\,-\,\hs d_{Fw\hh'}[\hs h(v,w)\hh]\, 
+\,h([Fw\hh'\nnh,v],w)\,+\,h([v,Fw],w\hh'\hh)\hs. 
\end{equation} 
On the other hand, by (\ref{dub}.a) and (\ref{trc}.ii) with $\,\varPhi=F\nnh$, for 
any $\,\mathcal{V}$-par\-al\-lel section $\,v\,$ of $\,\mathcal{V}\nh$, 
\begin{enumerate} 
  \def\theenumi{{\rm\roman{enumi}}} 
\item[{\rm(i)}] $d_{Fw}[\hs h(v,w\hh'\hh)\hh]\, 
-\,\hs d_{Fw\hh'}[\hs h(v,w)\hh]\, 
=\,\beta(w\hh'\hh)\hs\theta(Fw,v)\,-\,\beta(w)\hs\theta(Fw\hh'\nnh,v)$, 
\item[{\rm(ii)}] $\beta(w\hh'\hh)\hs\theta(Fw,v)\, 
-\,\beta(w)\hs\theta(Fw\hh'\nnh,v)\, 
=\,\theta(v,F\bw)\hs\zeta(w,w\hh'\hh)/\hh2$, 
\item[{\rm(iii)}] $h([Fw\hh'\nnh,v],w)\,+\,h([v,Fw],w\hh'\hh)\, 
=\,h(\mathrm{D}_v(Fw),w\hh'\hh)\,-\,h(\mathrm{D}_v(Fw\hh'\hh),w)$. 
\end{enumerate} 
Applying $\,d_v$ to (\ref{trc}.i) and using (\ref{dub}.a) along with 
(\ref{dub}.c) and the Leib\-niz rule, we see that the difference of the 
right-hand sides in (iii) and (i) is 
$\,2\hs\zeta(w,w\hh'\hh)\hs[\hs d_v\trf+\trf\hs\alpha(v)]$. Thus, 
by (i) -- (iii), the expression (\ref{plu}) is equal to twice the right-hand 
side of (ii), plus the difference just mentioned, that is, to 
$\,\zeta(w,w\hh'\hh)\hs[2\hskip1.4ptd_v\trf+2\hs\trf\hs\alpha(v) 
-\theta(F\bw,v)]$, and (\ref{ape}.c) follows. (Note that 
$\,\zeta(w,w\hh'\hh)\,$ does not change when the pair $\,(w,w\hh'\hh)$ is 
replaced with $\,(w+Fw,w\hh'\nnh+Fw\hh'\hh)$, since $\,Fw\,$ and $\,Fw\hh'$ 
are sections of $\,\mathcal{V}=\mathrm{Ker}\,\zeta$, cf.\ (\ref{mor}) and 
(\ref{dub}.h).)

Next, for sections $\,v\,$ of $\,\mathcal{V}$, we have $\,\ta(v)=\alpha(v)\,$ 
(see the definition of the $\,1$-form $\,\alpha\,$ in Section~\ref{hrzd}). 
This is consistent with (\ref{ape}.b) in view of (\ref{mor}) and 
(\ref{vim}.ii), since $\,\mathcal{V}\,$ is $\,h$-null. Similarly, 
$\,\ta(\tw)=2\hh\tga(\zeta w)$, where $\,\tw=w+Fw\,$ and $\,w\,$ is any 
section of $\,\mathcal{H}$. (By (\ref{dub}.h), $\,\zeta\tw=\zeta w$.) Now 
(\ref{ape}.c) for $\,v=\zeta w$, (\ref{dub}.g) and (\ref{dub}.f) give 
(\ref{ape}.b).

Finally, (\ref{ape}.d) is immediate from Lemma~\ref{gwzww}(ii) applied to both 
$\,\mathcal{H}\,$ and $\,\thor=\mathcal{H}+F$ (where, in the latter case, 
$\,w,w\hh'\nnh,w\hh''$ are replaced by $\,\tw=w+Fw$, 
$\,\tw\hh'\nnh=w\hh'\nnh+Fw\hh'$ and $\,\tw\hh''\nnh=w\hh'\nnh+Fw\hh''$), 
along with (a), (\ref{mor}), (\ref{gxh}) and Remark~\ref{liebr}. As before, 
$\,\zeta(\tw\hh'\nnh,\tw\hh''\hh)=\zeta(w\hh'\nnh,w\hh''\hh)$. 
\end{proof} 
In a basic octuple 
$\,(M,\mathcal{V}\nh,\hh\mathrm{D},h,\alpha,\beta,\theta,\zeta)$, the 
$\,2$-form $\,\zeta\,$ treated as a morphism $\,\tm\to\tm$ (cf.\ 
(\ref{vim})) is a section of $\,\mathcal{F}\hs$ according to (\ref{mor}) and 
(\ref{vim}). By (\ref{trc}.i), $\,\trz=1$. Any function 
$\,f:M\to\bbR\,$ thus gives rise to the section $\,F=f\zeta\,$ of 
$\,\mathcal{F}\nnh$, with $\,\trfz=f$.

Consequently, every section of $\,\mathcal{F}\hs$ can be uniquely written as 
$\,F+f\zeta$, where $\,f:M\to\bbR\,$ and $\,F\,$ is a section of 
$\,\mathcal{F}\hs$ with $\,\trf=0$. Given two horizontal distributions 
$\,\mathcal{H}\,$ and $\,\thor=\mathcal{H}+(F+f\zeta)$, where $\,\trf=0$, 
relations (\ref{ape}.b) and (\ref{ape}.c) now yield 
\begin{equation}\label{tag} 
\begin{array}{rl} 
\mathrm{i)}&\hskip1pt\ta\,=\,\alpha\hs- 
2\hskip1.4ptd_{\hs\zeta(\hs\cdot\hs)}\hh f-10f\beta+h(F\bw,\,\cdot\,)\hs, 
\\ 
\mathrm{ii)}&\hskip1pt\tga(v)\, 
=\,\gamma(v)-\hs d_vf-2f\alpha(v)+\theta(F\bw,v)/\hh2 
\end{array} 
\end{equation} 
for any section $\,v\,$ of $\,\mathcal{V}\nh$. In fact, (\ref{brw}) and 
(\ref{trc}.i) give 
$\,[F^*\nnh\alpha](w)=\alpha(Fw)=h(Fw,\bw)=h(F\bw,w)\,$ for any vector field 
$\,w$, if $\,\trf=0$. Thus, 
\begin{equation}\label{faa} 
F^*\nnh\alpha=h(F\bw,\,\cdot\,)\hskip8pt\mathrm{for\ sections}\hskip5ptF 
\hskip5pt\mathrm{of}\hskip7pt\mathcal{F}\hskip6pt\mathrm{with}\hskip6pt\trf=0 
\hs. 
\end{equation} 
On the other hand, by (\ref{bwz}.i), $\,\zeta^*\nnh\alpha=2\beta\,$ and 
$\,h(\zeta\bw,\,\cdot\,)=-2\beta$, while (\ref{dub}.f) and (\ref{brw}) imply 
that $\,\theta(\zeta\bw,v)/\hh2=-\hs h(v,\bw)=-\hs\alpha(v)$.

\section{The first three conditions in Theorem~\ref{crvtc}}\label{cabl} 
\setcounter{equation}{0} 
Let us fix a basic octuple 
$\,(M,\mathcal{V}\nh,\hh\mathrm{D},h,\alpha,\beta,\theta,\zeta)\,$ and a real 
constant $\,K$.

Given a horizontal distribution $\,\mathcal{H}\,$ for 
$\,(M,\mathcal{V}\nh,\hh\mathrm{D},h,\alpha,\beta,\theta,\zeta)$, we denote by 
$\,g,\nabla\nnh,R,\alpha,\theta$ and $\,\bw\,$ the corresponding objects 
described in Section~\ref{hrzd}. Setting 
$\,\varXi(u,v,w)=R(w,u)\hs v\hs-\hs Kh(v,w)\hs u\,$ for sections $\,u,v\,$ of 
$\,\mathcal{V}\,$ and a vector field $\,w$, we obtain a section $\,\varXi\,$ 
of the vector bundle 
$\,\hs\mathrm{Hom}\hs(\mathcal{V}\otimes\mathcal{V}\otimes\tm,\tm)\,$ over 
$\,M$. We also define sections $\,B\,$ of $\,[\tam]^{\wedge2}$ and $\,\vt\,$ 
of $\,\hs\mathcal{H}^*\nnh$, by 
$\,B=\hs d\beta+2\beta\wedge\alpha+K\zeta/\hh2\,$ and 
$\,\vt(w)\hh\theta(u,v)=[\nabla_{\!w}\theta](u,v) 
+2\hh\alpha(w)\hh\theta(u,v)$, for sections $\,w\,$ of $\,\mathcal{H}\,$ and 
$\,u,v\,$ of $\,\mathcal{V}\nh$. Obviously, $\,\varXi,B\,$ and $\,\vt\,$ 
depend on $\,\mathcal{H}$, and $\,\vt\,$ is well defined, since 
$\,\mathcal{V}^{\wedge2}$ is trivialized by $\,\theta$.

Conditions (a), (b) and (c) in Theorem~\ref{crvtc} amount, respectively, to 
$\,\varXi=0$, $\,B=0$ and $\,\vt=0$. (See Remark~\ref{hrstr} below.) The 
simultaneous vanishing of $\,\varXi,B\,$ and $\,\vt\,$ is a special property 
of $\,\mathcal{H}$. To determine which choices of $\,\mathcal{H}\,$ have this 
property, we first describe the transformations that $\,\varXi,B\,$ and 
$\,\vt\,$ undergo when $\,\mathcal{H}\,$ is replaced by another horizontal 
distribution $\,\thor$. As pointed out in Section~\ref{dfhd}, $\,\thor\,$ is 
always the result of adding to $\,\mathcal{H}\,$ a section of $\,\mathcal{F}$. 
Writing an arbitrary section of $\,\mathcal{F}\hs$ uniquely as $\,F+f\zeta\,$ 
with $\,\trf=0\,$ and $\,f:M\to\bbR$ (cf.\ the end of Section~\ref{dfhd}), and 
denoting by $\,\tX,\tB\,$ and $\,\tvt\,$ the analogues of $\,\varXi,B\,$ and 
$\,\vt\,$ for the new horizontal distribution 
$\,\thor=\mathcal{H}+(F+f\zeta)$, we have, as shown in the next section, 
\begin{equation}\label{tfr} 
\begin{array}{rl} 
\mathrm{a)}\hskip-2pt&\tX=\varXi\nh-\oy F\hs,\\ 
\mathrm{b)}\hskip-2pt&\tB=B-2\beta\wedge\hs[\hs2\hskip1.4ptd_{\hs\zeta(\hs 
\cdot\hs)}\hh f-h(F\bw,\,\cdot\,)\hh]\hs,\\ 
\mathrm{c)}\hskip-2pt&\tvt(\tw)=\hskip2pt\vt(w)+ 
\mathrm{div}^{\mathcal{V}}\nh(Fw)-4\hs d_{\hs\zeta w}f-24f\beta(w) 
+3h(F\bw,w)\hs. 
\end{array} 
\end{equation} 
Here $\,d_{\hs\zeta(\hs\cdot\hs)}\hh f\,$ is defined as in 
Lemma~\ref{hpehf}, $\,\oy F\,$ is given by 
\begin{equation}\label{pfu} 
(\oy F)(u,v,w)\,=\,\hs\mathrm{D}_u\mathrm{D}_v(Fw) 
-\beta(w)\hs\theta(F\bw,u)\hh v 
+\hs\mathrm{D}_u[\hh\theta(F\bw,v)\hs\zeta w]/2 
\end{equation} 
for any $\,\mathcal{V}$-pro\-ject\-a\-ble vector field $\,w\,$ and 
$\,\mathcal{V}$-par\-al\-lel sections $\,u,v\,$ of $\,\mathcal{V}$, while, in 
(\ref{tfr}.c), $\,w$ stands for an arbitrary section of 
$\,\mathcal{H}\,$ and $\,\tw=w+(F+f\zeta)w$. (Thus, $\,\tw\,$ is a section of 
$\,\thor$.) Finally, the $\,\mathcal{V}\nnh$-di\-ver\-gence 
$\,\mathrm{div}^{\mathcal{V}}\nh u:M\to\bbR\,$ of any section $\,u\,$ of 
$\,\mathcal{V}\,$ is the (pointwise) trace of the bundle morphism 
$\,\mathrm{D}\hh u:\mathcal{V}\to\mathcal{V}\,$ sending each section $\,v\,$ 
of $\,\mathcal{V}\,$ to $\,\mathrm{D}_vu$, cf.\ Section~\ref{pmaf}. 
\begin{remark}\label{hrstr}For 
$\,(M,\mathcal{V}\nh,\hh\mathrm{D},h,\alpha,\beta,\theta,\zeta)$, 
$\,\mathcal{H}\,$ and $\,K\,$ as above, the condition $\,\varXi=0\,$ is 
equivalent to vanishing of $\,\varXi(u,v,w)\,$ whenever $\,u,v\,$ are sections 
of $\,\mathcal{V}\nh$, while $\,w$, rather than being an arbitrary vector 
field on $\,M$, is assumed to be a section of $\,\mathcal{H}$. In fact, as 
$\,\mathcal{V}\,$ is $\,h$-null, Lemma~\ref{gwzww}(i) gives 
$\,\varXi(u,v,w)=0\,$ if $\,u,v,w\,$ are sections of $\,\mathcal{V}\nh$. 
\end{remark} 
\begin{example}\label{twopl}For the basic octuple 
$\,(M,\mathcal{V}\nh,\hh\mathrm{D},h,\alpha,\beta,\theta,\zeta)\,$ 
associated, as in Section~\ref{afsy}, with a fixed 
tw\hbox{o\hh-\nh}\hskip0ptplane system $\,(\bs,\xi,\ts,\plane,c,\varOmega)$, 
let $\,\mathcal{H}\,$ be the horizontal distribution appearing in 
Lemma~\ref{gflat}, and let $\,K\,$ be a real constant. The objects 
$\,\varXi,B\,$ and $\,\vt\,$ then are given by 
\begin{equation}\label{xuw} 
\varXi(u,v,w)=-\hs K\varOmega(\yw,v)\hs u\hs,\hskip17pt 
B=K\nh\phi^{-1}\xi\wedge\ts\hs,\hskip17pt 
\vt=0\hs, 
\end{equation} 
where $\,\yw=\xi(w)\hh\rd+\ts(w)\hh c\,$ and $\rd\,$ is the radial vector 
field on $\,\plane$, while (\ref{pfu}) takes the form 
\begin{equation}\label{pfp} 
\begin{array}{rl} 
\mathrm{i)}\hskip0pt&(\oy F)(u,v,w)\,=\,\hs\mathrm{D}_u(H_vw)\, 
-\,\xi(w)\hs\varOmega(F\bw,u)\hh v\hs,\hskip18pt\mathrm{with}\\ 
\mathrm{ii)}\hskip0pt&H_vw\,=\,\mathrm{D}_v(Fw)\, 
-\,\varOmega(F\bw,v)\hh \yw\hh. 
\end{array} 
\end{equation} 
To justify (\ref{xuw}) and (\ref{pfp}), first note that 
$\,\varXi(u,v,w)=-\hs Kh(v,w)\hs u=-\hs K\varOmega(\yw,v)\hs u\,$ by 
(\ref{hvv}). Next, as an immediate consequence of (\ref{peo}), if $\,u\,$ is a 
section of $\,\mathcal{V}\nh$, 
\begin{equation}\label{duf} 
\mathrm{i)}\hskip6pt\phi=\varOmega(\rd,c)\hs,\hskip14pt 
\mathrm{ii)}\hskip6ptd_u\phi=\varOmega(u,c)\hs,\hskip14pt 
\mathrm{iii)}\hskip6ptd_c\phi=0\hs. 
\end{equation} 
On the other hand, (\ref{hvv}), (\ref{duf}.i) and (\ref{bwa}.ii) give 
\begin{equation}\label{hzw} 
h(c,w)=\xi(w)\hs\phi\hs,\hskip12pt 
h(\rd,w)=-\hs\ts(w)\hs\phi\hs,\hskip12pt 
h(\yw,w\hh'\hh)=-\hs\phi\hs(\xi\wedge\ts)(w,w\hh'\hh) 
\end{equation} 
for sections $\,w\,$ of $\,\mathcal{H}$. Furthermore, for such $\,w$, 
\begin{equation}\label{bef} 
\mathrm{a)}\hskip7pt\beta=\phi^{-2}\xi\hs, 
\hskip13pt\mathrm{b)}\hskip7pt\theta=\phi^2\varOmega\hs, 
\hskip13pt\mathrm{c)}\hskip7pt\zeta=2\hs\phi^{-1}\xi\wedge\ts\hs, 
\hskip13pt\mathrm{d)}\hskip7pt\zeta w=-2\phi^{-2}\yw\hs, 
\end{equation} 
where the first three equalities are obvious from (\ref{peo}), 
(\ref{aeh}), and the definition of $\,\theta\,$ in Section~\ref{hrzd}, 
while the last one is easily verified by taking the $\,g$-in\-ner products 
of both sides with any section $\,w\hh'$ of $\,\mathcal{H}$, and using 
(\ref{bef}.c) along with the last formula in (\ref{hzw}).

Also, $\,\alpha=-\hskip1.7ptd\hskip1.2pt\log\hs\phi$, as both sides agree on 
$\,\mathcal{V}\hs$ (by (\ref{peo}) and (\ref{aeh})), and vanish on 
$\,\mathcal{H}\,$ (due to the definition of $\,\alpha\,$ in 
Section~\ref{hrzd}, where $\,\gamma=0$, cf.\ Lemma~\ref{gflat}). Since 
$\,d\hskip1pt\xi=0$, (\ref{bef}.a) thus yields 
$\,d\beta+2\beta\wedge\alpha=0$, and so 
$\,B=K\zeta/\hh2=K\nh\phi^{-1}\xi\wedge\ts\,$ (see (\ref{bef}.c)), as 
required in (\ref{xuw}). The relation $\,\gamma=0\,$ in Lemma~\ref{gflat}, 
combined with (\ref{nwt}), implies in turn the last equality in (\ref{xuw}). 
Now (\ref{pfu}) and (\ref{bef}) give (\ref{pfp}).

In addition, for $\,\bw\,$ and $\,\bu\,$ as in (\ref{brw}) and 
Remark~\ref{aetbu} we have, in this case, 
\begin{equation}\label{xwb} 
\mathrm{a)}\hskip7pt\xi(\bw)=0\hs,\hskip13pt\mathrm{b)}\hskip7pt\ts(\bw) 
=\phi^{-1},\hskip13pt\mathrm{c)}\hskip7pt\zeta\bw=-2\phi^{-3}c\hs,\hskip13pt 
\mathrm{d)}\hskip7pt\bu=\phi^{-3}\nh c\hs. 
\end{equation} 
In fact, $\,\xi(\bw)=0\,$ by (\ref{bef}.a) and (\ref{bwz}.ii), while 
(\ref{trc}.ii) with $\,\varPhi=\ts$, (\ref{bef}.a), (\ref{bef}.c) and 
(\ref{bwa}.ii) yield $\,\ts(\bw)=\phi^{-1}\nnh$, and, as 
$\,\yw=\xi(w)\hh\rd+\ts(w)\hh c$, the third equality is immediate from the 
first two and (\ref{bef}.d). Finally, $\,h(\phi^{-3}\nh c,\,\cdot\,)=\beta\,$ 
by (\ref{bef}.a), the first formula in (\ref{hzw}), and (\ref{hvv}.a), so that 
(\ref{xwb}.d) follows. 
\end{example}

\section{Proof of (\ref{tfr})}\label{ptwe} 
\setcounter{equation}{0} 
Equality (\ref{tfr}.b) is obvious from (\ref{tag}.i), since 
$\,\tB-B=2\beta\wedge(\ta-\alpha)\,$ due to the fact that $\,\beta$ and 
$\,\zeta\,$ do not depend on $\,\mathcal{H}$.

We now establish (\ref{tfr}.a), assuming that $\,u,v\,$ are 
$\,\mathcal{V}$-par\-al\-lel sections of $\,\mathcal{V}\,$ (cf.\ (\ref{trv})), 
while $\,w$, in addition to being $\,\mathcal{V}$-pro\-ject\-a\-ble, is a 
section of $\,\mathcal{V}\,$ or a section of $\,\mathcal{H}$. In the former 
case, both sides equal $\,0$. Namely, Remark~\ref{hrstr} then shows that 
$\,\varXi(u,v,w)=0\,$ for any choice of a horizontal distribution, including 
$\,\tX(u,v,w)=0\,$ for $\,\thor$, while $\,(\oy F)(u,v,w)=0\,$ by (\ref{pfu}) 
since, for sections $\,w\,$ of $\,\mathcal{V}\nh$, (\ref{mor}) and 
(\ref{vim}.ii) yield $\,Fw=\zeta w=0\,$ and $\,\beta(w)=0$.

In the latter case, where $\,w\,$ is a section of $\,\mathcal{H}$, the 
relation $\,\thor=\mathcal{H}+(F+f\zeta)\,$ implies that 
$\,\tw=w+(F+f\zeta)w\,$ is a $\,\mathcal{V}$-pro\-ject\-a\-ble section of 
$\,\thor$, and, according to the preceding paragraph, 
$\,\tX(u,v,w)=\tX(u,v,\tw)$. Let us now evaluate $\,\varXi(u,v,w)\,$ (or, 
$\,\tX(u,v,\tw)$) with the aid of the equality 
$\,R(w,u)\hh v=\mathrm{D}_u[w,v]+2\beta(w)\hs\gamma(u)\hs v 
-\mathrm{D}_u[\gamma(v)\hs\zeta w]\,$ in Lemma~\ref{gwzww}(iii) (or, 
respectively, its analogue for $\,\thor$). Since $\,\beta(\tw)=\beta(w)\,$ and 
$\,\zeta\tw=\zeta w$, cf.\ (\ref{vim}.ii), we may thus express 
$\,\varXi(u,v,w)-\tX(u,v,w)=R(w,u)\hh v-\tR(\tw,u)\hh v\,$ as a sum of some 
terms containing $\,F\,$ and some terms involving $\,f$. By (\ref{tag}.ii), 
the former terms add up to $\,(\oy F)(u,v,w)$. (Since $\,v\,$ is 
$\,\mathcal{V}$-par\-al\-lel, $\,[w,v]-[\tw,v]=\mathrm{D}_v[(F+f\zeta)w]$.) On 
the other hand, the sum $\,S\,$ of the latter terms is zero. Namely, 
(\ref{tag}.ii) gives $\,S=\mathrm{D}_u\mathrm{D}_v(f\zeta w) 
+2\beta(w)\hs[\hh d_uf+2f\alpha(u)\hh]\hh v-(d_ud_vf)\hs\zeta w 
-2\hh(d_uf)\hh\alpha(v)\hs\zeta w 
-2f\{\hskip-.4ptd_u\hs[\hh\alpha(v)]\}\hs\zeta w 
-[\hh d_vf+2f\alpha(v)\hh]\hs\mathrm{D}_u(\zeta w)$. (The last four terms 
arise when $\,\mathrm{D}_u\,$ is applied to 
$\,-\hs[\hh d_vf+2f\alpha(v)\hh]\hs\zeta w$.) However, by the Leib\-niz rule, 
$\,\mathrm{D}_u\mathrm{D}_v(f\zeta w) 
=(d_ud_vf)\hs\zeta w+(d_uf)\hs\mathrm{D}_v(\zeta w) 
+(d_vf)\hs\mathrm{D}_u(\zeta w)+f\mathrm{D}_u\mathrm{D}_v(\zeta w)$, while, 
according to Remark~\ref{duzwe} and (\ref{dub}.d), 
$\,\mathrm{D}_u\mathrm{D}_v(\zeta w) 
=2\hh\mathrm{D}_u[\hh\alpha(v)\hs\zeta w-\beta(w)\hs v\hh] 
=2\hs\{\hskip-.4ptd_u\hs[\hh\alpha(v)]\}\hs\zeta w 
+2\alpha(v)\hs\mathrm{D}_u(\zeta w) 
-4\alpha(u)\hs\beta(w)\hs v$. The resulting cancellations show that 
$\,S\,$ equals 
$\,[\hh\mathrm{D}_v(\zeta w)-2\hh\alpha(v)\hs\zeta w\nh 
+\nh2\beta(w)\hs v\hh]d_uf$, and so $\,S=0\,$ in view of Remark~\ref{duzwe}. 
This yields (\ref{tfr}.a).

To prove (\ref{tfr}.c), let us fix $\,\mathcal{V}$-par\-al\-lel sections 
$\,u,v\,$ of $\,\mathcal{V}\,$ and a $\,\mathcal{V}$-pro\-ject\-a\-ble section 
$\,w$ of $\,\mathcal{H}$. Since $\,\nabla\,$ is tor\-sion\-free, the 
Leib\-niz rule gives $\,[\vt(w)-2\hh\alpha(w)]\hs\theta(u,v) 
=[\nabla_{\!w}\theta](u,v)=d_w[\hs\theta(u,v)]-\theta(\nabla_{\!u}w,v) 
-\theta(u,\nabla_{\!v}w)-\theta([w,u],v)-\theta(u,[w,v])$. However, 
$\,\mathcal{H}=\mathrm{Ker}\,\theta$, cf.\ Section~\ref{hrzd}, so that 
$\,\nabla_{\!u}w\,$ and $\,\nabla_{\!v}w\,$ can be replaced here with their 
$\,\mathcal{V}\,$ components, equal, by (\ref{ned}.b), to 
$\,-\hs\gamma(u)\hskip1pt\zeta w\,$ and $\,-\hs\gamma(v)\hskip1pt\zeta w$. 
Consequently, $\,-\hs\theta(\nabla_{\!u}w,v)-\theta(u,\nabla_{\!v}w) 
=\gamma(u)\hs\theta(\zeta w,v)+\gamma(v)\hs\theta(u,\zeta w)\,$ which, by 
Remark~\ref{cclsm}, equals 
$\,\gamma(\zeta w)\hs\theta(u,v)=\alpha(w)\hs\theta(u,v)/2$, as 
$\,\alpha(w)=2\gamma(\zeta w)\,$ (see Section~\ref{hrzd}). Thus, 
$\,\vt(w)\hs\theta(u,v)=d_w[\hs\theta(u,v)]+\theta([u,w],v)+\theta(u,[v,w]) 
+5\hh\alpha(w)\hs\theta(u,v)/2$.

Suppose now that $\,\thor=\mathcal{H}+F\nnh$, where $\,F\,$ is an 
{\it arbitrary\/} section of $\,\mathcal{F}\nnh$, not necessarily one with 
$\,\trf=0$. For $\,\tw=w+Fw$, the preceding equality, applied to both 
$\,\mathcal{H}\,$ and $\,\thor$, yields 
$\,[\tvt(\tw)-\vt(w)]\hs\theta(u,v)=d_{Fw}[\hs\theta(u,v)] 
+\theta([u,Fw],v)+\theta(u,[v,Fw]) 
+5\hs[\ta(\tw)-\alpha(w)]\hs\theta(u,v)/2$. As $\,u\,$ and $\,v\,$ are 
$\,\mathcal{V}$-par\-al\-lel, 
$\,[u,Fw]=\mathrm{D}_u(Fw)\,$ and $\,[v,Fw]=\mathrm{D}_v(Fw)$, so that 
Remark~\ref{trvol} gives $\,\theta([u,Fw],v)+\theta(u,[v,Fw]) 
=[\hs\mathrm{div}^{\mathcal{V}}\nh(Fw)]\hs\theta(u,v)$. Hence, by 
(\ref{dub}.e), $\,\tvt(\tw)-\vt(w)=-2\hh\alpha(Fw) 
+\mathrm{div}^{\mathcal{V}}\nh(Fw)+5\hs[\ta(\tw)-\alpha(w)]/2$. Substituting 
for $\,F\nnh$, in this last equality, the sum $\,F+f\zeta\,$ with $\,\trf=0$, 
we obtain (\ref{tfr}.c), as $\,\mathrm{div}^{\mathcal{V}}\nh(\zeta w)=0\,$ by 
Remark~\ref{duzwe} and (\ref{dub}.g), and so 
$\,\mathrm{div}^{\mathcal{V}}\nh[(F+f\zeta)w] 
=\mathrm{div}^{\mathcal{V}}\nh(Fw)+d_{\hs\zeta w}f$, while (\ref{faa}) and 
(\ref{dub}.g) give $\,\alpha((F+f\zeta)w)=h(F\bw,w)+2f\beta(w)$, and 
(\ref{tag}.i) with $\,\tw=w+(F+f\zeta)w\,$ yields 
$\,\ta(\tw)-\alpha(w)=-2\hskip1.4ptd_{\hs\zeta w}f-10f\beta(w)+h(F\bw,w) 
+\alpha((F+f\zeta)w)$.

\section{Dimension of a solution space}\label{doas} 
\setcounter{equation}{0} 
Let $\,\mathcal{H}\,$ be a horizontal distribution for a basic octuple 
$\,(M,\mathcal{V}\nh,\hh\mathrm{D},h,\alpha,\beta,\theta,\zeta)$.

Formula (\ref{pfu}) makes sense also when $\,F\nnh$, rather than being a 
section of $\,\mathcal{F}\nnh$, is just a section of the restriction of 
$\,\mathcal{F}\,$ to a surface $\,\plane\hh'$ embedded in $\,M\,$ and 
contained in a leaf of $\,\mathcal{V}\nh$. (In fact, (\ref{pfu}) involves only 
covariant derivatives in directions tangent to $\,\mathcal{V}\nh$.) The kernel 
of the operator $\,\oy$, acting on those sections $\,F\,$ of the restriction 
of $\,\mathcal{F}\,$ to $\,\plane\hh'$ which satisfy the additional condition 
$\,\trf=0$, is at most \hbox{eight-}\hskip0ptdi\-men\-sion\-al.

To show this, we first observe that, if $\,F\,$ is a section of 
$\,\mathcal{F}\,$ with $\,\oy F=0\,$ and $\,\trf=0$, then 
\begin{equation}\label{pfz} 
\begin{array}{rl} 
\mathrm{i)}\hskip0pt&h(\mathrm{D}_v(Fw),w\hh'\hh)\, 
-\,h(\mathrm{D}_v(Fw\hh'\hh),w)\, 
=\,\theta(v,F\bw)\hs\zeta(w,w\hh'\hh)/\hh2\hs,\\ 
\mathrm{ii)}\hskip0pt&h(Fw,w\hh'\hh)\,=\,h(Fw\hh'\nnh,w)\,\hs, 
\hskip19pt\mathrm{iii)}\hskip9pt\theta(\mathrm{D}_u(F\bw),v)\, 
=\,\theta(\mathrm{D}_v(F\bw),u) 
\end{array} 
\end{equation} 
for any sections $\,u,v\,$ of $\,\mathcal{V}\,$ and any 
$\,\mathcal{V}$-pro\-ject\-a\-ble local sections $\,w,w\hh'$ of 
$\,\mathcal{H}$.

Namely, (\ref{pfz}.ii) is obvious from (\ref{trc}.i), while (\ref{pfz}.i) 
easily follows if we apply $\,d_v$ to (\ref{pfz}.ii), use (\ref{dub}.a) along 
with the Leib\-niz rule, and then set $\,\varPhi w=\theta(v,Fw)\,$ in 
(\ref{trc}.ii). Finally, (\ref{pfz}.iii) is immediate if one 
skew-sym\-met\-rizes the right-hand side of (\ref{pfu}) in $\,u,v$, assuming 
$\,u,v\,$ to be $\,\mathcal{V}$-par\-al\-lel (so that 
$\,\mathrm{D}_u\mathrm{D}_v=\mathrm{D}_v\mathrm{D}_u\,$ as $\,\mathrm{D}\,$ 
is flat), and then uses (\ref{dub}.e) along with the Leib\-niz rule and 
Remark~\ref{duzwe}.

Note that the assumption $\,\trf=0\,$ alone yields (\ref{pfz}.ii) and 
(\ref{pfz}.i), while Remark~\ref{trvol} allows us to rewrite (\ref{pfz}.iii) 
as $\,\mathrm{div}^{\mathcal{V}}\nh(F\bw)=0$, and (by (\ref{mor}) and 
(\ref{vim}.ii)) $\,w,w\hh'$ in (\ref{pfz} may equivalently be just any 
$\,\mathcal{V}$-pro\-ject\-a\-ble local vector fields.

Condition $\,\oy F=0\,$ implies that $\,F\nnh$, restricted to any 
$\,\mathrm{D}\hs$-ge\-o\-des\-ic contained in $\,\plane\hh'\nnh$, satisfies 
a system of sec\-ond-or\-der linear ordinary differential equations 
solved for the second derivatives. To make sense of $\,F\bw\,$ in this 
context, here and below we choose $\,\phi\,$ as at the end of 
Section~\ref{boct}, thus getting $\,F\bw=\phi^{-1}\nnh F\phi\hh\bw$, where 
$\,\phi\hh\bw\,$ is $\,\mathcal{V}$-pro\-ject\-a\-ble by Lemma~\ref{pctbl}. 
Let us fix a point $\,x\in\plane\hh'\nnh$. Due to uniqueness of solutions, 
$\,F\,$ is completely determined, on $\,\plane\hh'\nnh$, by the pair 
$\,(F,\bz)$ consisting of its value at $\,x$, still denoted by $\,F\nnh$, 
and its $\,\mathcal{V}\nh$-differential $\,\bz\,$ at $\,x$. More precisely, 
$\,F\,$ is a linear operator $\,\mathcal{H}_x\to\mathcal{V}_{\nh x}$, and 
$\,\bz\,$ may be treated as a bi\-lin\-e\-ar mapping 
$\,\mathcal{V}_{\nh x}\times\mathcal{H}_x\to\mathcal{V}_{\nh x}$ with 
$\,\bz(v,w)=\mathrm{D}_v(Fw)\,$ (at $\,x$, for $\,v,w\,$ as in 
(\ref{pfz})). Thus, solutions $\,F\,$ to the equations $\,\oy F=0\,$ and 
$\,\trf=0\,$ on $\,\plane\hh'$ lead, at $\,x$, to pairs $\,(F,\bz)\,$ which 
are vectors in a $\,12$\diml\ space $\,\fb$. Rather than being arbitrary 
vectors in $\,\fb$, such $\,(F,\bz)\,$ are subject to the additional 
constraints stemming from (\ref{pfz}), which state that, for all 
$\,v\in\mathcal{V}_{\nh x}$ and $\,w,w\hh'\hn\in\mathcal{H}_x$, 
\begin{equation}\label{cst} 
\begin{array}{rl} 
\mathrm{i)}\hskip0pt&h(\bz(v,w),w\hh'\hh)\, 
-\,h(\bz(v,w\hh'\hh),w)\, 
=\,\phi^{-1}\hn\theta(v,F\phi\hh\bw)\hs\zeta(w,w\hh'\hh)/\hh2\hs,\\ 
\mathrm{ii)}\hskip0pt&h(Fw,w\hh'\hh)\,=\,h(Fw\hh'\nnh,w)\hs,\\ 
\mathrm{iii)}\hskip0pt&\theta(\bz(u,\phi\hh\bw),v)\, 
-\,\theta(\bz(v,\phi\hh\bw),u)\,=\,-\hs h(F\phi\hh\bw,\bw)\hs\theta(u,v)\hs.\\ 
\end{array} 
\end{equation} 
We obtain (\ref{cst}.iii) from (\ref{pfz}.iii), the relation 
$\,\alpha=-\hskip1.7ptd\hskip1.2pt\log\hs\phi\,$ on $\,\mathcal{V}\,$ 
(cf.\ Section~\ref{boct}), combined with Remark~\ref{cclsm} for the expression 
$\,\alpha(u)\hs\theta(F\phi\hh\bw,v)\,$ (tri\-lin\-e\-ar in 
$\,u,F\phi\hh\bw,v$), and (\ref{faa}). Since $\,\phi\hh\bw\,$ is 
$\,\mathcal{V}$-pro\-ject\-a\-ble (see above), it is useful here to rewrite 
$\,F\bw\,$ as $\,\phi^{-1}\nnh F\phi\hh\bw$.

By assigning to a linear operator 
$\,\varPhi:\mathcal{H}_x\to\mathcal{V}_{\nh x}$ the bi\-lin\-e\-ar form 
on $\,\mathcal{H}_x$ that sends $\,w,w\hh'$ to $\,(\varPhi w,w\hh'\hh)$, 
we obtain an isomorphism between the space of operators 
$\,\mathcal{H}_x\to\mathcal{V}_{\nh x}$ and the space of 
bi\-lin\-e\-ar forms $\,\mathcal{H}_{\nh x}\times\mathcal{H}_x\to\bbR\hs$. 
(See (\ref{vbi}).) We will say that a linear operator 
$\,\varPhi:\mathcal{H}_x\to\mathcal{V}_{\nh x}$ is $\,h${\it-self-ad\-joint\/} 
if it corresponds under this isomorphism to a form which is symmetric, that 
is, if $\,(\varPhi w,w\hh'\hh)\,=\,h(\varPhi w\hh'\nnh,w)\,$ for all 
$\,w,w\hh'\hn\in\mathcal{H}_x$. The space of $\,h$-self-ad\-joint 
operators $\,\mathcal{H}_x\to\mathcal{V}_{\nh x}$ is, obviously, 
\hbox{three\hh-}\hskip0ptdi\-men\-sion\-al.

Denoting by $\,\fb\hs'$ the subspace of $\,\fb\,$ formed by all 
$\,(F,\bz)\in\fb\,$ satisfying (\ref{cst}.ii) and (\ref{cst}.i), 
we have $\,\dim\hs\fb\hs'\nh=9$. In fact, the assignment 
$\,(F,\bz)\mapsto(F,\hat\bz)$, given by 
$\,\hat\bz(v,w) 
=\bz(v,w)-\phi^{-1}\hn\theta(v,F\phi\hh\bw)\hs\zeta w/4$, is an 
isomorphism $\,\fb\to\fb\,$ sending $\,\fb\hs'$ onto the space of 
all $\,(F,\hat\bz)\in\fb\,$ such that $\,F\,$ and 
$\,\hat\bz(v,\,\cdot\,)\,$ are $\,h$-self-ad\-joint for every 
$\,v\in\mathcal{V}_{\nh x}$, while the latter space is 
\hbox{nine\hs-}\hskip0ptdi\-men\-sion\-al (cf.\ the last paragraph).

Finally, the three conditions (\ref{cst}) together define 
an \hbox{eight-}\hskip0ptdi\-men\-sion\-al 
subspace of $\,\fb$. In fact, the subspace in question is 
the kernel of a linear functional on the space $\,\fb\hs'$ with 
$\,\dim\hs\fb\hs'\nh=9$. (Note that (\ref{cst}.iii) amounts to 
a single scalar equation, due to its skew-sym\-me\-try in $\,u,v$.) 
The functional in question is nonzero, since condition (\ref{cst}.iii) for 
$\,(F,\bz)\in\fb\,$ is {\it not\/} a consequence of (\ref{cst}.i) and 
(\ref{cst}.ii). An example $\,(F,\bz)\,$ proving the last claim may be defined 
as follows. We choose $\,F\,$ which is both $\,h$-self-ad\-joint and such that 
$\,h(F\phi\hh\bw,\bw)\ne0$. In other words, $\,\phi\hh\bw\,$ is not null for 
the symmetric bi\-lin\-e\-ar form corresponding to $\,F$. Thus, 
(\ref{cst}.ii) holds. Then we set 
$\,\bz(v,w)=\phi^{-1}\hn\theta(v,F\phi\hh\bw)\hs\zeta w/4$, which clearly 
gives (\ref{cst}.i). However, the left-hand side of (\ref{cst}.iii) equals 
here $\,1/4\,$ times $\,\theta(u,F\phi\hh\bw)\hs\theta(\zeta\bw,v) 
-\theta(v,F\phi\hh\bw)\hs\theta(\zeta\bw,u)$, which coincides with 
$\,\theta(\zeta\bw,F\phi\hh\bw)\hs\theta(u,v)/4\,$ (from Remark~\ref{cclsm} 
applied to the expression $\,\theta(u,F\phi\hh\bw)\hs\theta(\zeta\bw,v)$, 
tri\-lin\-e\-ar in $\,u,\zeta\bw,v$). Now, by (\ref{ztz}.ii), the left-hand 
side of (\ref{cst}.iii) is equal to 
$\,-\hs h(F\phi\hh\bw,\bw)\hs\theta(u,v)/2$, 
and hence different from the 
right-hand side, due to our choice of $\,F$.

\section{Explicit solutions for a tw\hbox{o\hh-\nh}\hskip0ptplane 
system}\label{esol} 
\setcounter{equation}{0} 
Let $\,K\,$ be a fixed real constant. We will now describe the set of all 
horizontal distributions $\,\mathcal{H}$, for any given basic octuple 
$\,(M,\mathcal{V}\nh,\hh\mathrm{D},h,\alpha,\beta,\theta,\zeta)$, which have 
the properties (a), (b) and (c) in Theorem~\ref{crvtc}.

Our discussion is local. Theorem~\ref{unimo} thus allows us to fix a 
tw\hbox{o\hh-\nh}\hskip0ptplane system 
$\,(\bs,\xi,\ts,\plane,c,\varOmega)$ and assume, without loss of 
generality, that 
$\,(M,\mathcal{V}\nh,\hh\mathrm{D},h,\alpha,\beta,\theta,\zeta)\,$ is its 
associated basic octuple. We denote by $\,\mathcal{H}\,$ the horizontal 
distribution appearing in Lemma~\ref{gflat}. Horizontal distributions 
$\,\thor\,$ satisfying (a) -- (c) in Theorem~\ref{crvtc}, written as 
$\,\thor=\mathcal{H}+(F+f\zeta)$, where $\,F$ is a section of 
$\,\mathcal{F}\hs$ with $\,\trf=0\,$ and $\,f:M\to\bbR\,$ (see the end of 
Section~\ref{dfhd}), are characterized, according to Section~\ref{cabl}, by 
simultaneous vanishing of $\,\tX,\tB\,$ and $\,\tvt\,$ in (\ref{tfr}), which 
is a system of three (usually non\-ho\-mo\-ge\-ne\-ous) linear partial 
differential equations with the unknowns $\,F$ and $\,f$. Specifically, 
$\,\tX=0\,$ is a sec\-ond-or\-der equation involving $\,F\,$ only, $\,\tB=0\,$ 
is of first order in $\,f\,$ and of order zero in $\,F\nnh$, while 
$\,\tvt=0\,$ is of first order in both $\,f\,$ and $\,F$. By (\ref{tfr}.a), 
$\,\tX=0\,$ if and only if $\,\oy F=\varXi\,$ (under the assumption that 
$\,\trf=0$).

In all three equations, only derivatives in $\,\hs\plane\,$ directions occur, 
so that we may fix $\,y\in\bs$ and restrict the unknowns $\,f\,$ and 
$\,F\,$ to the subset $\,\{y\}\times\plane_+\hs\approx\,\plane_+$, thus 
treating them as a function $\,f:\plane_+\to\bbR\,$ and, respectively, 
\begin{equation}\label{alo} 
\mathrm{a\ linear\ operator\ }\,w\mapsto Fw\,\mathrm{\ from\ }\,\dbs\, 
\mathrm{\ into\ the\ space\ of\ vector\ fields\ on\ }\,\plane_+, 
\end{equation} 
vector fields being identified with mappings $\,\plane_+\to\plane$. Here 
$\,\dbs\,$ is the translation vector space of $\,\bs$, canonically isomorphic 
both to $\,\tyb\,$ and to the fibre of $\,\mathcal{H}\,$ at any point of 
$\,\{y\}\times\plane_+$. The three equations, phrased in terms of such 
identifications, with fixed $\,y$, are solved below. As before, $\,\rd\,$ 
denotes the restriction to $\,\plane_+$ of the radial vector field, 
$\,\yw=\xi(w)\hh\rd+\ts(w)\hh c$ and 
$\,\varXi(u,v,w)=-\hs K\varOmega(\yw,v)\hs u$, cf.\ (\ref{xuw}), while 
$\,\phi:\plane_+\to(0,\infty)\,$ is given by $\,\phi=\varOmega(\rd,c)$. The 
symbol $\,F\bw\,$ stands for $\,\phi^{-1}F(\phi\hh\bw)$. (By (\ref{xwb}.a) -- 
(\ref{xwb}.b), or Lemma~\ref{pctbl}, $\,\phi\hh\bw\,$ is a 
$\,\mathcal{V}$-pro\-ject\-a\-ble section of $\,\mathcal{H}$, and may be 
identified with a constant vector field on $\,\bs$, so that this convention 
about the meaning of $\,F\bw\,$ agrees with our previous usage, such as in 
(\ref{pfp}).) 
\begin{lemma}\label{eight}Among all operators\/ {\rm(\ref{alo})}, those with 
$\,\oy F=\varXi\,$ and\/ $\,\trf=0\,$ form an 
\hbox{eight-}\hskip0ptdi\-men\-sion\-al af\-fine space\/ $\,\as\,$ 
containing\/ $\,F^K$ given by\/ $\,F^Kw=K\nh\phi\hs\ts(w)\hs\rd/2$. The 
translation vector space $\,\das\,$ of\/ $\,\as\,$ is the direct sum of three 
subspaces, of dimensions $\,2,\hs3\,$ and\/ $\,3$, consisting, respectively, 
of the operators\/ $\,F^\vf\nnh,F^\lambda\nnh,F^\mu$ sending $\,w\,$ to 
\begin{enumerate} 
  \def\theenumi{{\rm\roman{enumi}}} 
\item[{\rm(a)}] $2\hs\xi(w)\hs\vf-(d_\vf\phi)\hs\phi^{-1}\yw$, for any 
constant vector field\/ $\,\vf\,$ on\/ $\,\plane$, 
\item[{\rm(b)}] $\phi^{-1}[\hs\xi(w)\hs\lambda(\rd,\rd)\hh c 
+\ts(w)\hs\lambda(c,c)\hh\rd]$, for any constant 
symmetric $\,2$-ten\-sor field\/ $\,\lambda\,$ on\/ $\,\plane$, 
\item[{\rm(c)}] $\phi^{-1}[\hh\mu(\rd,\rd)\hh\yw-2\hh\mu(\yw,\rd)\hh\rd]$, 
for any constant symmetric $\,2$-ten\-sor field\/ $\,\mu\,$ on\/ $\,\plane$. 
\end{enumerate} 
If\/ $\,F=F^K\nh+F^\vf\nh+F^\lambda\nh+F^\mu\nnh$, then\/ 
$\,\phi^2F\bw=\ex\rd+\lx c$, where\/ 
$\,\ex=\lambda(c,c)-2\hh\mu(c,\rd)+K\nh\phi^2\nnh/\hh2$ and\/ 
$\,\lx=\mu(\rd,\rd)-\varOmega(\vf,c)$. Finally, for any\/ 
$\,w,w\hh'\nh\in\dbs$, 
\begin{equation}\label{hfk} 
\begin{array}{l} 
h(F^K\nh w,w\hh'\hh)=-\hs K\nh\phi^2\ts(w)\hs\ts(w\hh'\hh)/2\hs,\\ 
h(F^\vf\nh w,w\hh'\hh)=2\hh\varOmega(\rd,\vf)\hs\xi(w)\hs\xi(w\hh'\hh) 
-\varOmega(\vf,c)\hs[\hs\xi(w)\hh\ts(w\hh'\hh)+\xi(w\hh'\hh)\hh\ts(w)]\hs,\\ 
h(F^\lambda\nh w,w\hh'\hh)=\lambda(\rd,\rd)\hs\xi(w)\hs\xi(w\hh'\hh) 
-\lambda(c,c)\hh\ts(w)\hh\ts(w\hh'\hh)\hs,\\ 
h(F^\mu\nh w,w\hh'\hh)=\mu(\rd,\rd)\hs[\hs\xi(w)\hh\ts(w\hh'\hh) 
+\xi(w\hh'\hh)\hh\ts(w)]+2\hh\mu(c,\rd)\hh\ts(w)\hh\ts(w\hh'\hh)\hs. 
\end{array} 
\end{equation} 
\end{lemma} 
\begin{proof}First, (\ref{hzw}), (\ref{bwa}.ii), (\ref{hvv}) and 
(\ref{duf}.ii) yield (\ref{hfk}). Hence 
$\,\trfk=[\hskip-1.7pt[\hskip-.4ptF^\vf\hskip0pt]\hskip-1.7pt] 
=[\hskip-1.7pt[\hskip-.4ptF^\lambda\hskip0pt]\hskip-1.7pt] 
=[\hskip-1.7pt[\hskip-.4ptF^\mu\hskip0pt]\hskip-1.7pt]=0\,$ due to 
(\ref{trc}.i) and symmetry in $\,w,w\hh'$ of the right-hand sides in 
(\ref{hfk}).

Since $\,\xi(\bw)=0$, $\,\ts(\bw)=\phi^{-1}$ and 
$\,Y\hskip-2.7pt_{\bw}=\hs\phi^{-1}\nh c\,$ (see (\ref{xwb}) and 
(\ref{hvv}.c)), we have, by (\ref{duf}.ii), 
\begin{equation}\label{fkb} 
\begin{array}{l} 
F^K\nh\bw=K\nh\rd/2\hs,\hskip7ptF^\vf\bw=\varOmega(c,\vf)\hs\phi^{-2}\nh c, 
\hskip7ptF^\lambda\bw=\lambda(c,c)\hs\phi^{-2}\nh\rd\hs,\\ 
F^\mu\bw=\phi^{-2}[\hh\mu(\rd,\rd)\hs c-2\mu(c,\rd)\hh\rd]\hs. 
\end{array} 
\end{equation} 
We denote by $\hs H_v^K\nnh w,H_v^\vf w,H_v^\lambda w\,$ or $\,H_v^\mu w\,$ the 
expression (\ref{pfp}.ii) with $\,F\hn\,$ replaced by 
$\hs F^K\nnh,F^\vf\nnh,F^\lambda$ or $\,F^\mu\nnh$. Let $\,\as\,$ be the 
af\-fine space of all operators (\ref{alo}) with $\,\oy F=\varXi\,$ and 
$\,\trf=0$.

Using (\ref{duf}) and noting that 
$\,\varOmega(\rd,c)\hs v-\varOmega(\rd,v)\hs c=\varOmega(v,c)\hs\rd\,$ 
(see Remark~\ref{cclsm}), we obtain 
$\,H_v^K\nnh w=K\hh[\hh\varOmega(v,c)\hs\ts(w) 
-\varOmega(\rd,v)\hs\xi(w)/2]\hh\rd$. Now (\ref{pfp}.i) and (\ref{fkb}) 
combined with the equality 
$\,\varOmega(u,v)\hs\rd+\varOmega(\rd,u)\hs v=\varOmega(\rd,v)\hs u\,$ 
(immediate from Remark~\ref{cclsm}) yield $\,\oy F^K\nnh=\varXi$, with 
$\,\varXi(u,v,w)=-\hs K\varOmega(\yw,v)\hs u$. As $\,\trfk=0$, it follows 
that $\,F^K\nh\in\as$.

Furthermore, (\ref{duf}) gives $\,H_v^\vf w 
=\varOmega(c,\vf)\hs\phi^{-1}\hn\xi(w)\hh v$, and so 
$\,\oy F^\vf\nh=0\,$ by (\ref{pfp}.i) and (\ref{fkb}).

Next, $\,H_v^\lambda w=2\phi^{-1}\hn\xi(w)\hs\lambda(v,\rd)\hs c 
-\phi^{-2}\hn\xi(w)\hs[\hh\lambda(c,c)\hs\varOmega(\rd,v)\hs\rd 
+\varOmega(v,c)\hs\lambda(\rd,\rd)\hs c\hh]\,$ by (\ref{duf}), as 
$\,\varOmega(\rd,c)\hs v-\varOmega(\rd,v)\hs c=\varOmega(v,c)\hs\rd\,$ 
(see above). The relation $\,\oy F^\lambda\nh=0\,$ is now easily 
verified using (\ref{pfp}.i) and (\ref{fkb}) along with the equalities 
$\,\varOmega(u,v)\hs\rd+\varOmega(\rd,u)\hs v=\varOmega(\rd,v)\hs u$, 
$\,(d_u\phi)\hh\rd-\phi\hs u=\varOmega(u,c)\hh\rd+\varOmega(c,\rd)\hh u 
=\varOmega(u,\rd)\hs c\,$ (and hence 
$\,\varOmega(u,c)\hs\lambda(\rd,\rd)-\phi\hs\lambda(u,\rd) 
=\varOmega(u,\rd)\hs\lambda(c,\rd)$), 
$\,\varOmega(u,c)\hs\lambda(\rd,v)-\phi\hh\lambda(u,v) 
=\varOmega(u,c)\hs\lambda(\rd,v)+\varOmega(c,\rd)\hs\lambda(u,v) 
=\lambda(c,v)\hs\varOmega(u,\rd)$ and 
$\,\lambda(c,c)\hs\varOmega(\rd,v)+\lambda(c,v)\hs\varOmega(c,\rd) 
+\lambda(c,c)\hs\varOmega(\rd,v)=0$, due to Remark~\ref{cclsm} and 
(\ref{duf}.i).

Similarly, $\,H_v^\mu w= 
2\hh\phi^{-1}\hn\ts(w)\hs[\hh\mu(v,\rd)\hs c-\mu(v,c)\hs\rd\hh] 
-\phi^{-1}\mu(\rd,\rd)\hs\xi(w)\hs v$, since 
$\,\varOmega(\rd,v)\hs c+\varOmega(v,c)\hs\rd=\varOmega(\rd,c)\hs v$, 
which also gives 
$\,\varOmega(\rd,v)\hs\mu(c,\rd)+\varOmega(v,c)\hs\mu(\rd,\rd) 
=\varOmega(\rd,c)\hs\mu(v,\rd)=\phi\hs\mu(v,\rd)$. Therefore, 
$\,\oy F^\mu\nh=0$, in view of the relation 
$\,\varOmega(u,c)\hh\rd-\phi\hs u=\varOmega(u,\rd)\hs c\,$ (see 
above) and two further equalities, which are its immediate consequences: 
$\,\varOmega(\rd,u)\hs\mu(c,\rd)+\varOmega(u,c)\hs\mu(\rd,\rd) 
=\phi\hs\mu(c,\rd)\,$ and $\,\varOmega(u,c)\hs\mu(v,\rd) 
=\phi\hs\mu(u,v)+\mu(c,v)\hs\varOmega(u,\rd)$.

As established above, $\,F^\vf\nnh,F^\lambda\nnh,F^\mu\nh\in\das$, 
where $\,\das\,$ is the translation vector space of $\,\as$. The 
operators $\,F^\vf\nnh,F^\lambda\nnh,F^\mu$ together span a vector space 
of dimension $\,8$. In fact, assuming that 
$\,[F^\vf\nh+F^\lambda\nh+F^\mu\hh]\hs\phi\,$ vanishes identically, we will 
show that $\,\vf=0\,$ and $\,\lambda=\mu=0$. Namely, 
$\,[F^\vf\nh+F^\lambda\nh+F^\mu\hh]\hs\phi$, as a function on $\,\plane_+$ 
valued in the space of operators $\,\dbs\to\plane$, is a polynomial of 
degree at most $\,3$, with some homogeneous components 
$\,H\cst\nnh,H\lin\nnh,H\qdr\nnh,H\cub$ of degrees $\,0,1,2,3$. 
Clearly, $\,H\cst\hn w\,$ equals $\,-\hs d_\vf\phi\,$ times the constant term 
in $\,\yw=\xi(w)\hh\rd+\ts(w)\hh c$, that is, 
$\,0=H\cst\hn w=(d_\vf\phi)\hs\ts(w)\hh c$. (Our assumption is that 
$\,H\cst\nh=H\lin\nh=H\qdr\nh=H\cub\nh=0$.) Hence $\,d_\vf\phi=0$. 
Next, $\,0=H\lin\hn w=2\hs\xi(w)\hs\phi\hs\vf+\ts(w)\hs\lambda(c,c)\hs\rd$, 
and so $\,\vf=0\,$ due to linear independence of $\,\xi\,$ and $\,\ts$. 
Similarly, $\,H\cub\hn w\,$ is the cubic term in $\,\phi\hh F^\mu w$, and, 
therefore, $\,0=H\cub\hn w=-\hs\xi(w)\hs\mu(\rd,\rd)\hh\rd$, so that 
$\,\mu=0$. Thus, $\,F^\lambda\nh=0\,$ and hence $\,\lambda=0$, as $\,\xi\,$ 
and $\,\ts$ are linearly independent.

Consequently, $\,\dim\das\ge8$, while $\,\dim\das\le8\,$ according to 
Section~\ref{doas}, which shows that $\,\das\,$ is both 
\hbox{eight-}\hskip0ptdi\-men\-sion\-al and spanned by all 
$\,F^\vf\nnh,F^\lambda\nnh,F^\mu\nnh$, completing the proof. 
\end{proof} 
\begin{lemma}\label{dcfeq}For\/ $\,F\,$ with $\,\oy F=\varXi\,$ and\/ 
$\,\trf=0$, written uniquely as\/ 
$\,F=F^K\nh+F^\vf\nh+F^\lambda\nh+F^\mu\nnh$, cf.\ Lemma\/~{\rm\ref{eight}}, 
and a function\/ $\,f:\plane_+\to\bbR\hs$, the section\/ $\,F+f\zeta\,$ of\/ 
$\,\mathcal{F}\hs$ satisfies the conditions\/ $\,\tB=0\,$ and\/ $\,\tvt=0\,$ 
if and only if\/ 
$\,f=r\phi^3+[\lambda(c,X)-\mu(\rd,\rd)-\hh d_\vf\phi]\hs\phi/4\,$ for 
some $\,r\in\bbR\hs$. 
\end{lemma} 
\begin{proof}By (\ref{tfr}.b) with $\,B=K\nh\phi^{-1}\xi\wedge\ts\,$ and 
$\,\beta=\phi^{-2}\xi\,$ (cf.\ (\ref{xuw}) and (\ref{bef}.a)), 
the condition $\,\tB=0\,$ is equivalent to 
$\,K\nh\phi^2\hs\xi\wedge\ts 
=2\hs\xi\wedge\hs[\hs2\hs\phi\hskip1ptd_{\hs\zeta(\hs\cdot\hs)}\hh f 
-h(\phi\hh F\bw,\,\cdot\,)\hh]$. The formula for $\,\phi^2F\bw$ in 
Lemma~\ref{eight} gives 
$\,h(\phi\hh F\bw,\,\cdot\,)=-\hs[\lambda(c,c)-2\hh\mu(c,\rd) 
+K\nh\phi^2\nnh/2\hh]\hs\ts+[\hs\mu(\rd,\rd)-d_\vf\phi\hh]\hs\xi$, since 
$\,h(\rd,\hs\cdot\hs)=-\phi\hs\ts\,$ and $\,h(c,\hs\cdot\hs)=\phi\hs\xi\,$ 
(see (\ref{hzw})). Next, 
$\,\xi\wedge\hskip1ptd_{\hs\zeta(\hs\cdot\hs)}\hh f= 
-2\phi^{-2}(d_cf)\hs\xi\wedge\ts$, as one sees evaluating both sides on 
$\,(\bw,w)$, for $\,w\in\dbs$, and using (\ref{xwb}). Thus, $\,\tB=0\,$ if and 
only if 
\begin{equation}\label{dcf} 
4\hskip1.4ptd_cf\,=\,[\lambda(c,c)-2\hh\mu(c,\rd)]\hs\phi\hs. 
\end{equation} 
Let us now assume (\ref{dcf}). As 
$\,F=F^K\nh+F^\vf\nh+F^\lambda\nh+F^\mu\nnh$, the formulae in 
Lemma~\ref{eight} easily yield 
$\,\phi\,\mathrm{div}^{\mathcal{V}}\nh(Fw) 
=[\lambda(c,c)-2\hh\mu(c,\rd)+3K\nh\phi^2\nnh/2\hh]\hs\ts(w) 
+[2\hh\lambda(c,\rd)-3\hh\mu(\rd,\rd)-d_\vf\phi\hh]\hs\xi(w)\,$ and 
$\,3\hs\phi\hh h(F\bw,w) 
=-\hs3\hh[\lambda(c,c)-2\hh\mu(c,\rd)+K\nh\phi^2\nnh/2\hh]\hs\ts(w) 
+3\hh[\hs\mu(\rd,\rd)-d_\vf\phi\hh]\hs\xi(w)$, while (\ref{bef}.d) and 
(\ref{bef}.a) give 
$\hs-\hs4\hs\phi\hskip1.4ptd_{\hs\zeta w}f= 
8\hs\phi^{-1}(d_cf\hh)\hs\ts(w)+8\hs\phi^{-1}(d_\rd f\hh)\hs\xi(w)\,$ and 
$\,-24\hs\phi f\beta(w)=-24f\phi^{-1}\hs\xi(w)$. Adding the last four 
equalities side-by-side, and using the relation $\,\vt=0\,$ (cf.\ 
(\ref{xuw})), we see that, by (\ref{tfr}.c), $\,\tvt=0\,$ if and only if 
$\,4(d_\rd f-3f)\hs\phi^{-1}\nh=2\hskip1.4ptd_\vf\phi-\lambda(c,\rd)$, or, 
equivalently, $\,4\hskip1.4ptd_\rd(f\phi^{-3}) 
=2\hs\phi^{-2}\hs d_\vf\phi-\phi^{-2}\lambda(c,\rd)$. (The terms involving 
$\,\ts(w)\,$ add up to zero as a consequence of (\ref{dcf}).) The system 
formed by this last equation and 
$\,4\hskip1.4ptd_c(f\phi^{-3})=\phi^{-2}\lambda(c,c) 
-2\hs\phi^{-2}\mu(c,\rd)\,$ (which is immediate from (\ref{dcf}) and 
(\ref{duf}.iii)) determines the solution $\,f\phi^{-3}:\plane_+\to\bbR\,$ 
uniquely up to an additive constant. Our assertion now follows, since a 
solution may be defined by $\,4f\phi^{-3}\nh=[\hh\lambda(c,\rd)-\mu(\rd,\rd) 
-d_\vf\phi\hs]\hs\phi^{-2}\nnh$. 
\end{proof}

\section{The remaining condition in Theorem~\ref{crvtc}}\label{trct} 
\setcounter{equation}{0} 
Suppose that 
$\,(M,\mathcal{V}\nh,\hh\mathrm{D},h,\alpha,\beta,\theta,\zeta)\,$ and 
$\,\mathcal{H}\,$ are chosen as in Lemma~\ref{gflat}, for a fixed 
tw\hbox{o\hh-\nh}\hskip0ptplane system 
$\,(\bs,\xi,\ts,\plane,c,\varOmega)$, while $\,K\,$ is a given real 
constant, $\,F\,$ is a section of $\,\mathcal{F}\nnh$, and the new horizontal 
distribution $\,\thor=\mathcal{H}+F\,$ satisfies conditions (a) -- (c) in 
Theorem~\ref{crvtc}. Thus, $\,F\,$ can be uniquely written as 
$\,F=F^K\nh+F^\vf\nh+F^\lambda\nh+F^\mu\nh+f\zeta$, with the summands defined 
as in Lemmas~\ref{eight} and~\ref{dcfeq}, and hence depending, for any fixed 
$\,y\in\varSigma$, on the quadruple 
$\,(\vf,\lambda,\mu,r)\in V\nh\times\bbR\hs$, where 
$\,V\nnh=\plane\times[\plane^*]^{\odot2}\nnh\times[\plane^*]^{\odot2}\nh$ 
(thus, $\,\dim\,V\nh=8$). As $\,(\vf,\lambda,\mu,r)\,$ varies with $\,y$, 
rather than being a single element of $\,V\nh\times\bbR\hs$, it constitutes 
$\,(V\nh\times\bbR)$-val\-ued function on some connected open 
set $\,\,U\,$ in the af\-fine plane $\,\bs$. 
\begin{theorem}\label{cndtd}Under the above assumptions, given\/ 
$\,a\in\plane\,$ with\/ $\,\varOmega(a,c)=1$, condition\/ {\rm(d)} in 
Theorem\/~{\rm\ref{crvtc}} is satisfied by\/ $\,\thor=\mathcal{H}+F\,$ if and 
only if, for some function\/ $\,\gs:U\to\bbR\hs$, 
\begin{equation}\label{eqo} 
\begin{array}{l} 
[\hh2\mu_1(\rd,\rd)-\lambda_2(\rd,\rd)]\hs\phi 
=4\lambda(c,\rd)\hs\mu(\rd,\rd)-4\mu(c,\rd)\hs\lambda(\rd,\rd) 
-2K\nh\phi^2\varOmega(\rd,\vf)\hs,\\ 
\lambda_1(c,\rd)+2\hh\varOmega(\rd,\vf_2)=4\hh\mu(\vf,\rd)+\varOmega(\rd,a) 
-\phi\hskip1.6pt\mathrm{det}_\varOmega\lambda-\gs\hh\phi\hs, 
\end{array} 
\end{equation} 
where\/ $\,\rd\,$ is the radial vector field on\/ $\,\plane\,$ and\/ 
$\,\phi=\varOmega(\rd,c)$, as well as 
\begin{equation}\label{oqc} 
K\varOmega(\vf_1,c)+\gs_2=2K\lambda(c,\vf)+8\hh r\hh. 
\end{equation} 
The subscripts denote here the partial derivatives with respect to the 
af\-fine coordinates\/ $\,y^{\hs j}$ in\/ $\,\bs\,$ such that\/ 
$\,d\hskip.2pty^{\nh1}\nh=\xi\,$ and\/ $\,d\hskip.2pty^2\nh=\ts$. 
\end{theorem} 
The symbol $\,\mathrm{det}_\varOmega\lambda\,$ in (\ref{eqo}) represents the 
function $\,\,U\nh\to\bbR\,$ assigning to $\,y\in U\,$ the ratio 
$\,(\mathrm{det}\hskip2.7pt\varOmega)^{-1}\mathrm{det}\hskip2.7pt\lambda$, in 
which $\,\lambda\,$ stands for the value of $\,\lambda\,$ at $\,y$, and the 
determinants of the bi\-lin\-e\-ar forms 
$\,\lambda,\varOmega\in[\plane^*]^{\otimes2}$ are evaluated in any basis of 
$\,\plane$. Thus, by (\ref{duf}.i), 
\begin{equation}\label{det} 
[\hh\lambda(c,\rd)]^2\nh-\lambda(c,c)\hs\lambda(\rd,\rd) 
=-\hs\phi^2\mathrm{det}_\varOmega\lambda\hs, 
\end{equation} 
since one may use the basis $\,c,\rd\,$ (even though it depends on a point of 
$\,\plane_+$).

A proof of Theorem~\ref{cndtd} will be given in Sections~\ref{ptfp} 
and~\ref{ptsp}. 
\begin{theorem}\label{equiv}Under the same hypotheses as above, the 
assignment\/ $\,(\vf,\lambda,\mu,r)\mapsto(\vf,\lambda,\mu)$ defines a 
bijective correspondence between functions\/ 
$\,(\vf,\lambda,\mu,r):U\nh\to V\times\bbR\,$ satisfying conditions\/ 
{\rm(\ref{eqo})} -- {\rm(\ref{oqc})} for some\/ $\,\gs:U\to\bbR\hs$, and 
functions\/ $\,(\vf,\lambda,\mu):U\nh\to V\hs$ with 
\begin{equation}\label{eqn} 
\begin{array}{l} 
[\hh2\mu_1(\rd,\rd)-\lambda_2(\rd,\rd)]\hs\phi 
=4\lambda(c,\rd)\hs\mu(\rd,\rd)-4\mu(c,\rd)\hs\lambda(\rd,\rd) 
-2K\nh\phi^2\varOmega(\rd,\vf)\hs,\\ 
\lambda_1(c,c)+2\hh\varOmega(c,\vf_2)=4\mu(\vf,c)-r1\hs. 
\end{array} 
\end{equation} 
\end{theorem} 
\begin{proof}As $\,d_c\phi=0\,$ by (\ref{duf}.iii), applying $\,d_c$ to the 
second equality in (\ref{eqo}) we see that (\ref{eqo}) implies (\ref{eqn}). 
Conversely, suppose that $\,\vf,\lambda,\mu\,$ satisfy (\ref{eqn}) and 
$\,a\in\plane\,$ is a fixed vector with $\,\varOmega(a,c)=1$. For any 
$\,y\in\bs$, if $\,\psi\,$ denotes the restriction of the function 
$\,\lambda_1(c,\rd)+2\hh\varOmega(\rd,\vf_2)-4\hh\mu(\vf,\rd) 
-\varOmega(\rd,a)+\phi\hskip1.6pt\mathrm{det}_\varOmega\lambda\,$ to the 
set $\,\{y\}\times\plane_+\subset M\,$ (which we may identify with 
$\,\plane_+$), the second equality in (\ref{eqn}) and (\ref{duf}.iii) give 
$\,d_c\psi=0$. Since $\,\psi\,$ is a linear functional, (\ref{duf}.iii) 
implies that $\,\psi=-\hs\gs\hh\phi\,$ for some scalar $\,\gs$, depending on 
$\,y$, which yields the second relation in (\ref{eqo}). Our assertion now 
follows, as $\,\vf,\lambda,\mu\,$ uniquely determine $\,\gs\,$ (and $\,r$) via 
the second equality in (\ref{eqo}) and then, respectively, via (\ref{oqc}).
\end{proof} 
\begin{example}\label{lccne}Let us fix $\,K\in\bbR\,$ and $\,a\in\plane\,$ 
with $\,\varOmega(a,c)=1$. A solution $\,(\vf,\lambda,\mu):\bs\to V$ to 
(\ref{eqn}) can obviously be defined by setting $\,\vf=0$, $\,\mu=0\,$ and 
choosing the real-val\-ued component functions 
$\,\lambda(c,c),\hs\lambda(a,a),\hs\lambda(a,c)\,$ of $\,\lambda\,$ relative 
to the basis $\,c,a\,$ of $\,\plane\,$ in such a way that $\,\lambda(c,c)\,$ 
equals a constant minus the coordinate function $\,y^{\nh1}\nnh$, while 
$\,\lambda(a,a)\,$ and $\,\lambda(a,c)\,$ are arbitrary functions of 
$\,y^{\nh1}\nnh$. 
\end{example} 
\begin{remark}\label{unkfc}For $\,\lambda,\mu:U\nh\to[\plane^*]^{\odot2}$ as 
above, let $\,\dv\,$ be the function on $\,\,U\,$ valued in 
en\-do\-mor\-phisms of $\,\plane\,$ which corresponds to $\,\lambda\,$ via 
$\,\varOmega$, so that $\,\lambda(u,v)=\varOmega(\dv u,v)\,$ for all 
$\,u,v\in\plane$. Thus, 
$\,\lambda(c,\rd)\hs\mu(\rd,\rd)-\mu(c,\rd)\hs\lambda(\rd,\rd) 
=\mu(c,\rd)\hs\varOmega(\rd,\dv\rd)+\mu(\rd,\rd)\hs\varOmega(\dv\rd,c)$. Hence
\begin{equation}\label{lcx} 
\lambda(c,\rd)\hs\mu(\rd,\rd)-\mu(c,\rd)\hs\lambda(\rd,\rd) 
=\phi\hs\mu(\dv\rd,\rd) 
\end{equation} 
by (\ref{lcv}) with $\,v=\dv\rd$. This allows us to rewrite the first equality 
in (\ref{eqn}) or (\ref{eqo}) as 
\begin{equation}\label{tmx} 
2\mu_1(\rd,\rd)-\lambda_2(\rd,\rd)=4\mu(\dv\rd,\rd) 
-2K\nh\phi\hs\varOmega(\rd,\vf)\hs. 
\end{equation} 
We also have $\,\dv=\varOmega^{-1}\nnh\lambda$, meaning that 
$\,\dv\,$ is, at each $\,y\in U\nh$, the composite of $\,\lambda\,$ and the 
inverse of $\,\varOmega$, where the values of $\,\lambda\,$ and 
$\,\varOmega\,$ at $\,y\,$ are treated as linear operators 
$\,\plane\to\plane^*$ sending $\,u$ to $\,\lambda(u,\,\cdot\,)\,$ and 
$\,\varOmega(u,\,\cdot\,)$. Consequently, by (\ref{duf}.i), equations 
(\ref{eqn}) now take the form 
\begin{equation}\label{qln} 
\begin{array}{l} 
2\mu_1(\rd,\rd)-\lambda_2(\rd,\rd)=4\mu((\varOmega^{-1}\nnh\lambda)\rd,\rd) 
-2K\varOmega(\rd,c)\hs\varOmega(\rd,\vf)\hs,\\ 
\lambda_1(c,c)+2\hh\varOmega(c,\vf_2)=4\mu(\vf,c)-1\hs. 
\end{array} 
\end{equation} 
Fixing $\,a\in\plane\,$ such that $\,\varOmega(a,c)=1\,$ and using it to 
represent $\,\vf,\lambda\,$ and $\,\mu\,$ by their components 
$\,\mu(c,c),\hs\mu(a,a),\hs\mu(a,c),\hs\lambda(c,c),\hs\lambda(a,a), 
\hs\lambda(a,c),\hs\varOmega(c,\vf),\hs\varOmega(a,\vf)\,$ relative to the 
basis $\,c,a\,$ of $\,\plane$, which form an octuple of functions 
$\,\,U\nh\to\bbR\hs$, we see that (\ref{qln}) (and, therefore, (\ref{eqn})) 
is equivalent to the following system of four first-or\-der qua\-si-lin\-e\-ar 
partial differential equations with eight unknown real-val\-ued functions of two 
real variables: 
\begin{equation}\label{loc} 
\begin{array}{l} 
2\mu_1(c,c)-\lambda_2(c,c)=4\lambda(c,c)\hs\mu(a,c) 
-4\mu(c,c)\hs\lambda(a,c)\hs,\\ 
2\mu_1(a,a)-\lambda_2(a,a)=4\lambda(a,c)\hs\mu(a,a)-4\mu(a,c)\hs\lambda(a,a) 
-2K\nh\varOmega(a,\vf)\hs,\\ 
2\mu_1(a,c)-\lambda_2(a,c)=4\lambda(a,c)\hs\mu(a,c)-4\mu(c,c)\hs\lambda(a,a) 
-2K\nh\varOmega(c,\vf)\hs,\\ 
\lambda_1(c,c)+2\hh\varOmega(c,\vf_2)=4\hh\varOmega(a,\vf)\hs\mu(c,c) 
-4\hh\varOmega(c,\vf)\hs\mu(a,c)-1\hs. 
\end{array} 
\end{equation}
The subscripts in $\,\mu_1(c,c),\hs\varOmega(c,\vf_2)$, etc., stand for 
partial derivatives of $\,\mu(c,c)\,$ and $\,\varOmega(c,\vf)$.

In fact, one obtains (\ref{loc}) by applying $\,d_cd_c,\hs d_ad_a$ and 
$\,d_cd_a$ to the first equality in (\ref{qln}) and noting that 
$\,\mu(\vf,c)=\varOmega(a,c)\hs\mu(\vf,c)=\varOmega(a,\vf)\hs\mu(c,c) 
-\varOmega(c,\vf)\hs\mu(a,c)\,$ (cf.\ Remark~\ref{cclsm}). 
\end{remark}

\section{Some lemmas}\label{sole} 
\setcounter{equation}{0} 
Throughout this section we make the same assumptions as in Section~\ref{trct}, 
and $\,w,w\hh'\nnh,w\hh''$ always stand for constant vector fields on the 
af\-fine plane $\,\bs$, treated also as $\,\mathcal{V}$-pro\-ject\-a\-ble 
sections of $\,\mathcal{H}$. As $\,\xi\,$ and $\,\ts\,$ are constant 
$\,1$-forms on $\,\bs$, it follows that $\,\xi(w),\ts(w),\ts(w\hh'\hh)$, 
etc., are constant functions on $\,M=\bs\times\plane_+$. We also set 
$\,\tw=w+Fw$.

For $\,\,U\subset\bs\,$ as in Section~\ref{trct}, let the functions 
$\,\qx,\ex,\lx,\lx^\pm:U\times\plane_+\to\bbR\,$ be given by 
\begin{equation}\label{abd} 
\begin{array}{rl} 
\mathrm{i)}\hskip0pt&\qx=\lambda(\rd,\rd)+2\hh\varOmega(\rd,\vf)\hs,\hskip22pt 
\ex=\lambda(c,c)-2\hh\mu(c,\rd)+K\nh\phi^2\nnh/\hh2\hs,\\ 
\mathrm{ii)}\hskip0pt&\lx=\mu(\rd,\rd)-\varOmega(\vf,c)\hs,\hskip34pt 
\lx^\pm\nh=\lx\pm2f\phi^{-1}\hs. 
\end{array} 
\end{equation} 
Since $\,F=F^K\nh+F^\vf\nh+F^\lambda\nh+F^\mu\nh+f\zeta$, Lemmas~\ref{eight} 
--~\ref{dcfeq}, (\ref{duf}) and (\ref{hvv}.c) yield 
\begin{equation}\label{ffw} 
\begin{array}{rl} 
\mathrm{i)}\hskip0pt&2\hh\lp\nh=\mu(\rd,\rd)+\lambda(c,\rd) 
-3\hh\varOmega(\vf,c)+4\hh r\phi^2,\\ 
\mathrm{ii)}\hskip0pt&2\hh\lm\nh=3\mu(\rd,\rd)-\lambda(c,\rd) 
-\varOmega(\vf,c)-4\hh r\phi^2,\\ 
\mathrm{iii)}\hskip0pt&\phi\hh Fw=[\ex\ts(w)-\lp\nh\xi(w)]\hs\rd 
+[\lm\nh\ts(w)+\qx\hs\xi(w)]\hs c\hs,\\ 
\mathrm{iv)}\hskip0pt&d_{Fw}\phi=\ex\ts(w)-\lp\nh\xi(w)\hs. 
\end{array} 
\end{equation} 
(By (\ref{duf}.i) and Remark~\ref{cclsm}, $\,\phi=\varOmega(\rd,c)\,$ and 
$\,\varOmega(\vf,c)\hh\rd+\varOmega(\rd,\vf)\hs c+\varOmega(c,\rd)\hs\vf=0$, 
so that $\,\vf=\phi^{-1}\varOmega(\vf,c)\hh\rd 
+\phi^{-1}\varOmega(\rd,\vf)\hs c$, while (\ref{bef}.d) gives 
$\,\zeta w=-2\phi^{-2}\yw$.) Now, from (\ref{hzw}), 
\begin{equation}\label{hfw} 
h(Fw,w\hh'\hh)=\qx\hs\xi(w)\hs\xi(w\hh'\hh) 
+\lp\nh\xi(w)\hh\ts(w\hh'\hh) 
+\lm\nh\ts(w)\hs\xi(w\hh'\hh) 
-\ex\hs\ts(w)\hh\ts(w\hh'\hh)\hs. 
\end{equation} 
\begin{remark}\label{polyn}Any degree $\,k\,$ homogeneous polynomial function 
on $\,\plane$, valued in an arbitrary 
fi\-\hbox{nite\hh-}\hskip0ptdi\-men\-sion\-al vector space, is 
\begin{enumerate} 
  \def\theenumi{{\rm\alph{enumi}}} 
\item[{\rm(a)}] an eigen\-vec\-tor of $\,d_X$ for the eigenvalue $\,k$, 
\item[{\rm(b)}] an eigen\-vec\-tor of $\,d_X\nh-\mathrm{Id}\,$ for the 
eigen\-value $\,k-1$. 
\end{enumerate} 
(In fact, (b) is obvious from (a).) Examples of such functions include 
$\,\rd\,$ (valued in $\,\plane$, with $\,k=1$), as well as the real-val\-ued 
functions such as $\,\phi\,$ or $\,\lambda(c,\rd)\,$ (with $\,k=1$) and 
$\,\mu(\rd,\rd)$ (with $\,k=2$). Thus, for instance, $\,d_\rd\phi=\phi$. 
\end{remark} 
\begin{lemma}\label{brckt}The Lie bracket\/ $\,[Fw,Fw\hh'\hh]\,$ equals 
$\,\phi^{-2}(\xi\wedge\ts)(w,w\hh'\hh)\,$ times 
\[ 
[(2\hh\qx-d_\rd\qx)\hh\ex+(2\lm\nnh-d_\rd\lm)\lp\nnh 
+\qx\hs d_c\lm\nnh-\lm\nh d_c\qx\hh]\hs c 
+(\ex\hs d_\rd\lp\nnh-\lp\nh d_\rd\ex+\lm\nh d_c\lp\nnh 
+\qx\hs d_c\ex)\hh\rd\hh. 
\] 
\end{lemma} 
\begin{proof}We have $\,\phi^2[Fw,Fw\hh'\hh]=[\phi Fw,\phi Fw\hh'\hh] 
-\phi\hs(d_{Fw}\phi)Fw\hh'+\phi\hs(d_{Fw\hh'}\phi)Fw$. By (\ref{ffw}.iv), 
(\ref{ffw}.iii) and (\ref{bwa}.ii), the last two terms add up to 
$\,(\lp\nh\lm\nh+\ex\qx)\hs[(\xi\wedge\ts)(w,w\hh'\hh)]\hs c$. 
As $\,[\hh c,\rd\hh]=d_c\rd=c$, our assertion follows if we evaluate 
$\,[\phi Fw,\phi Fw\hh'\hh]\,$ using (\ref{ffw}.iii) and (\ref{bwa}.ii). 
\end{proof} 
\begin{lemma}\label{iprlb}The $\,h$-in\-ner product\/ 
$\,h([Fw\hh'\nnh,Fw\hh''\hh],w)\,$ equals $\,\zeta(w\hh'\nnh,w\hh''\hh)/2\,$ 
times 
\[ 
\begin{array}{l} 
[(2\hh\qx-d_\rd\qx)\hh\ex+(2\lm\nnh-d_\rd\lm)\lp\nnh 
+\qx\hs d_c\lm\nnh-\lm\nh d_c\qx\hh]\hs\xi(w)\\ 
\phantom{[(2\hh\qx-d_\rd\qx)\hh\ex+(2\lm\nnh-d_\rd\lm) 
\lp\nnh}+\hskip2pt(\lp\nh d_\rd\ex 
-\ex\hs d_\rd\lp\nnh-\lm\nh d_c\lp\nnh 
-\qx\hs d_c\ex)\hs\ts(w)\hh. 
\end{array} 
\] 
\end{lemma} 
\begin{proof}This is obvious from Lemma~\ref{brckt}, (\ref{hzw}) and 
(\ref{bef}.c). 
\end{proof} 
As in the statement of Theorem~\ref{cndtd}, here and in the next section the 
subscripts $\,(\hskip2.5pt)_j$, $\,j=1,2$, denote the directional derivatives 
in the directions of the constant vector fields $\,\partial_j$ on $\,\bs\,$ 
form\-ing the basis of $\,\dbs\,$ dual to the basis $\,\xi,\ts\,$ of 
$\,\dbs^*\nnh$. In other words, 
$\,(\hskip2.5pt)_j=\partial/\partial y^{\hs j}$ for the 
af\-fine coordinates $\,y^{\hs j}$ with $\,d\hskip.2pty^{\nh1}\nh=\xi\,$ and 
$\,d\hskip.2pty^2\nh=\ts$. 
\begin{remark}\label{dpsie}For any function $\,\chi\,$ on $\,\bs\,$ one has 
$\,d_w[\chi\hs\xi(w\hh'\hh)]-d_{w\hh'}[\chi\hs\xi(w)] 
=-\hs\chi_2(\xi\wedge\ts)(w,w\hh'\hh)$ and 
$\,d_w[\chi\hs\ts(w\hh'\hh)]-d_{w\hh'}[\chi\hs\ts(w)] 
=\chi_1(\xi\wedge\ts)(w,w\hh'\hh)$, since $\,d\chi=\chi_1\xi+\chi_2\ts$. 
(Cf.\ (\ref{bwa}).) 
\end{remark} 
\begin{lemma}\label{hwpfw}The expression\/ 
$\,h([w\hh'\nnh,Fw\hh''\hh],w)+h([Fw\hh'\nnh,w\hh''\hh],w)\,$ equals\/ 
$\,(\xi\wedge\ts)(w,w\hh'\hh)\,$ times\/ 
$\,(\lmo-\qx_2)\hs\xi(w) 
-(\ex_1+\lpt)\hs\ts(w)$, that is, 
$\,\zeta(w\hh'\nnh,w\hh''\hh)/2\,$ times\/ 
$\,(\lmo-\qx_2)\hh\phi\hs\xi(w) 
-(\ex_1+\lpt)\hh\phi\hh\ts(w)$. 
\end{lemma} 
\begin{proof}Treating $\,Fw\hh'$ and $\,[w,Fw\hh'\hh]\,$ as functions 
$\,\,U\nh\to\plane$, we have $\,[w,Fw\hh'\hh]=d_w(Fw\hh'\hh)$, so that, 
by (\ref{ffw}.iv) and Remark~\ref{dpsie}, 
$\,\phi\hs[w,Fw\hh'\hh]+\phi\hs[Fw,w\hh'\hh] 
=d_w(\phi\hs Fw\hh'\hh)-d_{w\hh'}(\phi\hs Fw)\,$ 
equals $\,(\xi\wedge\ts)(w,w\hh'\hh)\,$ times 
$\,(\lmo-\qx_2)\hs c+(\ex_1+\lpt)\hs\rd$. Now (\ref{duf}) and (\ref{bef}.c) 
yield our claim. 
\end{proof} 
\begin{lemma}\label{tgatw}If\/ $\,\tw=w+Fw$, then 
\[ 
\begin{array}{l} 
8\hh\tga(\tw)=2\hs[\hs3\mu_1(\rd,\rd)-2\lambda_2(\rd,\rd) 
-4\hh\varOmega(\rd,\vf_2)-\varOmega(\vf_1,c)-\lambda_1(c,\rd) 
-4\hs\phi^2r_1]\hs\phi\hs\xi(w)\\ 
\phantom{2\tga(\tw)}+\hskip2pt2\hs[4\mu_1(c,\rd)-\mu_2(\rd,\rd) 
+3\hh\varOmega(\vf_2,c)-\lambda_2(c,\rd)-2\lambda_1(c,c)-4\hs\phi^2r_2]\hs 
\phi\hs\ts(w)\\ 
\phantom{2\tga(\tw)}+\hskip2pt[\hs12\mu(c,\rd)\hs\lambda(\rd,\rd) 
-13\lambda(c,\rd)\hs\mu(\rd,\rd)+4K\nh\phi^2\varOmega(\rd,\vf) 
+12\hs r\phi^2\lambda(c,\rd)]\hs\xi(w)\\ 
\phantom{2\tga(\tw)}+\hskip2pt[\hs10\hs\varOmega(\vf,c)\hs\mu(\rd,\rd) 
-2\lambda(c,c)\hs\lambda(\rd,\rd)+3(\lambda(c,\rd))^2\hh]\hs\xi(w)\\ 
\phantom{2\tga(\tw)}+\hskip2pt[\hs8\mu(c,\rd)\hs\varOmega(\rd,\vf) 
-24\hs r\hs\varOmega(\vf,c)\hs\phi^2\hh]\hs\xi(w)\\ 
\phantom{2\tga(\tw)}+\hskip2pt[\hs4\lambda(c,c)\hs\varOmega(\rd,\vf) 
+\varOmega(\vf,c)\hs\lambda(c,\rd)]\hs\xi(w) 
+2\hs[\hh\varOmega(\vf,c)]^2\hs\xi(w)\\ 
\phantom{2\tga(\tw)}+\hskip2pt[\hs24\hs r\phi^2\mu(c,\rd) 
-2\mu(c,\rd)\hs\mu(\rd,\rd)+K\nh\phi^2\lambda(c,\rd)]\hs\ts(w)\\ 
\phantom{2\tga(\tw)}+\hskip2pt[\hs8\mu(c,c)\hs\lambda(\rd,\rd) 
-7\lambda(c,c)\hs\mu(\rd,\rd)]\hs\ts(w)\\ 
\phantom{2\tga(\tw)}+\hskip2pt[\hh2\lambda(c,\rd)\hs\mu(c,\rd) 
-12\hs r\lambda(c,c)\hs\phi^2\nh-6K\varOmega(\vf,c)\hs\phi^2\hh]\hs 
\ts(w)\\ 
\phantom{2\tga(\tw)}+\hskip2pt[16\mu(c,c)\hs\varOmega(\rd,\vf) 
-\lambda(c,c)\hs\lambda(c,\rd)+14\hs\varOmega(\vf,c)\hs\mu(c,\rd)]\hs 
\ts(w)\\ 
\phantom{2\tga(\tw)}+\hskip2pt\lambda(c,c)\hs\varOmega(\vf,c)\hs\ts(w)\hs. 
\end{array} 
\] 
\end{lemma} 
\begin{proof}As $\,\gamma=0\,$ (see Lemma~\ref{gflat}) and 
$\,[w\hh'\nnh,w\hh''\hh]=0$, (\ref{ape}.d) implies that 
$\,\tga(\tw)\hs\zeta(w\hh'\nnh,w\hh''\hh)$ equals the sum of the 
in\-ner-prod\-uct expressions appearing in Lemmas~\ref{iprlb} and~\ref{hwpfw}. 
Our assertion now follows from these two lemmas, (\ref{abd}.i), 
(\ref{ffw}.i), (\ref{ffw}.ii) and Remark~\ref{polyn}. (By 
Remark~\ref{polyn}(b), $\,2\hh\qx-d_\rd\qx=2\hh\varOmega(\rd,\vf)\,$ and 
$\,2\lm\nnh-d_\rd\lm\nh=-\hs\varOmega(\vf,c)-\lambda(c,\rd)/2$.) 
\end{proof} 
For $\,f\,$ as in Lemma~\ref{dcfeq}, any fixed vector $\,a\in\plane\,$ with 
$\,\varOmega(a,c)=1$, and any section $\,v\,$ of $\,\mathcal{V}\nh$, 
\begin{equation}\label{fgv} 
\begin{array}{rl} 
\mathrm{a)}\hskip0pt&4\hh\tga(v)=(\ex+K\nh\phi^2\nnh/\hh2)\hs\varOmega(\rd,v) 
-(\lx+4\hh r\phi^2)\hs\varOmega(v,c)\hs,\\ 
\mathrm{b)}\hskip0pt&\phi\hs\ta(\tw)=(\lx+4\hh r\phi^2)\hs\xi(w) 
-(\ex+K\nh\phi^2\nnh/\hh2)\hs\ts(w)\hs,\\ 
\mathrm{c)}\hskip0pt&\ta(\tw)+\phi^{-1}d_{Fw}\phi 
=-2\ea\hh\xi(w)-K\nh\phi\hs\ts(w)/2\hs,\hskip8pt\mathrm{for}\hskip6pt\ea 
=f\phi^{-2}\nh-2\hh r\phi\hs,\\ 
\mathrm{d)}\hskip0pt&\phi^{-2}h(v,w)=d_v\hs[\hh\phi^{-1}\ts(w) 
-\phi^{-1}\varOmega(\rd,a)\hs\xi(w)]\hs,\\ 
\mathrm{e)}\hskip0pt&\phi^{-2}\hs\tga(v)=d_v\ea+K\varOmega(\rd,v)/4\hs, 
\hskip8pt\mathrm{where}\hskip6pt\ea\hs=f\phi^{-2}\nh-2\hh r\phi\hs. 
\end{array} 
\end{equation} 
In fact, (\ref{ape}.c) with with $\,\trf=f\,$ and $\,\gamma=0\,$ (cf.\ 
Section~\ref{esol} and Lemma~\ref{gflat}) yields 
$\,\tga(v)=-\hs\phi\hs d_v(f\phi^{-1})+\varOmega(\ex\rd+\lm\nnh c,v)/2$, 
since $\,F\bw=\phi^{-2}(\ex\rd+\lm\nnh c)\,$ (from (\ref{ffw}.iii) and 
(\ref{xwb})), while $\,\phi^{-2}\theta=\varOmega\,$ (by (\ref{bef}.b)), and 
$\,\alpha(v)=-\hs\phi^{-1}d_v\phi\,$ (see (\ref{aeh}) and (\ref{peo})). Thus, 
(\ref{duf}) and the formula for $\,f\,$ in Lemma~\ref{dcfeq} yield 
$\,4\hh\tga(v)=2\ex\hh\varOmega(\rd,v) 
-(2\hh\lm\nh+8\hh r\phi^2)\hs\varOmega(v,c) 
+[\hh2\mu(v,\rd)-\lambda(c,v)]\hs\phi$. Relation (\ref{fgv}.a), that is, 
vanishing of $\,4\hh\tga(v)-(\ex+K\nh\phi^2\nnh/\hh2)\hs\varOmega(\rd,v) 
+(\lx+4\hh r\phi^2)\hs\varOmega(v,c)$, is now immediate from (\ref{abd}.i), as 
(\ref{ffw}.ii) gives 
$\,2\hh\lm\nh+8\hh r\phi^2\nh=3\mu(\rd,\rd)-\lambda(c,\rd)-\varOmega(\vf,c) 
+4\hh r\phi^2\nnh$, while 
\begin{equation}\label{lcv} 
\begin{array}{l} 
\phi\hs\lambda(c,v)=\varOmega(v,c)\hs\lambda(c,\rd) 
+\lambda(c,c)\hs\varOmega(\rd,v)\hs,\\ 
\phi\hs\mu(v,\rd)=\mu(c,\rd)\hs\varOmega(\rd,v) 
+\mu(\rd,\rd)\hs\varOmega(v,c)\hs, 
\end{array} 
\end{equation} 
as a consequence of (\ref{duf}.i) and Remark~\ref{cclsm}.

Since $\,\ta(\tw)=2\hh\tga(\zeta\tw)=2\tga(\zeta w)\,$ (see 
Section~\ref{hrzd}, (\ref{mor}) and (\ref{dub}.h)), (\ref{fgv}.b) follows from 
(\ref{fgv}.a), (\ref{bef}.d), (\ref{hvv}.c) and (\ref{duf}.i).

Next, (\ref{fgv}.b) and (\ref{ffw}.iv) give 
$\,[\hh\ta(\tw)+\phi^{-1}d_{Fw}\phi\hh]\hs\phi 
=(\lx-\lp\nh+4\hh r\phi^2)\hs\xi(w) 
-K\nh\phi^2\ts(w)/2$. As $\,\lx-\lp\nh=-2f\phi^{-1}$ (see 
(\ref{abd}.ii)), this equals 
$\,-\hs[2\ea\hh\xi(w)+K\nh\phi\hs\ts(w)/2\hh]\hs\phi$, which yields 
(\ref{fgv}.c).

Subtracting the left-hand side of (\ref{fgv}.d) from its right-hand side, and 
evaluating the difference with the aid of (\ref{duf}.ii) and (\ref{hvv}), we 
obtain zero, since (\ref{duf}.i) and Remark~\ref{cclsm} give 
$\,\varOmega(v,c)\hs\varOmega(\rd,a)-\varOmega(v,a)\hs\phi 
=\varOmega(v,\rd)\hs\varOmega(c,a)=\varOmega(\rd,v)$. We thus obtain 
(\ref{fgv}.d).

Finally, using differentiation by parts and (\ref{duf}.ii), we obtain 
$\,\phi^2d_v\ea=\phi d_v(\phi\ea)-\varOmega(v,c)\hs\phi\ea$, for 
$\,\ea\,$ as in (\ref{fgv}.e). Since, by Lemma~\ref{dcfeq} and 
(\ref{duf}.ii), 
\begin{equation}\label{ffa} 
4\hs\phi\ea=\lambda(c,X)-\mu(\rd,\rd)-\varOmega(\vf,c)-4\hh r\phi^2, 
\end{equation} 
the preceding equality, (\ref{fgv}.a), 
(\ref{duf}.ii) and the formulae for $\,\ex\,$ and $\,\lx\,$ in (\ref{abd}) 
yield an expression for $\,4\hh\tga(v)-4\hs\phi^2d_v\ea\,$ involving 
$\,K,\vf,\lambda,\mu,r\,$ and $\,\varOmega,c,\rd,\phi,v\,$ (but not 
$\,f,\qx,\ex,\lx\,$ or $\,\lx^\pm$), which equals 
$\,K\nh\phi^2\varOmega(\rd,v)\,$ in view of (\ref{lcv}). This proves 
(\ref{fgv}.e).

\section{Proof of Theorem~\ref{cndtd}, first part}\label{ptfp} 
\setcounter{equation}{0} 
We use the same assumptions and notations as at the beginning of 
Section~\ref{sole}.

Condition (d) in Theorem~\ref{crvtc} may be naturally split into two parts, 
which read 
\begin{equation}\label{dth} 
\mathrm{i)}\hskip6pt(d\hs\tga+2\hh\ta\wedge\tga)(v,\tw)=h(v,w)/2\hs, 
\hskip16pt\mathrm{ii)}\hskip6pt(d\hs\tga+2\hh\ta\wedge\tga)(\tw,\tw\hh'\hh)=0 
\end{equation} 
for all sections $\,w,w\hh'$ of $\,\mathcal{H}\,$ and $\,v\,$ of 
$\,\mathcal{V}\nh$, where $\,\tw=w+Fw\,$ and 
$\,\tw\hh'\nnh=w\hh'\nnh+Fw\hh'\nnh$. By (\ref{bwa}), 
$\,(d\hs\tga+2\hh\ta\wedge\tga)(v,\tw)-h(v,w)/2 
=d_v[\hs\tga(\tw)]-2\phi^{-1}(d_v\phi)\hh\tga(\tw) 
-d_\tw[\hs\tga(v)]-2\hh\ta(\tw)\hs\tga(v)-\tga([v,\tw\hh])-h(v,w)/2$, 
since, according to Section~\ref{hrzd}, (\ref{aeh}) and (\ref{peo}), 
$\,\ta(v)=\alpha(v)=-\hs\phi^{-1}d_v\phi$. Multiplied by 
$\,\phi^{-2}\nnh$, this becomes 
\begin{equation}\label{fmt} 
d_v[\hs\phi^{-2}\hs\tga(\tw)]-d_\tw[\hs\phi^{-2}\hs\tga(v)] 
-\phi^{-2}\hs\tga([v,\tw\hh]) 
-2\hh[\hh\ta(\tw)+\phi^{-1}d_{Fw}\phi\hh]\hs\phi^{-2}\hs\tga(v) 
-\phi^{-2}h(v,w)/2\hs, 
\end{equation} 
where we have first rewritten $\,\phi^{-2}d_\tw[\hs\tga(v)]\,$ as 
$\,d_\tw[\hs\phi^{-2}\hs\tga(v)]+2\hs\phi^{-3}(d_\tw\phi)\hh\tga(v)$, 
using differentiation by parts, and then noted that 
$\,d_\tw\phi=d_{w+Fw}\phi=d_{Fw}\phi$. 
\begin{lemma}\label{fmtde}For\/ $\,\ea\,$ given by\/ {\rm(\ref{ffa})} 
and a fixed vector\/ $\,a\in\plane\,$ with\/ $\,\varOmega(a,c)=1$, 
the expression\/ {\rm(\ref{fmt})}, that is, 
$\,[(d\hs\tga+2\hh\ta\wedge\tga)(v,\tw)-h(v,w)/2\hh]\hs\phi^{-2}\nnh$, is 
the result of applying $\,d_v$ to 
\begin{equation}\label{ftg} 
\begin{array}{l} 
\phi^{-2}\hs\tga(\tw)\,+\,(2\hs\phi)^{-1}\varOmega(\rd,a)\hs\xi(w)\, 
-\,(2\hs\phi)^{-1}\ts(w)\,-\hs\,d_\tw\ea\,+\,2\hh\ea\nh^2\xi(w)\, 
+\, K\nh\phi\ea\hh\ts(w)\\ 
+\hskip1.7ptK[\varOmega(\rd,\vf)+\lambda(\rd,\rd)/4\hh]\hs\xi(w) 
-K\varOmega(\rd,Fw/4)+K[\hs\mu(\rd,\rd)-\lambda(c,X)]\hh\ts(w)/2\hs. 
\end{array} 
\end{equation} 
\end{lemma} 
\begin{proof}The first and last terms in (\ref{fmt}) are the $\,d_v$-im\-ages 
of the first three terms in (\ref{ftg}), cf.\ (\ref{fgv}.d). In the sum 
of the remaining three terms in (\ref{fmt}), let us replace 
$\,\ta(\tw)+\phi^{-1}d_{Fw}\phi\,$ by the right-hand side of (\ref{fgv}.c), 
$\,\phi^{-2}\hs\tga(v)\,$ by the right-hand side of (\ref{fgv}.e), and 
$\,\phi^{-2}\hs\tga([v,\tw\hh])$ by an analogous expression involving, 
instead of $\,v$, the section $\,[v,\tw\hh]\,$ of $\,\mathcal{V}\,$ 
(see Remark~\ref{liebr}). After rearranging terms and noting that 
$\,K\nh\phi\hh(d_v\ea)\hh\ts(w)= 
d_v[K\nh\phi\ea\hh\ts(w)]-K\nh\ea\hh\varOmega(v,c)\hs\ts(w)\,$ (cf.\ 
(\ref{duf}.ii)), the result is 
\begin{equation}\label{dwd} 
-\hs d_\tw d_v\ea-d_{[v,\tw\hh]}\ea+4\ea\hh(d_v\ea)\hs\xi(w) 
+d_v[K\nh\phi\ea\hh\ts(w)]-d_v[\hh K\varOmega(\rd,Fw/4)] 
\end{equation} 
plus, as explained next, $\,K/4\,$ times 
\begin{equation}\label{aox} 
4\ea\hs[\hh\varOmega(\rd,v)\hs\xi(w)-\varOmega(v,c)\hs\ts(w)] 
+K\nh\phi\hs\varOmega(\rd,v)\hs\ts(w) 
+2\hh\varOmega(v,Fw)\hs. 
\end{equation} 
Namely, the terms 
$\,-K[\hh\varOmega(d_\tw\rd,v)+\varOmega(\rd,[v,\tw\hh])]/4$, 
originally present in the resulting expression, can be rewritten as 
$\,-\hs d_v[\hh K\varOmega(\rd,Fw/4)]\,$ plus $\,K/4\,$ times 
$\,2\hh\varOmega(v,Fw)$, since $\,[v,\tw\hh]=[v,w+Fw]=[v,Fw]=d_v(Fw)$, while 
$\,d_v\rd=v\,$ and $\,d_\tw\rd=d_{w+Fw}\rd=d_{Fw}\rd=Fw$.

As $\,-\hs d_\tw d_v-d_{[v,\tw\hh]}=-\hs d_vd_\tw$, (\ref{dwd}) is the 
$\,d_v$-im\-age of the sum of the fourth, fifth and sixth terms in 
(\ref{ftg}) along with the middle term in the second line of (\ref{ftg}). 
Next, cancelling the $\,K\,$ factor in the second line in (\ref{ftg}), minus 
the middle term, and then applying $\,4\hs d_v$, we get 
\begin{equation}\label{fov} 
[\hh4\hh\varOmega(v,\vf)+2\lambda(v,\rd)]\hs\xi(w) 
+[\hh4\hh\mu(v,\rd)-2\lambda(c,v)]\hh\ts(w)\hs, 
\end{equation} 
Therefore, our claim will follow if we show that (\ref{aox}) equals 
(\ref{fov}). To this end, let us rewrite (\ref{aox}), noting that $\,4\ea 
=[\hh\lambda(c,X)-\mu(\rd,\rd)-\varOmega(\vf,c)]\hs\phi^{-1}-4\hh r\phi\,$ 
(see (\ref{ffa})), and $\,Fw$, in $\,\varOmega(v,Fw)$, may be replaced by 
$\,\phi^{-1}$ times the right-hand side of (\ref{ffw}.iii). Substituting for 
the resulting occurrences of $\,\qx,\ex\,$ and $\,\lx^\pm$ the expressions in 
(\ref{abd}.i), (\ref{ffw}.i) and (\ref{ffw}.ii), and using (\ref{lcv}) along 
with the equalities 
$\,\lambda(c,\rd)\hs\varOmega(\rd,v)+\lambda(\rd,\rd)\hs\varOmega(v,c) 
=\lambda(v,\rd)\hs\phi\,$ and 
$\,\varOmega(v,c)\hs\varOmega(\rd,\vf)+\varOmega(\vf,c)\hs\varOmega(v,\rd) 
=\varOmega(v,\vf)\hs\phi$, both of which immediate from (\ref{duf}.i) and 
Remark~\ref{cclsm}, we see that (\ref{aox}) in fact coincides with 
(\ref{fov}). 
\end{proof} 
In the whole discussion following (\ref{dth}), which includes 
Lemma~\ref{fmtde} and its proof, we never made use of the fact that 
$\,\tga(\tw)\,$ is a specific function $\,U\times\plane_+\to\bbR\hs$, given by 
the formula in Lemma~\ref{tgatw}. This now allows us to solve (\ref{dth}.i) as 
a system of differential equations imposed on $\,\tga(\tw)\,$ treated 
$\,\tga(\tw)\,$ as an arbitrary function. The only assumption made about 
$\,\tga(\tw)\,$ is that its dependence on $\,w\,$ (via the relation 
$\,\tw=w+Fw$) should be val\-ue\-wise and linear, or, equivalently, 
$\,\tga(\tw)\,$ is a combination of $\,\xi(w)\,$ and $\,\ts(w)\,$ with some 
coefficients which are functions on a connected open set in $\,M$. 
\begin{lemma}\label{solns}Solving the system\/ {\rm(\ref{dth}.i)} for the 
unknown function\/ $\,\tga(\tw)$, the dependence of which on\/ $\,w\,$ is 
val\-ue\-wise and linear, we obtain 
\[ 
\begin{array}{l} 
8\hh\tga(\tw)=2\hs[\hs\lambda_1(c,\rd)-\mu_1(\rd,\rd)-\varOmega(\vf_1,c) 
-4\hs\phi^2r_1]\hs\phi\hs\xi(w)\\ 
\phantom{2\tga(\tw)}+\hskip2pt2\hs[\lambda_2(c,\rd)-\mu_2(\rd,\rd) 
-\varOmega(\vf_2,c)-4\hs\phi^2r_2]\hs\phi\hs\ts(w)\\ 
\phantom{2\tga(\tw)}+\hskip2pt[\hs3\lambda(c,\rd)\hs\mu(\rd,\rd) 
-4\mu(c,\rd)\hs\lambda(\rd,\rd)-4K\nh\phi^2\varOmega(\rd,\vf) 
+12\hs r\phi^2\lambda(c,\rd)]\hs\xi(w)\\ 
\phantom{2\tga(\tw)}+\hskip2pt[\hs2\lambda(c,c)\hs\lambda(\rd,\rd) 
-(\lambda(c,\rd))^2\nh-6\hh\varOmega(\vf,c)\hs\mu(\rd,\rd)]\hs\xi(w)\\ 
\phantom{2\tga(\tw)}-\hskip2pt[\hs8\mu(c,\rd)\hs\varOmega(\rd,\vf) 
+24\hs r\hh\varOmega(\vf,c)\hs\phi^2\nh+4\hs\phi\hs\varOmega(\rd,a) 
-4\gs\phi^2\hh]\hs\xi(w)\\ 
\phantom{2\tga(\tw)}+\hskip2pt[\hs4\lambda(c,c)\hs\varOmega(\rd,\vf) 
+\varOmega(\vf,c)\hs\lambda(c,\rd)]\hs\xi(w) 
+2\hs[\hh\varOmega(\vf,c)]^2\hs\xi(w)\\ 
\phantom{2\tga(\tw)}+\hskip2pt[\hs24\hs r\hh\phi^2\mu(c,\rd) 
-2\mu(c,\rd)\hs\mu(\rd,\rd)+K\nh\phi^2\lambda(c,\rd)]\hs\ts(w)\\ 
\phantom{2\tga(\tw)}+\hskip2pt[\hh\lambda(c,c)\hs\mu(\rd,\rd) 
+2\lambda(c,\rd)\hs\mu(c,\rd) 
-12\hs r\lambda(c,c)\hs\phi^2\nh 
+2K\varOmega(\vf,c)\hs\phi^2\nh+4\hj\hh\phi^2\hh]\hs\ts(w)\\ 
\phantom{2\tga(\tw)}+\hskip2pt[\hs4\hs\phi-\lambda(c,c)\hs\lambda(c,\rd) 
-2\hh\varOmega(\vf,c)\hs\mu(c,\rd)]\hs\ts(w) 
+\lambda(c,c)\hs\varOmega(\vf,c)\hs\ts(w)\hs, 
\end{array} 
\] 
where\/ $\,a\in\plane\,$ is fixed, with\/ $\,\varOmega(a,c)=1$, and\/ 
$\,\gs,\hj\,$ are arbitrary functions defined on an open subset of\/ 
$\,\bs$, so that, as functions in\/ $\,M=\bs\times\plane_+$, they are constant 
in the\/ $\,\plane_+$ direction. 
\end{lemma} 
\begin{proof}The dependence of $\,\tga(\tw)\,$ on $\,w\,$ is assumed to be 
val\-ue\-wise and linear. Thus, by Lemma~\ref{fmtde}, a function 
$\,\tga(\tw)\,$ is a solution to (\ref{dth}.i) if and only if (\ref{ftg}) is 
constant in the $\,\plane_+$ direction, that is, equal to 
$\,[\gs\hs\xi(w)+\hj\hh\ts(w)]/2\,$ for some functions $\,\gs,\hj\,$ 
defined on an open set in $\,\bs$. On the other hand, $\,8\hs\phi^2$ times 
(\ref{ftg}) equals 
$\,4\hh[\gs\hs\xi(w)+\hj\hh\ts(w)]\hs\phi^2$ if and only if 
$\,8\hh\tga(\tw)$ is given by the formula displayed in the lemma, as one 
easily verifies using the last identity in Remark~\ref{dpsie}, for 
$\,\chi=4\hs\phi\ea$, and (\ref{ffa}), along with four equalities, justified 
below: 
\begin{equation}\label{fft} 
\begin{array}{rl} 
\mathrm{a)}\hskip0pt&4\hs\phi^2d_\tw\ea=\phi\hs d_w(4\hs\phi\ea) 
+[\hh\lambda(c,c)-2\hh\mu(c,\rd)]\hs[\lm\nh\ts(w)+\qx\hs\xi(w)]\\ 
&\phantom{4\hs\phi^2d_\tw\ea}+\hskip2pt[\hh\varOmega(\vf,c) 
-\mu(\rd,\rd)-4\hh r\phi^2\hh]\hs[\ex\ts(w)-\lp\nh\xi(w)]\hs,\\ 
\mathrm{b)}\hskip0pt&16\hs\phi^2\ea^2\nnh=[\lambda(c,X)-\mu(\rd,\rd) 
-\varOmega(\vf,c)-4\hh r\phi^2\hh]^2,\\ 
\mathrm{c)}\hskip0pt&4K\nh\phi\ea=K[\lambda(c,X)-\mu(\rd,\rd)-\varOmega(\vf,c) 
-4\hh r\phi^2\hh]\hs,\\ 
\mathrm{d)}\hskip0pt&K\nh\phi\hs\varOmega(\rd,Fw)=K[\lm\nh\ts(w) 
+\qx\hs\xi(w)]\hs\phi\hs, 
\end{array} 
\end{equation} 
and then replacing $\,\qx,\ex,\lx^\pm$ by right-hand sides of (\ref{abd}.i), 
(\ref{ffw}.i) and (\ref{ffw}.ii).

Equalities in (\ref{fft}.b,\hs c,\hs d) are obvious: (\ref{fft}.b,\hs c) 
from (\ref{ffa}), and (\ref{fft}.d) from (\ref{ffw}.iii) along with 
(\ref{duf}.i). As for (\ref{fft}.a), it follows since 
$\,\tw=w+Fw\,$ and $\,4\hs\phi^2d_{Fw}\ea=d_{\phi\hs Fw}(4\hs\phi\ea) 
-(4\hs\phi\ea)\hs d_{Fw}\phi\,$ (differentiation by parts), which, by 
(\ref{ffw}), equals 
$\,[\ex\ts(w)-\lp\nh\xi(w)]\hs(d_\rd\nh 
-\mathrm{Id})(4\hs\phi\ea)$ plus $\,[\lm\nh\ts(w) 
+\qx\hs\xi(w)]\hs d_c(4\hs\phi\ea)$, while (\ref{ffa}) 
gives $\,(d_\rd-\mathrm{Id})(4\hs\phi\ea) 
=\varOmega(\vf,c)-\mu(\rd,\rd)-4\hh r\phi^2$ (see Remark~\ref{polyn}(b)) and 
$\,d_c(4\hs\phi\ea)=\lambda(c,c)-2\mu(c,\rd)\,$ (cf.\ (\ref{duf}.iii)). 
\end{proof} 
\begin{remark}\label{eqrhs}Conditions necessary and sufficient for 
equality of the right-hand sides of the two formulae for $\,8\hh\tga(\tw)$, 
provided by Lemmas~\ref{tgatw} and~\ref{solns}, can be described as 
follows. The coefficients of both $\,\xi(w)\,$ and $\,\ts(w)\,$ 
in the two right-hand sides, when restricted to 
$\,\{y\}\times\plane_+\hs\approx\,\plane_+$ for any given $\,y\in\bs$, 
are cubic polynomial functions in $\,\plane$, and so one may proceed by 
equating their cubic, quadratic, linear and constant homogeneous components. 
\begin{enumerate} 
  \def\theenumi{{\rm\alph{enumi}}} 
\item[{\rm(a)}]Equality between the cubic or, respectively, quadratic 
homogeneous components of the coefficient of $\,\xi(w)\,$ is equivalent to 
the first or, respectively, second relation in (\ref{eqo}), as one easily sees 
using (\ref{det}) and noting that 
$\,\varOmega(\vf,c)\hs\mu(\rd,\rd)+\mu(c,\rd)\hs\varOmega(\rd,\vf) 
=\mu(\vf,\rd)\hs\phi$, in view of Remark~\ref{cclsm} and (\ref{duf}.i). 
\item[{\rm(b)}]The quadratic or, respectively, linear homogeneous components 
of the coefficient of $\,\ts(w)$ in the two right-hand sides are equal if 
and only if 
\begin{equation}\label{tmo} 
[\hh2\mu_1(c,\rd)-\lambda_2(c,\rd)]\hs\phi 
=2\lambda(c,c)\hs\mu(\rd,\rd)-2\mu(c,c)\hs\lambda(\rd,\rd) 
+2K\varOmega(\vf,c)\hs\phi^2\nh+\hj\hs\phi^2, 
\end{equation} 
or, respectively, $\,[\hh\lambda_1(c,c)+2\hh\varOmega(c,\vf_2)]\hs\phi 
=4\mu(c,c)\hs\varOmega(\rd,\vf)+4\hh\varOmega(\vf,c)\hh\mu(c,\rd)-\phi$. 
The last equation may also be rewritten as 
$\,\lambda_1(c,c)+2\hh\varOmega(c,\vf_2)=4\mu(c,\vf)-1$, since, according to 
Remark~\ref{cclsm} and (\ref{duf}.i), 
\begin{equation}\label{mcc} 
\mu(c,c)\hs\varOmega(\rd,\vf)+\varOmega(\vf,c)\hh\mu(c,\rd) 
=\mu(c,\vf)\hs\varOmega(\rd,c)=\mu(c,\vf)\hs\phi\hs. 
\end{equation} 
\item[{\rm(c)}]The corresponding equalities involving other homogeneous 
components of either coefficient function are always satisfied. 
\item[{\rm(d)}]Subtracting from (\ref{tmo}) one-half of the equation obtained 
by applying $\,d_c$ to the first equality in (\ref{eqo}), and using 
(\ref{duf}.iii), we get $\,\hj=K\varOmega(c,\vf)$. Similarly, $\,d_c$ applied 
to the second equality in (\ref{eqo}) yields 
$\,\lambda_1(c,c)+2\hh\varOmega(c,\vf_2)=4\mu(c,\vf)-1$. 
\end{enumerate} 
\end{remark} 
\begin{lemma}\label{dcomp}Under the hypotheses of Theorem\/~{\rm\ref{cndtd}}, 
condition\/ {\rm(\ref{dth}.i)} is satisfied by all sections\/ $\,w\,$ of\/ 
$\,\mathcal{H}\,$ and\/ $\,v\,$ of\/ $\,\mathcal{V}\nh$, with\/ $\,\tw=w+Fw$, 
if and only if\/ {\rm(\ref{eqo})} holds for some function\/ 
$\,\gs:U\to\bbR\hs$.

Furthermore, the relations\/ {\rm(\ref{eqo})}, for any specific function\/ 
$\,\gs$, imply that\/ $\,8\hh\tga(\tw)\,$ is given by the formula in 
Lemma\/~{\rm\ref{solns}} with this\/ $\,\gs\,$ and\/ 
$\,\hj=K\varOmega(c,\vf)$. In addition, then 
\begin{equation}\label{eqt} 
\begin{array}{l} 
[\hh2\mu_1(c,\rd)-\lambda_2(c,\rd)]\hs\phi 
=2\lambda(c,c)\hs\mu(\rd,\rd)-2\mu(c,c)\hs\lambda(\rd,\rd) 
+K\varOmega(\vf,c)\hs\phi^2,\\ 
\lambda_1(c,c)+2\hh\varOmega(c,\vf_2)=4\mu(c,\vf)-1\hs. 
\end{array} 
\end{equation} 
\end{lemma} 
\begin{proof}Condition (\ref{dth}.i) is clearly equivalent to equality of 
the right-hand sides in Lemmas~\ref{tgatw} and~\ref{solns}, for suitably 
chosen functions $\,\gs,\hj:U\to\bbR\hs$. Thus, according to 
Remark~\ref{eqrhs}, it implies (\ref{eqo}), as well as (\ref{tmo}), the second 
relation in (\ref{eqt}), and the equality $\,\hj=K\varOmega(c,\vf)$. Now 
(\ref{tmo}) with $\,\hj=K\varOmega(c,\vf)\,$ yields the first relation in 
(\ref{eqt}).

Conversely, if (\ref{eqo}) is satisfied, Remark~\ref{eqrhs}(d) shows that, 
setting $\,\hj=K\varOmega(c,\vf)\,$ we obtain all four equalities involving 
homogeneous components, mentioned in Remark~\ref{eqrhs}(a), \hs(b). Thus, the 
right-hand side in Lemma~\ref{tgatw} coincides with that in Lemma~\ref{solns}. 
\end{proof}

\section{Proof of Theorem~\ref{cndtd}, second part}\label{ptsp} 
\setcounter{equation}{0} 
In addition to the assumptions and notations adopted at the beginning of 
Section~\ref{sole}, we also assume, throughout this section, that 
condition (\ref{dth}.i) holds for all sections $\,w\,$ of 
$\,\mathcal{H}\,$ and $\,v\,$ of $\,\mathcal{V}\nh$. As before, 
$\,w,w\hh'$ will from now on stand for constant vector fields on $\,\bs$, 
while $\,\tw=w+Fw$ and $\,\tw\hh'\nnh=w\hh'\nnh+Fw\hh'\nnh$. The 
subscripts $\,(\hskip2.5pt)_j$ will again denote the partial derivatives 
relative to the af\-fine coordinates $\,y^{\hs j}$ with 
$\,d\hskip.2pty^{\nh1}\nh=\xi\,$ and $\,d\hskip.2pty^2\nh=\ts$.

By Lemma~\ref{dcomp}, $\,8\hh\tga(\tw)\,$ is given by the formula in 
Lemma~\ref{solns} with $\,\hj=K\varOmega(c,\vf)$. Explicitly, 
\begin{equation}\label{gwe} 
\tga(\tw)=d_w\ty+\gy\xi(w)+\ny\ts(w)\hs, 
\end{equation} 
where the functions $\,\gy,\ny\,$ and $\,\ty=\phi^2\ea\,$ (cf. (\ref{ffa})) 
are defined by 
\begin{equation}\label{tnp} 
\begin{array}{l} 
4\hs\phi^{-1}\ty=\lambda(c,\rd)-\mu(\rd,\rd)-\varOmega(\vf,c)-4\hs\phi^2r\hs,\\ 
8\gy=3\lambda(c,\rd)\hs\mu(\rd,\rd)-4\mu(c,\rd)\hs\lambda(\rd,\rd) 
-4K\nh\phi^2\varOmega(\rd,\vf)+12\hs r\phi^2\lambda(c,\rd)\\ 
\phantom{8\gy}+\hskip2pt2\lambda(c,c)\hs\lambda(\rd,\rd)-[\hh\lambda(c,\rd)]^2 
\nh-6\hh\varOmega(\vf,c)\hs\mu(\rd,\rd)\\ 
\phantom{8\gy}-\hskip2pt8\mu(c,\rd)\hs\varOmega(\rd,\vf) 
-24\hs r\hs\varOmega(\vf,c)\hs\phi^2-4\hs\phi\hs\varOmega(\rd,a) 
+4\gs\hh\phi^2\\ 
\phantom{8\gy}+\hskip2pt4\lambda(c,c)\hs\varOmega(\rd,\vf) 
+\varOmega(\vf,c)\hs\lambda(c,\rd) 
+2\hs[\varOmega(\vf,c)]^2,\\ 
8\ny=24\hs r\phi^2\mu(c,\rd)-2\mu(c,\rd)\hs\mu(\rd,\rd) 
+K\nh\phi^2\lambda(c,\rd)\\ 
\phantom{8\ny}+\hskip2pt\lambda(c,c)\hs\mu(\rd,\rd) 
+2\lambda(c,\rd)\hs\mu(c,\rd)-12\hs r\lambda(c,c)\hs\phi^2\nh 
-2K\varOmega(\vf,c)\hs\phi^2\\ 
\phantom{8\ny}+\hskip2pt4\hs\phi-\lambda(c,c)\hs\lambda(c,\rd) 
-2\hh\varOmega(\vf,c)\hs\mu(c,\rd)+\lambda(c,c)\hs\varOmega(\vf,c)\hs. 
\end{array} 
\end{equation} 
\begin{lemma}\label{gwwep}We have\/ 
$\,2\hs[(d\hs\tga+2\hh\ta\wedge\tga)(\tw,\tw\hh'\hh)]\hs\phi 
=[(\xi\wedge\ts)(w,w\hh'\hh)]\hh\psi\,$ for all\/ 
$\,w\,$ and\/ $\,w\hh'\nnh$, where\/ $\,\psi:U\times\plane_+\to\bbR\,$ is the 
function given by 
\begin{equation}\label{psi} 
\begin{array}{l} 
\psi\,\,=\,\,2(\ny_1-\gy_2)\hh\phi+3(\lx+4\hh r\phi^2)\hs(\ty_2+\ny\hh) 
+3(\ex+K\nh\phi^2\nnh/\hh2)\hs(\ty_1+\gy)\\ 
\phantom{\psi\,\,\hs}+\hskip2.4pt2\hh\qx\,d_c(\ty_2+\ny\hh) 
-2\hh\lm d_c(\ty_1+\gy)-2\hh\lp d_\rd(\ty_2+\ny\hh) 
-2\hh\ex\,d_\rd(\ty_1+\gy)\hs. 
\end{array} 
\end{equation} 
\end{lemma} 
\begin{proof}We begin by observing that 
\begin{equation}\label{igw} 
\begin{array}{rl} 
\mathrm{i)}\hskip0pt&\tga(\tw)=(\ty_1+\gy)\hs\xi(w)+(\ty_2+\ny\hh)\hh\ts(w)\hs, 
\\ 
\mathrm{ii)}\hskip0pt&2\hh\tga([\tw,\tw\hh'\hh])\hs\phi 
=[(\xi\wedge\ts)(w,w\hh'\hh)]\hh[(\lx+4\hh r\phi^2)\hs(\ty_2+\ny\hh) 
+(\ex+K\nh\phi^2\nnh/\hh2)\hs(\ty_1+\gy)]\hs. 
\end{array} 
\end{equation} 
Namely, (\ref{igw}.i) is immediate from (\ref{gwe}) and the last identity in 
Remark~\ref{dpsie}, for $\,\chi=\ty$. To prove (\ref{igw}.ii), let us note 
that, for any section $\,v\,$ of $\,\mathcal{V}\nh$, (\ref{fgv}.a) and 
(\ref{hvv}) yield $\,4\hh\tga(v)=h(v,\hw)$, for the unique section $\,\hw\,$ 
of $\,\mathcal{H}\,$ with $\,\xi(\hw)=\ex+K\nh\phi^2\nnh/\hh2\,$ and 
$\,\ts(\hw)=\lx+4\hh r\phi^2\nnh$. This may be applied to 
$\,v=[\tw,\tw\hh'\hh]$, as $\,[w,w\hh'\hh]=0$, and so Remark~\ref{liebr} 
implies that $\,[\tw,\tw\hh'\hh]=[w,Fw\hh'\hh]+[Fw,w\hh'\hh]+[Fw,Fw\hh'\hh]\,$ 
is a section of $\,\mathcal{V}\,$ (cf.\ (\ref{mor})). Consequently, 
$\,4\hh\tga([\tw,\tw\hh'\hh])=h([w,Fw\hh'\hh],\hw)+h([Fw,w\hh'\hh],\hw) 
+h([Fw,Fw\hh'\hh],\hw)$. Since $\,\gamma=0\,$ (see Lemma~\ref{gflat}) and 
$\,[w,w\hh'\hh]=0$, (\ref{ape}.d) for the triple $\,(\hw,w,w\hh'\nnh)\,$ 
instead of $\,(w,w\hh'\nnh,w\hh''\hh)\,$ now gives 
$\,4\hh\tga([\tw,\tw\hh'\hh])=\tga(\hw+F\hw)\hs\zeta(w,w\hh'\hh)$. Thus, 
(\ref{igw}.ii) follows from (\ref{igw}.i), with $\,\tw=w+Fw\,$ replaced by 
$\,\hw+F\hw$, and (\ref{bef}.c).

Next, $\,2\hh\phi\hs d_{Fw}[\hh\tga(\tw\hh'\hh)] 
-2\hh\phi\hs d_{Fw\hh'\hh}[\hh\tga(\tw)]\,$ is equal to 
$\,(\xi\wedge\ts)(w,w\hh'\hh)\,$ times the second line in (\ref{psi}). To see 
this, we express 
$\,\phi\hs d_{Fw}[\hh\tga(\tw\hh'\hh)]=d_{\phi\hs Fw}[\hh\tga(\tw\hh'\hh)]\,$ 
with the aid of (\ref{ffw}.iii) and (\ref{igw}.i), using the fact that 
$\,\xi(w),\ts(w),\ts(w\hh'\hh)$, etc., are constant on 
$\,M=\bs\times\plane_+$, while $\,(\xi\wedge\ts)(w,w\hh'\hh)\,$ is given by 
(\ref{bwa}.ii).

Also, $\,d_w[\hh\tga(\tw\hh'\hh)] 
-d_{w\hh'}[\hh\tga(\tw)]=(\ny_1-\gy_2)\hs[(\xi\wedge\ts)(w,w\hh'\hh)]$. In 
fact, setting $\,\chi=\tga(\tw\hh'\hh)$, we get 
$\,d_w\chi=\chi_1\xi(w\hh'\hh)+\chi_2\ts(w\hh'\hh)\,$ from Remark~\ref{dpsie}, 
while, by (\ref{igw}.i), 
$\,\chi_j\nh=(\ty_{1j}+\gy_j)\hs\xi(w)+(\ty_{2j}+\ny_j)\hh\ts(w)\,$ for 
$\,j=1,2$, and our claim follows from (\ref{bwa}.ii), as 
$\,\ty_{12}\nh=\ty_{21}$.

Finally, (\ref{fgv}.b), (\ref{igw}.i) and (\ref{bwa}.ii) imply that 
$\,[(\ta\wedge\tga)(\tw,\tw\hh'\hh)]\hs\phi\,$ equals 
$\,(\xi\wedge\ts)(w,w\hh'\hh)\,$ times 
$\,(\lx+4\hh r\phi^2)\hs(\ty_2+\ny\hh) 
+(\ex+K\nh\phi^2\nnh/\hh2)\hs(\ty_1+\gy)$.

Our assertion is now immediate from (\ref{bwa}.iii) and the above equalities, 
including (\ref{igw}.i). 
\end{proof} 
The conclusion of Theorem~\ref{cndtd} is a trivial consequence of the next 
result. In fact, of the two conditions (\ref{dth}.i) and (\ref{dth}.ii), 
together constituting (d) in Theorem~\ref{crvtc}, the first holds, according 
to Lemma~\ref{dcomp}, if and only if (\ref{eqo}) does. Under the assumption 
(\ref{dth}.i) (or, (\ref{eqo})), made throughout this section, (\ref{dth}.ii) 
amounts, by Lemma~\ref{gwwep}, to vanishing of $\,\psi\,$ in (\ref{psi}), 
while (\ref{pse}) shows that this is equivalent to (\ref{oqc}). 
\begin{theorem}\label{psecu}The function\/ $\,\psi\,$ in 
Lemma\/~{\rm\ref{gwwep}} can also be expressed as 
\begin{equation}\label{pse} 
\psi\,=\,[\hh2K\lambda(c,\vf)+8\hh r-K\varOmega(\vf_1,c)-\gs_2]\hs\phi^3. 
\end{equation} 
\end{theorem} 
\begin{proof}We need to evaluate the right-hand side of (\ref{psi}), replacing 
$\,\ty,\gy,\ny\,$ (or, $\,\qx,\ex,\lx$) with the expressions provided by 
(\ref{tnp}) (or, (\ref{abd})), and $\,\lx^\pm$ with the right-hand sides in 
(\ref{ffw}.i) -- (\ref{ffw}.ii). To make this task manageable, we note that 
the ingredients of (\ref{psi}), when restricted to 
$\,\{y\}\times\plane_+\hs\approx\,\plane_+$ for any given $\,y\in\bs$, are 
polynomial functions in $\,\plane$, of degrees which are, at most: one 
(for $\,\phi$), two (for $\,\qx,\ex,\lx,\lx^\pm$ as well as 
$\,d_c(\ty_2+\ny\hh)\,$ and $\,d_c(\ty_1+\gy)$), three (for $\,\ny_1-\gy_2$ 
along with $\,\ty_2+\ny$, $\,\ty_1+\gy\,$ and their $\,d_\rd$-im\-ages), and, 
consequently, five for $\,\psi$. Rather than describing $\,\psi\,$ directly, 
we will derive formulae for its homogeneous components 
$\,\psi\qnt\nnh,\psi\qrt\nnh,\psi\cub\nnh,\psi\qdr\nnh,\psi\lin\nnh,\psi\cst$ 
of degrees $\,5,4,3,2,1\,$ and $\,0$. The number or the resulting terms will 
be reduced, since some of them can be consolidated or eliminated due to the 
fact that 
\begin{equation}\label{trm} 
\qx\hh\cst\hs=\,(\lx+4\hh r\phi^2)\lin\hs=\,[\hs d_\rd(\hs...\hs)]\cst\hs 
=\,[\hs d_c(\hs...\hs)]\cub\hs=\,0\hs. 
\end{equation} 
where $\hs...\hs$ stands for $\,\ty_2+\ny\,$ or $\,\ty_1+\gy$, and each 
homogeneous component of $\,d_\rd(\ty_2+\ny)\,$ or $\,d_\rd(\ty_1+\gy)\,$ 
equals the degree times $\,\ty_2+\ny\,$ or $\,\ty_1+\gy\,$ (cf.\ 
Remark~\ref{polyn}(b)).

Proceeding as stated above, we will show that 
$\,\psi\qnt\nh=\psi\cst\nh=\psi\lin\nh=\psi\qrt\nh=\psi\qdr\nh=0$, while 
$\,\psi\cub$ equals the right-hand side of (\ref{pse}). Specifically, 
\[ 
\begin{array}{l} 
\psi\qnt\nh=3(\lx+4\hh r\phi^2\nh-2\hh\lp)\qdr(\ty_2+\ny\hh)\cub\nh 
-3(\ex-K\nh\phi^2\nnh/\hh2)\qdr(\ty_1+\gy)\cub,\\ 
\psi\cst\nh=3(\lx+4\hh r\phi^2)\cst(\ty_2+\ny\hh)\cst\nh 
+3(\ex+K\nh\phi^2\nnh/\hh2)\cst(\ty_1+\gy)\cst\nh 
-(2\hh\lm)\cst d_c[(\ty_1+\gy)\lin\hh]\hs,\\ 
\psi\lin\nh=2(\ny_1-\gy_2)\cst\nh\phi 
+(3\hh\lx+12\hh r\phi^2\nh-2\hh\lp)\cst(\ty_2+\ny\hh)\lin\nh 
+(\ex+3K\nh\phi^2\nnh/\hh2)\cst(\ty_1+\gy)\lin\\ 
\phantom{\psi\lin\nh}+\hskip2pt3(\ex 
+K\nh\phi^2\nnh/\hh2)\lin(\ty_1+\gy)\cst\nh 
+2\hh\qx\hh\lin d_c[(\ty_2+\ny\hh)\lin\hh]\\ 
\phantom{\psi\lin\nh}-\hskip2pt(2\hh\lm)\lin d_c[(\ty_1+\gy)\lin\hh] 
-(2\hh\lm)\cst d_c[(\ty_1+\gy)\qdr\hh]\hs,\\ 
\psi\qrt\nh=2(\ny_1-\gy_2)\cub\nnh\phi 
+(3\hh\lx+12\hh r\phi^2\nh-6\hh\lp)\lin(\ty_2+\ny\hh)\cub\\ 
\phantom{\psi\qrt\nh}+\hskip2pt 
(3\hh\lx+12\hh r\phi^2\nh-4\hh\lp)\qdr(\ty_2+\ny\hh)\qdr\\ 
\phantom{\psi\qrt\nh}-\hskip2pt3(\ex-K\nh\phi^2\nnh/\hh2)\lin(\ty_1+\gy)\cub\nh 
-(\ex-3K\nh\phi^2\nnh/\hh2)\qdr(\ty_1+\gy)\qdr\\ 
\phantom{\psi\qrt\nh}+\hskip2pt2\hh\qx\hh\qdr d_c[(\ty_2+\ny\hh)\cub\hh] 
-(2\hh\lm)\qdr d_c[(\ty_1+\gy)\cub\hh]\hs. 
\end{array} 
\] 
As $\,(\lx+4\hh r\phi^2\nh-2\hh\lp)\qdr\nh=(\ex-K\nh\phi^2\nnh/\hh2)\qdr\nh 
=0$, we thus have $\,\psi\qnt\nh=0$. To conclude that 
$\,\psi\cst\nnh,\psi\lin$ and $\,\psi\qrt$ vanish as well, we use the relations 
\[ 
\begin{array}{l} 
8\hh(\ty_2+\ny\hh)\cst\nh=\lambda(c,c)\hs\varOmega(\vf,c)\hs,\hskip12pt 
4\hh(\ty_1+\gy)\cst\nh=[\varOmega(\vf,c)]^2,\\ 
8\hh(\ty_2+\ny\hh)\lin\nh=2\hh\varOmega(c,\vf_2)\hs\phi+4\hs\phi 
-2\hh\varOmega(\vf,c)\hs\mu(c,\rd)-\lambda(c,c)\hs\lambda(c,\rd)\hs,\\ 
8\hh(\ty_1+\gy)\lin\nh=2\hh\varOmega(c,\vf_1)\hs\phi 
+\varOmega(\vf,c)\hs\lambda(c,\rd) 
+4\hh\lambda(c,c)\hs\varOmega(\rd,\vf)\hs,\\ 
8\hh(\ty_2+\ny\hh)\qdr\nnh=2\lambda_2(c,\rd)\hs\phi\nh 
+\nnh\lambda(c,c)\hs\mu(\rd,\rd)\nh+\nh2\lambda(c,\rd)\hs\mu(c,\rd) 
\nnh-\hskip-2pt12\hs r\lambda(c,c)\hs\phi^2\nnh 
-\nh2K\varOmega(\vf,c)\hs\phi^2\nh,\\ 
8\hh(\ty_1+\gy)\qdr\nh=2\lambda_1(c,\rd)\hs\phi 
+2\lambda(c,c)\hs\lambda(\rd,\rd)-[\hh\lambda(c,\rd)]^2 
\nh-6\hh\varOmega(\vf,c)\hs\mu(\rd,\rd)\\ 
\phantom{8\hh(\ty_1+\gy)\qdr\nh}-\hskip2pt8\mu(c,\rd)\hs\varOmega(\rd,\vf) 
-24\hs r\hs\varOmega(\vf,c)\hs\phi^2-4\hs\phi\hs\varOmega(\rd,a) 
+4\gs\hh\phi^2,\\ 
8\hh(\ty_2+\ny\hh)\cub\nh=-2\hh\phi\hh\mu_2(\rd,\rd)-8\hs r_2\phi^3\nh 
+24\hs r\phi^2\mu(c,\rd)-2\mu(c,\rd)\hs\mu(\rd,\rd) 
+K\nh\phi^2\lambda(c,\rd)\hs,\\ 
8\hh(\ty_1+\gy)\cub\nh=-2\hh\phi\hh\mu_1(\rd,\rd)-8\hs r_1\phi^3\nh 
+3\lambda(c,\rd)\hs\mu(\rd,\rd)-4\mu(c,\rd)\hs\lambda(\rd,\rd)\\ 
\phantom{8\hh(\ty_1+\gy)\cub\nh}-\hskip2pt4K\nh\phi^2\varOmega(\rd,\vf) 
+12\hs r\phi^2\lambda(c,\rd)\hs. 
\end{array} 
\] 
Since $\,\psi\cst\nh=-\hh3\hh\varOmega(\vf,c)\hs(\ty_2+\ny\hh)\cst\nh 
+3\lambda(c,c)\hs(\ty_1+\gy)\cst\nh 
+\varOmega(\vf,c)\hs d_c[(\ty_1+\gy)\lin\hh]$, and, by (\ref{duf}.iii), 
$\,8\hs d_c[(\ty_1+\gy)\lin\hh]=-\hh3\hh\varOmega(\vf,c)\hs\lambda(c,c)$, we 
get $\,\psi\cst\nh=0$. Also, one easily verifies that 
\[ 
\begin{array}{l} 
8\hh(\ny_1-\gy_2)\cst\nh=\varOmega(\vf,c)\hs\lambda_1(c,c) 
+\lambda(c,c)\hs\varOmega(\vf_1,c)-4\hh\varOmega(\vf,c)\hs\varOmega(\vf_2,c) 
\hs,\\ 
8\hh(\ny_1-\gy_2)\lin\nh=2\mu(c,\rd)\hs\varOmega(c,\vf_1) 
-\lambda(c,c)\hs\lambda_1(c,\rd)-\lambda(c,\rd)\hs\lambda_1(c,c)\\ 
\phantom{8\hh(\ny_1-\gy_2)\lin\nh}-\hskip2pt2\hh\varOmega(\vf,c)\hs\mu_1(c,\rd) 
-4\lambda(c,c)\hs\varOmega(\rd,\vf_2)-4\hh\varOmega(\rd,\vf)\hs\lambda_2(c,c)\\ 
\phantom{8\hh(\ny_1-\gy_2)\lin\nh}+\hskip2pt\lambda(c,\rd)\hs\varOmega(c,\vf_2) 
-\varOmega(\vf,c)\hs\lambda_2(c,\rd)\hs,\\ 
8\hh(\ny_1-\gy_2)\qdr\nh=\lambda(c,c)\hs\mu_1(\rd,\rd) 
+\mu(\rd,\rd)\hs\lambda_1(c,c)+2\lambda(c,\rd)\hs\mu_1(c,\rd) 
+2\mu(c,\rd)\hs\lambda_1(c,\rd)\\ 
\phantom{8\hh(\ny_1-\gy_2)\qdr\nh}-\hskip2pt12\lambda(c,c)\hs\phi^2r_1\nh 
-12\hs r\phi^2\lambda_1(c,c)-2K\phi^2\varOmega(\vf_1,c) 
-2\lambda(c,c)\hs\lambda_2(\rd,\rd)-4\hs\phi^2\gs_2\\ 
\phantom{8\hh(\ny_1-\gy_2)\qdr\nh}-\hskip2pt2\lambda(\rd\hn,\nh\rd)\hs\lambda_2(c,c) 
\nh+\nh2\lambda(c,\nnh\rd)\hs\lambda_2(c,\nnh\rd)\nh 
+\nh6\hh\varOmega(\vf,c)\hs\mu_2(\rd,\nnh\rd)\nh 
+\nh6\mu(\rd\hn,\nh\rd)\hs\varOmega(\vf_2,c)\\ 
\phantom{8\hh(\ny_1-\gy_2)\qdr\nh}+\hskip2pt 
8\mu(c,\rd)\hs\varOmega(\rd,\vf_2)+8\hh\varOmega(\rd,\vf)\hs\mu_2(c,\rd) 
+24\hh\varOmega(\vf,c)\hs\phi^2r_2\nh+24\hs r\phi^2\varOmega(\vf_2,c)\hs,\\ 
8\hh(\ny_1-\gy_2)\cub\nh=24\hs r\phi^2\mu_1(c,\rd) 
+24\hs\phi^2\mu(c,\rd)\hs r_1-2\mu(c,\rd)\hs\mu_1(\rd,\rd)\\ 
\phantom{8\hh(\ny_1-\gy_2)\cub\nh}-\hskip2pt2\mu(\rd,\rd)\hs\mu_1(c,\rd) 
+K\nh\phi^2\lambda_1(c,\rd)-3\lambda(c,\rd)\hs\mu_2(\rd,\rd)\\ 
\phantom{8\hh(\ny_1-\gy_2)\cub\nh}-\hskip2pt3\mu(\rd,\rd)\hs\lambda_2(c,\rd) 
+4\mu(c,\rd)\hs\lambda_2(\rd,\rd)+4\lambda(\rd,\rd)\hs\mu_2(c,\rd)\\ 
\phantom{8\hh(\ny_1-\gy_2)\cub\nh}+\hskip2pt4K\nh\phi^2\varOmega(\rd,\vf_2) 
-12\hs r\phi^2\lambda_2(c,\rd)-12\hs\phi^2\lambda(c,\rd)\hs r_2\hh. 
\\ 
\end{array} 
\] 
As $\,(3\hh\lx+12\hh r\phi^2\nh-2\hh\lp)\cst\nh=0$, the above formula for 
$\,\psi\lin$ gives 
\[ 
\begin{array}{l} 
\psi\lin\nh=2(\ny_1-\gy_2)\cst\nh\phi 
+\lambda(c,c)\hs(\ty_1+\gy)\lin\nh-6\mu(c,\rd)\hs(\ty_1+\gy)\cst\\ 
\phantom{\psi\lin\nh}+\hskip2pt4\hh\varOmega(\rd,\vf)\hs d_c[(\ty_2 
+\ny\hh)\lin\hh]+\lambda(c,\rd)\hs d_c[(\ty_1+\gy)\lin\hh] 
+\varOmega(\vf,c)\hs d_c[(\ty_1+\gy)\qdr\hh]\hs, 
\end{array} 
\] 
so that $\,4\hh\psi\lin\nh=\phi\hs\varOmega(\vf,c)\hs 
[\hh\lambda_1(c,c)+2\hh\varOmega(c,\vf_2)+1] 
-4\hh\varOmega(\vf,c)\hs 
[\hh\mu(c,c)\hs\varOmega(\rd,\vf)+\varOmega(\vf,c)\hh\mu(c,\rd)]$. 
Thus, $\,4\hh\psi\lin\nh=\phi\hs\varOmega(\vf,c)\hs 
[\hh\lambda_1(c,c)+2\hh\varOmega(c,\vf_2)-4\mu(c,\vf)+1\hh]$, since 
$\,\mu(c,c)\hs\varOmega(\rd,\vf)+\varOmega(\vf,c)\hh\mu(c,\rd)=\mu(c,\vf)\,$ 
(see (\ref{mcc})).The second equation in (\ref{eqt}) now yields 
$\,\psi\lin\nh=0$. Similarly, 
\[ 
\begin{array}{l} 
\psi\qrt\nh=2(\ny_1-\gy_2)\cub\nh\phi 
-3\lambda(c,\rd)\hs(\ty_2+\ny\hh)\cub\nh 
+[\hh\mu(\rd,\rd)+4\hh r\phi^2\hh]\hs(\ty_2+\ny\hh)\qdr\\ 
\phantom{\psi\qrt\nh}+\hskip2pt6\mu(c,\rd)\hs(\ty_1+\gy)\cub\nh 
+K\nh\phi^2(\ty_1+\gy)\qdr\\ 
\phantom{\psi\qrt\nh}+\hskip2pt2\lambda(\rd,\rd)\hs d_c[(\ty_2+\ny\hh)\cub\hh] 
+[\hh4\hh r\phi^2\nh-3\mu(\rd,\rd)]\hs d_c[(\ty_1+\gy)\cub\hh]\hs, 
\end{array} 
\] 
and so, by (\ref{det}) and (\ref{duf}.iii), $\,4\hh\psi\qrt$ equals 
\[ 
\begin{array}{l} 
[\hh2\mu(\rd\hn,\nh\rd)\nnh+\nnh8\hh r\phi^2\hh] 
\hs 
[\hh2\hh\phi\hh\mu_1(c,\nnh\rd)\nh-\nh\phi\hh\lambda_2(c,\nnh\rd)\nh 
-\nh2\lambda(c,\nh c)\hs\mu(\rd\hn,\nh\rd)\nh 
+\nh2\mu(c,\nh c)\hs\lambda(\rd\hn,\nh\rd)\nh 
-\nnh K\nh\varOmega(\vf,\nh c)\hs\phi^2\hh]\\ 
\hskip13.5pt-\hskip2pt4\mu(c,\nnh\rd)\hs 
[\hh2\hh\phi\hh\mu_1(\rd\hskip-1.7pt,\nnh\rd)\nh 
-\nh\phi\hh\lambda_2(\rd\hskip-1.7pt,\nnh\rd)\nh 
-\nh4\lambda(c,\nnh\rd)\hh\mu(\rd\hskip-1.7pt,\nnh\rd)\nh 
+\nh4\mu(c,\nnh\rd)\hh\lambda(\rd\hskip-1.7pt,\nnh\rd)\nh 
+\nh2K\nh\phi^2\nh\varOmega(\rd\hskip-1.7pt,\vf)]\\ 
\hskip13.5pt+\hskip2pt2K\nh\phi^3[\hh\lambda_1(c,\rd)+2\hh\varOmega(\rd,\vf_2) 
-\varOmega(\rd,a)+\phi\hskip1.6pt\mathrm{det}_\varOmega\lambda+\gs\hh\phi\hh]\\ 
\hskip13.5pt-\hskip2pt8K\nh\phi^2\varOmega(\rd,\vf)\hs\mu(c,\rd) 
-8K\nh\phi^2\varOmega(\vf,c)\hs\mu(\rd,\rd)\hs. 
\end{array} 
\] 
According to Lemma~\ref{dcomp}, we have (\ref{eqt}) and (\ref{eqo}). Thus, in 
the above expression, the first two lines vanish, while the third is equal to 
$\,8K\nh\phi^3\mu(\vf,\rd)$. As Remark~\ref{cclsm} and (\ref{duf}.i) give 
\begin{equation}\label{fmq} 
\phi\hh\mu(\vf,\rd)=\varOmega(\rd,\vf)\hs\mu(c,\rd) 
+\varOmega(\vf,c)\hs\mu(\rd,\rd)\hs, 
\end{equation} 
we see that $\,\psi\qrt\nh=0$. Next, as before, by (\ref{trm}), 
$\,\psi\qdr$ equals $\,2(\ny_1-\gy_2)\lin\nh\phi\,$ plus 
\[ 
\begin{array}{l} 
(3\hh\lx+12\hh r\phi^2\nh-4\hh\lp)\cst(\ty_2+\ny\hh)\qdr\nh 
-(2\hh\lp)\lin(\ty_2+\ny\hh)\lin\nh 
+3(\hh\lx+4\hh r\phi^2)\qdr(\ty_2+\ny\hh)\cst\\ 
\phantom{\psi\qdr\nh}-\hskip1.8pt(\ex-3K\nh\phi^2\nnh/\hh2)\cst(\ty_1+\gy)\qdr 
\nh+(\ex+3K\nh\phi^2\nnh/\hh2)\lin(\ty_1+\gy)\lin\\ 
\phantom{\psi\qdr\nh}+\hskip1.8pt3(\ex+K\nh\phi^2\nnh/\hh2)\qdr(\ty_1+\gy)\cst 
\nh+2\hh\qx\hh\lin d_c[(\ty_2+\ny\hh)\qdr\hh] 
+2\hh\qx\hh\qdr d_c[(\ty_2+\ny\hh)\lin\hh]\\ 
\phantom{\psi\qdr\nh}-\hskip1.8pt(2\hh\lm)\cst d_c[(\ty_1+\gy)\cub\hh] 
-(2\hh\lm)\lin d_c[(\ty_1+\gy)\qdr\hh]-(2\hh\lm)\qdr d_c[(\ty_1+\gy)\lin\hh] 
\hs, 
\end{array} 
\] 
which can easily be rewritten as 
\[ 
\begin{array}{l} 
\psi\qdr\nh=2(\ny_1-\gy_2)\lin\nh\phi 
+3\hh\varOmega(\vf,c)\hs(\ty_2+\ny\hh)\qdr\nh 
-\lambda(c,\rd)\hs(\ty_2+\ny\hh)\lin\\ 
\phantom{\psi\qdr\nh}+\hskip2pt 
[\hh3\mu(\rd,\rd)+12\hh r\phi^2\hh]\hs(\ty_2+\ny\hh)\cst\nh 
-\lambda(c,c)\hs(\ty_1+\gy)\qdr\nh-2\mu(c,\rd)\hs(\ty_1+\gy)\lin\\ 
\phantom{\psi\qdr\nh}+\hskip2pt3K\nh\phi^2(\ty_1+\gy)\cst\nh 
+4\hh\varOmega(\rd,\vf)\hs d_c[(\ty_2+\ny\hh)\qdr\hh] 
+2\lambda(\rd,\rd)\hs d_c[(\ty_2+\ny\hh)\lin\hh]\\ 
\phantom{\psi\qdr\nh}+\hskip2pt\varOmega(\vf,c)\hs d_c[(\ty_1+\gy)\cub\hh] 
+\lambda(c,\rd)\hs d_c[(\ty_1+\gy)\qdr\hh] 
+[\hh4\hh r\phi^2\nh-3\mu(\rd,\rd)]\hs d_c[(\ty_1+\gy)\lin\hh]\hs. 
\end{array} 
\] 
Since $\,\varOmega(a,c)=1$, we thus obtain 
\[ 
\begin{array}{l} 
4\hh\psi\qdr\nh=\varOmega(c,\vf)\hs 
[\hh2\hh\phi\hh\mu_1(c,\nnh\rd)\nh-\nh\phi\hh\lambda_2(c,\nnh\rd)\nh 
-\nh2\lambda(c,\nh c)\hs\mu(\rd\hn,\nh\rd)\nh 
+\nh2\mu(c,\nh c)\hs\lambda(\rd\hn,\nh\rd)\nh 
-\nnh K\nh\varOmega(\vf,\nh c)\hs\phi^2\hh]\\ 
\phantom{4\hh\psi\qdr\nh}-\hskip2pt2\hs\phi\hh\lambda(c,c)\hs 
[\hh\lambda_1(c,\rd)+2\hh\varOmega(\rd,\vf_2) 
-\varOmega(\rd,a)+\phi\hskip1.6pt\mathrm{det}_\varOmega\lambda+\gs\hh\phi\hh]\\ 
\phantom{4\hh\psi\qdr\nh}+\hskip2pt 
8\lambda(c,c)\hs[\hh\varOmega(\rd,\vf)\hs\mu(c,\rd) 
+\varOmega(\vf,c)\hs\mu(\rd,\rd)]\hs, 
\end{array} 
\] 
where we have used (\ref{det}) and formulae displayed earlier in the proof. 
Again, as Lemma~\ref{dcomp} yields (\ref{eqt}) and (\ref{eqo}), in the above 
equality the right-hand side of the first line vanishes, and the second line 
equals $\,-\hh8\hs\phi\hh\lambda(c,c)\hs\mu(\vf,\rd)$. Consequently, by 
(\ref{fmq}), $\,\psi\qdr\nh=0$.

Finally, repeating the same steps for $\,\psi\cub\nnh$, we verify that 
$\,\psi\cub$ equals 
\[ 
\begin{array}{l} 
2(\ny_1-\gy_2)\qdr\nh\phi 
+3(\lx+4\hh r\phi^2\nh-2\hh\lp)\cst(\ty_2+\ny\hh)\cub\nh 
+(3\hh\lx+12\hh r\phi^2\nh-4\hh\lp)\lin(\ty_2+\ny\hh)\qdr\\ 
\phantom{\psi\cub\nh}+\hskip2pt 
(3\hh\lx+12\hh r\phi^2\nh-2\hh\lp)\qdr(\ty_2+\ny\hh)\lin\nh 
-3(\ex-K\nh\phi^2\nnh/\hh2)\cst(\ty_1+\gy)\cub\\ 
\phantom{\psi\cub\nh}-\hskip2pt(\ex-3K\nh\phi^2\nnh/\hh2)\lin(\ty_1+\gy)\qdr 
\nh+(\ex+3K\nh\phi^2\nnh/\hh2)\qdr(\ty_1+\gy)\lin\nh 
+2\hh\qx\hh\lin d_c[(\ty_2+\ny\hh)\cub\hh]\\ 
\phantom{\psi\cub\nh}+\hskip2pt 
2\hh\qx\hh\qdr d_c[(\ty_2+\ny\hh)\qdr\hh] 
-(2\hh\lm)\lin d_c[(\ty_1+\gy)\cub\hh] 
-(2\hh\lm)\qdr d_c[(\ty_1+\gy)\qdr\hh]\hs, 
\end{array} 
\] 
which we can again rewrite as 
\[ 
\begin{array}{l} 
\psi\cub\nh=2(\ny_1-\gy_2)\qdr\nh\phi 
+6\hh\varOmega(\vf,c)\hs(\ty_2+\ny\hh)\cub\nh 
-2\lambda(c,\rd)\hs(\ty_2+\ny\hh)\qdr\\ 
\phantom{\psi\cub\nh}+\hskip2pt 
[\hh2\mu(\rd,\rd)+8\hh r\phi^2\hh]\hs(\ty_2+\ny\hh)\lin\nh 
-3\lambda(c,c)\hs(\ty_1+\gy)\cub\nh 
+2\mu(c,\rd)\hs(\ty_1+\gy)\qdr\\ 
\phantom{\psi\cub\nh}+\hskip2pt2K\nh\phi^2(\ty_1+\gy)\lin\nh 
+4\hh\varOmega(\rd,\vf)\hs d_c[(\ty_2+\ny\hh)\cub\hh] 
+2\lambda(\rd,\rd)\hs d_c[(\ty_2+\ny\hh)\qdr\hh]\\ 
\phantom{\psi\cub\nh}+\hskip2pt\lambda(c,\rd)\hs d_c[(\ty_1+\gy)\cub\hh] 
+[\hh4\hh r\phi^2\nh-3\mu(\rd,\rd)]\hs d_c[(\ty_1+\gy)\qdr\hh]\hs. 
\end{array} 
\] 
Consequently, $\,4\hh\psi\cub$ is equal to 
\[ 
\begin{array}{l} 
2\lambda(c,c)\hs 
[\hh2\hh\phi\hh\mu_1(\rd,\rd)-\phi\hh\lambda_2(\rd,\rd) 
-4\lambda(c,\nnh\rd)\hh\mu(\rd,\rd) 
+4\mu(c,\nnh\rd)\hh\lambda(\rd,\rd) 
+2K\nh\phi^2\nh\varOmega(\rd\hskip-1.7pt,\vf)]\\ 
\hskip13.5pt+\hskip2pt4\hs[\hh2\mu(\nh\rd\hn,\nh\rd)\nnh+\nnh8\hh r\phi^2\hh]\hs 
[\mu(c,c)\hs\varOmega(\rd,\vf)+\varOmega(\vf,c)\hh\mu(c,\rd)]\\ 
\hskip13.5pt-\hskip2pt[\hh2\phi\hh\mu(\nh\rd\hn,\nh\rd)\nnh+\nnh8\hh r\phi^3\hh] 
\hs[\hh\lambda_1(c,c)+2\hh\varOmega(c,\vf_2)+1\hh]\\ 
\hskip13.5pt+\hskip2pt4\hh\phi\hh\mu(c,\nnh\rd)\hs 
[\hh\lambda_1(c,\rd)+2\hh\varOmega(\rd,\vf_2) 
-\varOmega(\rd,a)+\phi\hskip1.6pt\mathrm{det}_\varOmega\lambda+\gs\hh\phi\hh]\\ 
\hskip13.5pt-\hskip2pt16\hh\mu(c,\rd)\hs[\hh\mu(c,\rd)\hs\varOmega(\rd,\vf) 
+\mu(\rd,\rd)\hs\varOmega(\vf,c)]\\ 
\hskip13.5pt+\hskip2pt4\hs[\hh8\hh r-K\varOmega(\vf_1,c)-\gs_2]\hs\phi^3\nh 
+8K\nh\phi^2[\hh\varOmega(\vf,c)\hs\lambda(c,\rd) 
+\lambda(c,c)\hs\varOmega(\rd,\vf)]\hs. 
\end{array} 
\] 
Each of the six lines forming the above expression can now be evaluated as 
follows. The first line vanishes due to (\ref{eqo}). The second and third 
lines add up to zero as a consequence of (\ref{mcc}) and (\ref{eqn}). 
Similarly, the fourth and fifth lines cancel each other in view of (\ref{eqo}) 
and (\ref{lcv}) for $\,v=\vf$. Next, again by (\ref{lcv}) for $\,v=\vf$, the 
last line equals 
$\,4\hs[\hh2K\lambda(c,\vf)+8\hh r-K\varOmega(\vf_1,c)-\gs_2]\hs\phi^3\nnh$, 
which completes the proof. 
\end{proof}

\section{The lo\-cal-struc\-ture theorem}\label{lsth} 
\setcounter{equation}{0} 
Given a tw\hbox{o\hh-\nh}\hskip0ptplane system 
$\,(\bs,\xi,\ts,\plane,c,\varOmega)\,$ (see Section~\ref{afsy}) and a real 
constant $\,K$, let us choose 
$\,(M,\mathcal{V}\nh,\hh\mathrm{D},h,\alpha,\beta,\theta,\zeta)\,$ and 
$\,\mathcal{H}\,$ as in Lemma~\ref{gflat}. If 
$\,(\vf,\lambda,\mu):U\nh\to V\hs$ is a solution to (\ref{eqn}), defined on 
a nonempty connected open set $\,\,U\subset\bs$, and $\,(\vf,\lambda,\mu,r)\,$ 
corresponds to $\,(\vf,\lambda,\mu)\,$ in the sense of Theorem~\ref{equiv}, 
then, setting $\,\thor=\mathcal{H}+F\,$ for the section 
$\,F=F^K\nh+F^\vf\nh+F^\lambda\nh+F^\mu\nh+f\zeta$ of $\,\mathcal{F}\nnh$, 
with $\,F^K\nnh,F^\vf\nnh,F^\lambda\nnh,F^\mu$ and $\,f\,$ 
given by the formulae in Lemmas~\ref{eight} --~\ref{dcfeq}, we obtain a new 
horizontal distribution $\,\thor\,$ on $\,\,U\times\plane_+$, which gives rise 
to the neu\-tral-sig\-na\-ture metric $\,\tg$ on the 
\hbox{four\hh-}\hskip0ptman\-i\-fold $\,\,U\times\plane_+$ characterized by 
(\ref{gxh}) (for $\,\thor\,$ rather than $\,\mathcal{H}$). 
\begin{theorem}\label{lcstr}Suppose that\/ 
$\,(\bs,\xi,\ts,\plane,c,\varOmega)\,$ is a fixed 
tw\hbox{o\hh-\nh}\hskip0ptplane system and\/ $\,K\in\bbR\hs$. 
\begin{enumerate} 
  \def\theenumi{{\rm\alph{enumi}}} 
\item[{\rm(a)}] For any nonempty connected open set\/ $\,\,U\subset\bs\,$ and 
any\/ $\,(\vf,\lambda,\mu):U\nh\to V\hs$ with\/ {\rm(\ref{eqn})}, a suitable 
orientation of\/ $\,\,U\times\plane_+$ makes\/ $\,\tg\,$ constructed above a 
strictly \itnw\ self-du\-al neutral Ein\-stein metric of Pe\-trov type\/ {\rm 
III} with the scalar curvature\/ $\,12\hh K$. 
\item[{\rm(b)}] Conversely, in any strictly \itnw\ self-du\-al neutral 
Ein\-stein \hbox{four\hh-}\hskip0ptman\-i\-fold\/ $\,(M,g)\,$ of Pe\-trov 
type\/ {\rm III} with the scalar curvature\/ $\,12\hh K$, every point has a 
neighborhood isometric to an open submanifold of\/ 
$\,(U\times\plane_+\hs,\hs\tg)$, where\/ $\,\tg\,$ is obtained as above using 
our fixed\/ $\,(\bs,\xi,\ts,\plane,c,\varOmega)\,$ and\/ $\,K\,$ along with 
some solution\/ $\,(\vf,\lambda,\mu):U\nh\to V\hs$ to\/ {\rm(\ref{eqn})}. 
\end{enumerate} 
\end{theorem} 
\begin{proof}As $\,(\vf,\lambda,\mu,r)\,$ satisfies (\ref{eqo}) -- (\ref{oqc}) 
(see Theorem~\ref{equiv}), with some $\,\gs:U\to\bbR$ and $\,a\in\plane\,$ 
such that $\,\varOmega(a,c)=1$, Theorem~\ref{cndtd} implies condition (d) in 
Theorem~\ref{crvtc} for $\,\thor=\mathcal{H}+F\nnh$, while conditions (a) -- 
(c) in Theorem~\ref{crvtc}, for $\,\thor=\mathcal{H}+F\nnh$, are guaranteed by 
our use of the formulae in Lemmas~\ref{eight} --~\ref{dcfeq} to define $\,F\,$ 
(cf.\ Section~\ref{cabl}). Assertion (a) now follows from Theorem~\ref{crvtc}. 
(That $\,\tg\,$ represents the strictly \nw\ case is clear as $\,\beta\ne0$ 
everywhere, according to (\ref{bef}.a).)

Conversely, for $\,(M,g)\,$ as in (b), 
$\,(M,\mathcal{V}\nh,\hh\mathrm{D},h,\alpha,\beta,\theta,\zeta)\,$ used in the 
above construction may be assumed, by Theorem~\ref{unimo}, to coincide with 
the basic octuple that $\,(M,g)\,$ gives rise to in Theorem~\ref{chboc}. Under 
this identification $\,g\,$ corresponds to a metric $\,\tg\,$ on a connected 
open subset of $\,\bs\times\plane_+$ defined as in (\ref{gxh}), that is, 
forming the unique to\-tal-met\-ric extension of $\,h\,$ such that $\,\thor\,$ 
is $\,\tg$-null, for $\,\thor\,$ associated with $\,\tg\,$ as in 
Lemma~\ref{sdept}. Denoting by $\,\mathcal{H}\,$ the horizontal distribution 
of Lemma~\ref{gflat}, we have $\,\thor=\mathcal{H}+F\,$ for some section 
$\,F\,$ of $\,\mathcal{F}$. (See Section~\ref{dfhd}.) Conditions (a) -- (d) in 
Theorem~\ref{crvtc} thus are satisfied by $\,\thor\,$ and $\,\tg$. The 
discussion at the beginning of Section~\ref{esol}, combined with 
Lemmas~\ref{eight} and~\ref{dcfeq}, now shows that 
$\,F=F^K\nh+F^\vf\nh+F^\lambda\nh+F^\mu\nh+f\zeta$, with 
$\,(\vf,\lambda,\mu,r)\,$ (depending on $\,y\in\bs$) and $\,f\,$ as in 
Lemmas~\ref{eight} and~\ref{dcfeq}. Theorems~\ref{cndtd} and~\ref{equiv} now 
yield (\ref{eqn}), completing the proof. 
\end{proof} 
Since the manifolds described in Theorem~\ref{lcstr} are all curvature 
homogeneous (Remark~\ref{cvhmg}), it is natural to ask if they must also be 
locally homogeneous. As the next result explicitly shows, this is not the 
case, and not all such manifolds are Ric\-\hbox{ci\hh-}\hskip0ptflat. The 
latter conclusion answers a question raised by 
\hbox{D\'\i az\hh-}\hskip0ptRamos, Gar\-\hbox{c\'\i a\hs-}\hskip0ptR\'\i o and 
V\nnh\'az\-\hbox{quez\hh-}\hskip0ptLo\-ren\-zo 
\cite[Remark~3.5]{diaz-ramos-garcia-rio-vazquez-lorenzo}. 
\begin{theorem}\label{rfnrf}For any given\/ $\,K\in\bbR\hs$, applying 
Theorem~\/{\rm\ref{lcstr}(a)} to $\,(\vf,\lambda,\mu):\bs\to V\hs$ chosen 
as in Example\/~{\rm\ref{lccne}} we obtain a strictly \itnw\ self-du\-al 
neutral Ein\-stein \hbox{four\hh-}\hskip0ptman\-i\-fold of Pe\-trov type\/ 
{\rm III} with the scalar curvature\/ $\,12\hh K$, which is not locally 
homogeneous. 
\end{theorem} 
\begin{proof}In view of Theorem~\ref{lcstr}(a), we only need to verify that 
the resulting metric $\,\tg\,$ on $\,M=\bs\times\plane_+$ is not locally 
homogeneous. This will clearly follow once we show that the function 
$\,\tga(\bu):M\to\bbR\hs$, constituting a local invariant of $\,\tg\,$ (see 
the end of Remark~\ref{tsmgm}), is nonconstant on $\,\{y\}\times\plane_+$ for 
some $\,y\in\bs$.

By (\ref{fgv}.a) for $\,v=\bu$, (\ref{xwb}.d) and (\ref{duf}.i), 
$\,\tga(\bu)=K+[\hh\lambda(c,c)-2\mu(c,\rd)]\hs\phi^{-2}\nnh$. 
As $\,\lambda(c,c),\mu(c,\rd)\,$ and $\,\phi^2$ restricted to 
$\,\{y\}\times\plane_+\hs\approx\,\plane_+$ are homogeneous polynomial 
functions of degrees $\,0,1\,$ and $\,2$, the restriction of $\,\tga(\bu)\,$ 
to $\,\{y\}\times\plane_+$ is nonconstant on $\,\{y\}\times\plane_+$ for every 
$\,y\in\bs\,$ at which $\,\lambda(c,c)\ne0$. However, our choice of 
$\,(\vf,\lambda,\mu)\,$ gives $\,\lambda(c,c)\ne0\,$ at all $\,y$ lying 
outside a specific line in the plane $\,\bs$, which completes the proof. 
\end{proof} 
In Theorem~\ref{rfnrf} we used some solutions of (\ref{eqn}) to obtain 
examples of metrics with interesting geometric properties. A description 
of all solutions to (\ref{eqn}) will be given in Section~\ref{stsy}.

\section{The method of characteristics}\label{tmoc} 
\setcounter{equation}{0} 
Given an open set $\,\,U\subset\rto$ with the Cartesian coordinates 
$\,y^{\nh1}\nnh,y^2$ and an open interval $\,I\subset\bbR\hs$, let 
$\,z:U\nh\to I\hs$ be the unknown function in a first-or\-der 
qua\-si-lin\-e\-ar equation 
\begin{equation}\label{rzo} 
\rho\hs z_1+\hs\sigma\nh z_2\hs=\,\chi\hh,\hskip5pt\mathrm{where}\hskip4pt 
z_j\nh=\partial z/\partial y^{\hs j}\hskip5pt\mathrm{and}\hskip5pt 
\rho,\sigma,\chi\hskip5pt\mathrm{are\ functions\ of}\hskip5pt 
(y^{\nh1}\nnh,y^2\nnh,z)\in U\times I\nh. 
\end{equation} 
One solves (\ref{rzo}) using the method of characteristics, based on the 
observation that a function $\,z\,$ of the variables 
$\,(y^{\nh1}\nnh,y^2)$, defined on $\,\,U\,$ and valued in $\,I\nh$, 
satisfies (\ref{rzo}) if and only if its graph surface in 
$\,\,U\times I\hs$ is a union of integral curves of the vector field 
$\,(\rho,\sigma,\chi)$.

Thus, if $\,y\in \Upsilon\,$ for a curve $\,\Upsilon\,$ embedded in 
$\,\,U\nh$, the operation of restricting functions to $\,\Upsilon$ is a 
bijective correspondence between germs at $\,y\,$ of solutions 
$\,z:U\nh\to\bbR\,$ to (\ref{rzo}) such that 
\begin{equation}\label{tvr} 
\mathrm{the\ vector\ 
}\,\,(\rho(y^{\nh1}\nnh,y^2\nnh,z),\sigma(y^{\nh1}\nnh,y^2\nnh,z))\,\, 
\mathrm{\ at\ }\,\,y,\mathrm{\ with\ }z=z(y),\mathrm{\ is\ not\ tangent\ 
to\ }\,\,\Upsilon, 
\end{equation} 
and germs at $\,y\,$ of functions $\,z:\Upsilon\to I\hs$ having the property 
(\ref{tvr}).

\section{Connections in plane bundles over surfaces}\label{copl} 
\setcounter{equation}{0} 
In this and the next sections, by a {\it plane bundle\/} (or, {\it line 
bundle\/}) we always mean a real vector bundle of fibre dimension $\,2\,$ 
(or, respectively, $\,1$).

For a vector sub\-bundle $\,\lb\,$ of a real vector bundle $\,\pb\,$ over a 
manifold $\,\bs$, we denote by $\,\mathrm{Hom}\hs(\lb,\pb\nh/\nh\lb)\,$ 
the vector bundle over $\,\bs\,$ whose sections are bundle mophisms from 
$\,\pb\,$ into the quotient bundle $\,\pb\nh/\nh\lb$. The {\it fundamental 
tensor of\/} $\,\lb\,$ relative to any given connection $\,\bna\,$ in 
$\,\pb\,$ is then defined to be the bundle mophism 
$\,\ft\nnh_\lb:\tb\to\mathrm{Hom}\hs(\lb,\pb\nh/\nh\lb)\,$ assigning to a 
vector field $\,w\,$ on $\,\bs\,$ the mophism $\,\lb\to\pb\nh/\nh\lb\,$ that 
sends any section $\,v\,$ of $\,\lb\,$ to the image of $\,\bna_{\!w}v$ under 
the quotient projection $\,\pb\to\pb\nh/\nh\lb$. Note that $\,\ft\nnh_\lb$ is 
well defined in view of the Leib\-niz rule, and $\,\ft\nnh_\lb$ vanishes 
identically zero if and only if $\,\lb\,$ is $\,\bna\nh$-par\-al\-lel. 
\begin{lemma}\label{fndte}If\/ $\,\lb\,$ is a line sub\-bun\-dle of a plane 
bundle\/ $\,\pb\,$ over a manifold\/ $\,\bs\,$ and\/ $\,\bna\,$ is a 
connection in\/ $\,\pb\,$ such that the fundamental tensor\/ $\,\ft\nnh_\lb$ 
is nonzero everywhere in\/ $\,\bs$, then\/ $\,\mathrm{Ker}\,\ft\nnh_\lb$ is a 
co\-di\-men\-sion-one distribution on\/ $\,\bs$. 
\end{lemma} 
\begin{proof}In fact, $\,\mathrm{Hom}\hs(\lb,\pb\nh/\nh\lb)\,$ then is a line 
bundle. 
\end{proof} 
\begin{lemma}\label{trcls}For a trace\-less en\-do\-mor\-phism $\,\nd\,$ of 
a \hbox{two\hh-}\hskip0ptdi\-men\-sion\-al real vector space, 
\begin{enumerate} 
  \def\theenumi{{\rm\alph{enumi}}} 
\item[{\rm(i)}] $\nd^2$ is a multiple of the identity, 
\item[{\rm(ii)}] $\nd\,$ is di\-ag\-o\-nal\-izable and nonzero if and only 
if\/ $\,\mathrm{tr}\hskip3pt\nd^2>0$, 
\item[{\rm(iii)}] the kernel of\/ $\,\nd\,$ equals the image of\/ $\,\nd\,$ if 
and only if\/ $\,\nd\ne0\,$ and\/ $\,\mathrm{tr}\hskip3pt\nd^2=0$. 
\end{enumerate} 
\end{lemma} 
This is clear since, in some basis, $\,\nd\,$ is represented by one of the 
matrices
\[ 
\left[\begin{matrix}\kappa&\hskip2pt0\hskip-2pt\cr 
0&-\hs\kappa\end{matrix}\right],\hskip7pt 
\left[\begin{matrix}0&-\hs\kappa\cr 
\kappa&0\end{matrix}\right],\hskip7pt 
\left[\begin{matrix}0&1\cr 
0&0\end{matrix}\right],\hskip12pt\mathrm{with}\hskip5pt\kappa\in\bbR\hs. 
\] 
Given a plane bundle $\,\pb\,$ over a surface $\,\,U\nh$, let $\,\varOmega\,$ 
be an $\,\mathrm{SL}\hh(2,\bbR)${\it-struc\-ture\/} in $\,\pb\nnh$, that is, a 
section of $\,[\pb^*]^{\wedge2}$ without zeros. If $\,\br\,$ is the curvature 
tensor of any connection $\,\bna\,$ in $\,\pb\,$ such that 
$\,\bna\nh\varOmega=0$, then at each point $\,y\in U\,$ we have 
$\,\mathrm{tr}\,\br=0\,$ and one of the three relations 
$\,\mathrm{tr}\,\br^2>\hs0$, $\,\mathrm{tr}\,\br^2=\hs0\,$ and 
$\,\mathrm{tr}\,\br^2<\hs0$, meaning that the trace (in)e\-qual\-i\-ty in 
question is satisfied by the en\-do\-mor\-phism $\,\br_y(w,w\hh'\hh)\,$ of 
$\,\pb\hskip-1.5pt_y$ (cf.\ (\ref{ope})) with some, or any, basis $\,w,w\hh'$ 
of $\,\tyb$. We are going to consider the following conditions: 
\begin{equation}\label{nao} 
\begin{array}{rl} 
\mathrm{i)}\hskip0pt&\bna\nh\varOmega=0\hskip8pt\mathrm{and}\hskip6.5pt 
\mathrm{tr}\,\br^2>\hs0\hskip6.5pt\mathrm{at\ each\ point\ of}\hskip6.5ptU\nh, 
\\ 
\mathrm{ii)}\hskip0pt&\bna\nh\varOmega=0\hs,\hskip5pt\mathrm{while}\hskip6.5pt 
\mathrm{tr}\,\br^2=\hs0\hskip6.5pt\mathrm{and}\hskip6.5pt\br\ne0\hskip6.5pt 
\mathrm{everywhere\ in}\hskip6.5ptU\nh,\\ 
\mathrm{iii)}\hskip0pt&\bna\nh\varOmega=0\hskip5pt\mathrm{and}\hskip6.5pt\bna 
\hskip6.5pt\mathrm{is\ flat,\ that\ is,}\hskip6.5pt\br\hskip6.5pt 
\mathrm{vanishes\ identically.} 
\end{array} 
\end{equation} 
\begin{lemma}\label{sbbdl}For\/ $\,U\nh,\pb\nnh,\varOmega,\bna\,$ and\/ 
$\,\br\,$ as above, vector fields\/ $\,w,w\hh'$ on $\,\,U\,$ linearly 
independent at every point, and the corresponding bundle morphism\/ 
$\,\br(w,w\hh'\hh):\pb\to\pb\nnh$, cf.\ {\rm(\ref{ope})}, 
\begin{enumerate} 
  \def\theenumi{{\rm\alph{enumi}}} 
\item[{\rm(a)}] condition\/ {\rm(\ref{nao}.i)} holds if and only if\/ 
$\,\pb\,$ is the direct sum of two line sub\-bun\-dles\/ $\,\lb^\pm\nnh$, 
which are the eigen\-space bundles of\/ $\,\br(w,w\hh'\hh)$, 
\item[{\rm(b)}] condition\/ {\rm(\ref{nao}.ii)} is equivalent to the existence 
of a line sub\-bun\-dle\/ $\,\lb\,$ of\/ $\,\pb\,$ which is both the kernel 
and the image of\/ $\,\br(w,w\hh'\hh)$. 
\end{enumerate} 
\end{lemma} 
\begin{proof}The assertion is immediate from Lemma~\ref{trcls}. 
\end{proof} 
In our subsequent discussion, (\ref{nao}.ii) will be coupled with the 
additional condition 
\begin{equation}\label{fdt} 
\mathrm{the\ fundamental\ tensor\ of\ }\,\lb\,\mathrm{\ relative\ to\ }\,\bna\, 
\mathrm{\ is\ nonzero\ everywhere\ in\ }\hs\,U\nh, 
\end{equation} 
where the fundamental tensor is defined at the beginning of this section, and 
$\,\lb\,$ stands for the line sub\-bun\-dle of $\,\pb\hs$ characterized by 
Lemma~\ref{sbbdl}(b).

In the next lemma, {\it general position\/} means the same as in the third 
paragraph of the Introduction. Given a connection $\,\bna\,$ in a plane bundle 
$\,\pb\,$ over a surface $\,\,U\nh$, with a fixed 
$\,\mathrm{SL}\hh(2,\bbR)$-struc\-ture $\,\varOmega\,$ in $\,\pb\,$ such that 
$\,\bna\nh\varOmega=0$, local coordinates $\,y^{\hs j}$ in $\,\,U\nh$, and 
local trivializing sections $\,e_k$ of $\,\pb\,$ satisfying the condition 
$\,\varOmega(e_1,e_2)=1$, we denote by $\,\Gamma_{\hskip-2.7ptjk}^{\hh l}$ the 
corresponding component functions of $\,\bna\nnh$, characterized by 
$\,\bna_{\!w}e_k=w^{\hs j}\Gamma_{\hskip-2.7ptjk}^{\hh l}e_l$ 
(summation over repeated indices), for any vector field $\,w\,$ on $\,\bs\,$ 
with the component functions $\,w^{\hs j}\nnh$. We continue using subscripts 
for partial derivatives, the only exceptions being the symbols 
$\,e_k,\hs\Gamma_{\hskip-2.7ptjk}^{\hh l}$ and $\,\br_{12\hskip1.3ptk}{}^l\nnh$. 
\begin{lemma}\label{gamgp}Suppose that a connection\/ $\,\bna\,$ in a plane 
bundle\/ $\,\pb\,$ with an\/ $\,\mathrm{SL}\hh(2,\bbR)$-struc\-ture\/ 
$\,\varOmega\,$ over a surface\/ $\,\,U\,$ satisfies\/ {\rm(\ref{nao}.i)}, 
or\/ {\rm(\ref{nao}.ii)} along with\/ {\rm(\ref{fdt})}, or\/ 
{\rm(\ref{nao}.iii)}. Then, at points in general position, locally, for some\/ 
$\,y^{\hs j}$ and\/ $\,e_k$ as above, $\,\bna\,$ has one of the following 
descriptions, in which $\,\psi,\chi\,$ are functions of\/ 
$\,(y^{\nh1}\nnh,y^2)\,$ and $\,p\,$ is a function of one real variable\/{\rm:} 
\begin{enumerate} 
  \def\theenumi{{\rm\roman{enumi}}} 
\item[{\rm(I)}] $\Gamma_{\hskip-2.7pt12}^{\hh1}=e^{2\chi}\nh,\hskip6pt 
\Gamma_{\hskip-2.7pt22}^{\hh2}=-\hs\Gamma_{\hskip-2.7pt21}^{\hh1} 
=\chi_2\hh,\hskip6pt\Gamma_{\hskip-2.7pt22}^{\hh1}=0\,$ and either 
\begin{enumerate} 
  \def\theenumi{{\rm\alph{enumi}}} 
\item[{\rm(a)}] $\Gamma_{\hskip-2.7pt11}^{\hh1} 
=-\hs\Gamma_{\hskip-2.7pt12}^{\hh2}=\psi_1\hh,\hskip6pt 
\Gamma_{\hskip-2.7pt11}^{\hh2}=0\,\,$ and\/ 
$\,\,\Gamma_{\hskip-2.7pt21}^{\hh2}=e^{2\psi}\nnh$, or 
\item[{\rm(b)}] $\Gamma_{\hskip-2.7pt11}^{\hh1} 
=-\hs\Gamma_{\hskip-2.7pt12}^{\hh2}$ is arbitrary, 
$\hskip3pt\Gamma_{\hskip-2.7pt11}^{\hh2}=p(y^{\nh1})\hh e^{-2\chi}\hs\,$ and\/ 
$\,\,\Gamma_{\hskip-2.7pt21}^{\hh2}=0$, or 
\item[{\rm(c)}] $\Gamma_{\hskip-2.7pt11}^{\hh1} 
=-\hs\Gamma_{\hskip-2.7pt12}^{\hh2} 
=p(y^{\nh1})-\chi_1,\hskip6pt\Gamma_{\hskip-2.7pt11}^{\hh2}$ is arbitrary, 
and\/ $\,\Gamma_{\hskip-2.7pt21}^{\hh2}=0\hh$, 
\end{enumerate} 
\item[{\rm(II)}] $\Gamma_{\hskip-2.7pt11}^{\hh2} 
=\Gamma_{\hskip-2.7pt12}^{\hh1}=\Gamma_{\hskip-2.7pt21}^{\hh2} 
=\Gamma_{\hskip-2.7pt22}^{\hh1}=0$, while\/ 
$\,\Gamma_{\hskip-2.7pt11}^{\hh1} 
=-\hs\Gamma_{\hskip-2.7pt12}^{\hh2}$ and\/ 
$\,\Gamma_{\hskip-2.7pt22}^{\hh2}=-\hs\Gamma_{\hskip-2.7pt21}^{\hh1}$ are 
arbitrary, 
\item[{\rm(III)}] $\Gamma_{\hskip-2.7ptjk}^{\hh l}=0\,$ for all\/ $\,j,k,l$. 
\end{enumerate} 
Conversely, {\rm(I)} -- {\rm(III)} easily imply\/ {\rm(\ref{nao}.i)}, or\/ 
{\rm(\ref{nao}.ii)} and\/ {\rm(\ref{fdt})}, or\/ {\rm(\ref{nao}.iii)}.

The descriptions in\/ {\rm(I)} -- {\rm(III)} and the three cases in\/ 
{\rm(\ref{nao})} are related as follows\/{\rm:} {\rm(\ref{nao}.i)} is realized 
by\/ {\rm(I-a)}, {\rm(I-b)} and\/ {\rm(II);} {\rm(\ref{nao}.ii)} 
with\/ {\rm(\ref{fdt})} by\/ {\rm(I-c);} and\/ {\rm(\ref{nao}.iii)} by\/ 
{\rm(III)}. 
\end{lemma} 
\begin{proof}Obviously, (III) corresponds to (\ref{nao}.iii). From now on we 
assume (\ref{nao}.i), or (\ref{nao}.ii) with (\ref{fdt}). Let $\,\ft^\pm$ (or, 
$\,\ft$) be the fundamental tensor, relative to $\,\bna\nnh$, of the line 
sub\-bun\-dle $\,\lb^\pm$ (or, $\,\lb$) appearing in Lemma~\ref{sbbdl}.

If (\ref{nao}.i) holds and $\,\ft^+\nnh=\ft^-\nnh=0\,$ everywhere in 
$\,\,U\nh$, we obtain (II) by choosing $\,e_1$ and $\,e_2$ to be sections of 
$\,\lb^+$ and $\,\lb^-\nnh$. Conversely, it is clear that (II) gives 
$\,\ft^\pm\nnh=0$.

The remaining case, in general position, means that we have either 
(\ref{nao}.i) and $\,\ft^\pm\nnh\ne0$ everywhere for some fixed sign 
$\,\pm\hs$, or (\ref{nao}.ii) and (\ref{fdt}). If we unify the notation by 
setting $\,\ft=\ft^\pm$ and $\,\lb=\lb^\pm$ for this sign $\,\pm\hs$, then 
in both cases we may, locally, choose $\,e_k$ and $\,y^{\hs j}$ such that 
$\,e_2$ is a section of $\,\lb\,$ and 
$\,\mathrm{Ker}\,d\hskip.2pty^{\nh1}\nnh=\mathrm{Ker}\,\ft$, cf.\ 
Lemma~\ref{fndte}. The coordinate vector field $\,\partial_2$ thus spans the 
distribution $\,\mathrm{Ker}\,\ft$, and so 
$\,\Gamma_{\hskip-2.7pt22}^{\hh1}=0$, while 
$\,\Gamma_{\hskip-2.7pt12}^{\hh1}\ne0$. Let us change the sign of $\,y^{\nh1}$ 
(or $\,e_1$, or $\,e_2$) if necessary, so as to make 
$\,\Gamma_{\hskip-2.7pt12}^{\hh1}$ positive. Hence 
$\,\Gamma_{\hskip-2.7pt12}^{\hh1}=e^{2\chi}$ for a function $\,\chi$. By 
(\ref{cur}) with $\,\Gamma_{\hskip-2.7pt22}^{\hh1}=0\,$ and 
$\,\Gamma_{\hskip-2.7pt21}^{\hh1}=-\hs\Gamma_{\hskip-2.7pt22}^{\hh2}$, we now 
get $\,0=e^{-2\chi}\br_{122}{}^1\nnh 
=2\chi_2+\Gamma_{\hskip-2.7pt21}^{\hh1}-\Gamma_{\hskip-2.7pt22}^{\hh2} 
=2(\chi_2-\Gamma_{\hskip-2.7pt22}^{\hh2})$, which yields the first line in (I).

So far, $\,e_1$ and $\,y^2$ have been completely arbitrary except 
for the relations $\,\varOmega(e_1,e_2)=1\,$ and 
$\,d\hskip.2pty^{\nh1}\wedge\,d\hskip.2pty^2\ne\,0$. Each of the following 
four paragraphs begins with a specific gen\-er\-al-po\-si\-tion assumption, 
and uses some special choice of one or both of $\,e_1$ and $\,y^2\nnh$.

If (\ref{nao}.ii) and (\ref{fdt}) hold, choosing $\,e_1$ so that it spans a 
line sub\-bun\-dle which is $\,\bna\nh$-par\-al\-lel along the $\,y^2$ 
coordinate direction, we get $\,\Gamma_{\hskip-2.7pt21}^{\hh2}=0$. Since 
$\,e_2$ spans both the kernel and the image of 
$\,\br(\partial_{\hh1},\partial_{\hh2})\,$ (see Lemma~\ref{sbbdl}(b)), 
one has 
$\,0=\br_{121}{}^1\nnh 
=(\Gamma_{\hskip-2.7pt11}^{\hh1})_2-(\Gamma_{\hskip-2.7pt21}^{\hh1})_1$, 
which yields (I-c) as 
$\,\Gamma_{\hskip-2.7pt21}^{\hh1}=-\hs\Gamma_{\hskip-2.7pt22}^{\hh2} 
=-\chi_2$.

For the remainder of the proof, we assume (\ref{nao}.i) and choose $\,e_1$ to 
be a section of $\,\lb^\mp\nnh$, with the sign $\,\mp\,$ opposite to $\,\pm\,$ 
that was fixed earlier. If $\,\ft^\mp\nnh=0\,$ identically, we obviously 
obtain (I-b) with $\,p(y^{\nh1})=0$.

If, on the other hand, $\,\ft^\mp\nnh\ne0\,$ everywhere and 
$\,\mathrm{Ker}\,\ft^\mp\nnh=\mathrm{Ker}\,\ft^\pm\nnh$, so that 
$\,\Gamma_{\hskip-2.7pt21}^{\hh2}=0$, we get, as before, 
$\,\Gamma_{\hskip-2.7pt11}^{\hh2}=e^{2\psi}$ for a function $\,\psi\,$ and 
$\,0=e^{-2\psi}\br_{121}{}^2\nnh 
=2\psi_2+\Gamma_{\hskip-2.7pt22}^{\hh2}-\Gamma_{\hskip-2.7pt21}^{\hh1} 
=2(\psi_2+\Gamma_{\hskip-2.7pt22}^{\hh2})$. Hence 
$\,\psi_2=-\hs\Gamma_{\hskip-2.7pt22}^{\hh2}=-\chi_2$, that is, 
$\,2\psi=-2\chi+\log\hskip1.2ptp(y^{\nh1})\,$ for some positive function 
$\,p\,$ of one variable, which again gives (I-b).

Finally, if $\,\ft^\mp\nnh\ne0\,$ and $\,\mathrm{Ker}\,\ft^\mp$ is transverse 
to $\,\mathrm{Ker}\,\ft^\pm$ everywhere in $\,\,U\nh$, we choose $\,y^2$ with 
$\,\mathrm{Ker}\,dy^2\nnh=\mathrm{Ker}\,\ft^\mp$, so that 
$\,\Gamma_{\hskip-2.7pt11}^{\hh2}=0$, while $\,\Gamma_{\hskip-2.7pt21}^{\hh2}$, 
being nonzero, may as before be assumed positive end hence equal to 
$\,e^{2\psi}$ for some function $\,\psi$. Now 
$\,0=e^{-2\psi}\br_{121}{}^2\nnh 
=\Gamma_{\hskip-2.7pt11}^{\hh1}-\Gamma_{\hskip-2.7pt12}^{\hh2}-2\psi_1 
=2(\Gamma_{\hskip-2.7pt11}^{\hh1}-\psi_1)$, and (I-a) follows. 
\end{proof}

\section{Solutions to the system (\ref{eqn})}\label{stsy} 
\setcounter{equation}{0} 
We use the same assumptions and notations as at the beginning of 
Section~\ref{trct} and in Theorem~\ref{cndtd}. The subscripts 
$\,(\hskip2.5pt)_j$, $\,j=1,2$, stand, again, for the directional derivatives 
in the directions of the constant vector fields $\,\partial_j$ on $\,\bs\,$ 
forming the basis of $\,\dbs\,$ dual to the basis $\,\xi,\ts$ of 
$\,\dbs^*\nnh$. In other words, 
$\,(\hskip2.5pt)_j=\partial/\partial y^{\hs j}$ for the af\-fine coordinates 
$\,y^{\hs j}$ with $\,d\hskip.2pty^{\nh1}\nh=\xi\,$ and 
$\,d\hskip.2pty^2\nh=\ts$. Here $\,\dbs\,$ is the translation vector space of 
the af\-fine plane $\,\bs$, and, as before, $\,\,U\,$ denotes a nonempty 
connected open subset of $\,\bs$, while $\,K\,$ is any given real constant.

The area form $\,\varOmega\in[\plane^*]^{\wedge2}\nnh\smallsetminus\{0\}\,$ 
gives rise to an isomorphic identification between the space 
$\,\mathfrak{sl}\hs(\plane)\,$ of all trace\-less en\-do\-mor\-phisms of 
$\,\plane\,$ and the space $\,[\plane^*]^{\odot2}$ of all symmetric 
bi\-lin\-e\-ar forms on $\,\plane$, which associates with 
$\,\nd\in\mathfrak{sl}\hs(\plane)\,$ the form $\,\bz\,$ sending 
$\,u,v\in\plane\,$ to $\,\bz(u,v)=\varOmega(\nd u,v)$. (Symmetry of $\,\bz\,$ 
is immediate from Remark~\ref{trvol}.) Under this identification any function 
$\,(\vf,\lambda,\mu):U\nh\to V\hs$ corresponds to a function 
$\,(\vf,\dv,\ve):U\nh\to\plane 
\times\mathfrak{sl}\hs(\plane)\nh\times\mathfrak{sl}\hs(\plane)$. 
(As in Section~\ref{trct}, we set 
$\,V\nnh=\plane\times[\plane^*]^{\odot2}\nnh\times[\plane^*]^{\odot2}\nnh$.) 
The functions $\,\dv\,$ and $\,\ve$, taking values in en\-do\-mor\-phisms of 
$\,\plane$, can be multiplied by each other and multiplied by functions 
$\,\,U\nh\to\plane$, in the sense of the val\-ue\-wise operations of composing 
en\-do\-mor\-phisms and evaluating them on vectors. Thus, for instance, the 
val\-ue\-wise commutator $\,[\hs\ve,\dv\hs]\,$ is a function 
$\,\,U\nh\to\mathfrak{sl}\hs(\plane)$. 
\begin{lemma}\label{eqvnd}A function\/ $\,(\vf,\lambda,\mu):U\nh\to V\hs$ 
satisfies\/ {\rm(\ref{eqn})} if and only if 
\begin{equation}\label{teo} 
\mathrm{i)}\hskip6pt2\hh\ve_1-\dv_2 
=2\hs[\hs\ve,\dv\hs]-K\varOmega(c,\,\cdot\,)\hs\vf 
-K\varOmega(\vf,\,\cdot\,)\hs c\hs,\hskip16pt\mathrm{ii)}\hskip6pt 
\varOmega(c,\dv_1c-2\vf_2-4\ve\hh\vf)=1\ 
\end{equation} 
for\/ $\,(\vf,\dv,\ve)\,$ corresponding to\/ $\,(\vf,\lambda,\mu)\,$ as above. 
\end{lemma} 
\begin{proof}Let $\,\nd_{\mathrm{L}}$ and $\,\nd_{\mathrm{R}}$ be the 
left-hand and right-hand sides of (\ref{teo}.i). Also, let 
$\,\nd=\varOmega(c,\,\cdot\,)\hs\vf+\varOmega(\vf,\,\cdot\,)\hs c$. Both 
$\,\nd_{\mathrm{L}}$ and $\,\nd_{\mathrm{R}}$ are functions 
$\,\,U\nh\to\mathfrak{sl}\hs(\plane)$. (In fact, $\,\nd\,$ takes values in 
$\,\mathfrak{sl}\hs(\plane)$, since $\,\varOmega(\nd u,v)\,$ is symmetric in 
$\,u,v\in\plane$.) For the same reasons of symmetry, we have 
$\,\nd_{\mathrm{L}}\nh=\nd_{\mathrm{R}}$ if and only if 
$\,\varOmega(\nd_{\mathrm{L}}\rd,\rd)=\varOmega(\nd_{\mathrm{R}}\rd,\rd)$, 
where $\,\rd\,$ denotes, as usual, the radial vector field on $\,\plane$. 
Clearly, 
$\,\varOmega(\nd_{\mathrm{L}}\rd,\rd)=2\mu_1(\rd,\rd)-\lambda_2(\rd,\rd)$, 
and, in view of (\ref{duf}.i), 
$\,\varOmega(\nd\rd,\rd)=2\hh\phi\hs\varOmega(\rd,\vf)$. On the other hand, by 
(\ref{lcx}), the expression 
$\,\lambda(c,\rd)\hs\mu(\rd,\rd)-\mu(c,\rd)\hs\lambda(\rd,\rd)$, depending 
skew-sym\-met\-ri\-cal\-ly on $\,\lambda\,$ and $\,\mu$, equals 
$\,\phi\hs\varOmega(\ve\hh\dv\rd,\rd)$. In view of skew-sym\-me\-try, this is 
further equal to $\,\varOmega([\hs\ve,\dv\hs]\rd,\rd)\hs\phi/2$. Hence 
(\ref{teo}.i) is equivalent the first equality in (\ref{eqn}). Equivalence 
between (\ref{teo}.ii) and the second equality in (\ref{eqn}) is in turn 
obvious from the definitions of $\,\dv\,$ and $\,\ve$. 
\end{proof} 
Treating $\,\,U\nh\times\plane\,$ as the total space of a product vector 
bundle $\,\pb\,$ over $\,\,U\nh$, we may view functions $\,\,U\nh\to\plane\,$ 
(including the constant $\,c$) as sections of $\,\pb\nnh$, while 
$\,\varOmega\,$ then becomes an $\,\mathrm{SL}\hh(2,\bbR)$-struc\-ture in 
$\,\pb\nnh$, cf.\ Section~\ref{copl}. 
\begin{lemma}\label{bijct}Functions\/ $\,(\vf,\lambda,\mu):U\nh\to V\hs$ are 
in a natural bijective correspondence with pairs\/ $\,(\vf,\bna\hh)$, in 
which\/ $\,\vf\,$ is a section of the plane bundle\/ 
$\,\pb=U\nh\times\plane\,$ with the\/ $\,\mathrm{SL}\hh(2,\bbR)$-struc\-ture\/ 
$\,\varOmega$ and\/ $\,\bna\,$ is a connection in\/ $\,\pb\,$ such that\/ 
$\,\bna\nh\varOmega=0$. The correspondence associates with\/ 
$\,(\vf,\lambda,\mu)$ the pair\/ $\,(\vf,\bna\hh)\,$ obtained by treating\/ 
$\,\vf:U\nh\to\plane\,$ as a section of\/ $\,\pb\,$ and setting 
\begin{equation}\label{nbo} 
\bna_{\!1}v=v_1+\dv\hs v\hs,\hskip26pt 
\bna_{\!2}v=v_2+2\hh\ve v 
\end{equation} 
for functions\/ $\,v:U\nh\to\plane$, also viewed as sections of\/ $\,\pb\nnh$, 
where\/ $\,(\vf,\dv,\ve)\,$ is related to\/ $\,(\vf,\lambda,\mu)\,$ as in the 
lines preceding Lemma\/~{\rm\ref{eqvnd}}, and\/ $\,\bna_{\!j}$ denotes the\/ 
$\,\bna\nnh$-co\-var\-i\-ant derivative in the direction of the coordinate 
vector field\/ $\,\partial_j$ on\/ $\,\,U\nh$. 
\end{lemma} 
\begin{proof}This is obvious, as the condition $\,\bna\nh\varOmega=0\,$ 
accounts for the fact that $\,\dv\,$ and $\,\ve\,$ take values in the space 
$\,\mathfrak{sl}\hs(\plane)\,$ of {\it trace\-less\/} en\-do\-mor\-phisms of 
$\,\plane$. 
\end{proof} 
\begin{lemma}\label{eqfcn}For a function\/ $\,(\vf,\lambda,\mu):U\nh\to V\nh$, 
condition\/ {\rm(\ref{eqn})} holds if and only if, in terms of the pair\/ 
$\,(\vf,\bna\hh)\,$ corresponding to\/ $\,(\vf,\lambda,\mu)\,$ in the sense of 
Lemma\/~{\rm\ref{bijct}}, the curvature tensor\/ $\,\br\,$ of\/ $\,\bna\nnh$, 
and\/ $\,c,\vf\,$ treated as sections of\/ $\,\pb\nnh$, 
\begin{equation}\label{brd} 
\mathrm{i)}\hskip6pt\br(\partial_{\hh1},\partial_{\hh2}) 
=K\varOmega(c,\,\cdot\,)\hs\vf+K\varOmega(\vf,\,\cdot\,)\hs c\hs,\hskip24pt 
\mathrm{ii)}\hskip6pt\varOmega(c,\bna_{\!1}\nnh\bna_{\!1}c-2\bna_{\!2}\vf) 
=1\hs, 
\end{equation} 
with\/ $\,\bna_{\!j}$ and\/ $\,\partial_j$ as in Lemma\/~{\rm\ref{bijct}}, so 
that\/ $\,\br(\partial_{\hh1},\partial_{\hh2})\,$ is a morphism\/ 
$\,\pb\to\pb$, cf.\/ {\rm(\ref{ope})}. 
\end{lemma} 
\begin{proof}By (\ref{cur}) and (\ref{nbo}), 
$\,\br(\partial_{\hh1},\partial_{\hh2}) 
=\dv_2-2\hh\ve_1-2\hs[\hs\ve,\dv\hs]$, Thus, (\ref{brd}.i) is equivalent to 
(\ref{teo}.i). Next, since $\,c\,$ is constant, (\ref{nbo}) gives 
$\,\bna_{\!1}c=\dv\hs c$, $\,\bna_{\!1}\nnh\bna_{\!1}c=\dv_1c+\dv^2c$, and 
$\,\bna_{\!2}\vf=\vf_2+2\ve\hh\vf$. As $\,\dv^2$ takes values in multiples of 
the identity (see Lemma~\ref{trcls}(i)), the left-hand side of (\ref{brd}.ii) 
coincides with that of (\ref{teo}.ii). 
\end{proof} 
\begin{remark}\label{rddek}Suppose that (\ref{brd}.i) is satisfied by a 
connection $\,\bna\,$ in a plane bundle $\,\pb$ over a surface $\,\,U\nh$, 
an $\,\mathrm{SL}\hh(2,\bbR)$-struc\-ture $\,\varOmega\,$ in $\,\pb\,$ with 
$\,\bna\nh\varOmega=0$, the curvature tensor $\,\br\,$ of $\,\bna\nnh$, 
sections $\,c\,$ and $\,\vf\,$ of $\,\pb\,$ such that $\,c\,$ has no zeros, 
some vector fields $\,\partial_j$ on $\,\,U\nh$, linearly independent at every 
point, and a real constant $\,K$. Let $\,y\in U\nh$. 
\begin{enumerate} 
  \def\theenumi{{\rm\alph{enumi}}} 
\item[{\rm(a)}]If $\,c\,$ and $\,K\vf\,$ are linearly independent at $\,y$, 
they are eigenvectors of $\,\br(\partial_{\hh1},\partial_{\hh2})$, acting in 
the fibre $\,\pb\hskip-1.5pt_y$, for the nonzero eigenvalues 
$\,K\varOmega(\vf,c)\,$ and $\,-\hs K\varOmega(\vf,c)\,$ (at $\,y$). 
\item[{\rm(b)}]If $\,K\vf$, at $\,y$, is a nonzero multiple of $\,c$, then 
$\,c\,$ spans, at $\,y$, the kernel of 
$\,\br(\partial_{\hh1},\partial_{\hh2})\,$ acting in 
$\,\pb\hskip-1.5pt_y$ (which is at the same time its image). 
\item[{\rm(c)}]If $\,K\vf=0$, at $\,y$, then 
$\,\br(\partial_{\hh1},\partial_{\hh2})=0$, at $\,y$. 
\item[{\rm(d)}]In case (a), or (b), or (c), we have at $\,y$, respectively, 
$\,\mathrm{tr}\,\br^2>\hs0$, or $\,\mathrm{tr}\,\br^2=\hs0\,$ and $\,\br\ne0$, 
or $\,\br=0$. 
\end{enumerate} 
In fact, (a) -- (c) are obvious, and (d) is immediate from Lemma~\ref{trcls}. 
\end{remark} 
\begin{lemma}\label{nprll}Given\/ 
$\,\,U\nh,\pb\nnh,\varOmega,\bna\nnh,\br,c,\vf,\hs\partial_j$ and\/ $\,K\,$ 
as in Remark\/~{\rm\ref{rddek}}, if conditions\/ 
{\rm(\ref{nao}.ii)} and\/ {\rm(\ref{brd})} are both satisfied, and\/ $\,\lb\,$ 
is the line sub\-bun\-dle of\/ $\,\pb\,$ appearing in 
Lemma\/~{\rm\ref{sbbdl}(b)}, then the fundamental tensor of\/ $\,\lb\,$ 
relative to\/ $\,\bna\nnh$, defined in Section\/~{\rm\ref{copl}}, is nonzero 
at all points of a dense open subset of\/ $\,\,U\nh$. 
\end{lemma} 
\begin{proof}By (\ref{brd}.i), $\,c\,$ and $\,\vf\,$ span the image of 
$\,\br(\partial_{\hh1},\partial_{\hh2})$, so that they are sections of 
$\,\lb\,$ (cf.\ Lemma~\ref{sbbdl}(b)). If $\,\lb\,$ were 
$\,\bna\nh$-par\-al\-lel on some nonempty open set $\,\,U'\nnh\subset U\nh$, 
then $\,\bna_{\!1}\nnh\bna_{\!1}c$ and $\,\bna_{\!2}\vf\,$ would be sections 
of $\,\lb\,$ as well, and so 
$\,\varOmega(c,\bna_{\!1}\nnh\bna_{\!1}c-2\bna_{\!2}\vf)\,$ would vanish on 
$\,\,U'\nnh$, contradicting (\ref{brd}.ii). 
\end{proof} 
We will now describe all 
$\,\,U\nh,\pb\nnh,\varOmega,\bna\nnh,\br,c,\vf,\hs\partial_j$ and $\,K\,$ with 
the properties listed in Remark~\ref{rddek}, for which (\ref{brd}) holds. 
According to Lemma~\ref{eqfcn}, this is equivalent to solving (\ref{eqn}). 
(Since a quadruple $\,(U\nh,\pb\nnh,\varOmega,c)$, having those of the 
properties just named that pertain just to $\,\,U\nh,\pb\nnh,\varOmega$ and 
$\,c$, is, locally, unique up to struc\-ture-pre\-serv\-ing bundle 
isomorphisms and dif\-feo\-mor\-phisms of the base, it makes no difference 
whether we treat $\,(U\nh,\pb\nnh,\varOmega,c)\,$ as fixed, or allow it to 
vary.)

First, let $\,K\nh=0$. Condition (\ref{brd}.i) now amounts to flat\-ness of 
$\,\bna$. If we fix an arbitrary flat connection $\,\bna\,$ in $\,\pb\,$ and 
choose, locally, a section $\,a\,$ of $\,\pb\,$ with $\,\varOmega(a,c)=1$, 
(\ref{brd}.ii) is equivalent to requiring that 
$\,2\bna_{\!2}\vf=a+\chi\hh c+\nh\bna_{\!1}\nnh\bna_{\!1}c\,$ for some 
function $\,\chi$. With prescribed $\,c$ (as well as $\,a\,$ and $\,\chi$), 
this is a system of linear ordinary differential equations with parameters, 
and our task has been reduced to solving it for $\,\vf$, which is well 
understood.

Therefore, from now on, $\,K\nh\ne0$. We also make a gen\-er\-al-po\-si\-tion 
assumption: 
\begin{equation}\label{gpo} 
\begin{array}{rl} 
\mathrm{i)}\hskip0pt&c\,\,\mathrm{\ and\ }\,\vf\,\mathrm{\ are\ linearly\ 
independent\ at\ each\ point\ of\ }\,\,U\nh,\mathrm{\ or,}\\ 
\mathrm{ii)}\hskip0pt&\vf\,\,\mathrm{\ is,\ at\ every\ point,\ a\ nonzero\ 
multiple\ of\ }\,c,\mathrm{\ or,\ finally,}\\ 
\mathrm{iii)}\hskip0pt&\vf\,\,\mathrm{\ vanishes\ identically\ on\ }\,\,U\nh. 
\end{array} 
\end{equation} 
In a dense open subset $\,\,U'$ of $\,\,U\nh$, one of (\ref{gpo}.i), 
(\ref{gpo}.ii) and (\ref{gpo}.iii) will clearly hold if we replace $\,\,U\,$ 
by a suitable neighborhood of any given point of $\,\,U'\nnh$.

Rather than dealing with the three possibilities in (\ref{gpo}), we will 
consider those listed in (\ref{nao}). Namely, if (\ref{brd}) holds (which is 
what we want to achieve), each of the three lines in (\ref{gpo}) 
implies, according to Remark~\ref{rddek} and Lemma~\ref{sbbdl}, the 
corresponding line in (\ref{nao}). In case (\ref{nao}.ii), we also assume 
(\ref{fdt}), which, by Lemma~\ref{nprll}, is a gen\-er\-al-po\-si\-tion 
version of a condition necessary for (\ref{brd}).

To summarize, we now begin with a connection $\,\bna\,$ in a plane bundle 
$\,\pb\,$ over a surface $\,\,U$ and an 
$\,\mathrm{SL}\hh(2,\bbR)$-struc\-ture $\,\varOmega\,$ in $\,\pb$, which 
satisfy (\ref{nao}.i), or (\ref{nao}.ii) and (\ref{fdt}), or (\ref{nao}.iii), 
but are otherwise arbitrary. Locally, at points in general position, all such 
$\,\bna\,$ and $\,\varOmega\,$ are described by Lemma~\ref{gamgp}. We may 
therefore proceed with the next step: finding all sections $\,c\,$ and 
$\,\vf\,$ of $\,\pb\,$ with (\ref{brd}), such that $\,c\ne0\,$ everywhere.

In the cases (\ref{nao}.i) and (\ref{nao}.ii), let $\,\lb^\pm$ (or $\,\lb$) be 
as in Lemma~\ref{sbbdl}. Choosing, locally, sections $\,c\,$ of $\,\lb^+$ and 
$\,\vf\,$ of $\,\lb^-$ (or, sections $\,c\,$ and $\,\vf\,$ of $\,\lb$) without 
zeros, we easily see that $\,\br(\partial_{\hh1},\partial_{\hh2})\,$ equals 
some function without zeros times the right-hand side of (\ref{brd}.i). 
Multiplying $\,c\,$ or $\,\vf\,$ by a suitable function, we obtain 
(\ref{brd}.i) for the new $\,c\,$ and $\,\vf$. The choice of $\,c\,$ and 
$\,\vf$ that satisfies (\ref{brd}.i) is not unique, as we are free to switch 
$\,c\,$ with $\,\vf$, and/or use, instead of them, $\,zc\,$ and 
$\,z^{-1}\nh\vf$, for any function $\,z:U\nh\to\bbR\,$ without zeros. (In our 
discussion, $\,\,U\,$ may be repeatedly replaced with smaller connected open 
sets, still denoted by $\,\,U\nh$.) As (\ref{brd}.i) is already satisfied, we 
just need to determine which functions $\,z:U\nh\to\bbR\smallsetminus\{0\}\,$ 
lead to (\ref{brd}.ii) for $\,zc\,$ and $\,z^{-1}\nh\vf\,$ rather than $\,c\,$ 
and $\,\vf$. Such $\,z\,$ are easily seen to be characterized by (\ref{rzo}) 
with $\,\rho=\varOmega(c,\bna_{\!1}c)\hh z^2\nnh$, 
$\,\sigma=\varOmega(c,\vf)$, 
$\,\chi=[z-\varOmega(c,\bna_{\!1}\nnh\bna_{\!1}c)]\hh z^3\nnh/\hh2$, 
$\,I\nh=(0,\infty)\,$ or $\,I\nh=(-\infty,0)$, and $\,\,U\,$ treated as a 
connected open set in $\,\rto\nnh$, for which $\,y^{\nh1}\nnh,y^2$ serve as 
the Cartesian coordinates.

Until now our discussion has covered both cases (\ref{nao}.i) and 
(\ref{nao}.ii). We now separate them, first assuming (\ref{nao}.i). As 
$\,\sigma=\varOmega(c,\vf)\,$ is nonzero everywhere, we may fix $\,y\in U\,$ 
and choose the curve $\,\Upsilon\,$ with (\ref{tvr}) to be a line segment in 
$\,\,U\nh\subset\rto$ containing $\,y\,$ and parallel to the $\,y^{\nh1}$ 
coordinate direction. As stated at the end of Section~\ref{tmoc}, germs at 
$\,y\,$ of functions $\,z:U\nh\to\bbR\smallsetminus\{0\}\,$ leading to 
(\ref{brd}.ii) are in a natural bijective correspondence with germs at $\,y\,$ 
of arbitrary nonzero functions on $\,\Upsilon$, which concludes our 
description of solutions to (\ref{eqn}) in the case (\ref{nao}.i).

Next, suppose that (\ref{nao}.ii) and (\ref{fdt}) are satisfied. Thus, 
$\,\sigma=\varOmega(c,\vf)\,$ vanishes identically, while 
$\,\rho=\varOmega(c,\bna_{\!1}c)\hh z^2\nnh$, so that $\,\rho\ne0\,$ by 
(\ref{fdt}), as the values of $\,z\,$ range over $\,\bbR\smallsetminus\{0\}$. 
The transversality condition (\ref{tvr}) then holds for $\,\Upsilon\,$ chosen 
to be a line segment through a given point $\,y$, parallel to the $\,y^2$ 
coordinate 
direction. Once again, germs at $\,y\,$ of functions 
$\,z:U\nh\to\bbR\smallsetminus\{0\}\,$ that lead to (\ref{brd}.ii) may be 
bijectively identified with germs at $\,y\,$ of arbitrary functions 
$\,\Upsilon\to\bbR\smallsetminus\{0\}$.

Finally, let us assume (\ref{nao}.iii), so that $\,\bna\,$ is flat. As 
$\,K\nh\ne0\,$ and we want $\,c\,$ to have no zeros, (\ref{brd}.i) can only be 
satisfied if $\,\vf=0\,$ identically. Condition (\ref{brd}.ii) with 
$\,\vf=0$ reads $\,\varOmega(c,\bna_{\!1}\nnh\bna_{\!1}c)=1$. Fixing, 
locally, $\,\bna\nh$-par\-al\-lel sections $\,e_1,e_2$ of $\,\pb\,$ with 
$\,\varOmega(e_1,e_2)=1$ and writing 
$\,c=(\mathrm{Re}\,\chi)\hh e_1+(\mathrm{Im}\,\chi)\hh e_2$ for some 
unspecified function $\,\chi:U\nh\to\bbC\smallsetminus\{0\}$, we get 
$\,\varOmega(c,\bna_{\!1}\nnh\bna_{\!1}c)=(\rho^2\sigma_1)_1$, where, locally, 
$\,\chi=\rho e^{i\sigma}$ for real-val\-ued functions $\,\rho>0\,$ and 
$\,\sigma$. The resulting equation $\,(\rho^2\sigma_1)_1=1\,$ can be solved 
by repeated integration in the $\,y^{\nh1}$ direction, with any prescribed 
$\,\rho$.

\section{Final remarks}\label{fire} 
\setcounter{equation}{0} 
Let $\,(M,g)\,$ be a strictly \nw\ self-du\-al neutral Ein\-stein 
\hbox{four\hh-}\hskip0ptman\-i\-fold of Pe\-trov type III. The geometry of 
$\,g\,$ gives rise, at least locally, to further structures in $\,M$, with 
properties arguably simpler than those assumed for $\,g$. One such structure 
is the basic octuple 
$\,(M,\mathcal{V}\nh,\hh\mathrm{D},h,\alpha,\beta,\theta,\zeta)\,$ determined 
by $\,(M,g)\,$ as in Theorem~\ref{chboc}. A part of it is the af\-fine 
foliation $\,(\mathcal{V}\nh,\hh\mathrm{D})\,$ on $\,M$, which allows us to 
identify $\,M$, locally, with the total space of an af\-fine plane bundle 
$\,\apb\,$ over a surface $\,\bs\,$ (see Section~\ref{pmaf}). As in the proof 
of Theorem~\ref{unimo}, setting $\,\alb=\phi^{-1}(0)$, with $\,\phi\,$ chosen 
as at the end of Section~\ref{boct}, we obtain an 
af\-\hbox{fine\hh-}\hskip0ptline sub\-bun\-dle $\,\alb\,$ of $\,\apb$, 
which is well defined, since $\,\phi$, although nonunique, is nevertheless 
unique up to multiplication by functions $\,\bs\to(0,\infty)$. Denoting by 
$\,\pb\,$ and $\,\lb\,$ the vector bundles over $\,\bs$ associated with 
$\,\apb\,$ and $\,\alb$, we may treat functions $\,\phi\,$ mentioned above as 
positive sections of the ca\-non\-i\-cal\-ly-orient\-ed line bundle 
$\,[\pb\nh/\nh\lb]^*\nnh$. The conclusions about constancy along 
$\,\mathcal{V}\,$ in (iii) at the end of Section~\ref{boct} now allow us to 
treat $\,\beta,\zeta\,$ and $\,\theta\,$ as sections of the vector bundles 
$\,\tab\otimes[\pb\nh/\nh\lb]^{\otimes2}\nnh$, 
$\,[\tab]^{\wedge2}\otimes[\pb\nh/\nh\lb]\,$ and 
$\,([\pb]^{\wedge2}\otimes[\pb\nh/\nh\lb]^{\otimes2})^*$ over $\,\bs\,$ (all 
three pulled back, via the bundle projection $\,M\to\bs$, to the corresponding 
pull\-back bundles over $\,M$). Thus, $\,\beta\,$ viewed as a 
$\,[\pb\nh/\nh\lb]^{\otimes2}\nnh$-val\-ued $\,1$-form on $\,\bs\,$ 
distinguishes the \hbox{one\hh-}\hskip0ptdi\-men\-sion\-al distribution 
$\,\mathrm{Ker}\,\beta\,$ on $\,\bs$.

Despite being natural geometric invariants of $\,g$, the objects just listed 
do not play a prominent role in our lo\-cal-struc\-ture result 
(Theorem~\ref{lcstr}). On the contrary, most ingredients of the construction 
appearing in Theorem~\ref{lcstr}(a), such as the 
tw\hbox{o\hh-\nh}\hskip0ptplane system $\,(\bs,\xi,\ts,\plane,c,\varOmega)$, 
cannot be canonically recovered from the resulting metric $\,g$. (The extent 
to which $\,(\bs,\xi,\ts,\plane,c,\varOmega)$ fails to be unique is made 
explicit in the proof of Theorem~\ref{unimo}.)

The ideal form of a structure theorem for a class of metrics would be one 
involving a construction that uses only natural invariants of the metrics in 
question. As explained above, for metrics discussed in this paper such a goal 
still appears elusive.

\end{document}